# Homogenization and Shape Differentiation of Quasilinear Elliptic Equations

## Homogeneización y Diferenciación de Formas de Ecuaciones Elípticas Cuasilineales

**David Gómez Castro**

Director: Prof. Jesús Ildefonso Díaz Díaz

Dpto. de Matemática Aplicada &
Instituto de Matemática Interdisciplinar
Facultad de Matemáticas
Universidad Complutense de Madrid

Esta tesis se presenta dentro del
*Programa de Doctorado en Ingeniería Matemática, Estadística e Investigación Operativa*

Diciembre 2017

Quisiera dedicar esta tesis a mi abuelo, Ángel Castro,
con el que tantos momentos compartí de pequeño.

# Acknowledgements / Agradecimientos

The next paragraphs will oscillate between Spanish and English, for the convenience of the different people involved.

Antes de empezar, me gustaría expresar mi agradecimiento a todas aquellas personas sin las cuales la realización de esta tesis no hubiese sido posible o, de haberlo sido, hubiese resultado en una experiencia traumática e indeseable. En particular, me gustaría señalar a algunas en concreto.

Lo primero agradecer a mis padres (y a mi familia en general) su infinito y constante apoyo, sin el cual, sin duda, no estaría donde estoy.

A Ildefonso Díaz, mi director, quién hace años me llamó a su despacho (entonces en el I.M.I.) y me dijo que, en lugar de la idea sobre aerodinámica que se me había ocurrido para el T.F.G., quizás unas cuestiones en las que él trabajaba sobre "Ingeniería Química" podían tener más proyección para nuestro trabajo conjunto. Nadie puede discutirle hoy cuánta razón tenía. Ya desde aquel primer momento he seguido su consejo, y creo que me ha ido bien. Además de su guía temática (y bibliográfica), le quedaré por siempre agradecido por haberme presentado a tantas (y tan célebres) personas, que han enriquecido mis conocimientos, algunas de las cuales se han convertido en colaboradores, y a las que me refiero a más adelante.

Al Departamento de Matemática de Aplicada de la Universidad Complutense de Madrid, que me acogió durante el desarrollo de esta tesis. En particular agradecer a su director, Aníbal Rodríguez Bernal su esfuerzo para hacer posible mi docencia. También, al Instituto de Matemática Interdisciplinar y al joven programa de doctorado IMEIO, en el cual esta tesis será la primera del Dpto. de Matemática Aplicada, la organización de productivos cursos de doctorados. Por último, agradecer muy especialmente a Antonio Brú su compañía y consejo en tantas comidas y largas charlas.

Now, allow me to thank, in English, my many collaborators. I hope our joint work continues for many years.

To Prof. Häim Brezis, who graciouslly accepted me as a visitor during the months of April to July 2017 at the Technion. I learnt from him not only a lot of Mathematics, but also about Jewish culture and history, which I found fascinating.

To Prof. Tatiana Shaposhnikova (and her students Alexander Podolskii and Maria Zubova), whose graduate course in Madrid opened for me the field of critical size homogenization, which has been tremendously fruitful. I would like to thank her specially for her patience, and her attention to detail.

To Prof. Claudia Timofte, who first introduced me to the theory of homogenization.

To Prof. Jean Michel Rakotoson, who has shown to be an insightful and generous collaborator, from whom I have learnt a great deal. Our several deep discussions have enlightened me. Also, to Prof. Roger Temam, whose "big picture" view is invaluable.

Por último, a Cheri. Por todo.

# Abstract

This thesis has been divided into two parts of different proportions. The first part is the main work of the candidate. It deals with the optimization of chemical reactors, and the study of the effectiveness, as it will explained in the next paragraphs. The second part is the result of the visit of the candidate to Prof. Häim Brezis at the Israel Institute of Technology (Technion) in Haifa, Israel. It deals with a particular question about optimal basis in $L^2$ of relevance in Image Proccesing, which was raised by Prof. Brezis.

The first part of the thesis, which deals with chemical reactors, has been divided into four chapters. It studies well-established models which have direct applications in Chemical Engineering, and the notion of “effectiveness of a chemical reactor”. One of the main difficulties we faced is the fact that, due to the Chemical Engineering applications, we were interested in dealing with root-type nonlinearities.

The first chapter focuses on modeling: obtaining a macroscopic (homogeneous) model from a prescribed microscopic behaviour. This method is known as *homogenization*. The idea is to consider periodically repeated particles of a fixed shape $G_0$, at a distance $\varepsilon$, which have been rescaled by a factor $a_\varepsilon$. This factor is usually of the form $a_\varepsilon = C_0\varepsilon^\alpha$, where $\alpha \geq 1$ and $C_0$ is a positive constant. The aim is to study the different behaviours as $\varepsilon \to 0$, when the particles are no longer considered. It was known that depending of this factor there are usually different behaviours as $\varepsilon \to 0$. First, the case of *big particles* and *small particles* are treated differently. The latter, which have been the main focus of this chapter, are divided into subcritical, critical and supercritical holes. Roughly speaking, there is a critical value $\alpha^*$ such that the behaviours $\alpha = 1$ (big particles), $1 < \alpha < \alpha^*$(subcritical particles), $\alpha = \alpha^*$ (critical particles) and $\alpha > \alpha^*$ (supercritical particles) are significantly different.

The main focus of the thesis has been in the cases $\alpha > 1$, although some new results for the case $\alpha = 1$ have been obtained (see [DGCT15; DGCT16]). In the subcritical cases we have significantly improved the regularity of the nonlinearities that are allowed, by applying uniform approximation arguments (see [DGCPS17d]). We also proved that, when the diffusion depends on the gradient ($p$-Laplacian type) with $p$ greater than the spatial dimension, no critical scales exists (see [DGCPS17b]). Also, this thesis includes some unpublished estimates which give a unified study of these cases, and provide some new insights. The newest and most relevant results in this sections are the ones obtained for the critical case. The state of the art in this field was dealing only with the case in which the shape of the particles (or holes), $G_0$, is a ball. In this direction we have shown that the case of maximal monotone graphs behaves as expected, providing a common roof for results with Neumann, Dirichlet and even Signorini boundary conditions (see [DGCPS17a; DGCPS17c]). In this case we have shown that the "strange" nonlinear term appearing in the homogeneous problem is always smooth, even when the microscopic problem is not. This behaviour can be linked to Nanotechnological properties of some materials. Furthermore, we have managed, for the time in the literature, to study the cases in which $G_0$ is not a ball (see [DGCSZ17]). In this last paper, the techniques of which are very new, gives some seemingly unexpected results, that answer the intuition of the experts. The results mentioned above were obtained by applying a modification of Tartar's oscillating test function method. The periodical unfolding method has also been applied by this candidate, in some unpublished work, and this was acknowledged by the authors of [CD16].

The contributions presented in this chapter improve many different works in the literature, and it has been presented in Table 1.1. The work in this chapter has been presented in the international congress ECMI 2016 (Spain, 2016) and Nanomath 2016 (France, 2016). At the time of presentation of this thesis a new paper dealing with the critical case and general shape of the particles has been submitted for publication.

The second chapter deals with a priori estimates for the effectiveness factor of a chemical reactor, which is a functional depending of the solution of the limit behaviour deduced in Chapter 1. This problems comes motivated, for example, by the application to waste water treatment tanks. Once we have obtained a homogeneous model, our aim is to decide which reactors are of optimizer in some classes, and also provide bounds for the effectiveness. In this direction, we have dealt with Steiner symmetrization, which allows us to compare the solution of any product type domain $\Omega_1 \times \Omega_2$, which would represent the chemical reactor, with a cylinder of circular basis, $B \times \Omega_2$, where $|B| = |\Omega_1|$. In this direction we have published two papers [DGC15b; DGC16] dealing with convex and concave kinetics, which extends the pioneering paper [ATDL96]. The work in this chapter has been presented in the following congresses: MathGeo 2013 (Spain), 10th AIMS Conference in Dynamical Systems, Differential Equations and Applications (Spain, 2013), Nanomath 2014 (Spain, 2014), Mini-workshop in honour of Prof. G. Hetzer (USA, 2016).

The third chapter deals with direct shape optimization techniques. This can be organized into two sections: shape differentiation and convex optimization. Shape differentiation is a technique that, given a initial shape $\Omega_0$ characterizes the infinitesimal change of the solution of our homogeneous problem when we consider deformations $(I+\theta)(\Omega_0)$. In this directions two papers have been published. First, we studied the Fréchet differentiable case, that requires the kinetic to be twice differentiable (see [DGC15a]). This was a first step to the problem in which we were interested, the case of non-smooth kinetic. In this setting, the solution may even develop a dead core (see [GC17]). Another of the techniques of direct optimization we applied was the convex optimization of the domain $G_0$. If we only allow the admissible set of shapes $G_0$ to be convex sets, then we have some compactness results, that guaranty that there exist optimal sets in this family (see [DGCT16; DGCT15]). The work in this chapter has been presented in the 11th AIMS Conference on Dynamical Systems, Differential Equations and Applications (USA, 2017).

The fourth chapter deals with linear elliptic equations with a potential, $-\Delta u + Vu = f$, where the potential, $V$, "blows up" near the boundary. This kind of equations appear as a result of the shape differentiation process, in the non-smooth case. The problems with a transport term, $\nabla u \cdot b$, $-\Delta u + \nabla u \cdot b + Vu$ is studied, in collaboration with Profs. Jean-Michel Rakotoson and Roger Temam. Different results of existence, uniqueness and regularity of solutions of this equations are presented (see [DGCRT17]). One of this results is the fact that, shall the blow up of the potential $V$ be fast enough, the condition $Vu \in L^1$ can act as a boundary condition for $u$. Some unpublished results are included in this chapter, which improve some of the results of [DGCRT17], in a limit case, by applying an extension of the argument in [DGCRT17] suggested recently by Prof. Brezis to this candidate, and which have not been published. The work in this chapter has been presented in the 11th AIMS Conference on Dynamical Systems, Differential Equations and Applications (USA, 2017). At the time of submission of this thesis a new paper improving the results of [DGCRT17] is under development.

The fifth chapter develops the second part of the thesis, and includes results obtained during the 2017 visit to Prof. Brezis (see [BGC17]). They improve some previous results by Brezis in collaboration with the group of Prof. Ron Kimmel. We showed that the basis of eigenvalues of $-\Delta$ with Dirichlet boundary conditions is the unique basis to approximate functions in $H_0^1$ in $L^2$ in an optimal way.

Besides the contributions in this thesis, this candidate has also developed other projects. The author published, jointly with Prof. Brú and Nuño, a paper [BGCN17] which involved the simulation of the fractional Laplacian in bounded domains, and their study by new statistical physics techniques. Besides this work, the candidate studied the modeling of Lithium-ion batteries, and is about to publish work on the well-posedness of the Newman model. This was presented in the 11th AIMS conference on Dynamical Systems (USA, 2017).

# Resumen en castellano

Esta tesis se ha divido en dos partes de tamaños desiguales. La primera parte es la componente central del trabajo del candidato. Se encarga de la optimización de reactores químicos de lecho fijo, y el estudio de su efectividad, como se expondrá en los siguientes párrafos. La segunda parte es el resultado de la visita del candidato al Prof. Häim Brezis en el Instituto Tecnológico de Israel (Technion) en Haifa, Israel. Se entra en una pregunta concreta sobre bases óptimas en $L^2$, que es de importancia en Tratamiento de Imágenes, y que fue formulado por el Prof. Brezis.

La primera parte de la tesis, que estudio reactores químicos, se ha dividido en 4 capítulos. Estudia un modelo establecido que tiene aplicaciones directas en Ingeniería Química, y la noción de efectividad. Una de las mayores dificultades con la que nos enfrentamos es el hecho que, por las aplicaciones en Ingeniería Química, estamos interesados en reacciones de orden menor que uni (de tipo raíz).

El primer capítulo se centra en la modelización: obtener un modelo macroscópico (homogéneo) a partir de un comportamiento microscópico prescrito. A este método se le conoce como *homogeneización*. La idea es considerar partículas periódicamente repetidas, de forma fija $G_0$, a una distancia $\varepsilon$, y que han sido reescaladas por un factor $a_\varepsilon$. La expresión habitual de este factor es $a_\varepsilon = C_0\varepsilon^\alpha$, donde $\alpha \geq 1$ y $C_0$ es una constante positiva. El objetivo es estudiar los diferentes comportamientos cuando $\varepsilon \to 0$, y ya no se consideran las partículas. Primero, los casos de *partículas grandes* y *partículas pequeños* se tratan de formas distintas. Este segundo, que ha sido el central en esta tesis, se divide en subcrítico, crítico y supercrítico. En términos generales, existe un valor $\alpha^*$ tal que los comportamientos de los casos $\alpha = 1$ (partículas grandes), $1 < \alpha < \alpha^*$ (partículas subcríticos), $\alpha = \alpha^*$ (partículas críticos) y $\alpha > \alpha^*$ (partículas supercríticos) son significativamente distintos. El objetivo central de la tesis han sido los casos $\alpha > 1$, aunque se han obtenido también algunos resultados para el caso $\alpha = 1$ (ver [DGCT15; DGCT16]). En el caso subcrítico hemos mejorado significativamente la regularidad de las no-linealidades permitidas, por argumentos de aproximación uniforme (ver [DGCPS17d]). También hemos demostrado que, cuando la difusión depende del gradiente (operadores de tipo $p$-Laplaciano) con $p$ mayor que la dimensión espacial, entonces no existen escalas críticas (ver [DGCPS17b]). Además, esta tesis incluye estimaciones no publicadas que dan un estudio unificado de estos casos, e introducen nuevas perspectivas. Los resultados más nuevos y más relevantes en estas secciones son los que se refieren al caso crítico. El estado del arte en este campo era lidiar sólo con el caso en que la forma de la partícula es una esfera. En esta dirección hemos mostrado que los grafos maximales monótonos se comporta como era de esperar, dando un techo común a los resultados con condiciones de frontera Neumann, Dirichlet e incluso Signorini (ver [DGCPS17a; DGCPS17c]). En este caso hemos demostrado que el término "extraño" en el problema homogéneo es siempre regular, incluso cuando la no-linealidad del problema microscópico no lo es. Este comportamiento se puede enlazar con propiedes Nanotecnológicas de algunos materiales. También hemos conseguido estudiar el caso en las partículas no son esferas, si no que tienen una forma más general (ver [DGCSZ17]). En este último artículo, que emplea las técnicas muy nuevas, se dan algunos resultados de aspecto aparentemente inesperado, pero que satisfacen la intuición de los expertos. Todo los resultados presentados en el texto precedente son obtenidos utilizando modificaciones del método de funciones test oscilantes de Tartar. El método de desdoble periódico (*periodical unfolding* en inglés) también ha sido usado por el candidato, en un trabajo sin publicar, y ha sido reconocido en los agradecimientos de [CD16]. Las contribuciones de este capítulo mejoran muchos trabajos previous, como se ha presentado en la Tabla 1.1. El trabajo de este capítulo se ha presentado en los congresos internacionales ECMI 2016 (Spain) y Nanomath 2016 (France).

El segundo capítulo trata sobre estimaciones a priori del factor de efectividad de las reacciones químicas: un funcional que depende de la solución del problema límite obtenido en el Capítulo 1. Este problema viene motivado, por ejemplo, por la aplicación en reactores de tratamiento de aguas residuales. Una vez que se ha obtenido el problema homogeneizado, nuestro objetivo es decidir qué reactores son optimizadores de este funcional, y dar cotas para la efectividad. En este sentido, hemos trabajado con la optimización de Steiner, que permite comparar reactores de la forma $\Omega_1 \times \Omega_2$ con reactores cilíndricos de la forma $B \times \Omega_2$, donde $|B| = |\Omega_1|$. En esta dirección se han publicado dos trabajos, [DGC15b; DGC16] lidiando con no-linealidades convexas y cóncavas. El trabajo de este capítulo se ha presentado en los siguientes congresos: MathGeo 2013 (España), 10th AIMS Conference in Dynamical Systems, Differential Equations and Applications (España, 2013), Nanomath 2014 (España), Mini-workshop in honour of Prof. G. Hetzer (USA, 2016).

El tercer capítulo trata con técnicas de optimización de formas directas. Se ha organizado en dos secciones: diferenciación de formas y optimización convexa. La diferenciación de formas es una técnica que, dado un una forma inicial $\Omega_0$, caracteriza el cambio infinitesimal de la solución de nuestro problema homogeneizado cuando se considera una deformación $(I+\theta)(\Omega_0)$. En esta dirección se han publicado dos artículos. Primero hemos estudiado la diferenciabilidad en el sentido de Fréchet, que requiere que la cinética sea dos veces derivable (ver [DGC15a]). Éste fue un primer paso hacia el problema en que estábamos interesados, el caso no-suave. En este contexto, la solución puede desarrollar un *dead core* (ver [GC17]). Otra de las técnicas que hemos usado es la optimización convexa directa del dominio $G_0$. Si solo consideramos el conjunto admisible de formas $G_0$ dentro de la familia convexa, entonces podemos obtener la existencia de extremos (ver [DGCT16; DGCT15]). Este trabajo ha sido presentado en el 11th AIMS Conference in Dynamical Systems, Differential Equations and Applications (USA, 2017).

El cuarto capítulo trata con ecuaciones elípticas con un potencial, $-\Delta u + Vu = f$, donde el potencial, $V$, "explota" cerca del borde. Este tipo de ecuaciones aparecen como resultado del proceso de diferenciación de formas del Capítulo 3, en el caso en que aparece un dead core. El problema con un término de transporte, $\vec{b} \cdot \nabla u$, fue también estudiado. Se obtuvieron diferentes resultados de existencia, regularidad y unicidad de soluciones (ver [DGCRT17]). Uno de los resultados más sorprendentes es que si $V$ explota suficientemente rápido, entonces la condición $Vu \in L^1$, que se suponía habitualmente como púramente técnica, se convierte en una condición de contorno Dirichlet homogénea. Los resultados expuestos en esta tesis se han presentado en el 11th AIMS Conference in Dynamical Systems, Differential Equations and Applications (USA, 2016). Se incluyen en esta tesis algunos resultados no publicados, sugeridos por Häim Brezis, que mejoran a los publicados en algunos casos. En el momento del déposito se está trabajando en un borrador que mejora, aún más, estos resultados.

El quinto capítulo desarrolla la segunda parte de la tesis, e incluye resultados obtenidos durante la visita en 2017 al Prof. Häim Brezis (ver [BGC17]). Se mejoran algunos resultados previos con el grupo de Ron Kimmel, sobre la existencia y unicidad de bases óptimas para representación de funciones $H^1$ en $L^2$.

Además de las contribuciones incluídas en esta tesis, el candidato ha desarrollado otros proyectos. El candidato ha publicado, conjuntamente con Antonio Brú y Juan Carlos Nuño, un artículo [BGCN17] que incluye la simulación numérica de un laplaciano fraccionario en un dominio acotado, y su estudio mediante técnicas de física estadística. Además de este trabajo se ha estudiado la modelización de baterías de ion-Litio, y se va a publicar un trabajo sobre la buena formulación del modelo de Newman (problema abierto desde 1971). Este último trabajo se presentó en el congreso 11th de la AIMS (USA, 2016).

# Table of contents

# List of figures

# Notation

### Roman Symbols

$a_\varepsilon$ The size of particles in the homogenization process. Usually $a_\varepsilon = C_0\varepsilon^\alpha$

$B(x,R)$ Ball centered at $x$ of radius $R$

$\mathscr{C}^n(\Omega)$ where $n \in \mathbb{N} \cup \{\infty\}$. Space of $n$ times differentiable functions with continuous derivative.

$\mathscr{C}_c^n(\Omega)$ Space of $n$ times differentiable functions with continuous derivative and compact support.

$\mathscr{C}_{per}^n(Q)$ where $Q$ is an n-dimensional cube, space of functions that can be extended by periodicity to a function in $C^n(\mathbb{R}^n)$

$\mathscr{E}$ Effectiveness factor

$F(A)$ Space of all functions $f : A \to \mathbb{R}$

$g(u)$ Kinetic of the homogeneous problem for $u$. Chapters 2, 3 and 4

$\mathscr{A}$ $= \left(\frac{n-p}{p-1}\right)^{p-1} C_0^{n-p}\omega_n$

$\mathscr{B}_0$ $= \left(\frac{n-p}{C_0(p-1)}\right)^{p-1}$

$n$ Spatial dimension

$W^{1,p}(\Omega)$ Space of functions such that $u \in L^p(\Omega)$ and $\nabla u \in (L^p(\Omega))^n$.

$W_0^{1,p}(\Omega)$ Closure in $W^{1,p}(\Omega)$ of the space $\mathscr{C}_c^\infty(\Omega)$

$W^{1,p}(\Omega,\Gamma)$ where $\Gamma \subset \partial\Omega$. Closure in $W^{1,p}(\Omega)$ of the space of functions $f \in \mathscr{C}^\infty(\bar{\Omega})$ such that $\operatorname{supp} f \cap \Gamma = \emptyset$ $W^{1,p}_{loc}(\Omega)$ space of functions such that they are in $W^{1,p}(K)$, for every $K \subset \Omega$ compact

$W^{1,p}_{per}(Q)$ where $Q$ is an n-dimensional, space of functions in $W^{1,p}$ that are periodic in $Q$, i.e. they can be extended by periodicity to a function in $W^{1,p}_{loc}(\mathbb{R}^n)$

$Y$ In Chapter 1 it will represent a cube of side 1. Either, $Y = [0,1]^n$ or $(-\frac{1}{2}, \frac{1}{2})^n$ depending on the case

**Greek Symbols**

$\alpha$ Power such that $a_\varepsilon = C_0 \varepsilon^\alpha$

$\beta(\varepsilon)$ Corrector order of the boundary term in homogenization. Usually $\beta(\varepsilon) = \varepsilon^{-\gamma}$. Chapter 1

$\beta(w)$ Kinetic of the homogeneous problem for $w$. Chapter 2, 3 and 4

$\beta^*(\varepsilon)$ Critical value of $\beta(\varepsilon)$. $\beta^*(\varepsilon) = a_\varepsilon^{1-n}\varepsilon^n$. Chapter 1.

$\partial\Omega$ Boundary of the set $\Omega$

$\emptyset$ Empty set

$\varepsilon$ Small parameter destined to go to zero

$\eta$ Ineffectiveness factor

$\Omega$ Generic open set. Its smoothness will be specified on a case by case basis

$\omega_n$ Volume of the $n$-dimensional ball

$\partial f$ Partial derivative of function $f$

$\sigma(u)$ Kinetic of the problem for $u$. Chapter 1.

**Superscripts**

$*$ In Chapter 1 critical value. In Chapter 2 and onwards decreasing rearrangement

$\star$ Schwarz rearrangement

# Steiner rearrangement

**Subscripts**

$\varepsilon$ In functions or sets will indicate that it refers to the nonhomogeneous problem.

$*$ Nondescript rearrangement

**Other Symbols**

$\overline{\Omega}$ Closure of the set of $\Omega$

$\Subset$ Compact subset

$\preceq$ Comparison of concentrations

$\emptyset$ Empty set

$\ll$ given sequences $(a_\varepsilon)_{\varepsilon>0}$, $(b_\varepsilon)_{\varepsilon>0}$. $a_\varepsilon \ll b_\varepsilon$ if $\lim_{\varepsilon\to 0} a_\varepsilon b_\varepsilon^{-1} = 0$

$\gg$ given sequences $(a_\varepsilon)_{\varepsilon>0}$, $(b_\varepsilon)_{\varepsilon>0}$. $a_\varepsilon \gg b_\varepsilon$ if $b_\varepsilon \ll a_\varepsilon$

$\sim$ given sequences $(a_\varepsilon)_{\varepsilon>0}$, $(b_\varepsilon)_{\varepsilon>0}$. $a_\varepsilon \sim b_\varepsilon$ if $\lim_{\varepsilon\to 0} a_\varepsilon b_\varepsilon^{-1} \in (0,+\infty)$

**Acronyms / Abbreviations**

CIF Cauchy's Integral Formula

DCT Dominated Convergence Theorem

div Divergence operator. For a differentiable function $f : \mathbb{R}^n \to \mathbb{R}$ its definition is $\operatorname{div} f = \sum_{i=1}^{n} \frac{\partial f}{\partial x_i}$

ODE Ordinary Differential Equation

PDE Partial Differential Equation

# Part I

# Optimization of chemical reactors

# Introduction

## Rutherford Aris and the definition of effectiveness.

The monographs of R. Aris [Ari75] and Aris and Strieder [AS73] model the behaviour of chemical reactions in terms of partial differential equations. Their works are amongst the first to consider the microscopic behaviour of the system, and derive from it the macroscopic properties.

In their book, Aris and Strieder model a Chemical Reactor by an open set $\Omega$. In it, they model the concentration of a Chemical Reactant by a spatial function $c = c(x)$. For the constituent equation they introduce an spatial diffusion term $\operatorname{div}(D_e \nabla c)$ and a reaction term $r(c)$ (the amount of reaction that is produced as a function of the amount of reactant). Their spatial model results:

$$\operatorname{div}(D_e \nabla c) = r(c) \text{ in } \Omega.$$

This equation alone is ill-posed, since there are many solutions of this problem. To fix a single one another equation needs to be considered. The author choose to allow a flux in the boundary of $\partial\Omega$. The full model results:

$$\begin{cases} \operatorname{div}(D_e \nabla c) = r(c) & \text{in } \Omega, \\ D_e \vec{n} \cdot \nabla c = k_c(c_f - c), & \text{on } \partial\Omega, \end{cases} \tag{1}$$

where $k_c$ is a permitivity constant of the boundary and $c_f$ is a maximum concentration of the reactant admitted by the solvent.

One of the novelties in the mentioned book is the they also propose a model, which we will call *non homogeneous*, in which the reactor contains many microscopic particles, which they model by an open domain $G$ of $\mathbb{R}^n$ (usually $n = 2$ or $3$). In this model that the reaction is taking place on the boundary of the particles. The following model comes up:

$$\begin{cases} \operatorname{div}(D_e \nabla c) = 0, & \hat{\Omega} \\ D_e \vec{n} \cdot \nabla c = k_c(c_f - c), & \partial\Omega \\ D_e \vec{n} \cdot \nabla c = \hat{r}(c), & \partial G \end{cases} \tag{2}$$

where $\hat{\Omega} = \Omega \setminus G$ and $G$ represents pellets and $\hat{r}$ is a reaction rate, possibly different from $r$. They represent the situation as Figure 1.

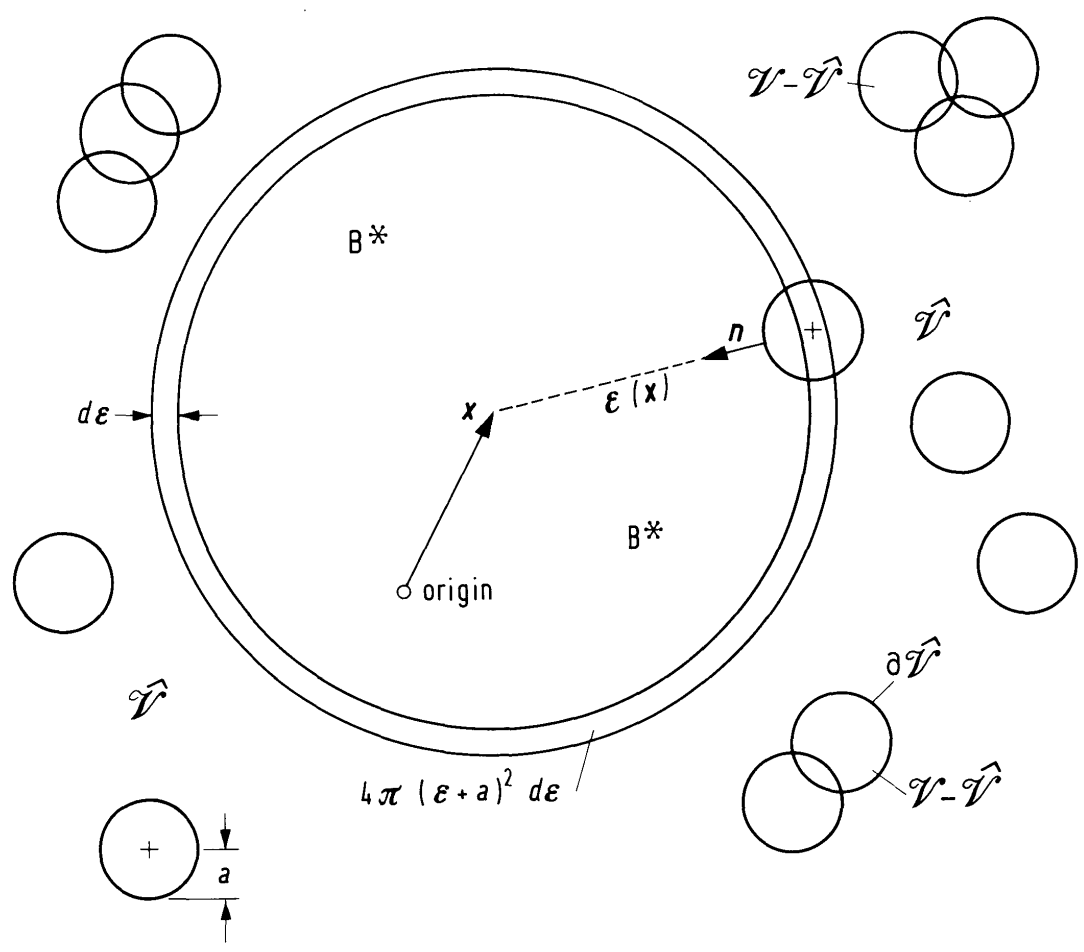

Fig. 1.5.1. Random spheres of radius $a$. The void region $\hat{\mathscr{V}}$, reaction zone $\mathscr{V} - \hat{\mathscr{V}}$, the interface $\partial\hat{\mathscr{V}}$, and the unit normal vector $\mathbf{n}$ to $\partial\hat{\mathscr{V}}$ are shown. The minimum distance from the point $x$ in the void to the reactive interface is $\varepsilon(x)$. The regions about the point $x$ of radius $\varepsilon + a$ from which spheres are excluded, and the adjacent shell of thickness $d\varepsilon$ in which at least one sphere center is found are also indicated.

Fig. 1 Domain as portrayed in [AS73]. In a small change of notation we have considered $\Omega$ instead of $V$ and $\hat{\Omega}$ instead of $\hat{V}$.

Aris and Strieder define the effectiveness factor of the chemical reactor as

$$\mathscr{E} = \frac{1}{|\Omega| r(c_f)} \int_\Omega r(c)$$

for the homogeneous problem, where $c$ is the solution of (1) and

$$\hat{\mathscr{E}} = \frac{1}{|\partial\hat{\Omega}| \hat{r}(c_f)} \int_{\partial\hat{\Omega}} \hat{r}(c)$$

for the non-homogeneous model, where $c$ is the solution of (2).

Another novelty in the work of Aris and Strieder is that, albeit by naive methods, they show which model of type (1) we must consider once we if we consider a constituent equation of form (2) (and viceversa).

This effectiveness factor is a very relevant quantity. There is a lot of mathematical literature dedicated to it (see, e.g., [BSS84; DS95]). It will be the main quantity under investigation in Chapters 2 and 3, through very different techniques.

Amongst other things, they were interested in the behaviour of the effectiveness in the different domains and, in particular, in choosing domains of optimal effectiveness.

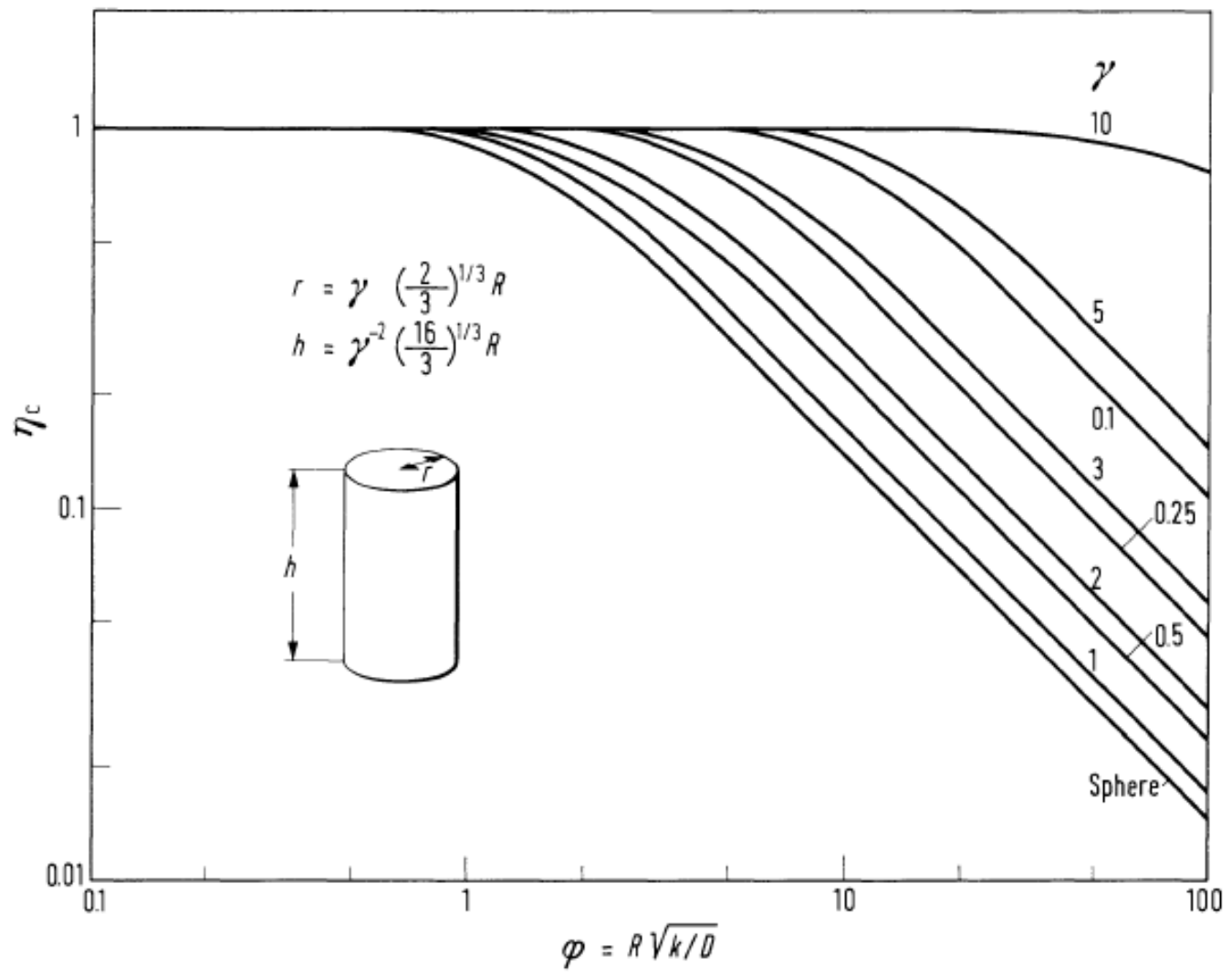

Fig. 4.5.1. Plot of effectiveness factor $\eta_c$ for a cylinder vs. $\phi$ where $R$ is the radius of a sphere having the same volume. (Note that this $\phi$ is $R\phi_c$).

Fig. 2 Image showing some (probably estimated) curves of the effectiveness factor $\mathscr{E}$ as a function of the shape parameters of the cylinder. Extract from [AS73]

Roughly speaking, the aim of Chapter 1 is to properly define the set $\hat{\Omega}$, to study in which sense we can pass from the equation over $\hat{\Omega}$ to the equation over $\Omega$ and in which sense we can pass from $\hat{\mathscr{E}}$ to $\mathscr{E}$.

# A comment on the notation

The first part of this thesis applies techniques for different problems from different fields of Applied Mathematics, which involve different communities. Whenever possible, we have tried to be consistent along the paper, but in some cases this would have made reading more inconvenient for some of the specialists. In particular the use of $\sigma, g$ and $\beta$ changes between Chapter 1 and the rest.

In Chapter 1 we define $w_\varepsilon$ as the solution of (1.9), and the $u_\varepsilon = 1 - w_\varepsilon$ as the solution of (1.12). Then, under some assumptions, we show that $\tilde{u}_\varepsilon$ to $u$, which is the solution of (1.180). This will be the relevant function studied in chapter 2 and onwards. In this setting we define the effectiveness as (1.14) and (1.15). Then, in Chapter 2 and onwards, $w$ is the solution of (2.1), whereas $u = 1 - w$ is the solution of (2.2), and the effectiveness is defined as (2.3) and (2.4).

# Chapter 1

# Deriving macroscopic equations from microscopic behaviour: Homogenization

## 1.1 Formulation of the microscopic problem

Let us present the precise mathematical formulation of the problem we will be interested in, and that is directly motivated by the problem proposed by Aris.

### 1.1.1 Open domain with particles

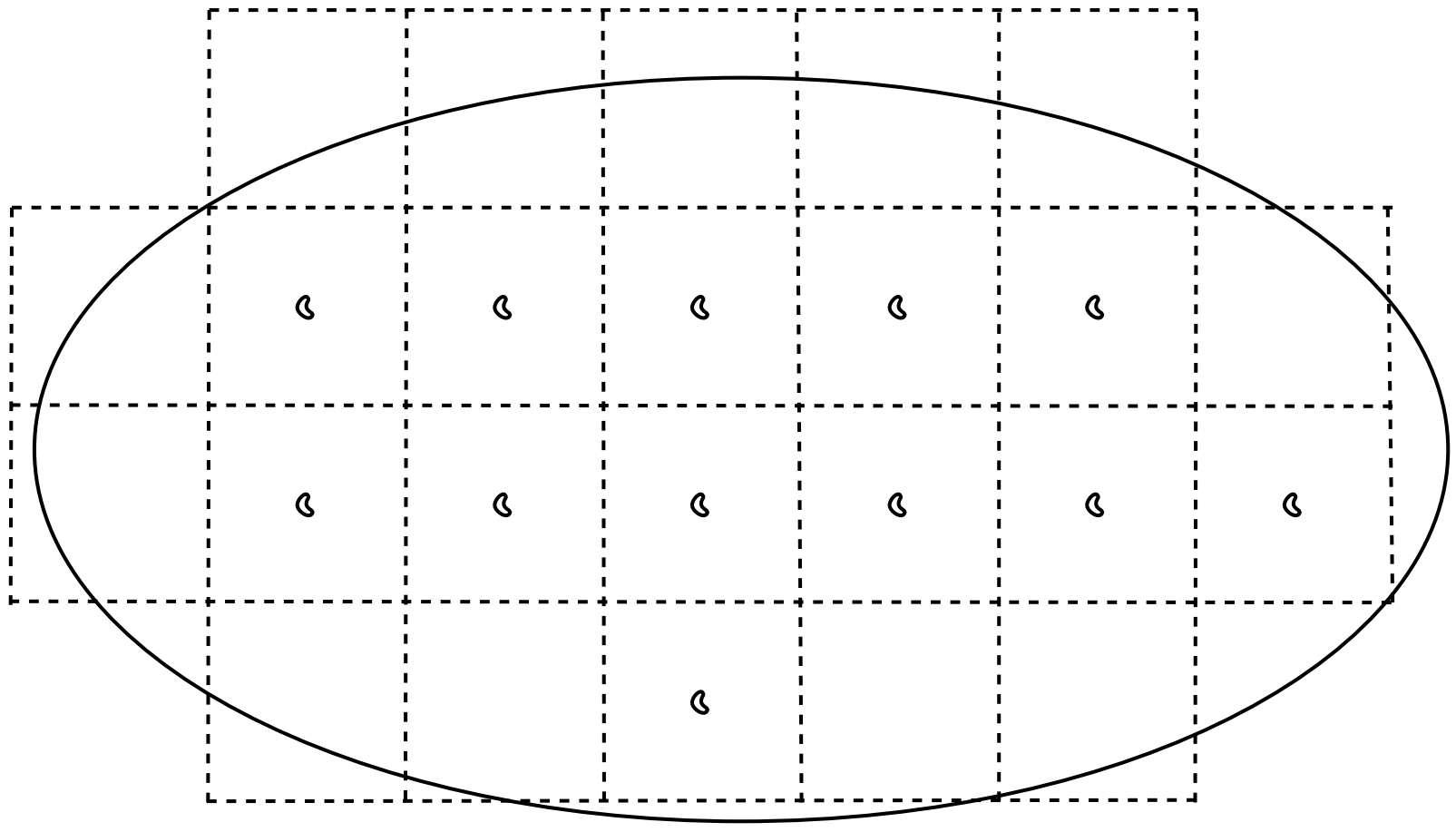

Fig. 1.1 The domain $\Omega_\varepsilon$.

Let us set up the geometrical framework. Let $\Omega \subset \mathbb{R}^n$ be an open set (bounded and regular, for simplicity), and let the shape of a generic inclusion (in our setting a particle, but it applies also to the case of a hole) be represented by a domain $G_0$ be an open set homeomorphic to a ball such that $\overline{G_0} \subset Y = (-\frac{1}{2}, \frac{1}{2})^n$ (i.e., there exists a invertible continuous map $\Psi : U \to V$ between open sets of $\mathbb{R}^n$, $U$ and $V$, where $G_0 \subset U$ and $V$ contains the open ball of radius one, $\Psi(G_0)$ is the ball and $\Psi^{-1}$ is continuous).

Considering the parameter $\varepsilon > 0$ the distance between the equispaced particles, we will typically be set in the following geometry

$$G_i^\varepsilon = \varepsilon i + a_\varepsilon G_0, \qquad i \in \mathbb{Z}^n \tag{1.1}$$

$$Y_i^\varepsilon = \varepsilon i + \varepsilon Y \tag{1.2}$$

$$\Upsilon_\varepsilon = \left\{ i \in \mathbb{Z}^n : \overline{Y_i^\varepsilon} \subset \Omega \right\}, \tag{1.3}$$

$$G^\varepsilon = \bigcup_{i \in \Upsilon_\varepsilon} G_i^\varepsilon, \tag{1.4}$$

$$S^\varepsilon = \bigcup_{i \in \Upsilon_\varepsilon} \partial G_i^\varepsilon, \tag{1.5}$$

$$\Omega_\varepsilon = \Omega \setminus \overline{G^\varepsilon}. \tag{1.6}$$

We will sometimes consider that

$$a_\varepsilon = C_0 \varepsilon^\alpha. \tag{1.7}$$

The parameter $\alpha \geq 1$ indicates the size of the particle relative to the repetition.

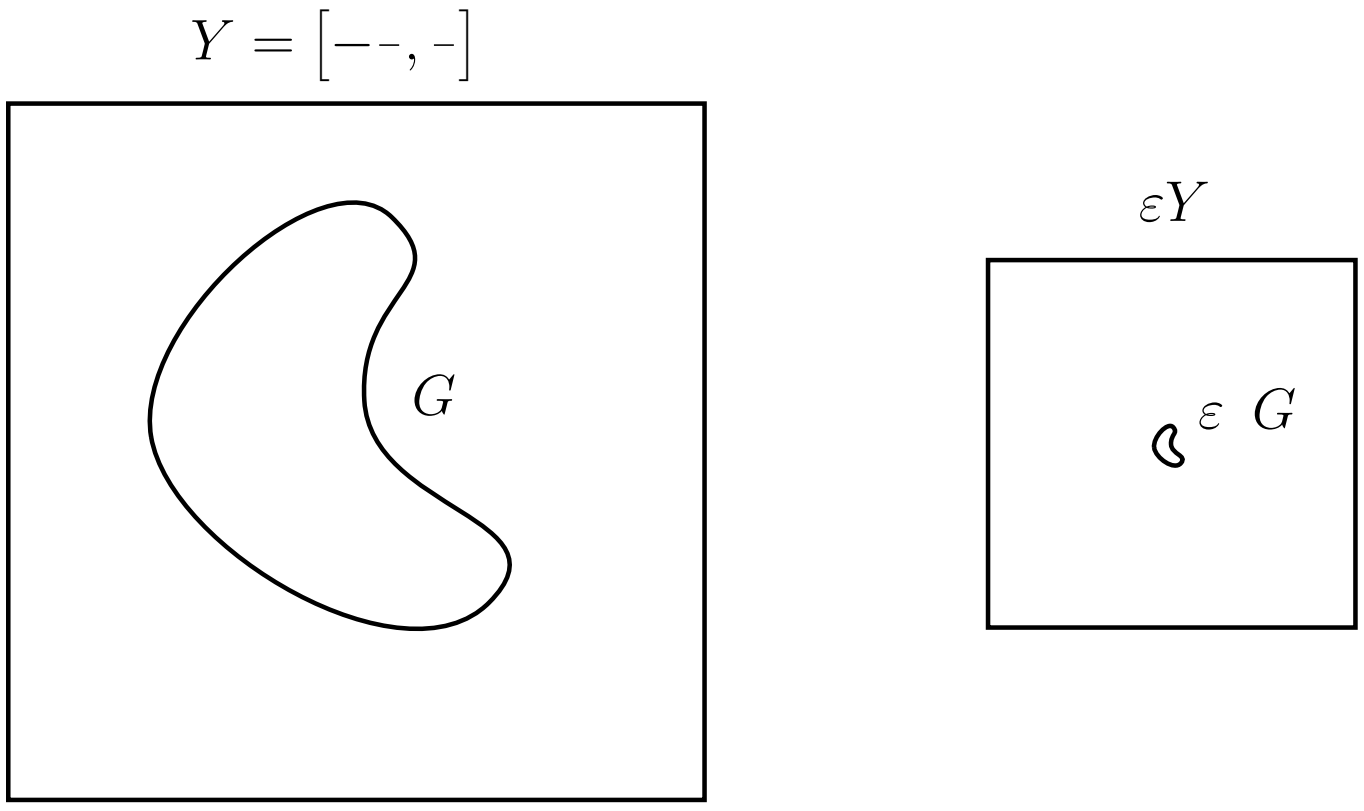

Fig. 1.2 The reference cell $Y$ and the scalings by $\varepsilon$ and $\varepsilon^\alpha$, for $\alpha > 1$. Notice that, for $\alpha > 1$, $\varepsilon^\alpha G_0$ (for a general particle shaped as $G_0$) becomes smaller relative to $\varepsilon Y$, which scales as the repetition. In most of our cases $G_0$ will be a ball $B_1(0)$.

### 1.1.2 Governing equation

We will consider that there is no diffusion inside $G^\varepsilon$, but it can be seen as a catalyst agent, producing a reaction on its boundary

$$\begin{cases} -\operatorname{div}(A^\varepsilon \nabla w^\varepsilon) = \hat{f}^\varepsilon & \Omega_\varepsilon, \\ A^\varepsilon \nabla w_\varepsilon \cdot n + \beta(\varepsilon)\hat{\sigma}(w_\varepsilon) = \beta(\varepsilon)\hat{g}^\varepsilon & S_\varepsilon, \\ w^\varepsilon = 1 & \partial\Omega. \end{cases} \tag{1.8}$$

Here $\hat{\sigma}$ is the reaction kinetics, typically a nondecreasing function. Notice that $\partial\Omega^\varepsilon = \partial\Omega \cup S^\varepsilon$. The parameter $\beta(\varepsilon)$ modulates the intensity of the reaction, and its value shall be precised shortly.

We will also deal, in this problem, with nonlinear diffusion

$$\begin{cases} -\Delta_p w_\varepsilon = \hat{f}^\varepsilon & \Omega_\varepsilon, \\ \dfrac{\partial w_\varepsilon}{\partial \nu_p} + \beta(\varepsilon)\hat{\sigma}(w^\varepsilon) = \beta(\varepsilon)\hat{g}^\varepsilon & S_\varepsilon, \\ w^\varepsilon = 1 & \partial\Omega, \end{cases} \tag{1.9}$$

where $p > 1$ and

$$\begin{aligned} -\Delta_p w &= \operatorname{div}(|\nabla w|^{p-2}\nabla w) \\ \frac{\partial w}{\partial \nu_p} &= |\nabla w|^{p-2}\nabla w \cdot n. \end{aligned}$$

The quasilinear diffusion operator $-\Delta_p$ represents the cases in which the diffusion coefficient depends of $|\nabla w|$ (see [Día85] and the references therein). Notice that for $p = 2$ we get the usual, linear, Laplacian operator. However, for $p > 2$ the operator becomes degenerate (the diffusion coefficient vanishes when $|\nabla w| = 0$) and for $1 < p < 2$ the operator becomes singular (the diffusion coefficient is unbounded as $|\nabla w| \to 0$).

**A change in variable** Boundary condition $w = 1$ is not nice in terms of functional spaces. We will thus choose the change in variable

$$u_\varepsilon = 1 - w_\varepsilon. \tag{1.10}$$

With the change of variable in mind we can rewrite (1.8) as

$$\begin{cases} -\operatorname{div}(A^\varepsilon \nabla u^\varepsilon) = f^\varepsilon & \Omega_\varepsilon, \\ A^\varepsilon \nabla u_\varepsilon \cdot n + \beta(\varepsilon)\sigma(u_\varepsilon) = \beta(\varepsilon) g^\varepsilon & S_\varepsilon, \\ u^\varepsilon = 0 & \partial\Omega, \end{cases} \tag{1.11}$$

and (1.9) as

$$\begin{cases} -\Delta_p u_\varepsilon = f^\varepsilon & \Omega_\varepsilon, \\ \dfrac{\partial u_\varepsilon}{\partial \nu_p} + \beta(\varepsilon)\sigma(u^\varepsilon) = \beta(\varepsilon) g^\varepsilon & S_\varepsilon, \\ u^\varepsilon = 0 & \partial\Omega, \end{cases} \tag{1.12}$$

where

$$\sigma(u) = \hat\sigma(1) - \hat\sigma(1-u), \qquad f^\varepsilon = \hat\sigma(1) - \hat f^\varepsilon, \qquad g^\varepsilon = \hat\sigma(1) - \hat g^\varepsilon. \tag{1.13}$$

### 1.1.3 Effectiveness and ineffectiveness

As motivated by the definition of Aris we define the *effectiveness* of the non-homogeneous problem as

$$\mathscr{E}_\varepsilon(\Omega, G_0) = \frac{1}{|S_\varepsilon|} \int_{S_\varepsilon} \hat\sigma(w_\varepsilon). \tag{1.14}$$

Since in this chapter we will deal with $u_\varepsilon$ rather than $w_\varepsilon$ let us define the *ineffectiveness functional*

$$\eta_\varepsilon(\Omega, G_0) = \frac{1}{|S_\varepsilon|} \int_{S_\varepsilon} \sigma(u_\varepsilon). \tag{1.15}$$

We have that

$$\eta_\varepsilon(\Omega, G_0) = \sigma(1) - \mathscr{E}_\varepsilon(\Omega, G_0). \tag{1.16}$$

Hence, in terms of convergence and optimization, analyzing one of the functionals is exactly the same as analyzing the other one.

### 1.1.4 Maximal monotone operators. *A common roof*

In some contexts, it is desirable to substitute the condition

$$\frac{\partial u_\varepsilon}{\partial \nu_p} + \beta(\varepsilon)\sigma(u^\varepsilon) = 0 \text{ on } S_\varepsilon \tag{1.17}$$

by the Dirichlet boundary condition

$$u_\varepsilon = 0 \text{ on } S_\varepsilon, \tag{1.18}$$

or even the case of Signorini type boundary condition (also known as boundary obstacle problem)

$$\begin{cases} u_\varepsilon \geq 0 & S_\varepsilon, \\ \partial_{\nu_p} u_\varepsilon + \beta(\varepsilon)\sigma_0(u_\varepsilon) \geq 0 & S_\varepsilon, \\ u_\varepsilon \left(\partial_{\nu_p} u_\varepsilon + \beta(\varepsilon)\sigma_0(u_\varepsilon)\right) = 0 & S_\varepsilon. \end{cases} \tag{1.19}$$

There is a unified presentation of this theory, using (1.9). The idea is to use *maximal monotone operators* (see [Bré73] or [Bro67])

**Definition 1.1.** Let $X$ be a Banach space and $\sigma : X \to \mathscr{P}(X')$. We say that $\sigma$ is a monotone operator if, for all $x, \hat{x} \in X$,

$$\langle x - \hat{x}, \xi - \hat{\xi} \rangle_{X \times X'} \geq 0 \qquad \forall \xi \in \sigma(x), \hat{\xi} \in \sigma(\hat{x}). \tag{1.20}$$

We define the domain of $\sigma$ as

$$D(\sigma) = \{x \in X : \sigma(x) \neq \emptyset\}. \tag{1.21}$$

Here $\emptyset$ is the empty set. We say that $\sigma$ is a maximal monotone operator if there is no other monotone operator $\widetilde{\sigma}$ such that $D(\sigma) \subset D(\widetilde{\sigma})$ and $\sigma(x) \subset \widetilde{\sigma}(x)$ for all $x \in X$.

It can be shown that any maximal monotone operator in $\mathbb{R}$ is given by a monotone functions, the jumps of which are filled by a vertical segment. It is immediate to prove the following:

**Proposition 1.1.** *Let $\sigma \in \mathscr{C}(\mathbb{R})$ be nondecreasing. Then $\sigma$ is a maximal monotone operator.*

Furthermore

**Proposition 1.2.** *Let $\sigma : \mathbb{R} \to \mathbb{R}$ be a nondecreasing function, and let $(x_n)_n$ be its set of discontinuities. Then, the function*

$$\widetilde{\sigma}(x) = \begin{cases} \sigma(x) & x \in \mathbb{R} \setminus \{x_n : n \in \mathbb{N}\}, \\ [\sigma(x_n^-), \sigma(x_n^+)] & x = x_n \text{ for some } n \in \mathbb{N}, \end{cases} \tag{1.22}$$

*is a maximal monotone operator.*

Boundary condition (1.18) can be written in terms of maximal monotone operators as (1.17) with

$$\sigma(x) = \begin{cases} \emptyset & x < 0, \\ \mathbb{R} & x = 0, \\ \emptyset & x > 0. \end{cases} \tag{1.23}$$

One the other hand, (1.19) can be written as (1.17) with

$$\sigma(x) = \begin{cases} \emptyset & x < 0, \\ (-\infty, 0] & x = 0, \\ \sigma_0(x) & x > 0. \end{cases} \tag{1.24}$$

Of course, the use of maximal monotone operators escapes the usual framework of classical solutions of PDEs. We will present the definition of weak solutions for this setting in Section 1.4.1.

Another advantage of maximal monotone operators is the simplicity to define inverses. For $\sigma : \mathbb{R} \to \mathscr{P}(\mathbb{R})$ we define its inverse in the sense of maximal monotone operators as the map $\sigma^{-1} : \mathbb{R} \to \mathscr{P}(\mathbb{R})$ given by

$$\sigma^{-1}(s) = \{x \in \mathbb{R} : s \in \sigma(x)\}. \tag{1.25}$$

It is a trivial exercise that $\sigma^{-1}$ is also a maximal monotone operator.

## 1.2 An introduction to homogenization

The main idea of this theory is to consider an inhomogeneous setting -be it due to some oscillating term in the equation or because of the domain itself- and decide which homogeneous equation can approximate the result in a “mean field approach” in order to “remove” these obstacles. Usually, studying the heterogeneous medium is not feasable, whereas the homogeneous equation can be easily undestood.

In order to fix notations, let us define some Sobolev spaces. For $\Omega$ a smooth set we define the space

$$W^{1,p}_{loc}(\Omega) = \{f : \Omega \to \mathbb{R} : \text{ for all } K \subset \Omega \text{ compact }, f \in W^{1,p}(K)\}. \tag{1.26}$$

For $\Gamma \subset \partial\Omega$ we define

$$W^{1,p}(\Omega,\Gamma) = \overline{\{f \in \mathscr{C}^\infty(\Omega) : f = 0 \text{ on } \Gamma\}}^{W^{1,p}}. \tag{1.27}$$

Some particular cases deserve their own notation:

$$\begin{aligned} W_0^{1,p}(\Omega) &= W^{1,p}(\Omega,\partial\Omega), \\ H^1(\Omega,\Gamma) &= W^{1,2}(\Omega,\Gamma), \\ H_0^1(\Omega) &= W_0^{1,2}(\Omega). \end{aligned}$$

Finally, for a cube $Q$, we define

$$W_{per}^{1,p}(Q) = \{f \in W^{1,p}(Q) : f \text{ can be extended by periodicity to } W_{loc}^{1,p}(\mathbb{R}^n)\}. \tag{1.28}$$

The theory of homogenezation is very broad and has been adapted to deal with many problem, from fluid mechanics to Lithium-ion batteries (see, e.g., [BC15]). Many books can be found which aim to give an introduction to this extensive and technical field (see, e.g., [BLP78; SP80; CD99; Tar10]).

In this Chapter we aim to give a comprehensive study of the problem in which the domain $\Omega$ contains some inclusions (or holes). This kind of problem, for which the literature is quite extensive, is known sometimes as the problem of "open domain with holes".

### 1.2.1 Some first results

To illustrate, in a very simple example, how some of the ideas work, let us go back to one the earliest results in homogenization. The idea behind the following example is a $G$-convergence argument (owed to Spagnolo [Spa68]).

**Example 1.1.** Let $a : \mathbb{R} \to \mathbb{R}$ be a $[0,1]$-periodic function such that $0 < \alpha \leq a \leq \beta$, $f \in L^2(0,1)$ and $a^\varepsilon(x) = a\left(\frac{x}{\varepsilon}\right)$ . We consider the one dimensional problem

$$\begin{cases} -\dfrac{d}{dx}\left(a^\varepsilon \dfrac{du_\varepsilon}{dx}\right) = f & x \in (0,1), \\ u_\varepsilon(0) = u_\varepsilon(1) = 0. \end{cases} \tag{1.29}$$

By multiplying by $u_\varepsilon$ and integrating, we have that the sequence $u_\varepsilon$ is bounded in $H_0^1(0,1)$, and therefore

$$u_\varepsilon \rightharpoonup u_0$$

in $H_0^1(0,1)$ and, by the same argument

$$a^\varepsilon \nabla u_\varepsilon = \xi_\varepsilon \rightharpoonup \xi_0$$

is convergent in $H^1(0,1)$ (since $f \in L^2$) and, in the limit

$$\begin{cases} -\dfrac{d}{dx}\left(\xi^0\right) = f & x \in (0,1), \\ u(0) = u(1) = 0. \end{cases} \tag{1.30}$$

holds. It is well-known that, for $h \in L^2(0,1)$, $h\left(\frac{\cdot}{\varepsilon}\right) \rightharpoonup \int_0^1 h$ in $L^2(0,1)$. Hence, up to a subsequence,

$$\nabla u_\varepsilon = \frac{1}{a^\varepsilon}\xi^\varepsilon \rightharpoonup \int_0^1 \frac{1}{a(x)} dx \cdot \xi^0 \tag{1.31}$$

in $L^2(0,1)$. Hence $\xi^0 = \frac{1}{\int_0^1 \frac{1}{a(x)} dx} \frac{du_0}{dx}$ and thus $u_0$ satisfies

$$\begin{cases} -\dfrac{d}{dx}\left(\dfrac{1}{\int_0^1 \frac{1}{a(x)} dx} \dfrac{du}{dx}\right) = f & x \in (0,1), \\ u(0) = u(1) = 0. \end{cases} \tag{1.32}$$

The term $a_0 = \frac{1}{\int_0^1 \frac{1}{a(x)} dx}$ is sometimes known as *effective diffussion coefficient*. This concludes this example.

One of the many works in homogenization in dimension higher than one is due to J.L Lions [Lio76], which contains a compendium of different references (e.g. [Bab76]). The focus of this work is the problem of oscillating coefficients

$$A^\varepsilon u_\varepsilon = f, \qquad A^\varepsilon v = \operatorname{div}\left(A\left(\frac{x}{\varepsilon}\right)\nabla v\right) \tag{1.33}$$

where $A = (a_{ij})$ is a matrix, $a_{ij} = a_{ji} \in L^\infty([0,1]^n)$ and are extended by periodicity. This models the behaviour of a periodical two phase composite (a material formed by the inclusion of two materials with different properties). This work is, no doubt, based on previous results, for example by Spagnolo (see, e.g., [Spa68]) on the limit behaviour of problems $-\operatorname{div}(A_k u_k)$ as $A_k \to A_\infty$.

The different approaches are very well presented in [BLP78] and [CD99].

### 1.2.2 Different techniques of homogenization

Here we will briefly present some of the most relevant methodologies applied in homogenization. Most of them have been applied to our problem, as we will see later.

**Multiple scales method** One of the possibilities in dealing with the limit consists on considering an expansion -which is known as asymptotical expansion- of the solutions as

$$u_\varepsilon(x) = u_0\left(x, \frac{x}{\varepsilon}\right) + \varepsilon u_1\left(x, \frac{x}{\varepsilon}\right) + \cdots, \tag{1.34}$$

and deriving the behaviour from there. This method, which is now known as *multiple scales method* is still very much in use to these days (see, e.g., [Día99; BC15]).

This kind of argument work in two steps. First, a formal deduction of the good approximation and a later rigorous proof. In particular, they use repeatedly the computation that, if $v = v(x, \xi)$ then

$$\frac{\partial}{\partial x_i}(v(x, \varepsilon^{-1}x)) = \frac{\partial v}{\partial x_i}(x, \varepsilon^{-1}x) + \varepsilon^{-1}\frac{\partial v}{\partial \xi_i}(x, \varepsilon^{-1}x). \tag{1.35}$$

Substituting (1.34) into $-\operatorname{div}(A^\varepsilon u^\varepsilon) = f$ and gathering terms the is seen that

$$u_\varepsilon(x) = u_0(x) + \varepsilon\hat{\xi}\left(\frac{x}{\varepsilon}\right)\cdot\nabla u_0 + \varepsilon^2\hat{\theta} : D^2 u_0 + \cdots, \tag{1.36}$$

and the equations for $u_0$, $\hat{\xi}$ and $\hat{\theta}$ can be found explicitly. The second part of this kind of argument is to estimate the convergence. It can be shown that

$$\left\| u_\varepsilon(x) - \left(u_0(x) + \varepsilon\hat{\xi}\left(\frac{x}{\varepsilon}\right)\cdot\nabla u_0 + \varepsilon^2\hat{\theta} : D^2 u_0\right)\right\|_{H^1(\Omega)} \le C\varepsilon^{\frac{1}{2}}. \tag{1.37}$$

Detailed examples can be found, e.g., in [CD99, Chapter 7], [BLP78] or, for the case of the elasticity equation, [OSY92].

**The Γ-convergence method** This method introduce by De Giorgi [DF75] and later developed in [DD83; Dal93]. The essential idea behind the Γ-convergence method is to study the problem in its energy form and study the conditions under which convergence of the energies implies convergence of their minimizers, the solutions of the elliptic problems. Here we present some results extracted from [Dal93].

**Definition 1.2.** Let $X$ be a topological space. The $\Gamma$-lower limit and $\Gamma$-upper limit of a sequence $(F_n)$ of functions $X \to [-\infty,\infty]$ are defined as follows

$$\left(\Gamma - \liminf_{n\to+\infty} F_n\right)(x) = \sup_{U\in\mathscr{N}(x)} \liminf_{n\to+\infty} \inf_{y\in U} F_n(y) \tag{1.38}$$

$$\left(\Gamma - \limsup_{n\to+\infty} F_n\right)(x) = \sup_{U\in\mathscr{N}(x)} \limsup_{n\to+\infty} \inf_{y\in U} F_n(y) \tag{1.39}$$

where $\mathscr{N}(x) = \{U \subset X \text{ open } : x \in U\}$ . If there exists $F : X \to [-\infty,+\infty]$ such that $F = \Gamma - \liminf_{n\to+\infty} F_n = \Gamma - \limsup_{n\to+\infty} F_n$ then $(F_n)$ we say that $F_n$ $\Gamma$-converges to $F$, and we denote it as

$$F = \Gamma - \lim_{n\to+\infty} F_n. \tag{1.40}$$

For the length of this section we will note

$$F' = \Gamma - \liminf_{n\to+\infty} F_n, \tag{1.41}$$

$$F'' = \Gamma - \limsup_{n\to+\infty} F_n. \tag{1.42}$$

The results that make this technique interesting for us are the following:

**Theorem 1.1.** *Suppose that $(F_n)$ are equi-coercive in $X$. Then $F'$ and $F''$ are coercive and*

$$\inf_{x\in X} F'(x) = \liminf_{n\to+\infty} \inf_{x\in X} F_n(x). \tag{1.43}$$

**Proposition 1.3.** *Let $x_n$ be a minimizer of $F_n$ in $X$ and assume that $x_n \to x$ in $X$. Then*

$$F'(x) = \liminf_{n\to\infty} F_n(x_n), \qquad F''(x) = \limsup_{n\to\infty} F_n(x_n). \tag{1.44}$$

In the context of homogenization we are mainly interested in the behavior of functionals

$$F_\varepsilon(u,A) = \begin{cases} \displaystyle\int_A f\left(\frac{x}{\varepsilon}, u(x), Du(x)\right) & u \in W^{1,p}(A), \\ +\infty & \text{otherwise,} \end{cases} \tag{1.45}$$

where $p > 1$.

Assume $f = f(y, Du)$ satisfies the following

i) For every $x \in \mathbb{R}^n$ the function $f(x,\cdot)$ is convex of class $\mathscr{C}^1$.

ii) For every $\xi \in \mathbb{R}^n$ is measurable and $Y$-periodic (where $Y$ is the unit cube).

iii) There exists $c_i \in \mathbb{R}$ $i = 0, \cdots, 4$ such that $c_1 \geq c_0 > 0$ and

$$c_0|\xi|^p \leq f(x,\xi) \leq c_1|\xi|^p + c_2, \tag{1.46}$$

$$\left|\frac{\partial f}{\partial \xi}(x,\xi)\right| \leq c_3|\xi|^{p-1} + c_4. \tag{1.47}$$

Then, let

$$f_0(\xi) = \inf_{v \in W^{1,p}_{per}(Y)} \int_Y f(y, \xi + Dv(y))dy. \tag{1.48}$$

Then for every sequence $\varepsilon_n \to 0$ we have that $F_{\varepsilon_n}$ $\Gamma$-converges to $F_0$ the functional defined by

$$F_0(u,A) = \begin{cases} \int_A f_0(Du) & u \in W^{1,p}(A), \\ +\infty & \text{otherwise.} \end{cases} \tag{1.49}$$

Function $f_0$ can be characterized in a more practical manner.

**Proposition 1.4.** *We have that*

$$f_0(\xi) = \int_Y f(y, Dv(y))dy, \tag{1.50}$$

*where u is the unique function that*

$$\begin{cases} v \in W^{1,p}_{loc}(\mathbb{R}^n), \\ Dv \textit{ is } Y-periodic, \\ \int_Y Dv = \xi, \\ \operatorname{div}(D_\xi f(y,Dv)) = 0 \textit{ in } \mathbb{R}^n \textit{ in the sense of distributions.} \end{cases} \tag{1.51}$$

Applying this method we can obtain the same result as in Example 1.1 in a way that can generalized to higher dimension

**Example 1.2.** Let $n = 1$ and $\Omega = Y = (0,1)$. Let us consider solutions of problem for the operator $-\operatorname{div}(a(\frac{x}{\varepsilon})u^\varepsilon)$. We consider the energy function

$$f(x,\xi) = a\left(\frac{x}{\varepsilon}\right)|\xi|^2. \tag{1.52}$$

Let us characterize $v$ given by Proposition 1.4. Since $(av')' = 0$ we have that that $v' = \frac{a_0}{a}$, where $k$ is a constant. Taking into account that $\int_0^1 v' = \xi$ we deduce that

$$a_0 = \frac{1}{\int_0^1 \frac{1}{a(x)}dx}\xi.$$

Therefore

$$\begin{aligned} f_0(\xi) &= \int_0^1 a(y)|v'(y)|^2 dy = \int_0^1 a(y)\frac{a_0^2}{a(y)^2}dy|\xi|^2 = a_0^2 \int \frac{1}{a(y)}dy|\xi|^2 \\ &= a_0|\xi|^2, \end{aligned}$$

where the effective diffusion is given by the coefficient

$$a_0 = \frac{1}{\int_0^1 \frac{1}{a(x)}dx}.$$

Hence, as we previously showed in Example 1.1, the limit of this kind $-\operatorname{div}(a(\frac{x}{\varepsilon})u^\varepsilon)$ is $-\operatorname{div}(a_0 u)$. Depending on the space $X$ we consider, we can fix one type of boundary condition or another. This completes this example.

In [Dal93] examples in higher dimensions are presented. This method was applied to our cases of interest, with some modification, by [Kai89] and [Gon97]. The details of this last paper (which is extremely synthetic and skips most of the computations) are given in Appendix 1.A.

**The two-scale convergence method** The two-scale method was introduced by Nguetseng [Ngu89] and later developed by Allaire [All92; All94]. The central definition of the theory is the following:

**Definition 1.3.** Let $(v_\varepsilon)$ be a sequence in $L^2(\Omega)$. We say that the sequence $v_\varepsilon$ *two scale converges* to a function $v_0 \in L^2(\Omega \times Y)$ if, for any function $\psi = \psi(x,y) \in \mathscr{D}(\Omega; \mathscr{C}^\infty_{\text{per}}(Y))$ one has

$$\lim_{\varepsilon \to 0} \int_\Omega v_\varepsilon(x)\psi\left(x, \frac{x}{\varepsilon}\right) = \frac{1}{|Y|}\int_\Omega \int_Y v_0(x,y)\psi(x,y). \tag{1.53}$$

By taking $\psi = \psi(x)$ in the previous definition it is immediate that

$$v_\varepsilon \rightharpoonup V^0 = \frac{1}{|Y|}\int_Y v_0(\cdot, y)dy \tag{1.54}$$

weakly in $L^2(\Omega)$. The key point of this theory is to study the convergence of functions of the type $\psi(x, \frac{x}{\varepsilon})$, and then apply them suitably to the weak formulation.

**Tartar's method of oscillating test functions** This method is due to L. Tartar (see [Tar77; Tar10; MT97]). The idea behind it is to consider the appropriate weak formulation and select suitable test functions $\varphi^\varepsilon$, with properties that, in the limit, reveal weak formulation of the homogeneous problem.

This is, basically, the general method applied to obtain the results of this thesis. As well shall see, its is not a straightforward recipy, and the choice of test function and their analysis can become a very hard task. Many detailed examples will be given in the following text. Perhaps the most illustrative, due to its simplicity, is Section 1.4.11.

One of the main difficulties that rise with this method in domains with particles or holes is the need of a common functional space, since $u_\varepsilon \in L^p(\Omega_\varepsilon)$. This leads to the construction of extension operators $P_\varepsilon : W^{1,p}(\Omega_\varepsilon) \to W^{1,p}(\Omega)$, that will be discussed in Section 1.4.5.

**The periodical unfolding method** The periodical unfolding method was introduced by Cioranescu, Damlamian and Griso in [CDG02; CDG08]. It consists on transforming the solution to a fixed domain $\Omega \times Y$. The case of particles (or holes) was considered in [CDGO08; CDDGZ12; CD16]. The latter paper acknowledges the contribution of the author of this thesis.

Let us present the reasoning in domains with particles (or holes). The idea is to decompose every point in $\Omega$ as a sum

$$x = [x]_Y + \{x\}_Y \tag{1.55}$$

where $[x]_Y$ is the unique element in $\mathbb{Z}^n$ such that $x - [x]_Y \in [0,1)^n$. That is, we have that $[\cdot]_Y$ is constant over $Y_\varepsilon^j$.

We define the operator

$$\mathscr{T}_{\varepsilon,\delta} : \varphi \in L^2(\Omega) \mapsto \mathscr{T}_{\varepsilon,\delta}(\varphi) \in L^p(\Omega \times \mathbb{R}^n)$$

as

$$\mathscr{T}_{\varepsilon,\delta}(\varphi)(x,z) = \begin{cases} \varphi\left(\varepsilon\left[\dfrac{x}{\varepsilon}\right]_Y + \varepsilon\delta z\right) & (x,z) \in \hat{\Omega}_\varepsilon \times \dfrac{1}{\delta}Y, \\ 0 & \text{otherwise,} \end{cases} \tag{1.56}$$

where

$$\hat{\Omega}_\varepsilon = \text{interior}\left( \bigcup_{\substack{\xi\in\mathbb{Z}^n:\\ \varepsilon(\xi+Y)\subset\Omega}} \varepsilon(\xi+\overline{Y}) \right). \tag{1.57}$$

Notice that $\mathscr{T}_{\varepsilon,\delta}(\varphi)(x,z)$ is piecewise constant in $x$. The boundary of $G_\varepsilon^j$ corresponds to $\hat{\Omega}_\varepsilon \times \partial G_0$.

The big advantage of this approach is that it removes the need to construct extension operators. Therefore, it allows to considers non-smooth shapes of $G_0$. This method as shown very good result, and the properties of $\mathscr{T}_{\varepsilon,\delta}(u_\varepsilon)$ are well understood, at least in the non-critical cases.

## 1.3 Literature review of our problem

Let us do a comprehensive review of the literature. Let us understand why. Recall the definition (1.3), (1.4) and (1.6). In a moderate abuse of notation, which will not lead to confusion let us use the notations:

- If $a \in \mathbb{R}$, $|a|$ will indicate its absolute value.
- If $A \subset \mathbb{R}^n$ is a set of dimension $m$ (i.e. $m = n$ if the domain is open, $m = n-1$ if $A = \partial U$ for $U$ open, etc.) then $|A|$ will indicate its $m$-dimensional Lebesgue measure.
- If $A$ is a finite set, then $|A|$ will indicate its cardinal (i.e. the number of the elements).

We also introduce two notations. Given $(a_\varepsilon)_{\varepsilon>0}, (b_\varepsilon)_{\varepsilon>0}$ sequence of real number we define

$$a_\varepsilon \sim b_\varepsilon \quad \equiv \quad \lim_{\varepsilon\to 0} \frac{a_\varepsilon}{b_\varepsilon} \in (0,+\infty) \tag{1.58}$$

$$a_\varepsilon \ll b_\varepsilon \quad \equiv \quad \lim_{\varepsilon\to 0} \frac{a_\varepsilon}{b_\varepsilon} = 0. \tag{1.59}$$

First, we estimate $|\Upsilon_\varepsilon|$. Since

$$\left| \Omega - \bigcup_{j\in\Upsilon_\varepsilon} (\varepsilon j + \varepsilon Y) \right| \to 0, \tag{1.60}$$

we have that

$$\frac{|\Omega|}{|\Upsilon_\varepsilon|\varepsilon^n} = \frac{|\Omega|}{|\Upsilon_\varepsilon||\varepsilon Y|} \to 1, \tag{1.61}$$

as $\varepsilon \to 0$. Hence

$$|\Upsilon_\varepsilon| \sim \varepsilon^{-n}. \tag{1.62}$$

We can therefore compute

$$\begin{aligned} |G_\varepsilon| &= \left| \bigcup_{j\in\Upsilon_\varepsilon} (\varepsilon j + a_\varepsilon G_0) \right| && (1.63) \\ &= |\Upsilon_\varepsilon||a_\varepsilon G_0| \\ &\sim \varepsilon^{-n} a_\varepsilon^n \\ &\sim (a_\varepsilon \varepsilon^{-1})^n. && (1.64) \end{aligned}$$

We see that the the case $a_\varepsilon \sim \varepsilon$ is different from the case $a_\varepsilon \ll \varepsilon$. In the first case, known, as case of *big particles*, $|G_\varepsilon|$ has a positive volume in the limit. We will show that this volumetric presence affects the kind of homogenized diffusion. The diffusion coefficient becomes a function of $G_0$. However, it has been shown that when $a_\varepsilon \ll \varepsilon$, case known as case of *small particles*, we have no volumetric contribution, and, as we will see, the diffusion is not affected.

Under some conditions, it has been reported that the case $a_\varepsilon \ll \varepsilon$ presents a critical scale $a_\varepsilon^*$ that separates different behaviour (the precise values and behaviours will be given later):

- In the subcritical case $a_\varepsilon^* \ll a_\varepsilon \ll \varepsilon$ the nature of the kinetic is preserved
- In the critical case $a_\varepsilon \sim a_\varepsilon^*$ the nature of the kinetic changes. This effect is known as the appearance of an *strange term.*
- in the supercritical case $a_\varepsilon \ll a_\varepsilon^*$ the problem behaves, in the limit, as the case $\sigma \equiv 0$. This case is not very relevant, as we will see in Section 1.4.11.

All of this problems have undergone extensive work, and Table 1.1 presents a detailed literature review, focusing on the case studied in terms of $a_\varepsilon$ and the regularity of $\sigma$. We hope this table puts the contributions of the author to this field in context.

### 1.3.1 Homogenization with big particles $a_\varepsilon \sim \varepsilon$

**Dirichlet boundary conditions on the particles.** The first work in this direction is [CS79]. It deals with problem (1.12) with $p = 2$ but fixing the value of $u_\varepsilon$ in $G_j^\varepsilon$ to a constant not necessarily zero rather than the Neumann boundary condition. This paper introduces the extension operator

$$P_\varepsilon : \{u \in H^1(\Omega_\varepsilon) : u = \text{ const. on } \partial G_i^\varepsilon,\ u = 0 \text{ on } \partial\Omega\} \to H_0^1(\Omega) \tag{1.65}$$

such that

$$\|P_\varepsilon v\|_{H^1(\Omega)} \leq C\|v\|_{H^1(\Omega_\varepsilon)}. \tag{1.66}$$

The advantage for the perforated domain problem is that $P_\varepsilon u_\varepsilon$ are all defined in the same space $H_0^1(\Omega)$, and therefore the convergence is easy to establish. In this paper the authors study several structures for matrix $A^\varepsilon$ in (1.11). The authors show that $P_\varepsilon u_\varepsilon \rightharpoonup u$, the solution of

$$\begin{cases} \operatorname{div}(a_0\nabla u) = f & \Omega, \\ u = 0 & \partial\Omega, \end{cases} \tag{1.67}$$

where $a_0$ is an effective diffusion matrix that depends of $G_0$. The nature of this result is similar to Example (1.1).

**Neumann boundary conditions on the particles** Later [CD88a] dealt with problem (1.12) with $p = 2$ (equivalently (1.12) with $A^\varepsilon = I$) in the case $\sigma(u) = au$ and $g^\varepsilon(x) = g\left(\frac{x}{\varepsilon} - j\right)$ in $\partial G_\varepsilon^j$. Their approach is an asymptotic expansion (this method is known as multiple scales method, see [SP80]).

If $\int_{\partial T} g(y)dS = 0$ they consider

$$u^\varepsilon(x) = u_0(x,y) + \varepsilon u_1(x,y) + \cdots \tag{1.68}$$

where $y = \frac{x}{\varepsilon}$. Otherwise they perform the expansion

$$u^\varepsilon(x) = \varepsilon^{-1}u_{-1}(x,y) + u_0(x,y) + \varepsilon u_1(x,y) + \cdots. \tag{1.69}$$

The result is that $P_\varepsilon u_\varepsilon \rightharpoonup u$, the solution of

$$\sum_{i,j} q_{ij}\frac{\partial^2 u}{\partial x_i \partial x_j} + \frac{|\partial T|}{|Y|}au = \begin{cases} f & \text{if } \int_{\partial T} g(y)dS = 0, \\ \int_{\partial T} g(y)dS & \text{otherwise,} \end{cases} \tag{1.70}$$

where $q_{ij}$ is given as

$$q_{ij} = \delta_{ij} + \frac{1}{|Y\setminus T|}\int_{Y\setminus T}\frac{\partial \chi_j}{\partial y_i}dy, \tag{1.71}$$

and $\chi_i$ are the solutions of the so-called *cell problems*:

$$\begin{cases} -\Delta\chi_i = 0 & \text{in } Y\setminus T, \\ \frac{\partial(\chi_i + y_i)}{\partial\nu} = 0 & \text{on } \partial T, \\ \chi_i & Y\text{-periodic}. \end{cases} \tag{1.72}$$

The surprising conclusion is that $u_\varepsilon \to u_0$ in the first case, but $\frac{1}{\varepsilon}u^\varepsilon \to u_{-1}$ in the second.

The nonlinear problem was later studied by Conca, J.I. Díaz, Timofte and Liñán in [CDT03; CDLT04]. Their technique, which involves oscillating test functions, requires the introduction of extension operators

$$P_\varepsilon : \{u \in H^1(\Omega_\varepsilon) : u = 0 \text{ on } \partial\Omega\} \to H^1_0(\Omega) \tag{1.73}$$

such that $(P_\varepsilon u)|_{\Omega_\varepsilon} = u$. This techniques allows for a fairly general class of nonlinearities $\sigma$, but not all can be considered. In particular, they consider the following cases $\sigma = \sigma(x,v)$

$$\left|\frac{\partial\sigma}{\partial v}(x,v)\right| \le C(1+|v|^q), \qquad 0 \le q < \frac{n}{n-2} \tag{1.74}$$

or

$$|\sigma(x,v)| \le C(1+|v|^q), \qquad 0 \le q < \frac{n}{n-2}. \tag{1.75}$$

The result is, naturally, that $P_\varepsilon u_\varepsilon \rightharpoonup u$ in $H^1_0(\Omega)$, where $u$ is the solution of

$$\begin{cases} -\operatorname{div}(a_0(G_0)\nabla u) + \frac{|\partial G_0|}{|Y\setminus G_0|}\sigma(u) = f & \text{in } \Omega, \\ u = 0 & \text{on } \partial\Omega, \end{cases} \tag{1.76}$$

and

$$a_0(G_0) = (q_{ij}) \tag{1.77}$$

is the effective diffusion matrix, where $q_{ij}$ are given by (1.71).

The same results as in [CDLT04] were obtained in [CDZ07] applying the unfolding method (developed for this case in [CDZ06]). The advantage of the unfolding method is that it reduces the regularity constraints on $\partial T$.

### 1.3.2 Homogenization with sub-critical small particles $a_\varepsilon^* \ll a_\varepsilon \ll \varepsilon$

The first result in this direction can be found in [CD88b]. The difficulty of this case is to understand the behaviour of the integrals $\beta(\varepsilon) \int_{S_\varepsilon} \cdot$ as $\varepsilon \to 0$. The theory in [CD88b] is rather comprehensive for the case $p = 2$. In their work the authors detect that, for $p = 2$, there exists values $a_\varepsilon^*$, depending on the spatial dimension, such that the behaviour changes. The term *critical case* already appears in this text. However, they are unable to specify what happens in this case.

The work in understanding the behavior of the boundary integral under the different situation of $p$ and regularity of $\sigma$ has been incremental over the last decades (see table 1.1). In Section 1.4.6 we present a unified approach that covers most of the results for the sub-critical cases. Through approximation techniques, developed in [DGCPS17d] and briefly presented in Section 1.4.9, the non smooth cases can be treated.

### 1.3.3 Homogenization with critical small particles

**Dirichlet boundary condition** First, Hruslov dealt with the Dirichlet homogeneous boundary condition on the holes [Hru72] (see the higher order case in [Hru77]) in a rather convoluted an functional way. In 1997 a measure theoretic analysis dealt the appearance of "strange term" in [CM97]. This later paper was much easier to understand.

**The Neumann boundary condition** The linear Neumann boundary condition was studied first in the homogeneous setting: [Hru79; Kai89; Kai90; Kai91]. The linear setting, $\sigma(u) = \lambda u$, was studied later in [OS96] (see also [OS95]).

In [Kai89; Kai91] (see also [Kai90]) an analysis of the nonlinear critical and subcritical cases is made by Kaizu, who is unable to properly characterize the change of the nonlinearity in the critical case. The first paper to properly study this case, and characterize the nature of the new ("strange") nonlinear function, is [Gon97], which applies the technique of $\Gamma$-convergence when $G_0$ is a ball, smooth $\sigma$ and $n = 3$. The very surprising result is the following

**Theorem 1.2.** *Let $G_0$ be a ball and assume*

*i)* $n = 3$, $a_\varepsilon = \varepsilon^\alpha$ $\alpha \le 3 = \alpha^*$, $\beta(\varepsilon) = \varepsilon^{-\gamma} \sim \beta^*(\varepsilon)$. *That is* $\gamma = 2\alpha - 3$

*ii)* $\frac{\partial \sigma}{\partial u} \ge C > 0$

*iii) There exists a nonnegative function $a \in \mathscr{C}_0^\infty(\Omega)$ such that*

$$\rho(x,u) = 2\int_0^u \sigma(x,s)ds + a(x) \geq 0. \tag{1.78}$$

*Then, $P_\varepsilon u_\varepsilon \rightharpoonup u$ in $H_0^1(\Omega)$ where $u$ is the solution of the problem*

$$\begin{cases} -\Delta u + 4\pi C(x,u) = f & \Omega, \\ u = 0 & \partial\Omega, \end{cases} \tag{1.79}$$

*and*

$$C(x,u) = \begin{cases} \sigma(x,u) & 2 < \alpha < 3, \\ H(x,u) & \alpha = \gamma = 3, \end{cases} \tag{1.80}$$

*being $H$ is the (unique) solution of functional equation*

$$H(s) = \sigma(s - H(s)). \tag{1.81}$$

However, there are a few steps that are not clearly justified in the paper. The computations of this paper have been detailed and explained in this thesis (as unpublished material) and correspond to Section 1.A.

The case of general $n$, non smooth $\sigma$ and $-\Delta_p$ were later studied in detail. This different problems introduce a different number of difficulties, and are a larger part of the work developed by the author of this thesis. For many years, the only case of $G_0$ that was understood was a ball:

- The usual Laplacian in $\mathbb{R}^n$ for $n = 3$: [Gon97]
- The usual Laplacian in $\mathbb{R}^n$ for $n \geq 2$: [ZS11] [ZS13]
- The $p$-Laplace operator and $2 < p < n$: [SP12] [Pod10] [Pod12] [Pod15].
- $n$-Laplacian for $n \geq 2$: The critical size of holes in the case $p < n$ is $a_\varepsilon = \varepsilon^{\frac{n}{n-p}}$. Naturally, this critical exponent $\alpha_p^*$ blows up as $p \to n$. As it turns out, a critical case also exists for the case $p = n$, and this was studied in [PS15].
- Roots and Heaviside type nonlinearity: [DGCPS16] .
- Signorini boundary conditions: [DGCPS17a]

- Maximal monotone operators and $1 < p < n$: [DGCPS17c]. This result covers all the previous cases under a common roof.

In 2017, the author of this thesis jointly with J.I. Díaz, T.A. Shaposhnikova and M.N. Zubova [DGCSZ17], considered (for the first time in the literature) the case of $G_0$ not a ball. The structure of the limit equation is unprecedented in the literature.

| | | | Dirichlet | $\sigma = 0$ | $\sigma = \lambda u$ | $0 < k_1 \leq \sigma \leq k_2$ <br> $\sigma(x,u) \in \mathscr{C}^1$ | $\lvert\sigma'(s)\rvert \leq C\lvert u\rvert^{p-1}$ <br> $(\sigma(x,u) - \sigma(x,v)) \geq C\lvert u-v\rvert^p$ | $\lvert\sigma'(u)\rvert \leq C(1+\lvert u\rvert^r)$ <br> or <br> $\lvert\sigma(u)\rvert \leq C(1+\lvert u\rvert^r)$ | Signorini | $\sigma$ m.m.g. |
|---|---|---|---|---|---|---|---|---|---|---|
| Big particles | $p = 2$ | $\alpha = 1$ | [CS79] | [CD88a] | | [CDLT04] <br> [CDZ07] (see also [CDZ06]) | | | | |
| Small non-critical particles | $p = 2$ | $1 < \alpha < \frac{n}{n-p}$ | | [CD88b] | [OS96] | [Gon97] ($N = 3$) <br> [ZS11][ZS13] | [SP12] | [Kai89],[DGCPS17d] | [JNRS14] | [Kai91] |
| | $2 < p < n$ | | | | | | | [DGCPS17d] | | |
| | $1 < p < 2$ | | | [Pod15] | | | | | | |
| | $p = n$ | $\frac{a_\varepsilon}{e^{-\varepsilon^{-\frac{n}{n-1}}}} \to 0$ | | [Podolskii and Shaposhnikova (to appear)] | | | | | | |
| | $p > n$ | $\alpha > 1$ | | [DGCPS17b] | | | | | | |
| Small critical particles | $p = 2$ | $\alpha = \frac{n}{n-p}$ <br> $G_0$ a ball | [Hru72] [CM97] | [OS96] | | [Gon97] ($N = 3$) <br> [ZS11][ZS13] | [SP12] | [Kai89] | [JNRS11] | [Kai91] , [DGCPS17c] |
| | $2 < p < n$ | | | | | | | $(\sigma = \lvert u\rvert^{q-1}u)$ | [GPPS15] | [DGCPS17c] |
| | $1 < p < 2$ | | [DGCPS17c] | | | | [DGCPS16] | [DGCPS17a] | | |
| | $p = 2$ | $\alpha = \frac{n}{n-p}$ <br> $G_0$ not a ball | | [DGCSZ17] | | | | | | |
| | $p = n$ | $\frac{a_\varepsilon}{e^{-\varepsilon^{-\frac{n}{n-1}}}} \to C \neq 0$ | | [PS15] | | | | | | |

Table 1.1 Schematic representation of bibliography for the homogenization problems (1.11), (1.12). Where $\alpha$ is present $a_\varepsilon = C_0\varepsilon^\alpha$. Gray background represents new results introduced by this thesis.

## 1.4 A unified theory of the case of small particles $a_\varepsilon \ll \varepsilon$

### 1.4.1 Weak formulation

When $\sigma : \mathbb{R} \to \mathbb{R}$ we usually define a weak solution of (1.12) as a function $u_\varepsilon \in W^{1,p}(\Omega_\varepsilon, \partial\Omega)$ such that

$$\int_{\Omega_\varepsilon} |\nabla u_\varepsilon|^{p-2} \nabla u_\varepsilon \cdot \nabla v + \beta(\varepsilon) \int_{S_\varepsilon} \sigma(u_\varepsilon) v = \int_{\Omega_\varepsilon} f^\varepsilon v + \beta(\varepsilon) \int_{S_\varepsilon} g^\varepsilon v \tag{1.82}$$

for all $v \in W^{1,p}(\Omega_\varepsilon, \partial\Omega)$. In Section 1.1.4 we introduced the concept of maximal monotone operator. When $\sigma$ is a maximal monotone operator this definition is no longer valid, since $\sigma(u_\varepsilon(x))$ may be multivalued. We change the previous equation by

**Definition 1.4.** We say that $u_\varepsilon \in W^{1,p}(\Omega_\varepsilon, \partial\Omega)$ is a weak solution of (1.12) if there exists $\xi \in L^p(S_\varepsilon)$ such that $\xi(x) \in \sigma(u_\varepsilon(x))$ for a.e. $x \in S_\varepsilon$ and

$$\int_{\Omega_\varepsilon} |\nabla u_\varepsilon|^{p-2} \nabla u_\varepsilon \cdot \nabla v + \beta(\varepsilon) \int_{S_\varepsilon} \xi v = \int_{\Omega_\varepsilon} f^\varepsilon v + \beta(\varepsilon) \int_{S_\varepsilon} g^\varepsilon v, \tag{1.83}$$

for all $v \in W^{1,p}(\Omega_\varepsilon, \partial\Omega)$.

Uniqueness of this kind of solution is a direct consequence of the monotonicity of $\sigma$ (see, e.g., [Día85]). However, it is not easy to show directly that there exist solutions of (1.12) in this sense. The energy formulation is much better for this task.

### 1.4.2 Energy formulation

Let us start by considering the usual case $\sigma : \mathbb{R} \to \mathbb{R}$. In this setting it is standard to define the energy functional over $W^{1,p}(\Omega_\varepsilon, \partial\Omega)$ as

$$J_\varepsilon(v) = \frac{1}{p} \int_{\Omega_\varepsilon} |\nabla v|^p + \beta(\varepsilon) \int_{S_\varepsilon} \Psi(v) - \int_{\Omega_\varepsilon} f^\varepsilon v - \beta(\varepsilon) \int_{S_\varepsilon} g^\varepsilon v, \tag{1.84}$$

where $\Psi(s) = \int_0^s \sigma(\tau)$. It is common to say that the energy formulation for (1.12) is

$$J_\varepsilon(u_\varepsilon) = \min_{v \in W^{1,p}(\Omega_\varepsilon, \partial\Omega)} J_\varepsilon(v). \tag{1.85}$$

For smooth $\sigma$ it can be shown that the unique solution $u_\varepsilon$ of (1.12) is the unique minimizer of this functional. One possible way to give a meaning to (1.12) is to use this formulation. For that we recall the concept of the subdifferential:

**Definition 1.5.** Let $X$ be a Banach space and $J : X \to (-\infty, +\infty]$ be a convex function. We usually define the domain of $J$ as

$$D(J) = \{x \in X : J(x) \neq +\infty\}. \tag{1.86}$$

We define the subdifferential of $J$ as the map $\partial J : X \to \mathscr{P}(X')$ given by

$$\partial J(x_0) = \{\xi \in X' : J(x) - J(x_0) \geq \langle \xi, x - x_0 \rangle \ \forall x \in X\}. \tag{1.87}$$

As is turns out, for every maximal monotone operator $\sigma$ defined in $\mathbb{R}$ there exists a convex function $\Psi$ such that $\sigma = \partial\Psi$. Furthermore, this $\Psi$ can be chosen uniquely under the extra condition $\Psi(0) = 0$.

This energy formulation connects directly with the concept of weak solution. The subdifferential $A_\varepsilon = \partial J_\varepsilon$ is given by the set of dual elements $\widehat{\xi}$ such that

$$\langle \widehat{\xi}, w \rangle = \int_{\Omega_\varepsilon} |\nabla v|^{p-2} \nabla v \cdot \nabla w + \beta(\varepsilon) \int_{S_\varepsilon} \xi w - \int_{\Omega_\varepsilon} f^\varepsilon w - \beta(\varepsilon) \int_{S_\varepsilon} g^\varepsilon w, \tag{1.88}$$

where $\xi(x) \in \sigma(v(x))$ for a.e. $x \in S_\varepsilon$ (see, e.g., [Lio69]).

The weak formulation of (1.12) is, precisely,

$$A_\varepsilon u_\varepsilon \ni 0. \tag{1.89}$$

### 1.4.3 Formulation as functional inequalities

Let us prove the equivalence between the weak and energy formulations:

**Lemma 1.4.1** (Chapter 1 in [ET99])**.** *Let $X$ be a reflexive Banach space, $J : X \to (-\infty, +\infty]$ be a convex functional $A = \partial J : X \to \mathscr{P}(X')$ be its subdifferential. Then the following are equivalent:*

*i) $u$ is a minimizer of $J$,*

*ii) $u \in D(A)$ and $0 \in Au$.*

*If either hold, then*

*iii) For every $v \in D(A)$ and $\xi \in Av$*

$$\langle \xi, v - u \rangle \geq 0. \tag{1.90}$$

*Furthermore, assume that $J$ is Gâteaux-differentiable on $X$ and $A$ is continuous on $X$ then* iii)) *is also equivalent to* i)).

**Remark 1.1.** Naturally, if there is uniqueness of iii) then the i)-iii) are also equivalent.

**Remark 1.2.** One should not confuse condition iii) with the -very similar- Stampacchia formulation (see e.g. [Día85]). For a bilinear form $a$ and a linear function $G$ the Stampacchia formulation is

$$a(u, v-u) \geq G(v-u) \tag{1.91}$$

for all $v$ in the correspondent space, whereas in formulation iii) we have $a(v, v-u)$.

From Lemma 1.4.1 we can extract some characterizing equations of the weak solution, which will be useful later.

**Proposition 1.5.** *Let $u_\varepsilon$ be a minimizer of $J_\varepsilon$. Then*

$$\int_{\Omega_\varepsilon} |\nabla u_\varepsilon|^{p-2} \nabla u_\varepsilon \cdot \nabla(v-u_\varepsilon) + \varepsilon^{-\gamma} \int_{S_\varepsilon} (\Psi(v) - \Psi(u_\varepsilon)) \geq \int_{\Omega_\varepsilon} f(v-u_\varepsilon), \tag{1.92}$$

*there exists $\xi \in \sigma(v(x))$ such that*

$$\int_{\Omega_\varepsilon} |\nabla v|^{p-2} \nabla v \cdot (v-u_\varepsilon) + \beta(\varepsilon) \int_{S_\varepsilon} \xi (v-u_\varepsilon) \geq \int_{\Omega_\varepsilon} f(v-u_\varepsilon), \tag{1.93}$$

*and*

$$\int_{\Omega_\varepsilon} |\nabla v|^{p-2} \nabla v \cdot (v-u_\varepsilon) + \beta(\varepsilon) \int_{S_\varepsilon} (\Psi(v) - \Psi(u_\varepsilon)) \geq \int_{\Omega_\varepsilon} f(v-u_\varepsilon), \tag{1.94}$$

*hold for all $v \in W^{1,p}(\Omega_\varepsilon, \partial\Omega)$.*

*Proof.* Let us assume that $u_\varepsilon$ is a minimizer of $J_\varepsilon$. Considering characterization iii) of Lemma 1.4.1 we have that

$$\int_{\Omega_\varepsilon} |\nabla v|^{p-2} \nabla v \cdot \nabla w + \beta(\varepsilon) \int_{S_\varepsilon} \xi w \geq \int_{\Omega_\varepsilon} f^\varepsilon w \tag{1.95}$$

for some $\xi$ such that $\xi(x) \in \sigma(u_\varepsilon(x))$. Since $\Psi$ is convex and $\sigma = \partial\Psi$ we have that

$$\Psi(v) - \Psi(u_\varepsilon) \geq \xi(v-u_\varepsilon). \tag{1.96}$$

Hence, (1.92) is proved.

Equation (1.93) can be obtained by considering the Brézis-Sibony characterization of the weak of (1.12) (see Lemme 1.1 of [BS71] or Theorem 2.2 of Chapter 2 in [Lio69]). Finally, let us prove (1.94). Consider the map $x \in \mathbb{R}^n \to |x|^p \in \mathbb{R}$. It is a convex map with derivative $D|x|^p = p|x|^{p-2}x$. Hence, for $a, b \in \mathbb{R}^n$ we have that

$$|a|^p - |b|^p \geq p|b|^{p-2} b \cdot (a-b). \tag{1.97}$$

Hence

$$|b|^p - |a|^p \leq p|b|^{p-2} b \cdot (b-a). \tag{1.98}$$

Considering $b = \nabla v$ and $a = \nabla u_\varepsilon$ we have that

$$|\nabla v|^p - |\nabla u_\varepsilon|^p \leq p|\nabla v|^{p-2} \nabla v \cdot (v - \nabla u_\varepsilon). \tag{1.99}$$

Taking into account this fact and that $u_\varepsilon$ is a minimizer of $J_\varepsilon$ we have that

$$0 \leq J(v) - J(u_\varepsilon) = \frac{1}{p}\int_{\Omega_\varepsilon} (|\nabla v|^p - |\nabla u_\varepsilon|^p) + \beta(\varepsilon)\int_{S_\varepsilon} (\Psi(v) - \Psi(u_\varepsilon)) - \int_{\Omega_\varepsilon} f(v - u_\varepsilon) \tag{1.100}$$

$$\leq \frac{1}{p}\int_{\Omega_\varepsilon} (|\nabla v|^p - |\nabla u_\varepsilon|^p) + \beta(\varepsilon)\int_{S_\varepsilon} (\Psi(v) - \Psi(u_\varepsilon)) - \int_{\Omega_\varepsilon} f(v - u_\varepsilon). \tag{1.101}$$

Thus, we have obtained (1.94). □

Under some conditions, one can show that these Variational Inequalities are, in fact, equivalent to the definition of weak and energy solutions. Since we will not need this, we give no further details here.

### 1.4.4 Existence and uniqueness of solutions

To prove the existence of solutions we can use Convex Analysis to prove the existence of minizers of $u_\varepsilon$, or consider a very strong theorem. To state in its broadest generality we introduce (following Brezis, see [Bre68]) the definition

**Definition 1.6.** Let $V$ be a reflexive Banach space. We say that $A : V \to V'$ is a *pseudo monotone operator* if it is bounded and it has following property: if $u_j \rightharpoonup u$ in $V$ and that

$$\limsup_{j \to +\infty} \langle T(u_j), u_j - u \rangle \leq 0,$$

then, for all $v \in X$,

$$\liminf_{j\to+\infty} \langle T(u_j), u_j - v\rangle \geq \langle T(u), u - v\rangle. \tag{1.102}$$

We can now state the theorem

**Theorem 1.3** ([Bre68], also Theorem 8.5 in [Lio69]). *Let $A : V \to V'$ be a pseudo-monotone operator and $\varphi$ a proper convex function lower semi-continuous such that*

$$\begin{cases} \text{there exist } v_0 \text{ such that } \varphi(v_0) < \infty \text{ and} \\ \dfrac{(Au, u - v_0) + \varphi(u)}{\|u\|} \to \infty, \text{ as } \|u\| \to \infty. \end{cases} \tag{1.103}$$

*Then, for $f \in V'$, there exists a unique solution of the problem*

$$(A(u) - f, v - u) + \varphi(v) - \varphi(u) \geq 0, \qquad \forall v \in V. \tag{1.104}$$

Uniqueness is a routine task. Let us give a sketch of proof, when $\sigma$ is a maximal monotone operator and $p \geq 2$. Assume that $u_\varepsilon^1, u_\varepsilon^2$ satisfy (1.83). Considering the difference between the two formulations

$$\int_{\Omega_\varepsilon} (|\nabla u_\varepsilon^1|^{p-2}\nabla u_\varepsilon^1 - |\nabla u_\varepsilon^2|^{p-2}\nabla u_\varepsilon^2)\cdot\nabla v + \beta(\varepsilon)\int_{S_\varepsilon}(\xi^1 - \xi^2)v = 0. \tag{1.105}$$

Taking $v = u_\varepsilon^1 - u_\varepsilon^2$, since $(\xi^1 - \xi^2)(u_\varepsilon^1 - u_\varepsilon^2)$ we have that

$$\int_{\Omega_\varepsilon} |\nabla(u_\varepsilon^1 - u_\varepsilon^2)|^p \leq 0. \tag{1.106}$$

There $u_\varepsilon^1 - u_\varepsilon^2$ is a constant. This constant is 0, due to the boundary condition. This concludes the proof.

We provide a complete proof of existence and uniqueness Considering the weak formulation the following result is immediate.

**Proposition 1.6** ([DGCPS17d]). *Let $p > 1$. Then, for every $\varepsilon > 0$ there exists a unique weak solution of (1.12) $u_\varepsilon \in W^{1,p}(\Omega_\varepsilon, \partial\Omega)$. Furthermore, there exists a constant $C$ independent of $\varepsilon$ such that*

$$\|\nabla u_\varepsilon\|_{L^p(\Omega_\varepsilon)}^{p-1} \leq C(\|f^\varepsilon\|_{L^{p'}(\Omega_\varepsilon)} + \beta(\varepsilon)\beta^*(\varepsilon)^{-1}\|g^\varepsilon\|_{L^\infty(S_\varepsilon)}). \tag{1.107}$$

Some extra information can be given about the pseudo-primitive $\Psi(u_\varepsilon)$.

**Proposition 1.7** ([DGCPS17c]). *There exists a unique $u_\varepsilon \in W^{1,p}(\Omega_\varepsilon, \partial\Omega)$ weak solution of* (1.92). *Besides, there exists $K > 0$ independent of $\varepsilon$ such that*

$$\|\nabla u_\varepsilon\|_{L^p(\Omega_\varepsilon)} + \varepsilon^{-\gamma} \|\Psi(u_\varepsilon)\|_{L^1(S_\varepsilon)} \leq K. \tag{1.108}$$

### 1.4.5 Extension operators

In order to introduce a definition of "convergence" we will need to construct an extension operator so that all solutions are extended to a common Sobolev space. If we do this correctly we will be able to take advantage of the compactness properties of this common space.

Let $A \subset B$. We say that $P$ is an extension operator if $P : F(A) = \{f : A \to \mathbb{R}\} \to F(B)$ and has the property that $P(f)|_A = f$. Let $p > 1$. We will say that a family of linear extension operator

$$P_\varepsilon : W^{1,p}(\Omega_\varepsilon) \to W^{1,p}(\Omega) \tag{1.109}$$

is uniformly bounded if there exists a constant $C > 0$, independent of $\varepsilon$, such that

$$\|P_\varepsilon u\|_{W_0^{1,p}(\Omega)} \leq C \|u\|_{W_0^{1,p}(\Omega_\varepsilon)} \quad \forall u \in W^{1,p}(\Omega_\varepsilon). \tag{1.110}$$

A family of operators with this property, for $1 \leq p < +\infty$, was constructed in [Pod15]. The idea is to apply the following theorem

**Theorem 1.4** (Theorem 7.25 in [GT01]). *Let $\Omega$ be a $C^{k-1,1}$ domain in $\mathbb{R}^n$, $k \geq 1$. Then (i) $\mathscr{C}^\infty(\bar{\Omega})$ is dense in $W^{k,p}(\Omega)$, $1 \leq p < +\infty$ and (ii) for any open set $\Omega' \supset\supset \Omega$ there exists a linear extension operator $E : W^{k,p}(\Omega) \to W_0^{k,p}(\Omega')$ such that $Eu = u$ in $\Omega$ and*

$$\|Eu\|_{W^{k,p}(\Omega')} \leq C \|u\|_{W^{k,p}(\Omega)} \tag{1.111}$$

*where $C = C(k, \Omega, \Omega')$.*

We consider a large ball $B$ such that $Y \Subset B$ and the linear extension operator

$$E : W^{1,p}(Y \setminus G_0) \to W^{1,p}(B) \tag{1.112}$$

such that

$$\|Eu\|_{W^{1,p}(B)} \leq C_0 \|u\|_{W^{1,p}(Y \setminus G_0)}. \tag{1.113}$$

In particular,

$$\|\nabla Eu\|_{L^p(B)} \leq C_1 \|\nabla u\|_{L^p(Y \setminus G_0)}. \tag{1.114}$$

Let us scale it down by $a_\varepsilon$:

$$E_{\varepsilon,j} : W^{1,p}((\varepsilon j + a_\varepsilon Y) \setminus G_\varepsilon^j) \to W^{1,p}(Y, G_0) \xrightarrow{E} W^{1,p}(B) \to W^{1,p}(\varepsilon j + a_\varepsilon B). \tag{1.115}$$

Notice that, rather than $Y_\varepsilon^j \setminus G_\varepsilon^j$ we are considering the $a_\varepsilon$-rescale of $Y$. By a simple change in variable we observe that

$$\|\nabla E_{\varepsilon,j} u\|_{L^p(\varepsilon j + B)} \leq C_1 \|u\|_{L^p((\varepsilon j + a_\varepsilon Y)\setminus G_\varepsilon^j)}. \tag{1.116}$$

Let $u \in W_0^{1,p}(\Omega_\varepsilon)$. Let us consider extend by 0 outside $\Omega$, i.e.

$$\widetilde{u}(x) = \begin{cases} u(x) & x \in \Omega_\varepsilon, \\ 0 & x \in \mathbb{R}^n \setminus \Omega. \end{cases} \tag{1.117}$$

We then define

$$P_\varepsilon u(x) = \begin{cases} E_{\varepsilon,j} \widetilde{u}(x) & x \in \varepsilon j + a_\varepsilon Y, j \in \Upsilon_\varepsilon, \\ u(x) & \text{otherwise}. \end{cases} \tag{1.118}$$

This is well defined, since the sets $\varepsilon j + a_\varepsilon Y$ does not overlap for $\varepsilon$ small. It is clear that $P_\varepsilon$ is linear, $P_\varepsilon u = u$ in $\Omega_\varepsilon$ and, by considering the sum over the space decomposition, we have the uniform bound (1.110), so $P_\varepsilon u \in W^{1,p}(\Omega)$. Since the boundary behaviour has not been modified, $P_\varepsilon u_\varepsilon \in W_0^{1,p}(\Omega)$. We can conclude

**Lemma 1.4.2.** *Let $G_0 \in \mathscr{C}^{0,1}$ such that $G_0 \Subset Y$. Then, there exists a uniformly bounded family of linear extension operators* (1.109).

#### 1.4.5.1 Extension operators and Poincaré constants

We will use the existence of a Poincaré constant for $W_0^{1,p}(\Omega)$, $C_{p,\Omega}$, such that

$$\|v\|_{L^p(\Omega)} \leq C_{p,\Omega} \|\nabla v\|_{L^p(\Omega)}, \qquad \forall v \in W_0^{1,p}(\Omega). \tag{1.119}$$

This constant $C_{p,\Omega}$ is known to exist for every domain $\Omega$ bounded. However, it is not trivial to show that all domain $\Omega_\varepsilon$ have a common constant. The following result is very often used in the literature but it is seldom stated. In [DGCPS17b] we took the time to prove it.

**Theorem 1.5** ([DGCPS17b])**.** *Let $p > 1$. If there exists a sequence of uniformly bounded extension operators in $W_0^{1,p}$ then there exists a uniform Poincaré constant for $W^{1,p}(\Omega_\varepsilon, \partial\Omega)$,*

*in the sense that*

$$\|u\|_{L^p(\Omega_\varepsilon)} \le C\|\nabla u\|_{L^p(\Omega_\varepsilon)} \qquad \forall u \in W_0^{1,p}(\Omega_\varepsilon) \text{ and } \varepsilon > 0, \tag{1.120}$$

*where C does not depend of $\varepsilon$. In particular, let*

$$\|\nabla P_\varepsilon u\|_{L^p(\Omega)} \le K_p\|\nabla u\|_{L^p(\Omega_\varepsilon)} \quad \forall u \in W_0^{1,p}(\Omega_\varepsilon), \tag{1.121}$$

*hold and $C_{p,\Omega}$ be a Poincaré constant for $W_0^{1,p}(\Omega)$. Then, $C = K_p C_{p,\Omega}$.*

#### 1.4.5.2 Convergence of the extension

Hence, the solution $u_\varepsilon$ can be extended, and $P_\varepsilon u_\varepsilon$ is a bounded sequence in $W^{1,p}(\Omega)$. Thus, it has a weak limit. The focus of the theory of homogenization is to characterize the equation satisfied by the limit function.

### 1.4.6 On treating the boundary measure and the appearance of a critical case

Treating the sequence of integrals $\int_{S_\varepsilon}$ is a delicate business. Before we begin their study rigorously, we will start by providing some intuitive (informal) computations.

#### 1.4.6.1 An informal approach

Let us focus first on (1.11). Its weak formulation reads

$$\int_{\Omega_\varepsilon} A^\varepsilon \nabla u_\varepsilon \nabla \varphi + \beta(\varepsilon)\int_{S_\varepsilon} \sigma(u_\varepsilon)\varphi = \int_{\Omega_\varepsilon} f^\varepsilon \varphi + \beta(\varepsilon)\int_{S_\varepsilon} g^\varepsilon \varphi \tag{1.122}$$

Let us see how the coefficient $\beta(\varepsilon)$ is decisive for the limit behaviour. First, we should keep in mind that (1.61). For a continuous function $\varphi \in \mathscr{C}^1(\Omega)$ we have, since $|\partial G_i^\varepsilon| = a_\varepsilon^{n-1}|\partial G_0|$ that

$$\beta(\varepsilon)\int_{S_\varepsilon} \varphi = \beta(\varepsilon)\sum_{i\in\Upsilon_\varepsilon}\int_{\partial G_\varepsilon^i} \varphi = \beta(\varepsilon)\sum_{i\in\Upsilon_\varepsilon}\left(\varphi(x_\varepsilon^i)|\partial G_\varepsilon^i| + \int_{\partial G_\varepsilon^i} \varphi'(\xi_\varepsilon^i(x))(x - x_\varepsilon^i)\right) \tag{1.123}$$

$$= \sum_{i\in\Upsilon_\varepsilon}\left(\beta(\varepsilon)\varphi(x_\varepsilon^i)a_\varepsilon^{n-1}|\partial G_0| + \alpha_\varepsilon^i\right). \tag{1.124}$$

If

$$\beta(\varepsilon)a_\varepsilon^{n-1} \sim \varepsilon^n, \tag{1.125}$$

we have (almost) a Riemann sum in the first term, except for the term $Y_\varepsilon^i \cap \partial\Omega \neq \emptyset$ (but this part has no contribution as $\varepsilon \to 0$). We check immediately that

$$\begin{aligned} |\alpha_\varepsilon^i| &\leq \|\varphi'\|_\infty 2 \operatorname{diam}(G_\varepsilon^i)|\partial G_i^\varepsilon|\beta(\varepsilon) \leq \|\varphi'\|_\infty 2 \operatorname{diam}(G_\varepsilon^i)|\partial G_0|\beta(\varepsilon)a_\varepsilon^{n-1} \\ &= 2\|\varphi'\|_\infty \operatorname{diam}(G_0)\beta(\varepsilon)a_\varepsilon^n \\ &= 2\|\varphi'\|_\infty \operatorname{diam}(G_0)a_\varepsilon\varepsilon^n. \end{aligned}$$

Hence, we can expect that

$$\beta(\varepsilon)\int_{S_\varepsilon}\varphi \to \begin{cases} C\int_\Omega \varphi & \beta(\varepsilon) \sim \beta^*(\varepsilon), \\ 0 & \beta(\varepsilon) \ll \beta^*(\varepsilon), \\ +\infty & \beta(\varepsilon) \gg \beta^*(\varepsilon), \end{cases} \tag{1.126}$$

where, recalling (1.125),

$$\beta^*(\varepsilon) = a_\varepsilon^{1-n}\varepsilon^n \tag{1.127}$$

as $\varepsilon \to 0$. Notice that

$$|S_\varepsilon| = \left|\bigcup_{j\in\Upsilon_\varepsilon}(\varepsilon j + \partial(a_\varepsilon G_0))\right| = |\Upsilon_\varepsilon||\partial(a_\varepsilon G_0)| \sim \varepsilon^{-n}a_\varepsilon^{n-1} \sim \frac{1}{\beta^*(\varepsilon)}. \tag{1.128}$$

Therefore, up to constants, this is an average

$$\beta^*(\varepsilon)\int_{S_\varepsilon} \text{ should behave like } \frac{1}{|S_\varepsilon|}\int_{S_\varepsilon}. \tag{1.129}$$

If there is any good behaviour, the only expectable result is that

$$\frac{1}{|S_\varepsilon|}\int_{S_\varepsilon} \to \frac{1}{|\Omega|}\int_\Omega. \tag{1.130}$$

This is true, at least, for constant functions.

**Remark 1.3.** In particular, if we consider the case $a_\varepsilon = C_0\varepsilon^\alpha$, $\beta(\varepsilon) = \varepsilon^{-\gamma}$ and $\beta^*(\varepsilon) = \varepsilon^{-\gamma^*}$ we can expect that

$$\gamma^* = \alpha(n-1) - n. \tag{1.131}$$

Thus, if $\beta(\varepsilon)$ is too small we cannot expect any reaction term in the limit equation, hence (in some sense) it becomes uninteresting. If $\beta(\varepsilon)$ is too large then there the reaction dominates the diffusion, we write

$$\beta(\varepsilon)^{-1}\beta^*(\varepsilon)\int_{\Omega_\varepsilon}\nabla u_\varepsilon\nabla\varphi+\beta^*(\varepsilon)\int_{S_\varepsilon}\sigma(u_\varepsilon)\varphi=\beta(\varepsilon)^{-1}\beta^*(\varepsilon)\int_{\Omega_\varepsilon}f\varphi+\beta^*(\varepsilon)\int_{S^\varepsilon}g^\varepsilon\varphi \tag{1.132}$$

and we see the diffusion term disappear in the effective equation.

This good intuitions are not always true. In Theorem 1.6 we will see some assumptions under which this intuitions hold.

#### 1.4.6.2 A trace theorem for $a_\varepsilon G_0$ in $\varepsilon Y$

The most difficult part of this analysis is the study of the boundary measure $\int_{S_\varepsilon}$, as well as the unexpected properties of the diffusion in the critical case. The following estimate will be fundamental to our study. The proof can be found for $p=2$ and $n\geq 2$ in [CD88b] and a different proof [OS96] in the case of balls. Here we extend the proof in [CD88b] to the case of $p>1$ and $n\geq 2$. Some of the following results were for $1<p<n$ were presented in [Pod15].

In the following pages we present a unified analysis of the different cases, similar to that of [CD88b], but including the cases $p\neq 2$.

**Lemma 1.4.3.** *Let $u\in W^{1,p}(Y_\varepsilon)$, $p>1$. Then*

$$\int_{a_\varepsilon G_0}|u|^p\leq Ca_\varepsilon^{n-1}\left(\varepsilon^{-n}\int_{Y_\varepsilon}|u|^p+\tau_\varepsilon\int_{Y_\varepsilon}|\nabla u|^p\right) \tag{1.133}$$

*where*

$$\tau_\varepsilon\sim\begin{cases}a_\varepsilon^{p-n} & p<n,\\ \ln\left(\dfrac{\varepsilon}{a_\varepsilon}\right)^{p-1} & p=n,\\ \varepsilon^{p-n} & p>n,\end{cases} \tag{1.134}$$

*and $C$ is a constant independent of $\varepsilon$ and $u$.*

*Proof.* Let

$$B_\varepsilon=B(0,\varepsilon)\setminus\overline{(a_\varepsilon G_0)} \tag{1.135}$$

and let $\varphi \in \mathscr{C}^\infty(\overline{B_\varepsilon})$. Since $G_0$ is star shaped then we can represent it in polar coordinates as

$$\partial G_0 = \{(\rho,\theta) : \rho = \Phi(\theta), \theta \in \Theta\}, \tag{1.136}$$

where $\Theta = [0,2\pi] \times [-\frac{\pi}{2},\frac{\pi}{2}]^{n-2}$.

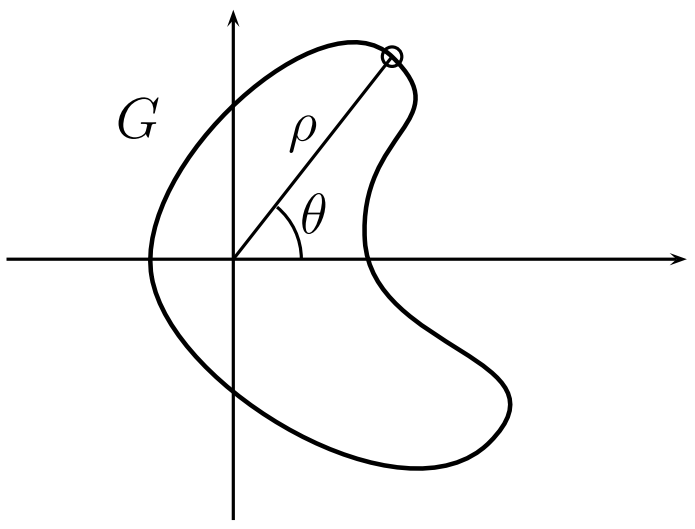

Fig. 1.3 The domain $G_0$ and its representation in polar coordinates.

The Jacobian can be written as $\rho^{n-1}J(\theta)$. Let us write $u$ in polar coordinates as $\chi(\rho,\theta) = u(x)$. Then, as in [CD88b],

$$\int_{a_\varepsilon G_0} |u|^p dx = a_\varepsilon^{n-1} \int_\Theta |\chi(a_\varepsilon \Phi(\theta),\theta)|^p J(\theta)F(\theta)d\theta, \tag{1.137}$$

where

$$F(\theta) = \prod_{i=1}^n \sqrt{\Phi(\theta)^2 + \left(\frac{\partial \Phi}{\partial \theta_i}\right)^2}. \tag{1.138}$$

We write, for any $\rho > a_\varepsilon \Phi(\theta)$ and $\theta \in \Theta$

$$\chi(a_\varepsilon \Phi(\theta),\theta) = \chi(\rho,\theta) - \int_{a_\varepsilon \Phi(\theta)}^\rho \frac{\partial \chi}{\partial t}(t,\theta)dt.$$

For $p > 1$, due to convexity

$$|\chi(a_\varepsilon \Phi(\theta),\theta)|^p \le 2^{p-1}|\chi(\rho,\theta)|^p + 2^{p-1}\left|\int_{a_\varepsilon \Phi(\theta)}^\rho \frac{\partial \chi}{\partial t}(t,\theta)dt\right|^p.$$

On the other hand

$$\left| \int\limits_{a_\varepsilon \Phi(\theta)}^{\rho} \frac{\partial \chi}{\partial t}(t,\theta)dt \right|^p \leq \left| \int\limits_{a_\varepsilon \Phi(\theta)}^{\rho} \frac{\partial \chi}{\partial t}(t,\theta) t^{\frac{n-1}{p}} t^{-\frac{n-1}{p}} dt \right|^p$$

$$\leq \left( \int\limits_{a_\varepsilon \Phi(\theta)}^{\rho} t^{-\frac{n-1}{p-1}} dt \right)^{p-1} \left( \int\limits_{a_\varepsilon \Phi(\theta)}^{\rho} \left| \frac{\partial \chi}{\partial t}(t,\theta) \right|^p t^{n-1} dt \right).$$

Taking

$$b_1 = \min_{\theta \in \Theta} \Phi(\theta) \qquad b_2 = \max_{\theta \in \Theta} \Phi(\theta)$$

we get

$$|\chi(a_\varepsilon \Phi(\theta),\theta)|^p \leq 2^{p-1} |\chi(\rho,\theta)|^p + 2^{p-1} \tau_\varepsilon \left( \int\limits_{a_\varepsilon \Phi(\theta)}^{\rho} \left| \frac{\partial \chi}{\partial t}(t,\theta) \right|^p t^{n-1} dt \right),$$

where

$$\tau_\varepsilon = \left( \int\limits_{a_\varepsilon b_1}^{\rho} t^{-\frac{n-1}{p-1}} dt \right)^{p-1}.$$

Integrating over $B_\varepsilon$ we obtain

$$\int_{\Theta}\int_{a_\varepsilon\Phi}^{\rho} |\chi(a_\varepsilon\Phi(\theta),\theta)|^p \rho^{n-1} J(\theta)F(\theta)d\rho d\theta$$

$$\leq 2^{p-1}\int_{\Theta}\int_{a_\varepsilon\Phi}^{\rho} |\chi(\rho,\theta)|^p \rho^{n-1} JFd\rho d\theta$$

$$+2^{p-1}\int_{\Theta}\int_{a_\varepsilon\Phi}^{\varepsilon} \tau_\varepsilon \left(\int_{a_\varepsilon\Phi(\theta)}^{\rho} \left|\frac{\partial\chi}{\partial t}(t,\theta)\right|^p t^{n-1}dt\right)\rho^{n-1}JFd\rho d\theta$$

$$\leq 2^{p-1}\int_{\Theta}\int_{a_\varepsilon\Phi}^{\rho} |\chi(\rho,\theta)|^p \rho^{n-1} JFd\rho d\theta$$

$$+2^{p-1}\tau_\varepsilon\tau_{2,\varepsilon}\int_{\Theta} \tau_\varepsilon \left(\int_{a_\varepsilon\Phi(\theta)}^{\rho} \left|\frac{\partial\chi}{\partial t}(t,\theta)\right|^p t^{n-1}dt\right)JFd\rho d\theta,$$

where

$$\tau_{2,\varepsilon} = \int_{b_1 a_\varepsilon}^{\varepsilon} \rho^{n-1}d\rho. \tag{1.139}$$

We can estimate the integral we wanted by

$$\int_{\Theta}\int_{a_\varepsilon\Phi}^{\rho} |\chi(a_\varepsilon\Phi(\theta),\theta)|^p \rho^{n-1} J(\theta)F(\theta)d\rho d\theta$$

$$\geq \tau_{3,\varepsilon}\int_{\Theta}\int_{a_\varepsilon\Phi}^{\rho} |\chi(a_\varepsilon\Phi(\theta),\theta)|^p J(\theta)F(\theta)d\rho d\theta$$

$$= \tau_{3,\varepsilon}\|\varphi\|^p_{L^p(\partial(a_\varepsilon G_0))},$$

where

$$\tau_{3,\varepsilon} = \int_{b_2 a_\varepsilon}^{\varepsilon} \rho^{n-1}d\rho. \tag{1.140}$$

Since $\tau_{3,\varepsilon}^{-1}\tau_{2,\varepsilon}$ is bounded we can conclude the estimate on $\|u\|_{L^p(\partial(a_\varepsilon G_0))}$. On the other hand

$$\tau_\varepsilon^{\frac{1}{p-1}} = \begin{cases} \dfrac{1}{1-\frac{n-1}{p-1}}\left(\varepsilon^{1-\frac{n-1}{p-1}} - (b_1 a_\varepsilon)^{1-\frac{n-1}{p-1}}\right) & p \neq n, \\ \ln\left(\dfrac{\varepsilon}{b_1 a_\varepsilon}\right) & p = n, \end{cases} \sim \begin{cases} a_\varepsilon^{1-\frac{n-1}{p-1}} & p < n, \\ \ln\left(\dfrac{\varepsilon}{a_\varepsilon}\right) & p = n, \\ \varepsilon^{1-\frac{n-1}{p-1}} & p > n, \end{cases} \tag{1.141}$$

which concludes the proof. □

**Remark 1.4.** It is not surprising that the $W^{1,n}(\Omega)$ for $\Omega \subset \mathbb{R}^n$ behaves differently. For example, radial solution solution of $\Delta_n u = 0$ in $\mathbb{R}^n$ includes $\ln|x|$, whereas for any other $p$-Laplacian radial solutions are of power type.

From this point forward we will assume that

$$G_0 \text{ is star-shaped.} \tag{1.142}$$

### 1.4.7 Behaviour of $\int_{S_\varepsilon}$ and appearance of $a_\varepsilon^*$

Define function $M_\varepsilon(x)$ as the unique $Y_\varepsilon$ - periodic built through the solution of the boundary value problem

$$\begin{cases} \Delta_p m_\varepsilon = \mu_\varepsilon & x \in Y_\varepsilon = \varepsilon Y \setminus \overline{a_\varepsilon G_0}; \\ \partial_{\nu_p} m_\varepsilon = 1 & x \in \partial(a_\varepsilon G_0) = S_\varepsilon^0; \\ \partial_{\nu_p} m_\varepsilon = 0 & x \in \partial Y_\varepsilon \setminus S_\varepsilon^0; \end{cases}, \qquad \int_{Y_\varepsilon} m_\varepsilon(x)dx = 0. \tag{1.143}$$

where $\mu_\varepsilon$ is a constant defined so as to satisfy the integrability condition

$$\mu_\varepsilon = \frac{\varepsilon^{-n} a_\varepsilon^{n-1} |\partial G_0|}{1 - (a_\varepsilon \varepsilon^{-1})^n |G_0|}. \tag{1.144}$$

That is

$$M_\varepsilon(x) = m_\varepsilon(x - P_\varepsilon^j), \qquad x \in Y_\varepsilon^j. \tag{1.145}$$

The aim of this section is to prove the following result

**Theorem 1.6.** *Assume that $a_\varepsilon \ll \varepsilon$,*

$$\beta(\varepsilon)\|M_\varepsilon\|_{L^p(\Omega_\varepsilon)}^{p-1} \to 0 \tag{1.146}$$

*as* $\varepsilon \to 0$ *and let*

$$\beta_0 = \lim_{\varepsilon \to 0} \mu_\varepsilon \beta(\varepsilon). \tag{1.147}$$

*Then, for all sequence* $v_\varepsilon \in L^p(\Omega)$ *such that* $v_\varepsilon \to v$ *in* $L^p(\Omega)$ *we have that*

$$\beta(\varepsilon) \int_{S_\varepsilon} v_\varepsilon dS \to \beta_0 \int_\Omega v dx \tag{1.148}$$

*as* $\varepsilon \to 0$.

Function $m_\varepsilon$ has the nice property of allowing us to write, for any test function $\varphi \in W^{1,p}(Y_\varepsilon)$,

$$-\int_{Y_\varepsilon} |\nabla m_\varepsilon|^{p-2} \nabla m_\varepsilon \nabla \varphi dx + \int_{S_\varepsilon^0} \varphi ds = \mu_\varepsilon \int_{Y_\varepsilon} \varphi dx. \tag{1.149}$$

We will use the following fact

**Lemma 1.4.4.** *Let* $p > 1$. *Then*

$$\|\nabla m_\varepsilon\|_{L^p(Y_\varepsilon)}^{p-1} \le C a_\varepsilon^{n-1} \left(\varepsilon^{-n+p} + \tau_\varepsilon\right)^{\frac{1}{p}}. \tag{1.150}$$

*Proof.* Setting in (1.149) $\varphi = m_\varepsilon$ and the definition of $m_\varepsilon(x)$, we obtain

$$\begin{aligned}
\|\nabla m_\varepsilon\|_{L^p(Y_\varepsilon)}^{p^2} &\le \left( \left| \int_{S_\varepsilon^0} m_\varepsilon ds \right| + \mu_\varepsilon \left| \int_{Y_\varepsilon} m_\varepsilon dx \right| \right)^p \\
&\le \left( \int_{S_\varepsilon^0} |m_\varepsilon| ds + \mu_\varepsilon \times 0 \right)^p \\
&\le \left( \left( \int_{S_\varepsilon^0} 1^{p'} \right)^{\frac{1}{p'}} \left( \int_{S_\varepsilon^0} |m_\varepsilon|^p \right)^{\frac{1}{p}} \right)^p \\
&\le \left( \int_{S_\varepsilon^0} 1 ds \right)^{p-1} \|m_\varepsilon\|_{L^p(S_\varepsilon^0)}^p \\
&\le C_1 a_\varepsilon^{(n-1)(p-1)} \|m_\varepsilon\|_{L^p(S_\varepsilon^0)}^p \le \\
&\le C_2 a_\varepsilon^{(n-1)(p-1)} a_\varepsilon^{n-1} \left( \varepsilon^{-n} \|m_\varepsilon\|_{L^p(Y_\varepsilon)}^p + \tau_\varepsilon \|\nabla m_\varepsilon\|_{L^p(Y_\varepsilon)}^p \right) \\
&\le C_3 a_\varepsilon^{p(n-1)} \left(\varepsilon^{-n+p} + \tau_\varepsilon\right) \|\nabla m_\varepsilon\|_{L^p(Y_\varepsilon)}^p
\end{aligned} \tag{1.151}$$

which concludes the proof. □

We have that

$$\frac{a_\varepsilon^{p(n-1)} \varepsilon^{-n+p}}{a_\varepsilon^{n(p-1)}} = (a_\varepsilon \varepsilon^{-1})^{p-n} \tag{1.152}$$

which in the case $a_\varepsilon = C_0 \varepsilon^\alpha$ results in

$$\frac{a_\varepsilon^{p(n-1)} \varepsilon^{-n+p}}{a_\varepsilon^{n(p-1)}} = C\varepsilon^{(p-n)(\alpha-1)}. \tag{1.153}$$

Using the previous estimates we get:

**Corollary 1.1.** *Let $p > 1$. Then*

$$\|\nabla m_\varepsilon\|_{L^p(Y_\varepsilon)} \leq \begin{cases} C a_\varepsilon^{\frac{n}{p}} & p < n, \\ C a_\varepsilon \ln \left( \dfrac{\varepsilon}{a_\varepsilon} \right)^{\frac{1}{n}} & p = n, \\ C a_\varepsilon^{\frac{n-1}{p-1}} \varepsilon^{\frac{p-n}{p(p-1)}} & p > n. \end{cases} \tag{1.154}$$

This allows us to write the following result:

**Corollary 1.2.** *Let $a_\varepsilon \ll \varepsilon$. Then, since $|\Upsilon_\varepsilon| \sim \varepsilon^{-n}$,*

$$\|\nabla M_\varepsilon\|_{L^p(\cup_j Y_\varepsilon^j)} \leq \begin{cases} C(a_\varepsilon \varepsilon^{-1})^{\frac{n}{p}} & 1 < p \leq n, \\ C(a_\varepsilon \varepsilon^{-1}) \ln(a_\varepsilon^{-1} \varepsilon)^{\frac{1}{n}} & p = n, \\ C(a_\varepsilon \varepsilon^{-1})^{\frac{n-1}{p-1}} & p > n. \end{cases} \tag{1.155}$$

**Corollary 1.3.** *Let $v_\varepsilon \in W^{1,p}(\Omega_\varepsilon, \partial\Omega)$. Then,*

$$\beta(\varepsilon) \int_{S_\varepsilon} v_\varepsilon = \rho_\varepsilon + \beta(\varepsilon)\mu_\varepsilon \sum_{j \in \Upsilon_\varepsilon} \int_{Y_\varepsilon^j} v_\varepsilon dx \tag{1.156}$$

*where*

$$0 \leq \rho_\varepsilon \leq C\beta(\varepsilon) \|M_\varepsilon\|_{L^p(\Omega_\varepsilon)}^{p-1}, \tag{1.157}$$

*and $C$ depends only on $\Omega$ and $\|v_\varepsilon\|_{W^{1,p}(\Omega_\varepsilon)}$.*

*Proof.*

$$\beta(\varepsilon)\int_{S_\varepsilon} v_\varepsilon = \beta(\varepsilon)\sum_{j\in\Upsilon_\varepsilon}\int_{Y_\varepsilon^j} div(|\nabla M_\varepsilon^j|^{p-2}\nabla M_\varepsilon^j v_\varepsilon) =$$

$$= \beta(\varepsilon)\sum_{j\in\Upsilon_\varepsilon}\int_{Y_\varepsilon^j} |\nabla M_\varepsilon^j|^{p-2}\nabla M_\varepsilon^j \nabla v_\varepsilon dx +$$

$$+ \beta(\varepsilon)\sum_{j\in\Upsilon_\varepsilon}\int_{Y_\varepsilon^j} (\Delta_p M_\varepsilon^j) v_\varepsilon dx =$$

$$= \beta(\varepsilon)\sum_{j\in\Upsilon_\varepsilon}\int_{Y_\varepsilon^j} |\nabla M_\varepsilon^j|^{p-2}\nabla M_\varepsilon^j \nabla v_\varepsilon dx +$$

$$+ \beta(\varepsilon)\sum_{j\in\Upsilon_\varepsilon}\mu_\varepsilon\int_{Y_\varepsilon^j} v_\varepsilon dx \tag{1.158}$$

Using Hölder's inequality

$$\beta(\varepsilon)\int_{\Omega_\varepsilon} |\nabla M_\varepsilon|^{p-1}|\nabla v_\varepsilon| dx \le C\beta(\varepsilon)\|M_\varepsilon\|_{L^p(\Omega_\varepsilon)}^{p-1}. \tag{1.159}$$

which concludes the proof. □

This is the reason why critical scales appear in the homogenization process for $p \le n$ and none can appear when $p > n$. The critical case occurs when $\rho_\varepsilon \not\to 0$. In particular, if $\rho_\varepsilon \to C \neq 0$ (where $\rho_\varepsilon$ is the quantity given by (1.156)) then the critical case rises, as we will see in Section 1.4.12.

*Proof of Theorem 1.6.* Due to Corollary 1.3 and

$$\left|\sum_{j\in\Upsilon_\varepsilon}\int_{Y_\varepsilon^j} v_\varepsilon - \int_\Omega v_\varepsilon\right| \le \|v_\varepsilon\|_{L^p(\Omega)}\left|\Omega\setminus\bigcup_{j\in\Upsilon_\varepsilon} Y_\varepsilon^j\right| \to 0, \tag{1.160}$$

which completes the proof. □

Also, this explicit computation explains the *a priori* strange formula for the critical scales. Consider the good scaling $\beta^*(\varepsilon)$.

**Corollary 1.4.** *We have that*

$$\beta^*(\varepsilon)\|M_\varepsilon\|^{p-1}_{L^p(\Omega_\varepsilon)} \leq \begin{cases} Ca_\varepsilon^{\frac{p-n}{p}} \varepsilon^{\frac{n}{p}} & p < n, \\ C\varepsilon \ln\left(a_\varepsilon^{-1}\varepsilon\right)^{\frac{n-1}{n}} & p = n, \\ C\varepsilon & p > n. \end{cases} \tag{1.161}$$

**Remark 1.5.** The right hand side of (1.161) is rather significant. As will see immediately in Theorem 1.6, the fact that this right hand goes to 0 as $\varepsilon \to 0$ is a sufficient condition for the integrals to behave nicely in the limit, and thus will see later that we are in the subcritical case. A priori, these estimates need to be sharp. However, as will see in Section 1.4.8, it is sharp, and $a_\varepsilon^*$ is the value such that the RHS of (1.161) converges to a constant.

**Remark 1.6.** If $a_\varepsilon = C_0\varepsilon^\alpha$ and $\beta(\varepsilon) = \varepsilon^{-(\alpha(n-1)-n)}$, then the result implies that

$$\varepsilon^{-\gamma}\int\limits_{S_\varepsilon} v_\varepsilon dS \to C_0^{n-1}|\partial G_0|\int_\Omega v dx \tag{1.162}$$

as $\varepsilon \to 0$ if $\alpha < \frac{n}{n-p}$. If $\alpha > \frac{n}{n-p}$ then

$$\varepsilon^{-\gamma}\int\limits_{S_\varepsilon} v_\varepsilon dS \to 0. \tag{1.163}$$

#### 1.4.7.1 $L^p - L^q$ estimates for $S_\varepsilon$

It is obvious that there are $L^p - L^q$ estimates for $S_\varepsilon$, in the sense that, if $0 < p < q$, for every $\varepsilon > 0$ there exists a constant $C_\varepsilon$ such that

$$\left(\int_{S_\varepsilon} |v|^p\right)^{\frac{1}{p}} \leq C_\varepsilon \left(\int_{S_\varepsilon} |v|^q\right)^{\frac{1}{q}} \qquad \forall v \in L^q(S_\varepsilon). \tag{1.164}$$

The interesting question is whether we can do this with uniform constant $C_\varepsilon$. The fact is that such results are true, but we have to be careful with the choice of constants. We will use this in the following sections.

**Lemma 1.4.5.** *Let* $0 < p < q$*. Then, there exists* $C$*, independent of* $\varepsilon$*, such that*

$$\left(\beta^*(\varepsilon)\int_{S_\varepsilon} |v|^p\right)^{\frac{1}{p}} \leq C\left(\beta^*(\varepsilon)\int_{S_\varepsilon} |v|^q\right)^{\frac{1}{q}} \qquad \forall v \in L^q(S_\varepsilon). \tag{1.165}$$

### 1.4.8 The critical scales $a_\varepsilon^*$

It has been long present in the literature that the critical size of holes in the case $1 < p < n$ is

$$a_\varepsilon^* = \varepsilon^{\frac{n}{n-p}}. \tag{1.166}$$

We will see in Sections 1.4.10, 1.4.11 and 1.4.12 that the three situations are entirely different. This critical value aligns precisely with estimate (1.161). As indicated in Remark 1.5, $a_\varepsilon^*$ is the value such that the RHS of (1.161) converges to a constant. If we write $a_\varepsilon = C_0 \varepsilon^\alpha$, this critical exponent $\alpha^* = \frac{n}{n-p}$ blows up as $p \to n$.

As it turns out, a critical case also exists for the case $p = n$, and this was studied in [PS15]. The critical choice, as presented in that paper, is the one that satisfies

$$\beta(\varepsilon) a_\varepsilon^{n-1} \varepsilon^{-n} \to C_1^2, \tag{1.167}$$

$$\frac{1}{\beta(\varepsilon)^{\frac{1}{n-1}} a_\varepsilon \ln \frac{4a_\varepsilon}{\varepsilon}} \to -C_2^2, \tag{1.168}$$

where $C_1, C_2 \neq 0$. Again, estimate (1.161) is sharp. Although this a bit more convoluted. Equation (1.167) only indicates $\beta(\varepsilon) \sim \beta^*(\varepsilon)$. Let us read (1.168) carefully

$$\begin{aligned} 1 &\sim -\frac{1}{\beta(\varepsilon)^{\frac{1}{n-1}} a_\varepsilon \ln \frac{4a_\varepsilon}{\varepsilon}} \\ &\sim \frac{1}{\beta(\varepsilon)^{\frac{1}{n-1}} a_\varepsilon \ln \frac{\varepsilon}{4a_\varepsilon}} \end{aligned}$$

Since $D_\varepsilon \sim 1$ is equivalent to $\frac{1}{D_\varepsilon} \sim 1$ we have that

$$\begin{aligned} 1 &\sim \beta(\varepsilon)^{\frac{1}{n-1}} a_\varepsilon \ln \frac{\varepsilon}{a_\varepsilon} \\ &\sim a_\varepsilon^{-1} \varepsilon^{\frac{n}{n-1}} a_\varepsilon \ln \frac{\varepsilon}{a_\varepsilon} \\ &\sim \varepsilon^{\frac{n}{n-1}} \ln \frac{\varepsilon}{a_\varepsilon}. \end{aligned}$$

This is exactly what we anticipated in (1.161). Again, $a_\varepsilon^*$ is the value such that the RHS of (1.161) converges to a constant. We can give the critical scale explictly

$$a_\varepsilon^* = \varepsilon e^{-\varepsilon^{-\frac{n}{n-1}}}. \tag{1.169}$$

We point out that this critical scale is not of the form $a_\varepsilon^* = \varepsilon^\alpha$, but rather the convergence to 0 of $a_\varepsilon^*$ is much faster as $\varepsilon \to 0$.

For $p > n$ we deduce that there exists no critical scale. In [DGCPS17b] the authors first noted that no critical scale in the realm $a_\varepsilon = \varepsilon^\alpha$ may exist. The estimates in this thesis go much further. Since, for $p > n$, the RHS of (1.161) always converges to 0 we can guaranty that no $a_\varepsilon$ critical may exist. In this sense we can say that, for $p > n$, $a_\varepsilon^* = 0$. With this notation, the cases $p > n$ in Theorem 1.7 are a direct improvement of the results in [DGCPS17b].

To summarize, going forward we will keep in mind that the critical scale is precisely

$$a_\varepsilon^* = \begin{cases} \varepsilon^{\frac{n}{n-p}} & 1 < p < n, \\ \varepsilon e^{-\varepsilon^{-\frac{n}{n-1}}} & p = n, \\ 0 & p > n. \end{cases} \tag{1.170}$$

### 1.4.9 Double approximation arguments

In the process of homogenization is typically more convenient to work with a function $\sigma$ that is as smooth as possible. Many authors have only stated their results for such smooth $\sigma$. Since the central theme of this thesis deals with root-type $\sigma$, it was one of our aims to develop a framework to extend the result to general $\sigma$. A natural way to do this, which has been successful in the past, is to consider uniform approximations. This is the subject of this section.

The following comparison results are obtained in [DGCPS17d]. They allow us to extend the results proved for $\sigma$ smooth to the case of $\sigma$ non as smooth.

**Lemma 1.4.6** ([DGCPS17d]). *Let $\sigma, \hat\sigma$ be continuous nondecreasing functions such that $\sigma(0) = 0$ and $u, \hat u$ be their respective solutions of* (1.12) *with $f^\varepsilon = f \in L^{p'}(\Omega)$ and $g^\varepsilon = 0$. Then, there exist constants C depending on p, but independent of $\varepsilon$, such that*

*i) If $1 < p < 2$*

$$\|\nabla(u_\varepsilon - \hat u_\varepsilon)\|_{L^p(\Omega_\varepsilon)} \le C\beta(\varepsilon)\beta^*(\varepsilon)^{-1}\|\sigma - \hat\sigma\|_{\mathscr{C}(\mathbb{R})}\left(\|\nabla u_\varepsilon\|^{2-p}_{L^p(\Omega_\varepsilon)} + \|\nabla \hat u_\varepsilon\|^{2-p}_{L^p(\Omega_\varepsilon)}\right)^{\frac{2}{p}}. \tag{1.171}$$

*ii) If $p \ge 2$ then*

$$\|\nabla(u_\varepsilon - \hat u_\varepsilon)\|^{p-1}_{L^p(\Omega_\varepsilon)} \le C\beta(\varepsilon)\beta^*(\varepsilon)^{-1}\|\sigma - \hat\sigma\|_{\mathscr{C}(\mathbb{R})}. \tag{1.172}$$

We need to study a sufficiently large family of functions $\sigma$ so that the uniform $\mathbb{R}$ approximation by smooth functions is possible. The following condition seems to fit our purposes:

$$|\sigma(t)-\sigma(s)| \leq C(|t-s|^{\alpha}+|t-s|^{p}) \qquad \forall t,s \in \mathbb{R}, \tag{1.173}$$

for some $0<\alpha\leq 1$ and $p\geq 1$. It represents "local Hölder" continuity, in the sense that there is no need for the function to be differentiable. On the other hand, as $|s-t|\to+\infty$, the function $\sigma$ behaves like a power, and then $\sigma$ can be non sublinear at infinity.

We have the following approximation result:

**Lemma 1.4.7.** *Let $\sigma \in \mathscr{C}(\mathbb{R})$, nondecreasing and there exists $0<\alpha\leq 1$, $p>1$ such that (1.173) holds. Then, for every $0<\delta<\frac{1}{4C}$ there exists $\sigma_\delta \in \mathscr{C}(\mathbb{R})$ (piecewise linear) such that*

$$\|\sigma_\delta-\sigma\|_{\mathscr{C}(\mathbb{R})} \leq \delta, \tag{1.174}$$

$$0 \leq \sigma'_\delta \leq D\delta^{1-\frac{1}{\alpha}}, \tag{1.175}$$

*where $D$ depends only on $C,\alpha,p$.*

The idea now is to consider the solution $u_{\varepsilon,\delta}$ of problems

$$\begin{cases} -\Delta_p u_{\varepsilon,\delta} = f^{\varepsilon} & \Omega_\varepsilon, \\ \partial_{\nu_p} u_{\varepsilon,\delta} + \beta(\varepsilon)\sigma_\delta(u_{\varepsilon,\delta}) = \beta(\varepsilon) g^{\varepsilon} & S_\varepsilon, \\ u_{\varepsilon,\delta} = 0 & \partial\Omega. \end{cases} \tag{1.176}$$

The process is the following. Pass to the limit as $\varepsilon\to 0$ for $\delta$ fixed, and characterize the limit of the solution $P_\varepsilon u_{\varepsilon,\delta}$ as $\varepsilon \to 0$ to some function $u_\delta$. Then study the limit of function $u_\delta$ as $\delta\to 0$ to a certain function $\widehat{u}$, the equation for which can be obtained through standard theory. The idea is to construct uniform bounds, that allow us to show $\widehat{u}$ is the limit of $P_\varepsilon u_\varepsilon$. We will see an example of application of this reasoning in the following section.

### 1.4.10 Homogenization of the subcritical cases $a_\varepsilon^* \ll a_\varepsilon \ll \varepsilon$

In this section we will study the limit behaviour for

$$1<p<+\infty \qquad a_\varepsilon^* \ll a_\varepsilon \ll \varepsilon. \tag{1.177}$$

Due to the definitions of $\beta_0$ (see (1.147)) and $\mu_\varepsilon$ (see (1.144)) we have that

$$\beta_0 = |\partial G_0| \lim_{\varepsilon \to 0} \beta(\varepsilon) \beta^*(\varepsilon)^{-1}. \tag{1.178}$$

The aim of this subsection will be to prove the following result:

**Theorem 1.7.** *Let* $1 < p < n$, $f^\varepsilon = f \in L^{p'}(\Omega)$,$g^\varepsilon = g \in W^{1,\infty}(\Omega)$, $a_\varepsilon^* \ll a_\varepsilon \ll \varepsilon$, $\sigma \in \mathscr{C}(\mathbb{R})$ *nondecreasing such that* $\sigma(0) = 0$ *and*

$$|\sigma(v)| \le C(1 + |u|^{p-1}). \tag{1.179}$$

*Then the following results hold:*

*i) Let* $\beta_0 < +\infty$. *Then, up to a subsequence* $P_\varepsilon u_\varepsilon \rightharpoonup u$ *in* $W_0^{1,p}(\Omega)$, *where u is the unique solution of*

$$\begin{cases} -\Delta_p u + \beta_0 \sigma(u) = f + \beta_0 g & \Omega, \\ u = 0 & \partial\Omega. \end{cases} \tag{1.180}$$

*ii) Let* $\beta_0 = +\infty$, $g = 0$ *and* $\sigma \in \mathscr{C}^1$. *Then, up to a subsequence* $P_\varepsilon u_\varepsilon \rightharpoonup u$ *in* $W_0^{1,p}(\Omega)$ *and u satisfies*

$$u(x) \in \sigma^{-1}(0) \tag{1.181}$$

*a.e. in* $\Omega$. *In other words,* $\sigma(u(x)) = 0$ *for a.e.* $x \in \Omega$. *In particular, if* $\sigma$ *is strictly increasing then* $u = 0$.

The regularity of $\sigma$ will be the key difficulty of our approach. As mentioned in the previous section, let us first study the smooth case.

**Remark 1.7.** When $a_\varepsilon = C_0 \varepsilon^\alpha$ and $\beta_\varepsilon = \varepsilon^{-\gamma^*}$ then it is easy to compute that

$$\beta_0 = |\partial G_0| C_0^{n-1}. \tag{1.182}$$

This coefficient is obtained in all cases.

**Smooth kinetic** Just the estimates on the normal derivatives allows to homogenize the noncritical case directly if $\sigma$ is a uniformly Lipschitz continuous function, since in that case the sequence $\sigma(u_\varepsilon)$ in $W^{1,p}(\Omega_\varepsilon, \partial\Omega)$ is bounded, and we can pass to the limit in the standard weak formulation. However, a further analysis allows to do so in many different settings.

Even the case of $\sigma$ monotone nondecreasing such that $\sigma(0) = 0$ and $\sigma'$ locally bounded is easy to understand. Then, we can pass to the limit in formulation (1.93). If we consider test functions $v \in W_0^{1,\infty}(\Omega)$, then $\sigma(v) \in W_0^{1,\infty}(\Omega)$. From the definition of (1.12) it is immediate that $\|\nabla u_\varepsilon\|_{W_0^{1,p}(\Omega)}$ is bounded and hence that $P_\varepsilon u_\varepsilon \rightharpoonup u$ in $W^{1,p}(\Omega_\varepsilon)$. Then, if $\beta(\varepsilon) \sim \beta^*(\varepsilon)$

$$\int_\Omega |\nabla v|^{p-2}\nabla v \cdot (v-u)dx + \beta_0 \int_\Omega \sigma(v)(v-u) \geq \int_{\Omega_\varepsilon} fvdx. \tag{1.183}$$

The case $\beta(\varepsilon) \gg \beta^*(\varepsilon)$ can be studied in a even easier way. In the weak formulation we get

$$\beta^*(\varepsilon)\beta(\varepsilon)^{-1} \int_{\Omega_\varepsilon} |\nabla u_\varepsilon|^p \nabla u_\varepsilon \cdot \nabla v + \beta^*(\varepsilon) \int_{S_\varepsilon} \sigma(u_\varepsilon)v = \beta^*(\varepsilon)\beta(\varepsilon)^{-1} \int_{S_\varepsilon} fv \tag{1.184}$$

for any $v \in W^{1,p}(\Omega_\varepsilon, \delta\Omega)$. Then, at least for $\sigma \in W^{1,\infty}(\mathbb{R})$ monotone nondecreasing such that $\sigma(0) = 0$, as $\varepsilon \to 0$

$$\int_\Omega \sigma(u_\varepsilon)v = 0. \tag{1.185}$$

Hence $\sigma(u_\varepsilon) = 0$. It is important to remark that in the previous literature the limits were identified to $u_\varepsilon = 0$, but this is only because $\sigma$ is required to be strictly increasing (see Table 1.1).

**Non smooth kinetic** The case of $\sigma \in \mathscr{C}(\mathbb{R})$, nondecreasing and $\sigma(0) = 0$ and the case $\beta(\varepsilon)\beta^*(\varepsilon)^{-1} \to 0$ can be treated thanks to the approximation Lemma 1.4.6. In essence

$$\|u_\varepsilon - u_{\varepsilon,\delta}\|_{W^{1,p}} \leq C\|\sigma - \sigma_\delta\|_\infty^\alpha \tag{1.186}$$

for some power $\alpha > 0$, where $\sigma_\delta$ is smooth. Hence, as $\varepsilon \to 0$, $P_\varepsilon u_{\varepsilon,\delta} \rightharpoonup u_\delta$ in $W^{1,p}(\Omega)$, where $u_\delta$ is the solution of

$$\begin{cases} -\Delta_p u_\delta + \beta_0 \sigma_\delta(u_\delta) = f + \beta_0 g & \Omega, \\ u_\delta = 0 & \partial\Omega. \end{cases} \tag{1.187}$$

Furthermore, $P_\varepsilon u_\varepsilon \rightharpoonup u$ in $W^{1,p}(\Omega)$, and the uniform comparison holds in the limit

$$\|u - u_\delta\|_{W^{1,p}} \leq C\|\sigma - \sigma_\delta\|_\infty^\alpha. \tag{1.188}$$

It is easy to show that, as $\delta \to 0$, we have that $u_\delta \rightharpoonup \widehat{u}$ in $W^{1,p}(\Omega)$, the solution of (1.180). As we pass $\delta \to 0$ in (1.188) we deduce that $u = \widehat{u}$.

### 1.4.11 Homogenization of the supercritical case $a_\varepsilon \ll a_\varepsilon^*$

As mentioned before this case is not very relevant. The proof is very simple. Here we present briefly the proof by Shaposhnikova and Zubova in [ZS13].

**Theorem 1.8.** *Let $\sigma \in W^{1,\infty}(\mathbb{R})$ and let us us consider $u_\varepsilon$ the solution of* (1.12) *for $p = 2$, where $f^\varepsilon = f \in L^2(\Omega)$ and $g^\varepsilon = 0$. Let $a_\varepsilon \ll a_\varepsilon^*$. Then, $P_\varepsilon u_\varepsilon \rightharpoonup u$ in $H_0^1(\Omega)$, where $u$ is the unique solution of*

$$\begin{cases} -\Delta u = f & \Omega, \\ u = 0 & \partial\Omega. \end{cases} \tag{1.189}$$

**Remark 1.8.** Notice that the previous result is independent of $\beta(\varepsilon)$. If $g \not\equiv 0$ then, due to Proposition 1.6 we should consider only the cases $\beta(\varepsilon) \ll \beta^*(\varepsilon)$ and $\beta(\varepsilon) \sim \beta^*(\varepsilon)$.

*Proof.* We already know that $P_\varepsilon u_\varepsilon \rightharpoonup u$ in $H_0^1(\Omega)$ independently of $\beta(\varepsilon)$, due to Proposition 1.6, since $g^\varepsilon = 0$. Let

$$K_0 = \max_{y \in G_0} |y|.$$

Let us define, for $j \in \Upsilon_\varepsilon$, functions $\psi_\varepsilon^j \in C_0^\infty(\Omega)$ such that $0 \le \psi_\varepsilon^j \le 1$ and

$$\psi_\varepsilon^j(x) = \begin{cases} 0 & \text{if } |x - P_\varepsilon^j| \ge 2K_0 a_\varepsilon, \\ 1 & \text{if } |x - P_\varepsilon^j| \le K_0 a_\varepsilon, \end{cases} \qquad |\nabla \psi_\varepsilon^j| \le K a_\varepsilon^{-1}, \tag{1.190}$$

and let

$$\psi_\varepsilon = \sum_{j \in \Upsilon_\varepsilon} \psi_\varepsilon^j. \tag{1.191}$$

It is clear that $\psi_\varepsilon = 1$ in $G_\varepsilon^j$. Due to the size of the support, it is also easy to check that

$$\psi_\varepsilon \to 0 \quad \text{in } H^1(\Omega). \tag{1.192}$$

Let $\varphi \in H_0^1(\Omega)$. Taking $\varphi(1-\psi_\varepsilon)$ as a test function in the weak formulation of (1.12) for $p = 2$, we have that

$$\int_{\Omega_\varepsilon} \nabla u_\varepsilon \nabla(\varphi(1-\psi_\varepsilon)) + \beta(\varepsilon) \int_{S_\varepsilon} \sigma(u_\varepsilon)\varphi(1-\psi_\varepsilon) = \int_{\Omega_\varepsilon} f\varphi(1-\psi_\varepsilon). \tag{1.193}$$

Since $(1-\psi_\varepsilon) = 0$ on $S_\varepsilon$ we have that

$$\beta(\varepsilon) \int_{S_\varepsilon} \sigma(u_\varepsilon)\varphi(1-\psi_\varepsilon) = 0. \tag{1.194}$$

On the other hand $\varphi(1-\psi_\varepsilon) \to \varphi$ in $H^1$ we have that the limit as $\varepsilon \to 0$ of equation (1.193) is

$$\int_\Omega \nabla u \nabla \varphi = \int_\Omega f\varphi. \tag{1.195}$$

This completes the proof. □

**Remark 1.9.** As it is clearly seen in the proof, the information of the limit weak formulation is revealed by the choice of a specific sequence of test function. The auxiliary function $\psi_\varepsilon$ oscillate, by construction, with the repetition of particle. This is precisely why this method is known as *oscillating test function*.

In the following sections we will present the results obtained by the author in the critical cases.

### 1.4.12 Homogenization of the critical case when $G_0$ is a ball and $1 < p < n$

In this section we will study the behaviour for

$$1 < p < n \qquad a_\varepsilon = C_0 \varepsilon^\alpha \qquad \beta(\varepsilon) = \varepsilon^{-\gamma} \qquad \alpha = \frac{n}{n-p} \qquad \gamma = \alpha(n-1)-n.$$

In this cases, the limit behaviour is the solution of the following problem:

$$\begin{cases} -\Delta_p u + \mathscr{A}|H(u)|^{p-2}H(u) = f & \Omega \\ u = 0 & \partial\Omega \end{cases} \tag{1.196}$$

where

$$\mathscr{A} = \left(\frac{n-p}{p-1}\right)^{p-1} C_0^{n-p}\omega_n \tag{1.197}$$

and $H$ is the solution of the functional equation

$$\mathscr{B}_0|H(x,s)|^{p-2}H(x,s) = \sigma(x,s-H(x,s)) - g(x) \tag{1.198}$$

$$\mathscr{B}_0 = \left(\frac{n-p}{C_0(p-1)}\right)^{p-1} \tag{1.199}$$

where $g^\varepsilon(x) = g(x)$.

As it can be seen in Table 1.1 there are many previous works in this direction. The term $|H(u)|^{p-2}H(u)$ is usually refered to as *strange term* in homogenization. Since $\sigma$ and $H$ are different functions it can be said that the nature of the reaction changes. This change of behaviour between the critical and subcritical cases has driven some researchers to make a connection between this critical case and the unexpected properties of well-studied elements when presented in nanoparticle form. For example, while a presentation as of gold as microparticles is inert (behaviour at microscale and macroscale coincide), some studies have shown that gold nanoparticles are, in fact, catalysts (see [Sch+02]).

#### 1.4.12.1 Weak convergence

This case is the trickiest. In this direction we first studied the case of power type reactions $\sigma(u) = |u|^{q-1}u$, where $0 \le q < 1$. The case $q = 0$ corresponds to the case of the Heaviside functions (which needs to be understood in the sense of maximal monotone operators). In this direction we published [DGCPS16]. The results and techniques applied in these cases were later generalized in [DGCPS17c], that deals with a general maximal monotone operator and $1 < p < n$. We dealt firstly with the case $g = 0$, $\sigma = \sigma(u)$ and $G_0$ is a ball.

The good setting for this equation is the energy setting, and we consider the definition of solution (1.94).

As its turns out, equation (1.198), which can be rewritten for maximal monotone operators as

$$\mathscr{B}_0 |H(s)|^{p-2} H(s) \subset \sigma(s - H(s)), \tag{1.200}$$

has the following nice property

**Lemma 1.4.8.** *Let $\sigma$ be a maximal monotone operator. Then* (1.200) *has a unique solution $H$, a nondecreasing nonexpansive function $\mathbb{R} \to \mathbb{R}$ (i.e. $0 < H' \le 1$ a.e.).*

In fact, function $H$ can be written in the following way

$$H(r) = (I + \sigma^{-1} \circ \Theta_{n,p})^{-1}(r), \tag{1.201}$$

where

$$\Theta_{n,p}(s) = \mathscr{B}_0 |s|^{p-2} s. \tag{1.202}$$

and $\mathscr{B}_0$ is given by (1.199)

We proved the following

**Theorem 1.9** ([DGCPS17c]). *Let $n \ge 3$, $1 < p < n$, $\alpha = \frac{n}{n-p}$, $\gamma = \alpha(p-1)$ and $G_0$ be a ball. Let $\sigma$ be any maximal monotone operator of $\mathbb{R}^2$ with $0 \in \sigma(0)$ and let $f \in L^{p'}(\Omega)$. Let $u_\varepsilon \in W_0^{1,p}(\Omega_\varepsilon, \partial\Omega)$ be the (unique) weak solution of problem* (1.12) *where $f^\varepsilon = f$ and $g^\varepsilon = 0$. Then there exists an extension $\tilde{u}_\varepsilon$ of $u_\varepsilon$ such that $\tilde{u}_\varepsilon \rightharpoonup u$ in $W_0^{1,p}(\Omega)$ as $\varepsilon \to 0$ where* $u \in W_0^{1,p}(\Omega)$ *is the (unique) weak solution of the problem* (1.196) *associated to the function H, defined by* (1.201).

We seek to apply oscillating test functions of the form $v_\varepsilon = v - H(v)W_\varepsilon$, where $v$ is a test function of the limit problem. For this we define the auxiliary problem

$$W_\varepsilon = \begin{cases} w_\varepsilon(x - P_\varepsilon^j) & x \in B_\varepsilon^j \setminus \overline{G_\varepsilon^j}, \\ 0 & x \notin \cup_j B_\varepsilon^j, \\ 1 & x \in \cup_j G_\varepsilon^j, \end{cases} \tag{1.203}$$

where $P_\varepsilon^j$ is the center of $Y_\varepsilon^j = \varepsilon j + \varepsilon[-\frac{1}{2}, \frac{1}{2}]^n$, $B_\varepsilon^j = B\left(P_\varepsilon^j, \frac{\varepsilon}{4}\right)$ , $G_\varepsilon^j = B\left(P_\varepsilon^j, a_\varepsilon\right)$ and $w_\varepsilon$ is the solution of

$$\begin{cases} -\Delta_p w_\varepsilon = 0 & a_\varepsilon < |x| < \dfrac{\varepsilon}{4}, \\ w_\varepsilon = 0 & |x| = \dfrac{\varepsilon}{4}, \\ w_\varepsilon = 1 & |x| = a_\varepsilon. \end{cases} \tag{1.204}$$

This function can be computed explicitly

$$w_\varepsilon(x) = \frac{|x|^{-\frac{n-p}{p-1}} - \left(\frac{\varepsilon}{4}\right)^{-\frac{n-p}{p-1}}}{a_\varepsilon^{-\frac{n-p}{p-1}} - \left(\frac{\varepsilon}{4}\right)^{-\frac{n-p}{p-1}}}. \tag{1.205}$$

The profile of this radial function can be seen in Figure 1.4

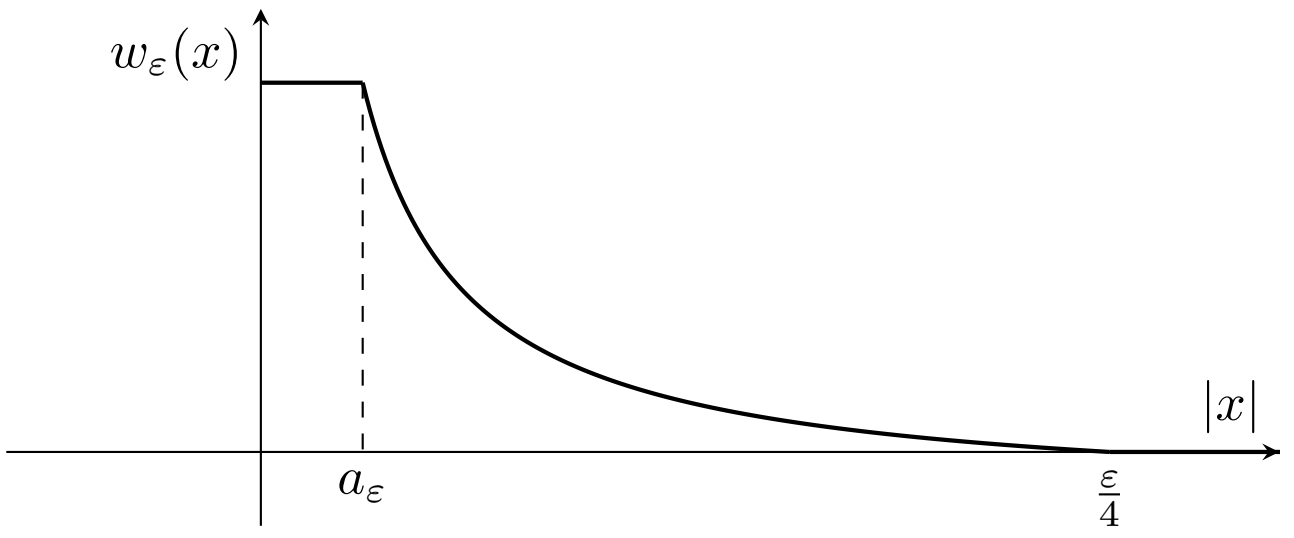

Fig. 1.4 Function $w_\varepsilon$

For $1 \le q \le p$

$$\int_\Omega |\nabla W_\varepsilon|^q dx \le K\varepsilon^{\frac{n(p-q)}{n-p}}, \tag{1.206}$$

hence

$$W_\varepsilon \to 0 \qquad \begin{cases} \text{strongly in } W_0^{1,q}(\Omega) \text{ if } 1 \le q < p, \\ \text{weakly in } W_0^{1,p}(\Omega). \end{cases} \tag{1.207}$$

The second statement may not seem obvious. However, in $W^{1,p}$ the norm is bounded, and hence there must exist a weak limit. This limit must coincide with the $W^{1,q}$ limits for $q < p$, and therefore must be 0.

**Theorem 1.10** (Theorem 1.2 in [DGCPS17c]). *Let $1 < p < n$ and $u_\varepsilon \in W_0^{1,p}(\Omega_\varepsilon, \partial\Omega)$ be a sequence of uniformly bounded norm, $v \in \mathscr{C}_c^\infty(\Omega)$, $h \in W^{1,\infty}(\Omega)$ and let*

$$v_\varepsilon = v - hW_\varepsilon. \tag{1.208}$$

*Then*

$$\lim_{\varepsilon\to 0}\left(\int_{\Omega_\varepsilon} |\nabla v_\varepsilon|^{p-2}\nabla v_\varepsilon \cdot \nabla(v_\varepsilon - u_\varepsilon)dx\right) = \lim_{\varepsilon\to 0}\left(I_{1,\varepsilon} + I_{2,\varepsilon} + I_{3,\varepsilon}\right), \tag{1.209}$$

*where*

$$I_{1,\varepsilon} = \int_{\Omega_\varepsilon} |\nabla v|^{p-2}\nabla v \cdot \nabla(v - u_\varepsilon)dx \tag{1.210}$$

$$I_{2,\varepsilon} = -\varepsilon^{-\gamma}\mathscr{B}_0 \int_{S_\varepsilon} |h|^{p-2}h(v - h - u_\varepsilon)dS \tag{1.211}$$

$$I_{3,\varepsilon} = -A_\varepsilon \varepsilon \sum_{j\in\Upsilon_\varepsilon} \int_{\partial T_\varepsilon^j} |h|^{p-2}h(v - u_\varepsilon)dS, \tag{1.212}$$

*and $A_\varepsilon$ is a bounded sequence. Besides, if $\tilde{u}_\varepsilon$ is an extension of $u_\varepsilon$ and $\tilde{u}_\varepsilon \rightharpoonup u$ in $W_0^{1,p}(\Omega)$ then, for any $v \in W_0^{1,p}(\Omega)$*

$$\lim_{\varepsilon\to 0}\int_{\Omega_\varepsilon} |\nabla v|^{p-2}\nabla v \cdot \nabla(v - hW_\varepsilon - u_\varepsilon)dx = \int_{\Omega} |\nabla v|^{p-2}\nabla v \cdot \nabla(v - u)dx. \tag{1.213}$$

Applying the convexity inequality

$$\Psi(v - H(v)) - \Psi(u_\varepsilon) \le \mathscr{B}_0|H(v)|^{p-2}H(v) \tag{1.214}$$

(this is the reason why $\mathscr{B}_0|H(v)|^{p-2}H(v) \in \sigma(v - H(v))$) and hence the good choice is $h = H(v)$. Thus, we show that we have $P_\varepsilon u_\varepsilon \rightharpoonup u$ in $W^{1,p}(\Omega)$ and $u$ satisfies that, at least for

$v \in W^{1,\infty}(\Omega)$,

$$\int_\Omega |\nabla v|^{p-2}\nabla v \cdot f(v-u) + \mathscr{A} \int_\Omega |H(v)|^{p-2} H(v)(v-u) \geq \int_\Omega f v, \tag{1.215}$$

with $\mathscr{A}$ given by (1.197) which is enough to conclude the result.

#### 1.4.12.2 Some examples

A relevant case in the applications corresponds to the Signorini type boundary condition (1.19), which can be written with the maximal monotone operator (1.24), given

$$H(s) = \begin{cases} H_0(s) & s \geq 0, \\ s & s < 0, \end{cases} \tag{1.216}$$

where

$$\mathscr{B}_0 |H_0(s)|^{p-2} H_0(s) = \sigma_0(s - H_0(s)), \quad s > 0. \tag{1.217}$$

This result was obtained previously in [JNRS14] by *ad hoc* techniques. In [DGCPS17c], we provide it as a corollary of a more general theory.

#### 1.4.12.3 Strong convergence with correctors

It was known in the literature that, at least for smooth $\sigma$,

$$\|u - (u_\varepsilon - H(u_\varepsilon)W_\varepsilon)\|_{W^{1,p}(\Omega_\varepsilon, \partial\Omega)} \to 0 \tag{1.218}$$

as $\varepsilon \to 0$. Nonetheless, it seems that no one had noticed that $W_\varepsilon$ converges strongly to 0 in $W^{1,q}$ for $q < p$. From this fact, we deduce immediately that

$$\|u - u_\varepsilon\|_{W^{1,q}(\Omega_\varepsilon, \partial\Omega)} \to 0, \qquad \text{for } q < p. \tag{1.219}$$

In the case of Signorini boundary conditions we proved the strong convergence (with the corrector term for $q = p$ and without it when $q < p$), which, for the case $1 < p < 2$, was published in [DGCPS17a].

### 1.4.13 Homogenization of the critical case when $G_0$ is not a ball and $p = 2$

For many years, several authors have tried to find a functional equation similar to (1.198) for the case in which $G_0$ is not a ball. In [DGCSZ17] we proved that such equation does not exist in a strict sense. Nonetheless, a equation of form (1.196) still holds, but with a more complicated function $H$.

Let $G_0$ be diffeomorphic to a ball, $p = 2$ and $a_\varepsilon = C_0 \varepsilon^{\alpha^*}$. Then, for any given constant $u \in \mathbb{R}$, we define $\widehat{w}(y; G_0, u)$, for $y \in \mathbb{R}^n \setminus G_0$, as the solution of the following one-parametric family of auxiliary external problems associated to the prescribed asymmetric geometry $G_0$ and the nonlinear microscopic boundary reaction $\sigma(s)$:

$$\begin{cases} -\Delta_y \widehat{w} = 0 & \text{if } y \in \mathbb{R}^n \setminus \overline{G_0}, \\ \partial_{\nu_y} \widehat{w} - C_0 \sigma(u - \widehat{w}) = 0, & \text{if } y \in \partial G_0, \\ \widehat{w} \to 0 & \text{as } |y| \to \infty. \end{cases} \tag{1.220}$$

We will prove in Section 4 that the above auxiliary external problems are well defined and, in particular, there exists a unique solution $\widehat{w}(y; G_0, u) \in H^1(\mathbb{R}^n \setminus \overline{G_0})$, for any $u \in \mathbb{R}$.

**Definition 1.7.** Given $G_0$, we define $H_{G_0} : \mathbb{R} \to \mathbb{R}$ by means of the identity

$$\begin{aligned} H_{G_0}(u) &:= \int_{\partial G_0} \partial_{\nu_y} \widehat{w}(y; G_0, u) \, \mathrm{d}S_y \\ &= C_0 \int_{\partial G_0} \sigma(u - \widehat{w}(y; G_0, u)) \, \mathrm{d}S_y, \quad \text{for any } u \in \mathbb{R}. \end{aligned} \tag{1.221}$$

**Remark 1.10.** Let $G_0 = B_1(0) := \{x \in \mathbb{R}^n : |x| < 1\}$ be the unit ball in $\mathbb{R}^n$. We can find the solution of problem (1.220) in the form $\widehat{w}(y; G_0, u) = \frac{\mathscr{H}(u)}{|y|^{n-2}}$, where, in this case, $\mathscr{H}(u)$ is proportional to $H_{B_1(0)}(u)$. We can compute that

$$\begin{aligned} H_{G_0}(u) &= \int_{\partial G_0} \partial_\nu \widehat{w}(u, y) \, \mathrm{d}S_y \\ &= \int_{\partial G_0} (n-2) H_{G_0}(u) \, \mathrm{d}S_y \\ &= (n-2) \mathscr{H}(u) \omega(n), \end{aligned}$$

where $\omega(n)$ is the area of the unit sphere. Hence, due to (1.221), $\mathscr{H}(u)$ is the unique solution of the following functional equation

$$(n-2)\mathscr{H}(u) = C_0\sigma(u-\mathscr{H}(u)). \tag{1.222}$$

In this case, it is easy to prove that $H$ is nonexpansive (Lipschitz continuous with constant 1). As mentioned before, this equation has been considered in many papers (see [DGCPS17c] and the references therein).

In [DGCSZ17] we proved several results on the regularity and monotonicity of the homogenized reaction $H_{G_0}(u)$ below. Concerning the convergence as $\varepsilon \to 0$ the following statement collects some of the more relevant aspects of this process:

**Theorem 1.11.** *Let $n \geq 3$, $a_\varepsilon = C_0\varepsilon^{-\gamma}$, $\gamma = \frac{n}{n-2}$, $\sigma$ a nondecreasing function such that $\sigma(0) = 0$ and that satisfies* (1.223).

$$|\sigma(s)-\sigma(t)| \leq k_1|s-t|^\alpha + k_2|s-t| \quad \forall s,t \in \mathbb{R}, \quad \textit{for some } 0 < \alpha \leq 1. \tag{1.223}$$

*Let $u_\varepsilon$ be the weak solution of* (1.12) *with $p = 2$, $f^\varepsilon = f \in L^2(\Omega)$ and $g^\varepsilon = 0$. Then there exists an extension to $H_0^1(\Omega)$, still denoted by $u_\varepsilon$, such that $u_\varepsilon \rightharpoonup u_0$ in $H^1(\Omega)$ as $\varepsilon \to 0$, where $u_0 \in H_0^1(\Omega)$ is the unique weak solution of*

$$\begin{cases} -\Delta u_0 + C_0^{n-2} H_{G_0}(u_0) = f & \textit{in } \Omega, \\ u_0 = 0 & \textit{on } \partial\Omega. \end{cases} \tag{1.224}$$

**Remark 1.11.** Since $|H_{G_0}(u)| \leq C(1+|u|)$ it is clear that $H_{G_0}(u_0) \in L^2(\Omega)$.

**Lemma 1.4.9.** *$H_{G_0}$ is a nondecreasing function. Furthermore:*

*i) If $\sigma$ satisfies* (1.223)*, then so does $H_{G_0}$.*

*ii) If $\sigma \in \mathscr{C}^{0,\alpha}(\mathbb{R})$, then so is $H_{G_0}$.*

*iii) If $\sigma \in \mathscr{C}^1(\mathbb{R})$, then $H_{G_0}$ is locally Lipschitz continuous.*

*iv) If $\sigma \in W^{1,\infty}(\mathbb{R})$, then so is $H_{G_0}$.*

## 1.5 Homogenization of the effectiveness factor

We conlude the theoretical results in this chapter by presenting some results on the convergence of the effectiveness factor. They can be found in [DGCT16].

Here, as in [CDLT04], we consider the following regularity assumptions:

$$|g'(v)| \le C(1+|v|^q), \qquad 0 \le q < \frac{n}{n-2}, \tag{1.225}$$

and we consider the strictly increasing and uniformly Lipschitz condition:

$$0 < k_1 \le g'(u) \le k_2. \tag{1.226}$$

We proved the following result, which seems to fulfill the intuitions expressed by Aris in his many works on the subject:

**Theorem 1.12** ([DGCT16]). *Assume that* $p = 2$, $a_\varepsilon = C_0 \varepsilon^\alpha$, $1 \le \alpha < \frac{n}{n-2}$ *and*

- *If* $\alpha = 1$, (1.225),
- *If* $1 < \alpha < \frac{n}{n-2}$, (1.226).

*Then*

$$\frac{1}{|S_\varepsilon|} \int_{S_\varepsilon} \sigma(u_\varepsilon) \to \frac{1}{|\Omega|} \int_\Omega \sigma(u) \qquad \text{as } \varepsilon \to 0. \tag{1.227}$$

This result was later improved in

**Theorem 1.13** ([DGCPS17d]). *Let* $p > 1$, $a_\varepsilon^* \ll a_\varepsilon \ll \varepsilon$, $\beta \sim \beta^*$ *and* $\sigma$ *be continuous such that* $\sigma(0) = 0$. *Let* $u_\varepsilon$ *and* $u$ *be the solutions of* (1.12) *and* (1.180). *Lastly, assume either:*

*i)* $\sigma$ *is uniformly Lipschitz continuous (*$\sigma' \in L^\infty$*), or*

*ii)* $\sigma \in \mathscr{C}(\mathbb{R})$ *and there exists* $0 < \alpha \le 1$ *and* $q > 1$ *such that we have* (1.173) *and*

$$(\sigma(t) - \sigma(s))(t-s) \ge C|t-s|^q, \qquad \forall t, s \in \mathbb{R}. \tag{1.228}$$

*Then* (1.227) *holds.*

**Remark 1.12.** As we have seen, the behaviour of $\int_{S_\varepsilon}$ in the critical case is more convoluted. Thus, a convergence of type (1.227) should not be expected. However, results of similar nature, applying the strong convergence with corrector (1.218) are work under development.

## 1.6 Pointwise comparison of solutions of critical and non-critical solutions

Since we do not have a natural definition of effectiveness in the critical case, the claim that it is "more effective" than the non critical case -a claim that is often made in the Nanotechnology community- is difficult to test. However, we know that, for the non critical cases the effectiveness is increasing with the value of $w$. Thus, we have studied whether we can find a pointwise comparison of the critical and noncritical limits.

Assume in (1.9) that $\hat{g}_\varepsilon = 0$. Then (1.12) becomes

$$\begin{cases} -\Delta_p u_\varepsilon = 0 & x \in \Omega_\varepsilon, \\ \partial_{\nu_p} u_\varepsilon + \varepsilon^{-\gamma}\sigma(u_\varepsilon) = \varepsilon^{-\gamma}\hat{\sigma}(1) & x \in S_\varepsilon, \\ u_\varepsilon = 0 & x \in \partial\Omega. \end{cases} \tag{1.229}$$

Notice that the presence of $w_\varepsilon = 1$ on $\partial\Omega$ is translated to a source in $S_\varepsilon$ for $u_\varepsilon$. When $a_\varepsilon \sim a_\varepsilon^*$, the strange term $H$ is the solution of

$$\mathscr{B}_0 |H(s)|^{p-2} H(s) = \sigma(s - H(s)) - \hat{\sigma}(1). \tag{1.230}$$

Then $w_\varepsilon$ converges weakly in $W^{1,p}(\Omega)$ to $w_{\text{crit}}$, the solution of

$$\begin{cases} -\Delta_p w_{\text{crit}} + \mathscr{A} |h(w_{\text{crit}})|^{p-2} h(w_{\text{crit}}) = 0 & \Omega, \\ w_{\text{crit}} = 1 & \partial\Omega, \end{cases} \tag{1.231}$$

and $h$ is given by

$$|h(w)|^{p-2} h(w) = \hat{\sigma}(1) - |H(1-w)|^{p-2} H(1-w). \tag{1.232}$$

In the noncritical cases, $a_\varepsilon^* \ll a_\varepsilon = C_0 \varepsilon^\alpha \ll \varepsilon$, we know that an extension of $w_\varepsilon$ converges weakly in $W^{1,p}(\Omega)$ to $w_{\text{non-crit}}$, the solution of

$$\begin{cases} -\Delta_p w_{\text{non-crit}} + \hat{\mathscr{A}} \hat{\sigma}(w_{\text{non-crit}}) = 0 & \Omega, \\ w_{\text{non-crit}} = 1 & \partial\Omega, \end{cases} \tag{1.233}$$

with $\hat{\mathscr{A}} = C_0^{n-1} |\partial G_0|$.

We showed, first in [DGC17] under restricted assumptions and later in [DGCPS17c], that a pointwise comparison holds. We stated the following theorem:

**Theorem 1.14** ([DGCPS17c])**.** *Let* $n \geq 3, p \in [2,n)$*,* $a_\varepsilon \sim a_\varepsilon^*$*,* $f^\varepsilon = 0$ *and* $\hat{\sigma} \in \mathscr{C}(\mathbb{R})$ *non decreasing such that* $\sigma(0) = 0$ *Then, we have that*

$$w_{crit} \geq w_{non\text{-}crit}. \tag{1.234}$$

The critical case produces a pointwise " better" reaction.

## 1.7 Some numerical work for the case $\alpha = 1$

To obtain explicit numerical solutions of the different homogeneous and nonhomogeneous problems COMSOL Multiphysics was applied[1]. Also, using the LiveLink tool, it allows to create a Matlab code that we have used to generalize the construction of the obstacles in our domains.

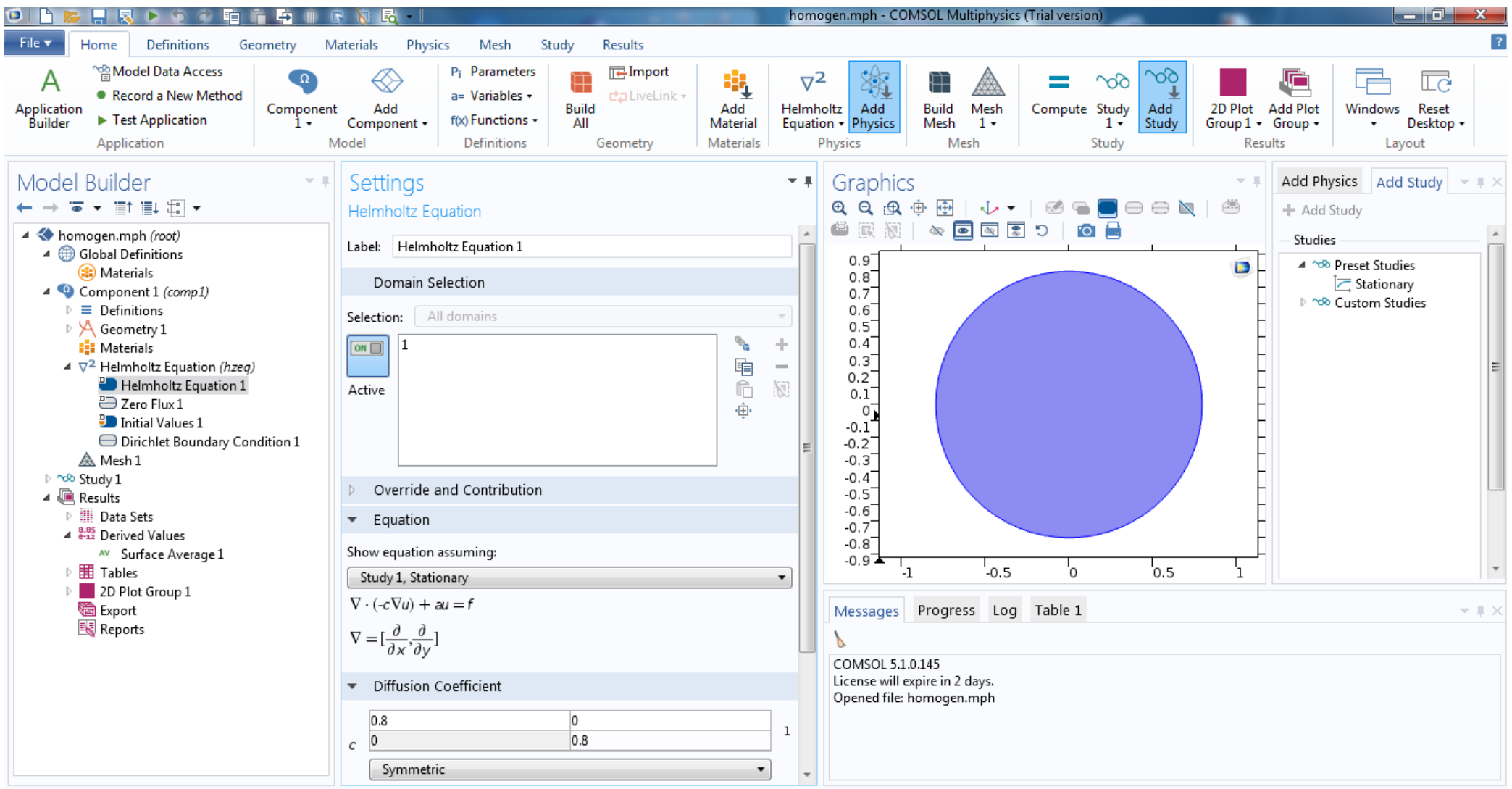

Fig. 1.5 Interfase of the COMSOL software

[1]The author wishes to thank Carlos Arechalde Pérez, Pablo Cañones Martín, Denis Coccolo Góngora, Nadia Loy and Amarpreet Kaur for their work during the IX Modelling Week UCM 2015, where the images were produced under the guidance of this author.

### 1.7.1 Numerical solutions of the non-homogeneous problem

We have simulated that each pellet is inside a periodicity cell. The input parameter of the function is $\varepsilon$, the side of this periodicity cell, that has four times the area of the pellet. Figure 1.6 shows what happened if we change the value of $\varepsilon$

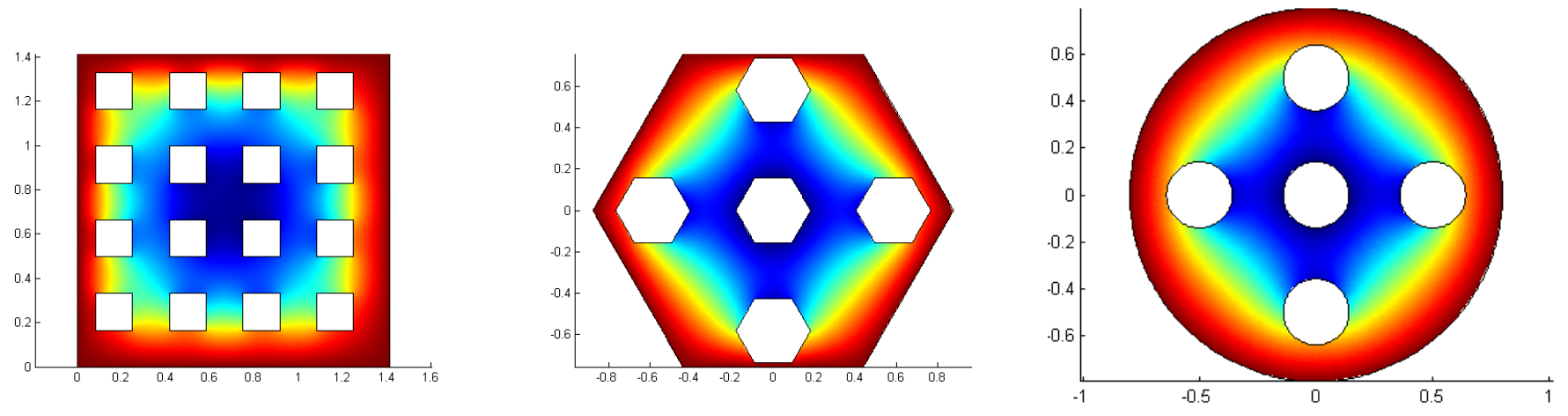

(a) $\varepsilon = \frac{1}{3}$

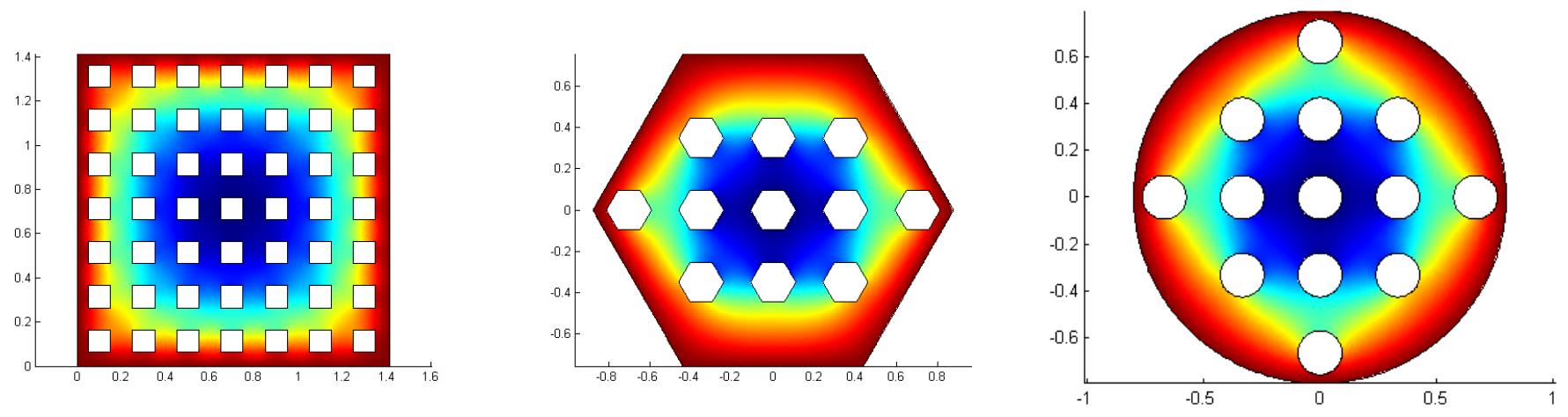

(b) $\varepsilon = \frac{1}{5}$

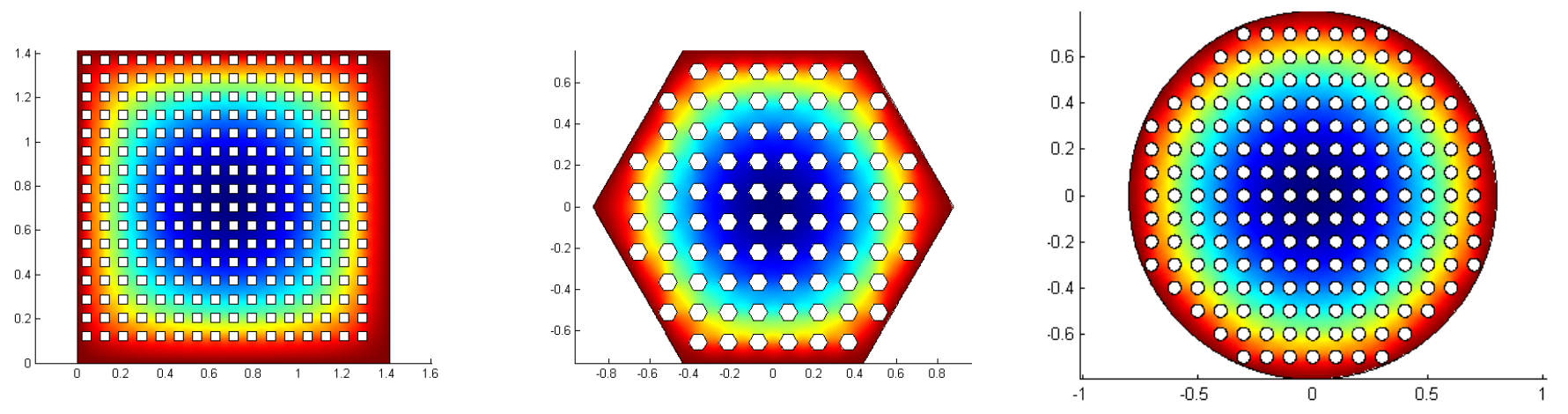

(c) $\varepsilon = \frac{1}{12}$

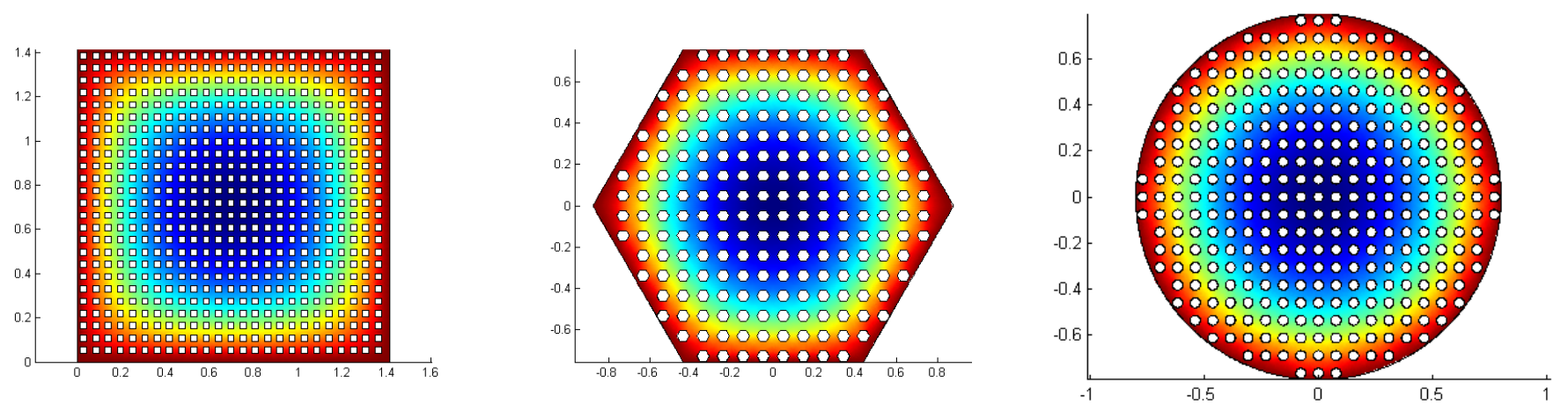

(d) $\varepsilon = \frac{1}{18}$

Fig. 1.6 Level set of the solution of (1.8) for $A^{\varepsilon} = I$, $\sigma(u) = u$ and $a_{\varepsilon} = \varepsilon$. Different values of $\varepsilon$, of domain $\Omega$ and $G_0$ are presented

### 1.7.2 Numerical solutions of the cell problem

The cell problem has solutions which are rather characteristic. If the domain $G_0$ is symmetric with respect to the axis so is the solution. Due to the periodic boundary conditions, is usually not easy to simulate the solution with "black box" software.

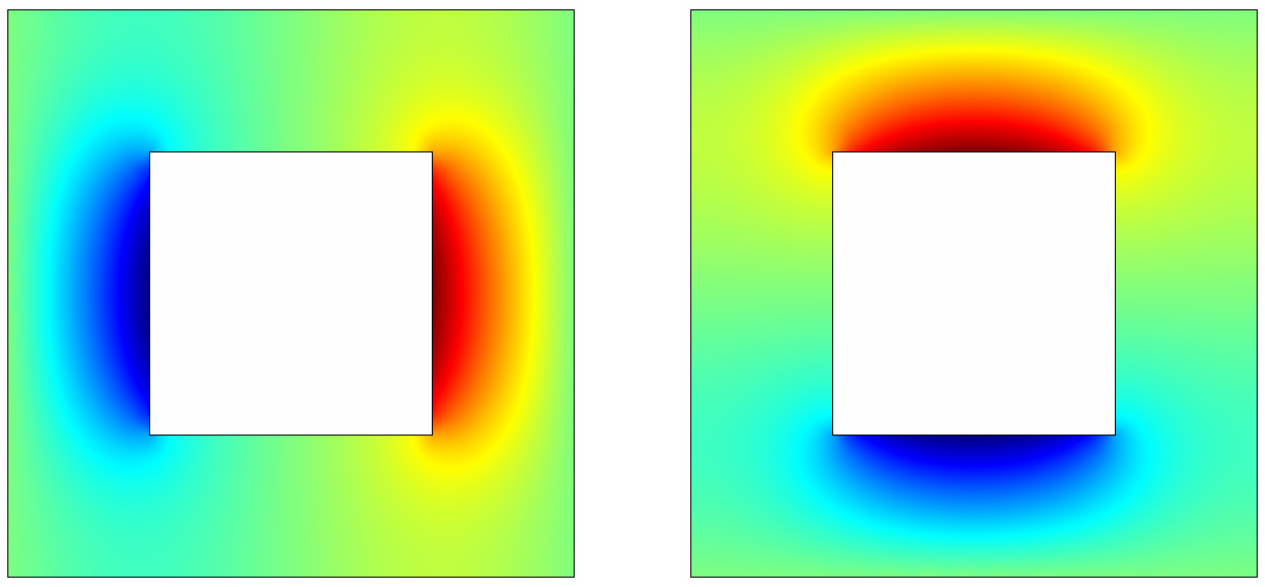

(a) $\chi_1$ solution (b) $\chi_2$ solution

Fig. 1.7 Level set of the solutions of (1.72) for $G_0$ a square

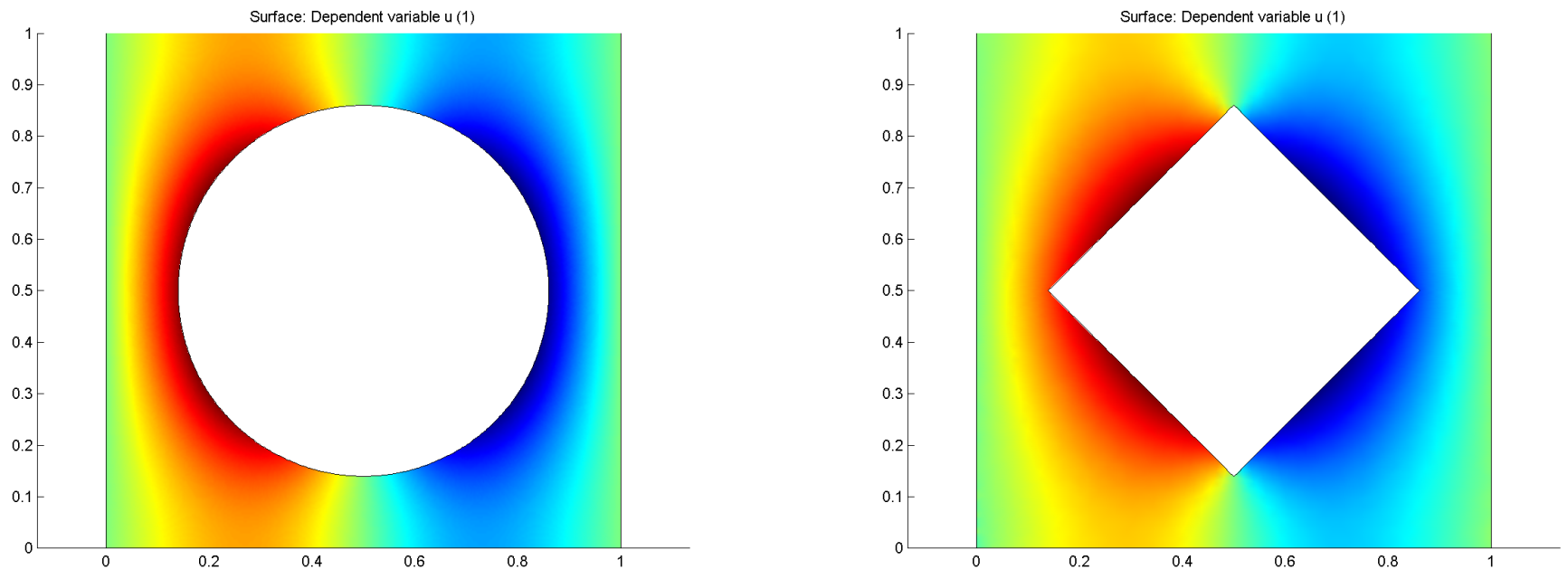

(a) Circular inclusions (b) Square inclusions

Fig. 1.8 Two obstacles $T$, and the level sets of the solution of the cell problem (1.72)

### 1.7.3 Numerical solutions of the homogeneous problem

From all of the nonhomogeneous simulations, the most interesting results are obtained for the smallest $\varepsilon$, because we can see that diffussion in this case is pretty similar to the homogenized problem, as we expected, because of the theoretical results.

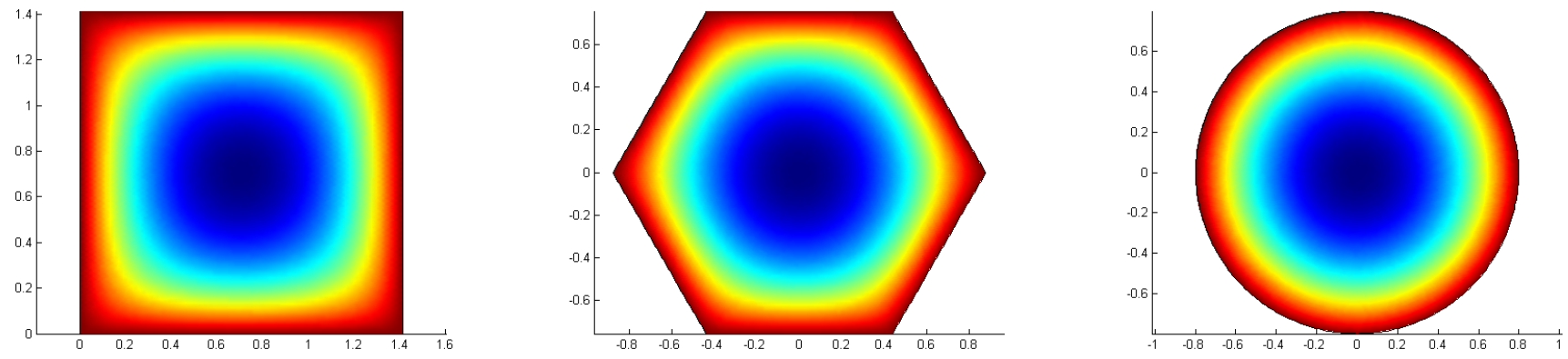

Fig. 1.9 Level sets of the solution of the homogenized problem (1.76), corresponding to the different cases in Figure 1.6

### 1.7.4 Approximation of the numerical solutions

The $L^2$ convergence is guarantied by the theoretical results.

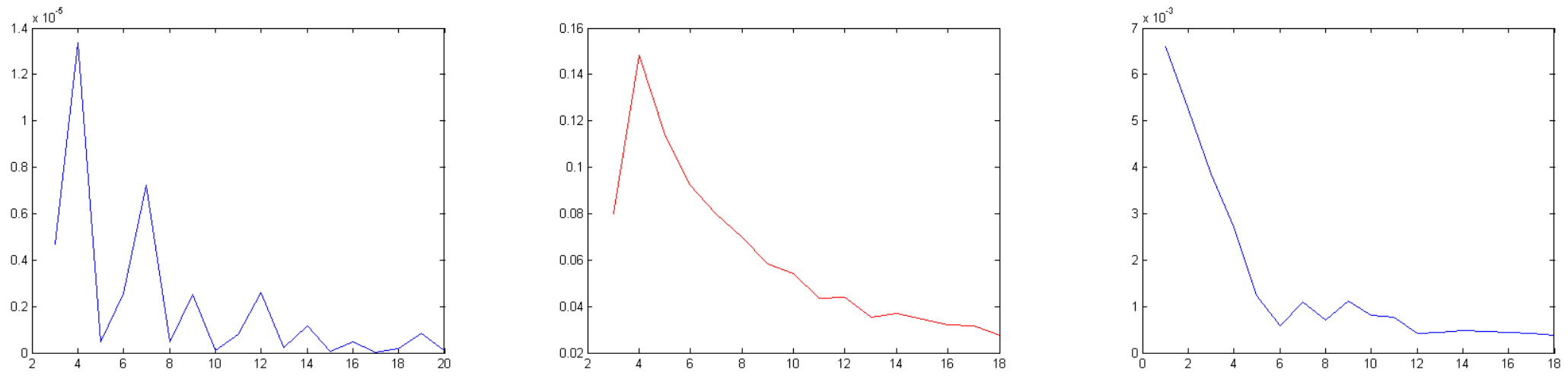

(a) Squares inside square (b) Hexagons inside hexagon (c) Circles inside circle

Fig. 1.10 $L^2$ norm convergence of $\widetilde{u}_\varepsilon \to u$

Nonetheless, the $L^\infty$ convergence has never been proven in the theorical setting. However, the numerical solutions seem to converge.

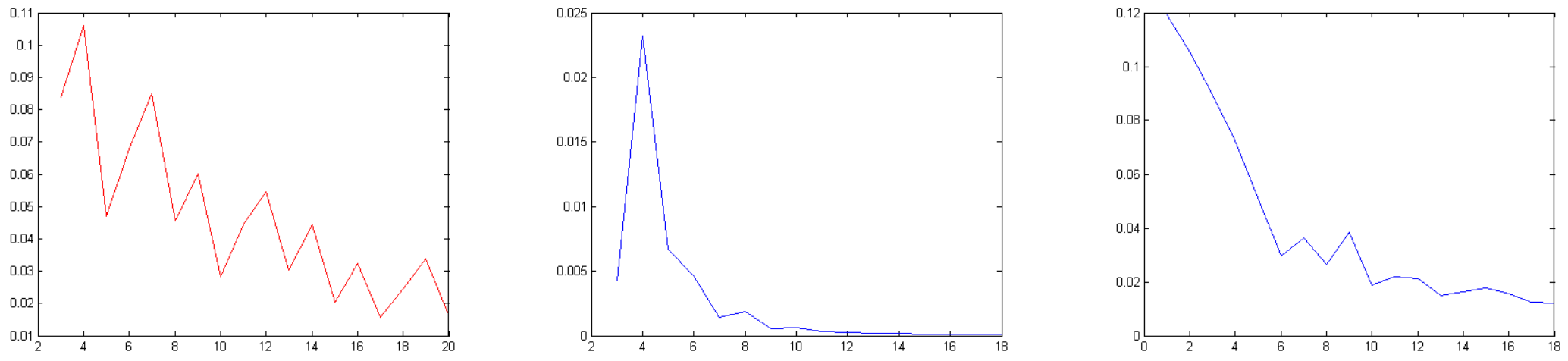

(a) Squares inside square (b) Hexagons inside hexagon (c) Circles inside circle

Fig. 1.11 $L^\infty$ norm convergence of $\widetilde{u}_\varepsilon \to u$

## 1.7.5 Convergence of the effectiveness

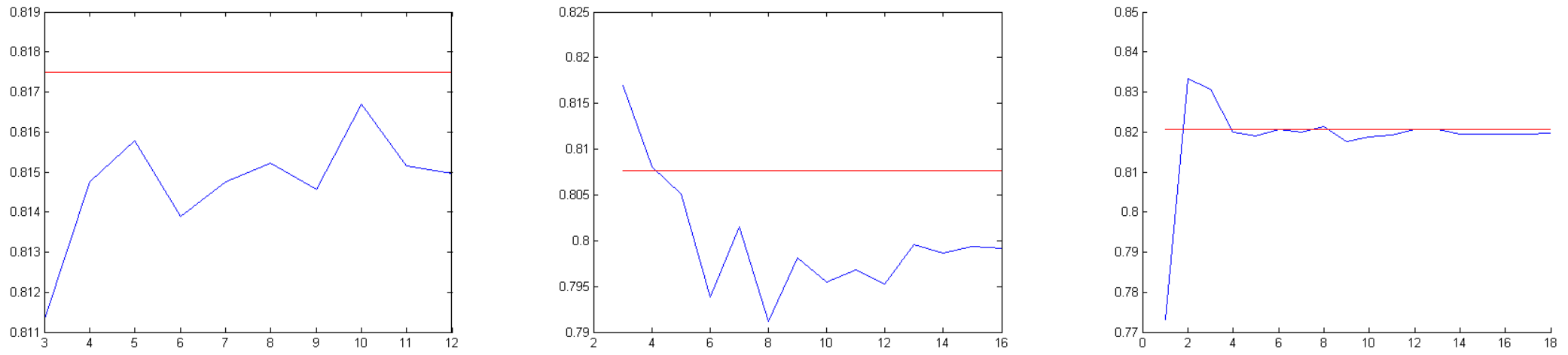

(a) Squares inside square (b) Hexagons inside hexagon (c) Circles inside circle

Fig. 1.12 Convergence efectiveness result: Red line shows the value of non homogeneous problem. Blue line shows the convergence of the homogeneous problem as a function of the value $n = \frac{1}{\varepsilon}$. Notice the order of magnitude in the graphs.

# Appendix 1.A Explanation of [Gon97]

The first paper to properly characterize the change of nonlinear kinetic is [Gon97], which applies the technique of Γ-convergence when $G_0$ is a ball. However, there are a few steps that are not clear (at least to a part of the community). We will try to clarify them in this section.

## 1.A.1 A Γ-convergence theorem

First, we introduce a Γ-convergence theorem proved in [Gon97] under *ad-hoc* assumptions. More general statements of similar nature can be found in [Dal93].

**Theorem 1.15** ([Gon97]). *Let $X_\varepsilon, X$ be Hilbert spaces and $\Phi_\varepsilon$ and $\Phi$ be functionals in these spaces. Let us assume that $\Phi_\varepsilon$ satisfy the following conditions:*

*i) There exists $\theta > 0$ such that*

$$\Phi_\varepsilon(u+v) \geq \Phi_\varepsilon(u) + L_\varepsilon(u;v) + C\|v\|_\varepsilon^\theta \qquad \forall u, v \in X_\varepsilon \tag{1.235}$$

*holds, where $L_\varepsilon$ is the linear functional with respect to v given by the Fréchet differential of $\Phi_\varepsilon$ at a point u,*

*ii) $u_\varepsilon \in X_\varepsilon$ is a minimizer of $\Phi_\varepsilon$ such that*

$$\Phi_\varepsilon(0) \geq \Phi_\varepsilon(u^\varepsilon) \geq C_1\|u_\varepsilon\|_\varepsilon^2 - C_2 \tag{1.236}$$

*where $C_1$ and $C_2$ are constants.*

*Suppose there exists a set $M \subset X$ that is everywhere dense in X and operators $P_\varepsilon : X_\varepsilon \to X, R_\varepsilon : M \to X_\varepsilon$ satisfying the conditions*

*i) $\|P_\varepsilon w^\varepsilon\| \leq C\|w^\varepsilon\|_\varepsilon, \forall w^\varepsilon \in X_\varepsilon$*

*ii) $P_\varepsilon R_\varepsilon w \to w$ weakly in X as $\varepsilon \to 0$, for every $w \in M$*

*iii) $\lim_{\varepsilon\to 0} \Phi_\varepsilon R_\varepsilon w = \Phi w$, for all $w \in M$*

*iv) For any $\gamma^\varepsilon \in X_\varepsilon$ such that $P_\varepsilon\gamma_\varepsilon \to \gamma$ weakly in X and any $w \in M$*

$$\overline{\lim_{\varepsilon\to 0}} |L_\varepsilon(R_\varepsilon w; \gamma_\varepsilon)| \leq \Psi(\|w\|)\|\gamma\| \tag{1.237}$$

*where $\Psi(t)$ is a continuous function of $t \geq 0$.*

*Then* $P_\varepsilon u^\varepsilon \to u$ *weakly in* $X$ *as* $\varepsilon \to 0$ *and* $u$ *is a minimizer of* $\Phi$.

With Goncharenko's notation $P_\varepsilon$ is a sort of extension operator (at least asymptotically), where $R_\varepsilon$ is an adaptation operator.

### 1.A.2 Proof of Theorem 1.2

The result stated by Goncharenko is Theorem 1.2. Let us prove this result.

The choice of spaces is, naturally,

$$X = H_0^1(\Omega) \qquad X_\varepsilon = H^1(\Omega_\varepsilon, \partial\Omega) \qquad M = C_c^2(\Omega). \tag{1.238}$$

Let us define

$$V_\varepsilon = \left\{ v \in X_\varepsilon : \frac{\partial v}{\partial n} + \sigma^\varepsilon(v) = 0, S_\varepsilon \right\}.$$

The main argument of the paper is to show the $\Gamma$-convergence of the energy functional. We define the energy function to be minimized by the solutions:

$$\Phi_\varepsilon(v^\varepsilon) = \int_{\Omega^\varepsilon} (|\nabla v^\varepsilon|^2 + 2 f^\varepsilon v^\varepsilon) \mathrm{d}x + \int_{S_\varepsilon} \varepsilon^{-\gamma} \rho(v^\varepsilon) \mathrm{d}\Gamma.$$

It is easy to check that $\Phi_\varepsilon$ satifies the hypothesis of Theorem 1.15. The technique of the proof by Goncharenko passes by the construction of the following operator

$$R_\varepsilon : M = \mathscr{C}^2(\Omega) \to V^\varepsilon$$

which imposes the boundary condition to any $C^2$ functions. Let $\varphi(t)$ be continuous ($0 \leq \varphi \leq 1$) and such that

$$\varphi(t) = \begin{cases} 1 & t \leq \frac{3}{2}. \\ 0 & t \geq 2. \end{cases}$$

Let us suppose there is only one particle, a ball of center 0. Let us say we want a behaviour of the type

$$R_\varepsilon w \sim \begin{cases} w(x_i) & \varepsilon^\alpha - scale, \\ w(x) + F_\varepsilon w(x) & \varepsilon - scale, \\ w(x) & 1 - scale, \end{cases}$$

for some operator $F_\varepsilon$, to be defined later. Then we can define $R_\varepsilon$ in three parts

$$R_\varepsilon w = w(x) + (w(0) - w(x))\varphi\left(\frac{|x|}{\varepsilon^\alpha}\right) + F_\varepsilon w(x)\varphi\left(\frac{4|x|}{\varepsilon}\right),$$

each one of them captures the behaviour in three different zones. Since the particle is a ball we can try with a operator $F_\varepsilon$ that yields radial functions: $F_\varepsilon w(x) = F_\varepsilon w(|x|)$.

Let us impose the condition $R_\varepsilon w \in V^\varepsilon$. The definition of this function in the paper comes out of the blue to serve its purpose. We wanted here to explain the rationale behind this choice. First, in a neighbourhood of $\{|x| = \varepsilon^\alpha\} = G_\varepsilon^0$ we find that

$$\varphi\left(\frac{|x|}{\varepsilon^\alpha}\right) = 1,$$

so in this neighborhood

$$R_\varepsilon w(x) = w(x_0) + F_\varepsilon w(x), \qquad \varepsilon^\alpha < |x| < \frac{3}{8}\varepsilon^\alpha.$$

Therefore

$$\nabla R_\varepsilon w(x) = (F_\varepsilon w)'(|x|)\frac{x}{|x|}$$

and, on $\{|x| = \varepsilon^\alpha\}$, we have

$$\frac{\partial R_\varepsilon w}{\partial n} + \varepsilon^{-\gamma}\sigma(R_\varepsilon w(x)) = (F_\varepsilon w)'(\varepsilon^\alpha) + \varepsilon^{-\gamma}\sigma(w(0) + F_\varepsilon w(\varepsilon^\alpha)).$$

where, we remind that $\gamma = 2\alpha - 3$. Hence, it will be useful to take $F_\varepsilon w$ a function such that $\phi' = -\phi^2$. This is, precisely $\phi(s) = \frac{A}{s}$. If $F_\varepsilon w(|x|) = -\frac{A^\varepsilon}{|\varepsilon|}$ (so that the derivative later on has a nice sign) we have that

$$A^\varepsilon \varepsilon^{-2\alpha} = \varepsilon^{-2\alpha+3}\sigma(w(0) - A^\varepsilon \varepsilon^{-\alpha}).$$

Hence

$$A^\varepsilon \varepsilon^{-3} = \sigma(w(0) - A^\varepsilon \varepsilon^{-\alpha}).$$

Now take $B^\varepsilon = A^\varepsilon \varepsilon^{-3}$ then

$$B^\varepsilon = \sigma(w(0) + B^\varepsilon \varepsilon^{3-\varepsilon}),$$

in the limit the equation would become

$$\begin{cases} B^0 = \sigma(w(0)) & \alpha < 3, \\ B^0 = \sigma(w(0) - B^0) & \alpha = 3. \end{cases}$$

The previous reasoning explains the following choices. When $\alpha = \gamma = 3$ let us take $A^\varepsilon = A\varepsilon^3$ where $A$ is the solution of the implicit equation

$$A = \sigma(w(0) - A).$$

Since $H$ is the solution of (1.81) we have that

$$A = H(w(0)).$$

On the other hand, for $\alpha < 3$ let us choose simply $A = \sigma(w(0))$.

Eventually

$$R_\varepsilon w = w(x) + (w(0) - w(x))\varphi\left(\frac{|x|}{\varepsilon^\alpha}\right) - \frac{A^\varepsilon}{|x|}\varphi\left(\frac{4|x|}{\varepsilon}\right).$$

It is important that we do not loose the notion that $0 \le A^\varepsilon \le \varepsilon^3$. Eventually, since the particles are balls located at the points $x_i$

$$R_\varepsilon w(x) = w(x) + \sum_{i=1}^{N(\varepsilon)} (w(x_i) - w(x))\varphi\left(\frac{|x - x_i|}{\varepsilon^\alpha}\right) - \sum_{i=1}^{N(\varepsilon)} \frac{A_i^\varepsilon}{|x - x_i|}\varphi\left(\frac{4|x - x_i|}{\varepsilon}\right),$$

where

$$A_i^\varepsilon = A_i\varepsilon^3, \qquad A_i = \begin{cases} \sigma(w(x_i)) & \alpha < 3, \\ H(w(x_i)) & \alpha = 3. \end{cases} \tag{1.239}$$

The deduction of the estimates of the convergence were also not detailed in Goncharenko's paper. We give the details in the following lines. We have that

$$\nabla R_\varepsilon w = \nabla w + \sum_{i=1}^{N(\varepsilon)} (-\nabla w(x))\varphi\left(\frac{|x - x_i|}{\varepsilon^\alpha}\right) + \sum_{i=1}^{N(\varepsilon)} (w(x_i) - w(x))\varphi'\left(\frac{|x - x_i|}{\varepsilon^\alpha}\right)\frac{x - x_i}{|x - x_i|\varepsilon^\alpha} \tag{1.240}$$

$$+ \sum_{i=1}^{N(\varepsilon)} \frac{A_i^\varepsilon}{|x - x_i|^2}\varphi\left(\frac{4|x - x_i|}{\varepsilon}\right) - \sum_{i=1}^{N(\varepsilon)} \frac{A_i^\varepsilon}{|x - x_i|}\varphi'\left(\frac{4|x - x_i|}{\varepsilon}\right)\frac{x - x_i}{|x - x_i|\varepsilon}. \tag{1.241}$$

When integrating the squares of the above (which is easily done is spherical coordinates) only a couple of terms survive: for sure the terms without $\varphi$ and, from the ones with $\varphi$ only the most singular one, the term

$$\frac{A_i^\varepsilon}{|x-x_i|^2}\varphi\left(\frac{4|x-x_i|}{\varepsilon}\right)$$

is to be significant. Let us see how. First we integrate

$$\int_{\Omega_\varepsilon}\left|\frac{A_i^\varepsilon}{|x-x_i|^2}\varphi\left(\frac{4|x-x_i|}{\varepsilon}\right)\right|^2 \mathrm{d}x = A_i^2\varepsilon^6\int_{\varepsilon^\alpha}^{C\varepsilon}\frac{1}{r^4}\varphi\left(\frac{4r}{\varepsilon}\right)r^2\mathrm{d}r$$
$$= A_i^2\varepsilon^5\int_{\varepsilon^{\alpha-1}}^{1}\frac{1}{s^2}\varphi(4s)\mathrm{d}s \sim A^2\varepsilon^{6-\alpha}.$$

The other terms go to 0 as $\varepsilon^p, p>3$ and so

$$\int_{\Omega_\varepsilon}|R_\varepsilon w|^2\mathrm{d}x = \int_{\Omega_\varepsilon}|w|^2+\sum_{i=1}^{N(\varepsilon)}C\varepsilon^p = \int_{\Omega_\varepsilon}|w|^2+C\varepsilon^{p-3}, \qquad p>3,$$

and

$$\int_{\Omega_\varepsilon}|\nabla R_\varepsilon w|^2\mathrm{d}x = \int_{\Omega_\varepsilon}|\nabla w|^2+\sum_{i=1}^{N(\varepsilon)}A_i\varepsilon^{3-\alpha}+C\varepsilon^{p-3}, \qquad p>3,$$

since $A_i = H(w(x_i))$ as the author says

$$\Phi_\varepsilon(R_\varepsilon w) = \int_{\Omega^\varepsilon}(|\nabla w|^2++2f^\varepsilon R_\varepsilon w)\mathrm{d}x+\sum_{i=1}^{N(\varepsilon)}A_i^2\varepsilon^{6-\alpha}$$
$$+\int_{S_\varepsilon}\varepsilon^{-\gamma}\rho(R_\varepsilon w)\mathrm{d}\Gamma+E(\varepsilon,w)$$
$$= \int_{\Omega^\varepsilon}(|\nabla w|^2+2f^\varepsilon R_\varepsilon w)\mathrm{d}x+\sum_{i=1}^{N(\varepsilon)}A_i^2\varepsilon^3$$
$$+\varepsilon^{-\gamma}\sum_{i\in\Upsilon_\varepsilon}\int_{G_i^\varepsilon}\rho(w(x_i)-A_i\varepsilon^{\varepsilon-\alpha})\mathrm{d}\Gamma+E(\varepsilon,w).$$

where $E(\varepsilon,w)$ goes to 0. Classical integration results guaranty the convergence of the Riemann sum. When $\alpha=3$

$$\sum_{i=1}^{N(\varepsilon)}A_i^2\varepsilon^{6-\alpha} = \sum_{i=1}^{N(\varepsilon)}H(w(x_i))^2\varepsilon^3 \to \int_\Omega H(w(x))^2.$$

Notice that, if $\alpha < 3$, then

$$\sum_{i=1}^{N(\varepsilon)} A_i^2 \varepsilon^{6-\alpha} = \varepsilon^{3-\alpha} \sum_{i=1}^{N(\varepsilon)} \sigma(w(x_i))^2 \varepsilon^3 \to 0.$$

On the other hand

$$\varepsilon^{-\gamma} \sum_{i \in \Upsilon_\varepsilon} \int_{G_i^\varepsilon} \rho(w(x_i) - A_i \varepsilon^{\varepsilon-\alpha}) \mathrm{d}\Gamma \to 4\pi \begin{cases} \int_\Omega \rho(x) & \alpha < 3, \\ \int_\Omega \rho(x - H(x)) & \alpha = 3. \end{cases} \tag{1.242}$$

Thus, the $\Gamma$-limit is

$$\Phi(w) = \int_\Omega |\nabla w|^2 + 2fw + \begin{cases} \rho(x) & 2 < \alpha < 3, \\ H(w(x))^2 + \rho(x - w(x)) & \alpha = 3. \end{cases} \tag{1.243}$$

Hence, we see the appearance of the "strange term" for $\alpha = 3$. If $\alpha < 3$ then we do not have this term, as it is noted on the paper. Surprisingly, note that the strange comes out of the diffusion operator. Applying Theorem 1.15 we conclude the proof of Theorem 1.2.

# Chapter 2

# Optimizing the effectiveness: symmetrization techniques

From this chapter on, due to the common practice of notation in this fields, which does not coincide with the practice the homogenization community will use the notation that follows for the homogeneous problem derived in the previous section:

$$\begin{cases} -\Delta w + \beta(w) = \hat{f} & \text{in } \Omega, \\ w = 1 & \text{on } \partial\Omega. \end{cases} \tag{2.1}$$

As we will play now with different domains $\Omega$ we will denote this solution $w_\Omega$. By introducing the change in variable $u = 1 - w$ the problem can be reformulated as

$$\begin{cases} -\Delta u + g(u) = f & \text{in } \Omega, \\ u = 0 & \text{on } \partial\Omega, \end{cases} \tag{2.2}$$

where $g(u) = \beta(1) - \beta(1-u)$ and $f = \beta(1) - \hat{f}$. In this case we write the effectiveness factor as:

$$\mathscr{E}(\Omega) = \frac{1}{|\Omega|}\int_\Omega \beta(w_\Omega)dx \tag{2.3}$$

and the ineffectiveness $\eta(\Omega) = \beta(1) - \mathscr{E}(\Omega)$ as

$$\eta(\Omega) = \frac{1}{|\Omega|}\int_\Omega g(u_\Omega)dx. \tag{2.4}$$

Roughly, we aim to find extremal sets $\Omega$ which maximize and minimize this functional, applying rearrangement techniques.

## 2.1 Geometric rearrangement: Steiner and Schwarz

As stated by Polya and Szegö [PS51] *symmetrization* is a geometric operation invented by Jakob Steiner[1] (see [Ste38] for the original reference). The original idea of Steiner symmetrization, as presented in [PS51] is purely geometrical: Considering a body $B$ and a plane construct another body $B^*$ such that:

- it is symmetrical with respect to the plane and
- for every line perpendicular to the plane, the intersections between it and the bodies $B$ and $B^*$ have the same lengths.

This process is shown to not increase the surface area, and to maintain the volume unchanged (this is simply a consequence of Fubini-Tonelli's theorem). By taking different planes we can deduce that, for fixed volume, a convex domain is a minimizer of surface area.

Later H. Schwarz[2] applied a similar method, but in which symmetrization was taken with respect to a line. As Steiner did with is symmetrization, Schwarz proved that Schwarz symmetrization leaves the volumen unchanged but diminishes (in the sense that it never increases) the surface area. In particular, Schwarz rearrangement can be obtained as a limit of Steiner symmetrizations. This was done for convex bodies in [PS51, p. 190], [Lei80, p. 226], and for non convex bodies in [BLL74].

For some reason, the original definition of Schwarz symmetrization was diffused in the literature, as noted by Kawohl in [Kaw85, p. 16]:

> Polya and Szegö distinguish between Schwarz and point symmetrization. Their definition of "symmetrization of a set with respect to a point" coincides with our [in his book] definition of Schwarz symmetrization and is commonly refered to as Schwarz symmetrization [Ban80b; Lio80; Mos84].

We will present later the definition as it commonly used nowadays. Let us start by saying a few words about isoperimetric inequalities.

[1] Jakob Steiner (18 March 1796 – 1 April 1863) was a Swiss mathematician who worked primarily in geometry.

[2] Karl Hermann Amandus Schwarz (25 January 1843 – 30 November 1921) was a German mathematician, known for his work in complex analysis. Do not confuse with Laurent-Moïse Schwartz (5 March 1915 – 4 July 2002), a French mathematician. The later pioneered the theory of distributions, which gives a well-defined meaning to objects such as the Dirac delta function.

## 2.2 Isoperimetric inequalities

As mentioned in the previous section, one of the most classical result obtained via symmetrization techniques is the isoperimetric inequality:

> From all $n$-dimensional bodies of a given volume, the $n$-ball is the one of least surface.

In the plane the solution was believed to be the circle from the time of Kepler. However the first succesful attempt towards proving this result mathematically in dimension 2 was made by Steiner in 1838 (see [Ste38]). This isoperimetic inequality in the plane can be written as follows. Let $\Omega$ be a smooth domain in $\mathbb{R}^2$, let $A$ be its surface area and $L$ its perimeter. Then

$$4\pi A \leq L^2. \tag{2.5}$$

Of course, equality holds for the circle. The surprising fact is that holds *only* when the domain is a circle.

Since this initial result there has been substancial research in this direction. For example, Hurwitz in 1902 applied Fourier series (see [Hur01]) and, in 1938, E. Schmidt made a proof using the arc length formula, Green's theorem and the Cauchy-Schwarz inequality (see [Sch39]).

A generalization of the isoperimetric inequalities is already well known. It can be written in the following terms:

**Theorem 2.1** (Federer, 1969 [Fed69]). *Let $S \subset \mathbb{R}^n$ be such that $\overline{S}$ has finite Lebesgue measure. Then*

$$n\omega_n^{\frac{1}{n}} L^n(\overline{S})^{\frac{n-1}{n}} \leq M_*^{n-1}(\partial S) \tag{2.6}$$

*where $\omega_n$ is the volume of the n-ball, $M_*^n$ is the Minkowski content and $L^n$ is the Lebesgue measure.*

In [Tal16; Ban80a; Rak08; BK06] the reader can find a survey on the study isoperimetric inequalities.

This geometrical inequality is equivalent to a result which is of interest to the specialist in Partial Differential Equations: the Sobolev inequality

$$n\omega_n^{\frac{1}{n}} \left( \int_{\mathbb{R}^n} |u|^{\frac{n}{n-1}} \right)^{\frac{n-1}{n}} \leq \int_{\mathbb{R}^n} |\nabla u|, \qquad \forall u \in W_c^{1,1}(\mathbb{R}^n). \tag{2.7}$$

## 2.3 From a geometrical viewpoint to rearrangement of functions

Once a particular type of rearrangement $\Omega_*$ of a set $\Omega$ is understood there is a natural way to define the rearrangement of a function $u : \Omega \to \mathbb{R}$. Consider the level sets:

$$\Omega_c = \{x \in \Omega : u(x) \geq c\}, \qquad c \in \mathbb{R}. \tag{2.8}$$

One can define the rearrangement $u_*$ of $u$ as:

$$u_* : \Omega_* \to \mathbb{R}, \qquad u_*(x) = \sup\{c \in \mathbb{R} : x \in (\Omega_c)_*\}. \tag{2.9}$$

Over the years, several other different types of rearrangements have been developed, an applied with success to different types of problems, with particularly good results in geometry and function theory, specially in PDEs. A catalogue of this techniques can be found in [Kaw85] (although there are many others excelent references, e.g., [Rak08; Ban80b]).

As noted in [PS51] both Steiner and Schwarz symmetrization reduce the Dirichlet integral of functions vanishing in the boundary. Informally

$$\int_{\Omega_*} |\nabla u_*|^2 \leq \int_{\Omega} |\nabla u|^2, \qquad \text{if } u = 0 \text{ on } \partial\Omega \tag{2.10}$$

This immediately appeals to the imagination of the PDE specialist. What seemed as a purely geometrical tool becomes a functional one.

## 2.4 The coarea formula

Symmetrization is the art of understanding the level set. The following result, known as the coarea formula, allows us to make consider level sets as domain of integration in Fubini-Tonelli theorem-like change of variable. For smooth functions it follows directly as a change of variables. A more general form it was stated by Federer in [Fed59] for Lipschitz functions and for bounded variation functions by Fleming and Rishel in [FR60]. We present the result as it appears in Federer's book [Fed69].

**Theorem 2.2.** *Let u be a Lipschitz function. Then for all $g \in L^1(\Omega)$ then*

$$\int_\Omega g(x)|\nabla u(x)|dx = \int_{-\infty}^{+\infty} \left( \int_{u^{-1}(t)} g(x) dH_{n-1}(x) \right) dt \tag{2.11}$$

*where $H_n$ is n-dimensional Hausdorff measure.*

The usual formulation is the particular case $g \equiv 1$. This formula, jointly with the isoperimetric inequality, gives a proof of the Sobolev inequality for $W^{1,1}(\mathbb{R}^n)$ given by (2.7).

## 2.5 Schwarz rearrangement

### 2.5.1 Decreasing rearrangement

For the purpose of this thesis we will focus mainly on two types of rearrangements: Schwarz and Steiner rearrangements. In particular we will be interested in studying this rearrangement as a tool in studying the Laplace operator, and other operators in divergence form as a first step for the consideration of problem (2.1). First, we introduce the (modern) definition of Schwarz symmetrization

**Definition 2.1.** Let $\Omega \subset \mathbb{R}^n$. We define the Schwarz rearrangement of $\Omega$ as

$$\Omega^\star = B(0,R), \qquad \text{such that } |\Omega^\star| = |\Omega|. \tag{2.12}$$

where $B(0,R)$ as ball centered at 0 of radius $R$.

The process of symmetrization for this kind of problems was introduced by Faber [Fab23] and Krahn [Kra25; Kra26] in their proof of the Rayleigh's conjecture, which can be stated in the following terms

**Theorem 2.3** (Rayleigh-Faber-Krahn). *Let*

$$\lambda(\Omega) = \min_{u \in H_0^1(\Omega)} \frac{\int_\Omega |\nabla u|^2}{\int_\Omega u^2}. \tag{2.13}$$

*Then*

$$\lambda(\Omega) \geq \lambda(\Omega^\star). \tag{2.14}$$

In modern terms $\lambda(\Omega)$ is, of course, known as the first eigenvalue for the Laplace operator, and $\lambda(\Omega)$ can be though as the smallest real number such that

$$\begin{cases} -\Delta u = \lambda u, & \Omega, \\ u = 0, & \partial\Omega, \end{cases} \tag{2.15}$$

has a nontrivial solution. The nontrivial solutions of this problem are known as eigenfunction. They will be used extensively in Part II.

**Definition 2.2.** Let $u : \Omega \to \mathbb{R}$ be a measurable function. We define the *distribution function of* $u$, $\mu : [0,+\infty) \to [0,+\infty)$, as

$$\mu(t) = |\{x \in \Omega : |u(x)| > t\}| \tag{2.16}$$

and the *decreasing rearrangement of* $u$, $u^* : [0,+\infty) \to \mathbb{R}$ as

$$u^*(s) = \sup\{t \geq 0 : \mu(t) > s\}. \tag{2.17}$$

**Definition 2.3.** We introduce the Schwarz rearrangement $u^\star$ of $u$ as

$$u^\star(x) = u^*(\omega_n |x|^n) \tag{2.18}$$

where $\omega_n$ represents the volume of the $n$-dimensional unit ball.

### 2.5.2 The three big inequalities and one big equation

There are several inequalities involving these rearrangements which will be of a great importance for us.

The Hardy-Littlewood-Polya inequality is very import since it allows to bound products in $L^2$. It can be stated as follows

**Theorem 2.4** (Hardy-Littlewood-Polya, 1929 [HLP29])**.** *Let $\Omega$ be a measurable subset of $\mathbb{R}^n$ and $f, g$ be non negative measurable functions. Then*

$$\int_\Omega fg \leq \int_0^{|G|} f^* g^*. \tag{2.19}$$

The second remarkable inequality is Riesz's inequality . It is very useful in order to make a priori comparisons with Green's kernel

**Theorem 2.5** (Riesz, 1930 [Rie30]). *Let $\Omega$ be a measurable subset of $\mathbb{R}^n$ and $f, g$ be non negative measurable functions. Then*

$$\int_{\mathbb{R}^n\times\mathbb{R}^n} f(x)g(x-y)h(y)dxdy \leq \int_{\mathbb{R}^n\times\mathbb{R}^n} f^\star(x)g^\star(x-y)h^\star(y)dxdy. \tag{2.20}$$

Both results, and further techniques, were compiled in one of the references text in the subject [HLP52].

Only a few years after the first appearance of these two results, in 1945, Polya and Szegö publish [PS45], were they introduce the following inequality, to prove that the capacity of a condenser diminishes or remains unchanged by applying the process of Schwarz symmetrization.

**Theorem 2.6** (Polya-Szego [PS45]). *Let $u : \mathbb{R}^n \to \mathbb{R}^+$ in $W^{1,p}(\mathbb{R}^n)$ where $1 \leq p < \infty$. Then*

$$\int_{\mathbb{R}^n} |\nabla u^\star|^p \leq \int_{\mathbb{R}^n} |\nabla u|^p. \tag{2.21}$$

This result was also useful in the proof of the Choquard conjecture by Lieb [Lie77].

The collaboration between Polya and Szegö continued in time, and they updated [PS51] over several editions. This is a reference text in isoperimetric inequalities and the use of different rearrangements in Mathematical Physics.

The final inequality could be one the most important in the theory, because it is used to convert the original PDE for $u$ to a PDE for $\mu$.

**Theorem 2.7** ([BZ87]). *Let $u \in W^{1,p}$ for some $1 \leq p < \infty$. Then the following holds:*

*i)* $\mu$ *is one-to-one.*

*ii)* $u^* \circ \left(\frac{\mu(t)}{\omega_n}\right)^{\frac{1}{n}} = Id.$

*iii) We can decompose* $\mu$ *as*

$$\mu(t) = |\{x \in \Omega : |\nabla u| = 0, u(x) > t\}| + \int_t^{+\infty} \int_{u^{-1}(t)} |\nabla u|^{-1} dH^{n-1} ds. \tag{2.22}$$

*iv) For almost all t,*

$$-\infty < \frac{d\mu}{dt} \leq \int_{\{x\in\Omega : u(x)=t\}} \frac{-1}{|\nabla u|} dH^{n-1}, \qquad a.e.\ t \in \mathbb{R}. \tag{2.23}$$

*Equality holds in the previous item if* $|\nabla u| \neq 0$.

*v) For almost all t,* $\frac{d\mu}{dt} < 0$.

### 2.5.3 Concentration and rearrangement

Even though the stronger results show that we have a pointwise comparison $u^\star \leq v$ this is not the case in general. However, there is a property much nicer in terms of rearrangements: the *concentration*. This term, which appears frequently in the mathematical literature, must not be confused with the chemical concept of concentration. As will see from the following definition they are entirely different.

**Definition 2.4.** Let $\Omega$ be an open set of $\mathbb{R}^n$ and let $\psi \in L^1(\Omega_1), \phi \in L^1(\Omega_2), |\Omega_1| = |\Omega_2|$. We say that the concentration of $\phi$ is less or equal than the concentration of $\psi$, and we denote this by $\phi \preceq \psi$ if

$$\int_0^t \phi^*(s)ds \leq \int_0^t \psi^*(s)ds, \qquad \forall t \in [0, |\Omega|]. \tag{2.24}$$

Equivalently,

$$\int_{B_r(0)} \phi^\star(x)dx \leq \int_{B_r(0)} \psi^\star(x)dx. \tag{2.25}$$

The following lemma is a very important result. It allows to understand the importance of convex functions in symmetrization.

**Lemma 2.5.1.** *Let* $y, z \in L^1(0, M)$*, y and z nonnegative. Suppose that y is nonincreasing and*

$$\int_0^t y(s)ds \leq \int_0^t z(s)ds, \qquad \forall t \in [0, M]. \tag{2.26}$$

*Then, for every continuous non decreasing convex function* $\Phi$ *we have*

$$\int_0^t \Phi(y(s))ds \leq \int_0^t \Phi(z(s))ds, \qquad \forall t \in [0, M]. \tag{2.27}$$

Applying this result it is possible to obtain the following properties (see [HLP52; HLP29; RF88; CR71; ATL89]).

**Proposition 2.1.** *Let $\Omega_1, \Omega_2$ be Borel sets in $\mathbb{R}^{N_1}, \mathbb{R}^{N_2}$ respectively such that $|\Omega_1| = |\Omega_2|$ and $\phi_i \in L^1_+(\Omega_i)$ (i.e. it is is in $L^1$ and non-negative). Then, the following are equivalent:*

*i)* $\phi_1 \preceq \phi_2$

*ii)* $F(\phi_1) \preceq F(\phi_2)$ *for every nondecreasing, convex function* $F$ *on* $[0,+\infty)$ *such that* $F(0) = 0$.

*iii) For all* $\psi \in L^1 \cap L^\infty(\Omega_1)$

$$\int_{\Omega_1} \varphi_1 \psi \le \int_0^{|\Omega_2|} \varphi_2^* \psi^* = \int_{\Omega_2^\star} \varphi_2^\star \psi^\star \tag{2.28}$$

*iv) For all* $\psi$ *nonincreasing on* $(0, |\Omega_1|)$, $\psi \in L^1 \cap L^\infty(0, |\Omega_1|)$

$$\int_0^{|\Omega_1|} \varphi_1^* \psi \le \int_0^{|\Omega_2|} \varphi_2^* \psi. \tag{2.29}$$

### 2.5.4 Schwarz symmetrization of elliptic problems

It is precisely in 1962, in a book in honor of Polya edited by Szegö [Wei62] that Weinberger extends the results by Faber and Krahn to obtain isoperimetric results for the Dirichlet problem with general elliptic self adjoint operator

$$L = \sum_{i,j=1}^N \partial_{x_i}\left(a_{ij}(x)\partial_{x_j}\right). \tag{2.30}$$

In 1976 Bandle [Ban76a] gives a pointwise estimates of the decreasing rearrangements of the solution of $-\Delta u + \alpha u + 1$ with Dirichlet boundary condition. In 1978, in [Ban78], she gives estimates on the Green kernel. In the same year Alvino and Trombetti [AT78] present result similar to Weinberger's for degenerate (non elliptic) equations.

In 1979 Talenti [Tal79] apply Schwarz symmetrization techniques, to improve upon the results of Weinberger and Bandle. He focuses on non linear elliptic equations of the form

$$\begin{cases} -\operatorname{div}(a(x,u,\nabla u)) + g(x,u) = 0 & \Omega, \\ u = h & \partial\Omega, \end{cases} \tag{2.31}$$

under the hypothesis

i) There exists $A : [0,+\infty) \to \mathbb{R}$ convex such that :

i) $a(x,u,\xi)\cdot\xi \geq A(|\xi|)$

ii) $A(r)/r \to 0$ as $r \to 0$

ii) $g$ is measurable and

$$(g(x,u)-g(x,0))u \geq 0 \tag{2.32}$$

iii) $h \in L^\infty(\Omega)$ and

$$\int_\Omega A(|\nabla h|) < \infty. \tag{2.33}$$

In the case $h \equiv b$ a constant Talenti introduces the "rearranged problem"

$$\begin{cases} -\operatorname{div}(\frac{A(|\nabla v|)}{|\nabla v|^2}\nabla v) = f^\star & \Omega^\star, \\ v = b & \partial\Omega^\star. \end{cases} \tag{2.34}$$

He manages to prove, for the first time in literature, that we can compare pointwise $u^\star$ with a different solution (which is easier to compute analytically). In fact

$$u^\star \leq v, \qquad a.e.\ \Omega^\star \tag{2.35}$$

Notice that, if $a(x,u,\xi) = \xi$, the operator is the usual Laplace operator and $A(\xi) = \xi^2$.

In 1980 P.L. Lions [Lio80] provides a simpler proof of this result in the linear case with $h \equiv 0$ and extends the estimates to operator in the form $A = Lu + c$ were $L$ is second order elliptic and $c : \Omega \to \mathbb{R}$. He compares the problems

$$\begin{cases} -\Delta u + cu = f & \Omega, \\ u = 0 & \partial\Omega, \end{cases} \quad \text{and} \quad \begin{cases} -\Delta v + (c^+)^{\star\star}v - (c^-)^\star v = f^\star & \Omega^\star \\ v = 0 & \partial\Omega^\star \end{cases} \tag{2.36}$$

where $\phi^{\star\star}(x) = \phi^{**}(\omega_n|x|)$ and $\phi^{**}(s) = \phi^*(|\Omega| - s)$. This last function is known as the increasing rearrangement of $\phi$. Analogously to the Schwarz rearrangement $\phi^{\star\star}$ represent represents the unique radial increasing function with the same distribution function as $|\phi|$ (it can be defined by analogy to the Schwarz rearrangement).

Lions shows that,

$$\int_\Omega F(|u|) \leq \int_{\Omega^\star} F(v) \tag{2.37}$$

for all function $F : \mathbb{R}_+ \to \mathbb{R}_+$ positive, increasing and convex. Besides, if $c \le 0$, we have that

$$u^\star \le v, \qquad a.e.\ \Omega^\star. \tag{2.38}$$

At the end of this paper Lions indicates that some nonlinear cases are immediately covered by the linear case, by simply taking a sequence of solutions $-\Delta u_{n+1} = g(u_n)$ and applying the comparison, since $g(u_n)^\star = g(u_n^\star)$. Around the same time, in 1979, Chiti in [CP79] (see [Chi79] in Italian) proved a similar result, by using a limit of simple functions. However, his result was presented as a Orlicz norm result, rather than as a PDE result.

The ideas behind both [Tal79] and [Lio80] is explained very elegantly in 1990 by Talenti in [Tal90]. By applying the Hardy-Littlewood-Polya inequality, the isoperimetric inequality in $\mathbb{R}^n$ and the coarea formula one the boundary value problem can be rewritten in terms of a ODE for the distribution function

$$n\omega_n^{\frac{1}{n}}\mu(t)^{1-\frac{1}{n}}B\left(n\omega_n^{\frac{1}{n}}\frac{\mu(t)^{1-\frac{1}{n}}}{\mu'(t)}\right) \le \int_0^{\mu(t)} g(s)ds. \tag{2.39}$$

Returning to a chronological order, in 1984, in the more general context of Mathematical Physics, Mossino publishes a book [Mos84] which contains a number of interesting statements on elliptic problems. However, none of the results are necessary to the interest of this Thesis.

In 1985 Díaz in [Día85] polishes some of the previously presented rearrangement techniques, in order to obtain estimates for free boundaries that rise in problem (2.40), when $g$ is not a Lipschitz function. The results are extended to the p-Laplace operator and the regularity hypothesis are simplified. The following theorem is stated:

**Theorem 2.8.** *Let $u, v$ be the solutions of problems*

$$\begin{cases} -\operatorname{div}(A(x,u,\nabla u)) + g(u) = f_1 & \Omega, \\ u = 0 & \partial\Omega, \end{cases} \tag{2.40}$$

*and*

$$\begin{cases} -\Delta_p v + g(v) = f_2 & \Omega, \\ v = 0 & \partial\Omega, \end{cases} \tag{2.41}$$

*where*

*i)* *A is a Caratheodory function such that* $A(x,u,\xi) \geq |\xi|^p$

*ii)* *g is continuous non decreasing such that* $g(0) = 0$

*iii)* $f_1 \in L^{p'}(\Omega)$ *such that* $f_1 \geq 0$

*iv)* $f_2 \in L^{p'}(\Omega^\star)$ *such that* $f_2 = f_2^\star$

*v)* $\int_0^t f_1^* \leq \int_0^t f_2^*$ *for all* $t \in [0, |\Omega|]$

*Then*

$$\int_0^t g(u^*) \leq \int_0^t g(v^*), \qquad t \in [0, |\Omega|]. \tag{2.42}$$

*In particular, for any convex nondecreasing real function* $\Phi$

$$\int_\Omega \Phi(g(u)) \leq \int_{\Omega^\star} \Phi(g(v)). \tag{2.43}$$

Notice that, in this more general setting, the result is not as strong as in the case withouth the absortion. We do not have a pointwise comparison of $u^\star$ and $v$. Behind the proof is the following lemma

**Lemma 2.5.2.** *Let u be the solution of* (2.40)*. Then* $u^*$ *is absolutely continuous in* $[0, |\Omega|]$ *and*

$$-\frac{du^*}{ds}(s) \leq \left( \frac{1}{n\omega_n^{\frac{1}{n}} s^{\frac{n-1}{n}}} \right)^{\frac{p}{p-1}} \left[ \int_0^s f_1^*(\theta)d\theta - \int_0^s g(u^*(\theta))d\theta \right]. \tag{2.44}$$

*Let v be the solution of* (2.41)*, then equality holds if* $g_2$ *is radial.*

In the proof of this lemma there are two main ingredients. The first is a general statement on the distribution function.

**Lemma 2.5.3.** *Let* $z \in W_0^{1,p}(\Omega)$*,* $z \geq 0$*. Then if* $\mu(t) = |\{x \in \Omega : z(x) > t\}|$ *one has*

$$n\omega_n^{\frac{1}{n}} \mu(t)^{\frac{n-1}{n}} \leq \left( -\frac{d\mu}{dt}(t) \right)^{\frac{1}{p}} \left( -\frac{d}{dt} \int_{\{z(x)>t\}} |\nabla z(x)|^p dx \right)^{\frac{1}{p}}. \tag{2.45}$$

The second one is a particular computation on an integral of the solution of problem (2.40).

**Lemma 2.5.4.** *Let $u \in W_0^{1,p}(\Omega)$ be a nonnegative solution of* (2.40). *Then the function*

$$\Psi(t) = \int_{\{u(x)>t\}} A(x,u,\nabla u)\cdot \nabla u dx \tag{2.46}$$

*is a decreasing Lipschitz continuous function of $t \in [0,+\infty)$, and the inequality*

$$0 \le -\frac{d\Psi}{dt}(t) \le \int_0^{\mu(t)} f_1^*(s)ds - \int_0^{\mu(t)} g(u^*(s))ds. \tag{2.47}$$

The last assertion of the theorem is a consequence of the Lemma 2.5.1.

Finally in 1990 Alvino, Trombetti and Lions, in [ATL90], prove the result for a general operator in the form

$$Lu = -\operatorname{div}(A(x)\nabla u) + \nabla(b(x)u) + d(x)\cdot\nabla u + c(x)u \tag{2.48}$$

under weak restriction on the operators. For completeness, and under the general operator $-\operatorname{div}(a(x,u,\nabla u))$, jointly with Ferone in [AFTL97], they extend the Polya-Szegö inequality and some of the previous lemmas to a more general context. As a corollary of their analysis they manage to obtain an comparison of the form $u^\star \le C_n v$.

### 2.5.5 Schwarz rearrangement of parabolic problems

The application of rearrangement techniques to parabolic equations was considered for the first time in 1976 by C. Bandle in [Ban76b]. This result was announced in 1975 in a *Comptes Rendus* note (see [Ban75]). The compared problems in this case are

$$\begin{cases} \dfrac{\partial u}{\partial t} - \Delta u = f(x) & \Omega\times(0,\infty),\\ u = 0 & \partial\Omega\times[0,+\infty),\\ u(\cdot,0)=u_0 & \overline{\Omega} \end{cases}, \quad \text{and} \quad \begin{cases} \dfrac{\partial v}{\partial t} - \Delta v = f^\star(x), & \Omega^\star\times(0,\infty),\\ v = 0 & \partial\Omega^\star\times[0,+\infty),\\ v(\cdot,0)=v_0^\star & \overline{\Omega^\star}, \end{cases} \tag{2.49}$$

where $\Omega$ has to be piecewise analytic. As described in the text:

> The proof is based on a differential inequality and uses very much a system of curvilinear coordinates defined by the level surfaces of $u(x,t)$. The introduction of those coordinates requires a strong assumption on the regularity of $u(x,t)$.

Let $\mu, \tilde{\mu}$ be the distribution functions of $u$ and $v$ respectively. By defining

$$H(a,t) = \int_0^t \mu(t,s)ds, \qquad \tilde{H}(a,t) = \int_0^a \tilde{\mu}(t,s)ds. \tag{2.50}$$

Then

$$\begin{cases} -\dfrac{\partial H}{\partial t} + p(a)\dfrac{\partial^2 H}{\partial a^2} + \displaystyle\int_0^{\cdot} f^*(s)ds \geq 0 \\ H(0,t) = 0, \\ \dfrac{\partial H}{\partial a}(0,t) = \max\limits_{\overline{\Omega}} u(x,t), \\ \dfrac{\partial H}{\partial a}(|\Omega|,t) = 0, \end{cases} \qquad \begin{cases} -\dfrac{\partial \tilde{H}}{\partial t} + p(a)\dfrac{\partial^2 \tilde{H}}{\partial a^2} + \displaystyle\int_0^{\cdot} f^*(s)ds = 0 \\ \tilde{H}(0,t) = 0, \\ \dfrac{\partial \tilde{H}}{\partial a}(0,t) = \max\limits_{x\in\overline{\Omega^\star}} v(x,t), \\ \dfrac{\partial \tilde{H}}{\partial a}(|\Omega|,t) = 0. \end{cases} \tag{2.51}$$

From this we can conclude that $u \preceq v$.

After this initial result, that worked only for smooth classical solutions some generalizations appeared. In 1982 J. L. Vazquez showed similar results for the porous medium equation: $u_t - \Delta\varphi(u) = f$ (see [Váz82]). Mossino and Rakotoson in 1986 (see [MR86]) obtained a similar result under weaker regularity by considering the directional derivative of the rearrangement, a technique that first appeared in [MT81].

A generalization comes in 1992, by Alvino, Trombetti and P.L. Lions in [ALT92], by applying the techniques in [MR86], simple properties of the fundamental solution and semigroup theory (in particular the Trotter formula which will be introduced in Section 2.7.3), which allows to reduce the regularity condition on the data and the solution. They obtain the expect comparison between the problem:

$$\begin{cases} \frac{\partial u}{\partial t} - \operatorname{div}(A(x,t)\nabla u) + a(x,t)u = f(x,t) & \Omega \times (0,T), \\ u = 0 & \partial\Omega \times (0,T), \\ u = u_0 & \Omega \times \{0\}, \end{cases} \tag{2.52}$$

under the assumptions

$$\begin{cases} a_{ij}, a \in L^\infty(\Omega \times (0,T)), f \in L^2(\Omega \times (0,T)), u_0 \in L^2(\Omega), \\ \exists \nu \in L^\infty(0,T), \alpha > 0, \forall \xi \in \mathbb{R}^n : \qquad \xi^t A(x,t)\xi \geq \nu(t)|\xi|^2 \geq \alpha|\xi|^2, \end{cases} \tag{2.53}$$

and the problem

$$\begin{cases} \frac{\partial v}{\partial t} - \nu(t)\Delta v + \{a_1 - a_2\}v = g(x,t) & \Omega^\star \times (0,T), \\ v = 0 & \partial\Omega^\star \times (0,T), \\ v = u_0^\star & \Omega^\star \times \{0\}. \end{cases} \tag{2.54}$$

under the assumptions

$$\begin{cases} a_i \in L^\infty_+(\Omega^\star \times (0,T)), v_0 \in L^2_+(\Omega^\star), g \in L^2_+(\Omega^\star \times (0,T)), \\ a_2, -a_1, g \text{ are spherically symmetric, nonincreasing} \\ \text{ with respect to } |x|, \text{ for almost all } t \in (0,T), \end{cases} \tag{2.55}$$

As stated in [ALT92] the result reads:

**Theorem 2.9.** *Let $u$ be the solution of* (2.52) *under hypothesis* (2.53) *and $v$ the solution* (2.54) *under hypothesis* (2.55). *Assume further that,*

$$\begin{cases} u_0 \preceq v_0, \quad a^-(t) \preceq a_2(t), \quad f(t) \prec q(t), \\ \int_0^r (a_1)^{**} \leq \int_0^r (a^{**})_+, \textit{ for all } r \in [0, |\Omega|], \textit{ for a.e. } t \in (0,T) \end{cases} \tag{2.56}$$

*Then, for all $t \in [0,T]$,*

$$u(t) \preceq v(t) \tag{2.57}$$

*(in the sense of Definition 2.4).*

## 2.6 A differentiation formula

Most results in PDEs using symmetrizaqtion techniques pass by the consideration of differentiation formulas of the following function

$$H(s,y) = \int_{\{u(x,y)>t\}} f(x,y)dx. \tag{2.58}$$

In 1998 Ferone and Mercaldo in [FM98] state a second order differentiation formula for rearrangements (citing works in Steiner rearrangement, [ATDL96], which we will present in the following section)

**Theorem 2.10.** *Let $\Omega = \Omega' \times (0,h)$ u be a nonnegative function in $W^{2,p}(\Omega)$, were $p > n+1$ and let $f$ be Lipschitz in $\overline{\Omega}$. Assume that*

$$|\{x \in \Omega : |\nabla_x u| = 0, u(x,y) \in (0, \sup u(\cdot,y))\}| = 0, \qquad \forall y \in (0,h). \tag{2.59}$$

*Then we have*

*i) For any $y \in (0,h)$, $H$ is differentiable with respect to $s$ for a.e. $s \geq 0$ and*

$$\frac{\partial H}{\partial t}(t,y) = - \int\limits_{u(x,y)=t} \frac{f(x,y)}{|\nabla_x u|} dH^{n-1}(x). \tag{2.60}$$

*ii) For every fixed $s$, $H$ is differentiable with respect to $y$ and*

$$\frac{\partial H}{\partial y}(s,y) = \int\limits_{u(x,y)>t} \frac{\partial f}{\partial y}(x,y)dx + \int\limits_{u(x,y)=t} \frac{\partial u}{\partial y}(x,y)\frac{f(x,y)}{|\nabla_x u|} dH^{n-1}(x). \tag{2.61}$$

From this we can extract the following corollary, which is of capital importance in symmetrization.

**Corollary 2.1.** *Let $u \in W^{1,\infty}(\Omega \times (0,h))$ be nonnegative. Then*

$$\int\limits_{u(x,y)>u^*(s,y)} \frac{\partial^2 u}{\partial y^2}(x,y)dx = \frac{\partial^2}{\partial y^2}\int_0^s u^*(\sigma)d\sigma - \int_{u(x,y)=u^*(s,y)} \frac{\left(\frac{\partial u}{\partial y}\right)^2}{|\nabla_x u|} dH^{n-1}(x) \tag{2.62}$$

$$+ \left( \int\limits_{u(x,y)=u^*(s,y)} \frac{\frac{\partial u}{\partial y}(x,y)}{|\nabla_x u|} dH^{n-1}(x) \right)^2 \left( \int\limits_{u(x,y)=u^*(s,y)} \frac{1}{|\nabla_x u|} dH^{n-1}(x) \right)^{-1}. \tag{2.63}$$

## 2.7 Steiner rearrangement

The idea of the Schwarz rearrangement (in the modern definition) is to consider radially decreasing functions. A smart analysis of the pros and cons of performing this symmetrization is presented in [ATDL96]:

> On the one hand, these results [on the Schwarz symmetrization] make the problem of determining a priori estimates easier by turning it into a one-dimensional problem; on the other hand, by this symmetrization process, the differential problem may lose properties that arise from the symmetry of the data with respect to a group of variables. In order to preserve this kind of symmetry, it is usefull to check whether comparison results hold when a partial symmetrization such as Steiner symmetrization is used.

In 1992 Alvino, Díaz, Lions and Trombetti introduce a new definition of Steiner symmetrization, which differs slightly from the one in [PS51]. We will follow this new definition. We point out that the new definition (which we give in precise terms below) can be obtained, as it was the case in of Schwartz symmetrization, as a limit of Steiner symmetrization perpendicular to a hyperplane, in the sense presented in [PS51].

The idea behind this (new) Steiner rearrangement is to symmetrize radially only in some variables, and therefore only works in product domains $\Omega = \Omega' \times \Omega''$.

**Definition 2.5.** Let

$$\Omega = \Omega' \times \Omega'' \subset \mathbb{R}^{n_1} \times \mathbb{R}^{n_2}. \tag{2.64}$$

We usually refer to the variables in $\Omega'$ as $x$, and to the variables in $\Omega''$ as $y$. We define the Steiner rearrangement of $\Omega$ with respect to the variables $x$ as

$$\Omega^{\#} = B(0,R) \times \Omega'' \quad \text{where } |B(0,R)| = |\Omega'|, \tag{2.65}$$

where $B(0,R) \subset \mathbb{R}^{n_1}$ is the ball centered at 0 of radius $R$.

**Remark 2.1.** Notice that

$$\Omega^{\#} = (\Omega')^{\star} \times \Omega''. \tag{2.66}$$

We can define the functional rearrangement as follows:

**Definition 2.6.** Let $\Omega = \Omega' \times \Omega'' \subset \mathbb{R}^{n_1} \times \mathbb{R}^{n_2}$, $u : \Omega' \to \mathbb{R}$ be a measurable function. We define the *distribution function of* $u$, $\mu : [0,+\infty) \times \Omega'' \to [0,+\infty)$, as

$$\mu(t,y) = |\{x \in \Omega : |u(x,y)| > t\}|, \tag{2.67}$$

and the *decreasing rearrangement of* $u$, $u^* : [0,+\infty) \times \Omega'' \to \mathbb{R}$ as

$$u^*(s,y) = \sup\{t \geq 0 : \mu(t,y) > s\}. \tag{2.68}$$

Finally we introduce the *Steiner rearrangement of* $u$ as

$$u^{\#}(x,y) = u^*(\omega_{n_1}|x|^{n_1}, y), \qquad (x,y) \in \Omega^{\#}, \tag{2.69}$$

where $\omega_n$ represents the volume of the $n$-dimensional ball.

**Remark 2.2.** Notice that

$$u^{\#}(x,y) = (u(\cdot,y))^{\star}(x), \qquad (x,y) \in (\Omega')^{\star} \times \Omega'' = \Omega^{\#}. \tag{2.70}$$

Naturally, the Steiner rearrangement shares, for every $y$, the same properties as the Schwarz rearrangement.

### 2.7.1 Steiner rearrangement of elliptic equations

As announced in [ADLT92], in [ATDL96] Alvino, Trombetti, Díaz and P.-L. Lions prove the following result

**Theorem 2.11.** *Let*

$$\begin{aligned} Lu = &- \sum_{i,j=1}^{n} \frac{\partial}{\partial x_j}\left(a_{ij}(x,y)\frac{\partial u}{\partial x_i}\right) - \sum_{h,k=1}^{m} \frac{\partial}{\partial y_k}\left(b_{hk}(x,y)\frac{\partial u}{\partial y_h}\right) \\ &- \sum_{i=1}^{n}\sum_{h=1}^{m} \frac{\partial}{\partial y_h}\left(c_{ih}(y)\frac{\partial u}{\partial x_i}\right) - \sum_{i=1}^{n}\sum_{h=1}^{m} \frac{\partial}{\partial x_i}\left(d_{hi}(y)\frac{\partial u}{\partial y_h}\right) \end{aligned} \tag{2.71}$$

*and let $u$ be a weak solution of the Dirichlet problem*

$$\begin{cases} Lu = f & \Omega, \\ u = 0 & \partial\Omega. \end{cases} \tag{2.72}$$

*We assume the following:*

*i) Coefficients $a_{ij}, b_{hk}, c_{ih}, d_{hl}$ and $f$ belong to $L^{\infty}(\Omega)$,*

*ii) (ellipticity condition) there exists $\nu > 0$ such that, for every $(\xi,\eta) \in \mathbb{R}^n \times \mathbb{R}^m$ and a.e. $(x,y) \in \Omega$*

$$\begin{aligned} &\sum_{i,j=1}^{n} a_{ij}(x,y)\xi_i\xi_j + \sum_{h,k=1}^{m} b_{hk}(y)\eta_h\eta_k \\ &\qquad + \sum_{i=1}^{n}\sum_{h=1}^{m} c_{ih}(y)\xi_i\eta_h + \sum_{i=1}^{n}\sum_{h=1}^{m} d_{hi}(y)\xi_i\eta_h \geq |\xi|^2 + \nu|\eta|^2, \end{aligned} \tag{2.73}$$

*iii) $\Omega = \Omega' \times \Omega''$ open, bounded subset of $\mathbb{R}^n$ and $\mathbb{R}^m$, respectively.*

*Let v be the weak solution of the problem*

$$\begin{cases} -\Delta_x v - \sum_{h,k=1}^m \frac{\partial}{\partial y_k}\left(a_{hk}(x,y)\frac{\partial u}{\partial y_h}\right) = f & \Omega^\#, \\ v = 0 & \partial\Omega^\#. \end{cases} \tag{2.74}$$

*Then we have, for any* $y \in \mathbb{R}^m$

$$\int_0^s u^*(\sigma,y)d\sigma \le \int_0^s v^*(\sigma,y)d\sigma. \tag{2.75}$$

The proof of this result is highly technical. It uses, as it was the cases in previous results in the literature, the differential geometry behind the level sets, considering specially the case of $C^1$ solutions. Besides, in the effort of taking about the most general elliptic operator, the presence of subscripts $i, j, h, k$ makes the work quite baroque.

In 2001 Chiacchio and Monetti in [CM01] (see also [Chi04]) introduce lower order terms to same equation. They deal with operators in the form:

$$Lu = -\Delta u - \sum_{i=1}^n \frac{\partial}{\partial x_i}(b_i(y)u) - \sum_{j=1}^m \frac{\partial}{\partial y_j}(\tilde{b}_j(y)u) + \sum_{i=1}^n d_i(y)\frac{\partial u}{\partial x_i} + \sum_{i=1}^m \tilde{d}_j(y)\frac{\partial u}{\partial y_j} + c(y)u \tag{2.76}$$

Later, in 2009, Chiacchio studies the eigenvalue problem (see [Chi09]).

## 2.7.2 Steiner rearrangement of linear parabolic problems

By applying [FM98] the following result can be deduced immediately. It is not written formally in any known paper, however, it is mentioned in [Chi04] and [Chi09].

**Proposition 2.2.** *Let u and v be the weak solutions of*

$$\begin{cases} \frac{\partial u}{\partial t} - \Delta u = f & \Omega \times (0,T), \\ u = 0 & \partial\Omega \times (0,T), \\ u(0) = u_0 & \Omega, \end{cases} \quad \textit{and} \quad \begin{cases} \frac{\partial v}{\partial t} - \Delta v = f^\# & \Omega^\# \times (0,T), \\ v = 0 & \partial\Omega^\# \times (0,T), \\ v(0) = v_0 & \Omega^\#, \end{cases} \tag{2.77}$$

*and let*

$$U(t,s,y) = \int_0^s u^*(t,\sigma,y)d\sigma, \qquad V(t,s,y) = \int_0^s v^*(t,\sigma,y)d\sigma. \tag{2.78}$$

*Then, there exists* $g(s) \geq 0$ *such that*

$$U_t - g(s)U_{ss} - \Delta_y U \leq \int_0^s f^*(\sigma, y)d\sigma, \qquad V_t - g(s)V_{ss} - \Delta_y V = \int_0^s f^*(\sigma, y)d\sigma. \quad (2.79)$$

Hence we deduce that

**Proposition 2.3.** *For a.e.* $y \in \Omega''$ *and a.e.* $t > 0$

$$u_0(\cdot, y) \preceq v_0(\cdot, y) \implies u(t, \cdot, y) \preceq v.(t, \cdot, y). \quad (2.80)$$

### 2.7.3 The Trotter-Kato formula

In order to treat our problem we will apply the Neveu-Trotter-Kato theorem, that characterizes the convergence of the semigroup in terms of the convergence of its generators. The abstract statement can be found in [Bré73].

**Theorem 2.12.** *Let* $(A^n)$ *and* $A$ *be maximal monotone operators such that* $\overline{D(A)} \subset \cap_n \overline{D(A^n)}$. *Let* $S_n$ *and* $S$ *be the semigroups generated by* $-A^n$ *and* $-A$ *respectively. The following properties are equivalent:*

*i) For every* $x \in \overline{D(A)}$, $S_n(\cdot)x \to S(\cdot)x$ *uniformly in compact subsets of* $[0, +\infty)$.

*ii) For every* $x \in \overline{D(A)}$ *and every* $\lambda > 0$, $(I + \lambda A^n)^{-1}x \to (I + \lambda A)^{-1}x$.

There is an important corollary to this theorem, that allows us to study the semigroup of an operator given by a sum of operators as the sequential application of the semigroup of each of these operators.

**Proposition 2.4** ([Bré73, Proposition 4.4 (p. 128)]). *Let* $A, B$ *be univoque maximal monotone operators such that* $A + B$ *is maximal monotone. Let* $S_A, S_B, S_{A+B}$ *be the semigroups associated to* $-A, -B, -(A+B)$. *Let* $C$ *be a closed convex subset of* $\overline{D(A)} \cap \overline{D(B)}$ *such that* $(I + \lambda A)^{-1}(C) \subset C$ *and* $(I + \lambda B)^{-1}C \subset C$. *Then, for every* $x \in C \cap \overline{D(A) \cap D(B)}$,

$$\left[ S_A\left(\frac{\cdot}{n}\right) S_B\left(\frac{\cdot}{n}\right) \right]^n x \to S_{A+B}(\cdot)x \quad (2.81)$$

*uniformly in every compact subset of* $[0, +\infty)$.

### 2.7.4 Steiner rearrangement of semilinear parabolic problems

In [DGC15b] and [DGC16] J.I. Díaz and myself apply Proposition 2.3, the Trotter-Kato theorem (in the form of Proposition 2.4) and explicit comparisons of the pointwise ODE

$$u_t + g(u) = f \tag{2.82}$$

to obtain the following symmetrization results. The proofs can be found in the corresponding papers collected in the Appendix.

**Theorem 2.13** ([DGC15b]). *Let $g$ be concave, verifying*

$$\int_0^\tau \frac{d\sigma}{g(\sigma)} < \infty, \qquad \forall \tau > 0. \tag{2.83}$$

*Let $h \in W^{1,\infty}(0,T)$, such that $h(t) \geq 0$ for all $t \in (0,T)$, $f \in L^2(0,T : L^2(\Omega))$ with $f \geq 0$ in $(0,T)$ and let $u_0 \in L^2(\Omega)$ be such that $u_0 \geq 0$. Let, $u \in C([0,T] : L^2(\Omega)) \cap L^2(\delta, T : H_0^1(\Omega))$ and $v \in C([0,T] : L^2(\Omega^\#)) \cap L^2(\delta, T : H_0^1(\Omega^\#))$, for any $\delta \in (0,T)$, be the unique solutions of*

$$(P) \begin{cases} \frac{\partial u}{\partial t} - \Delta u + h(t) g(u) = f(t), & \text{in } \Omega \times (0,T), \\ u = 0, & \text{on } \partial\Omega \times (0,T), \\ u(0) = u_0, & \text{on } \Omega, \end{cases}$$

$$(P^\#) \begin{cases} \frac{\partial v}{\partial t} - \Delta v + h(t) g(v) = f^\#(t), & \text{in } \Omega^\# \times (0,T), \\ v = 0, & \text{on } \partial\Omega^\# \times (0,T), \\ v(0) = v_0, & \text{on } \Omega^\#, \end{cases}$$

*where $v_0 \in L^2(\Omega^\#)$, $v_0 \geq 0$ is such that*

$$\int_0^s u_0^*(\sigma, y) d\sigma \leq \int_0^s v_0^*(\sigma, y) d\sigma, \qquad \forall s \in [0, |\Omega'|].$$

*Then, for any $t \in [0,T]$ and $s \in [0, |\Omega'|]$*

$$\int_0^s u^*(t, \sigma, y) d\sigma \leq \int_0^s v^*(t, \sigma, y) d\sigma.$$

**Theorem 2.14** ([DGC16]). *Let $\beta$ be a concave continuous nondecreasing function such that $\beta(0) = 0$. Give $T > 0$ arbitrary and let $f \in L^2(0,T : L^2(\Omega))$ with $f \geq 0$ in $(0,T)$ and let $w_0 \in L^2(\Omega)$ be such that $0 \leq w_0 \leq 1$. Let $w \in C([0,T] : L^2(\Omega)) \cap L^2(\delta, T : H_0^1(\Omega))$ and*

*$z \in C([0,T]:L^2(\Omega^\#)) \cap L^2(\delta,T:H_0^1(\Omega^\#))$, for any $\delta \in (0,T)$, be the unique solutions of*

$$(P)\begin{cases}\frac{\partial w}{\partial t} - \Delta w + \lambda \beta(w) = f(t) & \text{in } \Omega \times (0,T),\\ w = 1 & \text{on } \partial\Omega \times (0,T),\\ w(0) = w_0 & \text{on } \Omega,\end{cases}$$

$$(P^\#)\begin{cases}\frac{\partial z}{\partial t} - \Delta z + \lambda \beta(z) = f^\#(t), & \text{in } \Omega^\# \times (0,T),\\ z = 1, & \text{on } \partial\Omega^\# \times (0,T),\\ z(0) = z_0, & \text{on } \Omega^\#,\end{cases}$$

*where $z_0 \in L^2(\Omega^\#)$, $0 \le z_0 \le 1$ is such that*

$$\int_s^{|\Omega'|} z_0^*(\sigma,y)d\sigma \le \int_s^{|\Omega'|} w_0^*(\sigma,y)d\sigma, \qquad \forall s \in [0,|\Omega'|] \text{ and a.e. } y \in \Omega''.$$

*Then, for any $t \in [0,T]$, $s \in [0,|\Omega'|]$ and a.e. $y \in \Omega''$*

$$\int_s^{|\Omega'|} z^*(t,\sigma,y)d\sigma \le \int_s^{|\Omega'|} w^*(t,\sigma,y)d\sigma. \tag{2.84}$$

In terms of the comparison of the effectiveness we have the following consequence:

**Corollary 2.2.** *In the assumptions of Theorem 2.16, for any $t \in [0,+\infty)$ we have*

$$\int_{\Omega^\#} \beta(z(t,x))dx \le \int_\Omega \beta(w(t,x))dx. \tag{2.85}$$

We include now an unpublished alternative proof to one given in [DGC15b]. Let us define $S_A$ as the solution of

$$\begin{cases}u_t + Au = 0, & (0,T)\times\Omega\\ u = 0, & x \in \partial\Omega\\ u = u_0, & t = 0\end{cases} \tag{2.86}$$

and $S_A(t)u_0 = u(t)$. If we solve in $\Omega^\#$ the semigroup operator will be called $S_A^\#$. For the remainder of the text $\beta : \mathbb{R} \to \mathbb{R}$ is a nondecreasing function, and if $A = \beta$ then $A$ represents the Nemitskij operator associated to $\beta$ in the sense that $Au = \beta \circ u$.

**Proposition 2.5.** *Let $f \in L^2(\Omega)$ and $u_0 \preceq v_0$ . Then, for a.e. $t \in (0,T)$*

$$S_{-\Delta - f}(t)u_0 \preceq S^{\#}_{-\Delta - f^{\#}}(t)v_0. \tag{2.87}$$

*Proof.* We proceed as in [Ban76b]. We write a parabolic inequality for $U(t,\sigma,y)$, whereas for $V$ the equality holds and the result follows from the comparison principle. □

**Proposition 2.6.** *Let $\beta$ be convex and $u_0 \preceq v_0$. Then, for a.e. $t \in (0,T)$*

$$S_{\beta}(t)u_0 \preceq S^{\#}_{\beta}(t)v_0. \tag{2.88}$$

**Proposition 2.7.** *Let $f \in L^2(\Omega)$ and $u_0 \preceq v_0$. Then, for a.e. $t \in (0,T)$*

$$S_{-\Delta+\beta-f}(t)u_0 \preceq S^{\#}_{-\Delta+\beta-f^{\#}}(t)v_0 \tag{2.89}$$

*Proof.* The Trotter-Kato formula applies

$$\left( S_{-\Delta - f}\left(\frac{t}{n}\right) S_{\beta}\left(\frac{t}{n}\right)\right)^n u_0 \to S_{-\Delta+\beta-f}(t)u_0, \tag{2.90}$$

$$\left( S^{\#}_{-\Delta - f^{\#}}\left(\frac{t}{n}\right) S^{\#}_{\beta}\left(\frac{t}{n}\right)\right)^n v_0 \to S^{\#}_{-\Delta+\beta-f^{\#}}(t)v_0. \tag{2.91}$$

Thus, by applying the previous propositions

$$\left( S_{-\Delta - f}\left(\frac{t}{n}\right) S_{\beta}\left(\frac{t}{n}\right)\right)^n u_0 \preceq \left( S^{\#}_{-\Delta - f^{\#}}\left(\frac{t}{n}\right) S^{\#}_{\beta}\left(\frac{t}{n}\right)\right)^n v_0. \tag{2.92}$$

as $n \to +\infty$ in $L^2(\Omega)$ and $L^2(\Omega^{\#})$ uniformly in $t$ for $t \in [0,T]$. Hence, the comparison holds in the limit. □

*Proof of Theorem 2.16.* Let us first assume that $f \in \mathscr{C}([0,T];L^2(\Omega))$. Let us split $[0,T]$ in $n$ parts $t_k^n = \frac{k}{n}T$ and let $f^n$ be a piecewise constant function, given by

$$f^n(t) = f(t_k^n) \text{ for } t \in [t_k^n, t_{k+1}^n).$$

Take the mild solution $u^n$. Since $f^n$ is piecewise constant, a simple induction argument shows that we can take the semigroups piecewise, for $t \in [t_k^n, t_{k+1}^n)$

$$u^n(t) = S_{-\Delta+\beta-f(t_k^n)}(t-t_k^n) S_{-\Delta+\beta-f(t_{k-1}^n)}\left(\frac{T}{n}\right) \cdots S_{-\Delta+\beta-f(0)}\left(\frac{T}{n}\right) u_0 \tag{2.93}$$

Applying the same reasoning for $v_n$ we get, for $t \in [t_k^n, t_{k+1}^n)$

$$v^n(t) = S^{\#}_{-\Delta+\beta-f(t_k^n)^{\#}}(t-t_k^n)\, S^{\#}_{-\Delta+\beta-f(t_{k-1}^n)^{\#}}\left(\frac{T}{n}\right)\cdots S^{\#}_{-\Delta+\beta-f(0)^{\#}}\left(\frac{T}{n}\right) v_0 \tag{2.94}$$

Applying the formulation (2.93) and (2.94) and Proposition 2.7 finite times we have that, for a.e. $t \in [0,T]$

$$u^n(t) \preceq v^n(t). \tag{2.95}$$

It is easy to check that $f^n \to f$ in $L^2((0,T)\times\Omega)$ and $(f^n)^{\#} = (f^{\#})^n \to f^{\#}$ in $L^2((0,T)\times\Omega^{\#})$ (see, for example, [Rak08]). Due to the properties of the equations we have that $u_n \to u$ in $L^2((0,T)\times\Omega)$ and $v_n \to v$ in $L^2((0,T)\times\Omega^{\#})$. Therefore, the comparison (2.95) holds also in the limit, which concludes the proof for $f \in \mathscr{C}([0,T];L^2(\Omega))$. For $f \in L^2((0,T)\times\Omega)$ we take a sequence of functions $(f^n) \in \mathscr{C}([0,T]\times\Omega)$, $f^n \to f$ in $L^2((0,T)\times\Omega)$. Due to the continuity of $\cdot^{\#} : L^2(\Omega\times(0,T)) \mapsto L^2(\Omega^{\#}\times(0,T))$ the result follows. □

### 2.7.5 Steiner symmetrization of semilinear elliptic problems

In [DGC15b] and [DGC16] the proof of the semilinear elliptic problem is done by passing to the limit in the parabolic problem. If the comparison holds for every time it holds for the limit elliptic problem. The details can be found in the indicated papers.

**Theorem 2.15** ([DGC15b]). *Let g be concave, verifying*

$$\int_0^\tau \frac{d\sigma}{g(\sigma)} < \infty, \qquad \forall \tau > 0. \tag{2.96}$$

*Let $f \in L^2(\Omega)$ with $f \geq 0$. Let $u \in H_0^1(\Omega)$ and $v \in H_0^1(\Omega^{\#})$ be the unique solutions of*

$$(P)\begin{cases} -\Delta u + g(u) = f & \text{in } \Omega, \\ u = 0, & \text{on } \partial\Omega \end{cases}$$

$$(P^{\#})\begin{cases} -\Delta v + g(v) = f^{\#}, & \text{in } \Omega^{\#}, \\ v = 0, & \text{on } \partial\Omega^{\#}. \end{cases}$$

*Then, for any $s \in [0,|\Omega'|]$ and a.e. $y \in \Omega''$.*

$$\int_0^s u^*(\sigma,y)d\sigma \leq \int_0^s v^*(\sigma,y)d\sigma.$$

**Theorem 2.16** ([DGC16]). *Let $\beta$ be a concave continuous nondecreasing function such that $\beta(0) = 0$. Let $f \in L^2(\Omega)$ with $f \geq 0$, $0 \leq w_0 \leq 1$. Let $w \in H^1(\Omega)$ and $z \in H^1(\Omega^\#)$ be the unique solutions of*

$$(P) \begin{cases} -\Delta w + \lambda \beta(w) = f & \text{in } \Omega, \\ w = 1 & \text{on } \partial\Omega, \end{cases}$$

$$(P^\#) \begin{cases} -\Delta z + \lambda \beta(z) = f^\# & \text{in } \Omega^\#, \\ z = 1 & \text{on } \partial\Omega^\#, \end{cases}$$

*Then, for any $s \in [0, |\Omega'|]$ and a.e. $y \in \Omega''$*

$$\int_s^{|\Omega'|} z^*(\sigma, y) d\sigma \leq \int_s^{|\Omega'|} w^*(\sigma, y) d\sigma. \tag{2.97}$$

# 2.8 Other kinds of rearrangements

## 2.8.1 Relative rearrangement

Another rearrangement technique which has brought a lot of results in the last century is known as relative rearrangement. We will use its properties in Chapter 4.

Let us define this rearrangement

**Definition 2.7.** Let $u$ be a measurable function, $v \in L^1(\Omega)$ and, for $s \in [0, |\Omega|]$

$$w(s) = \begin{cases} \displaystyle\int_{Q(s)} v(x)dx, & \text{if } |P(s)| = 0, \\ \displaystyle\int_{Q(s)} v(x)dx + \int_0^{s-|Q(s)|} (v|_{P(s)})^*(\sigma)d\sigma, & \text{if } |P(s)| \neq 0. \end{cases} \tag{2.98}$$

where

$$P(s) = \{x \in \Omega : u(x) = u^*(s)\}, \tag{2.99}$$

$$Q(s) = \{x \in \Omega : u(x) > u^*(s)\}. \tag{2.100}$$

The relative rearrangement of $v$ with respect to $u$ as

$$v_u^*(s) = \frac{dw}{ds}(s). \tag{2.101}$$

It is known that if $v \in L^p(\Omega)$ then $\|v_u^*\|_{L^p(0,|\Omega|)} \leq \|v\|_{L^p(\Omega)}$. Many other properties can be found in [Rak88].

One of such results is the following:

**Theorem 2.17** (Rakotoson and Temam [RT90]). *Let $\Omega$ be a bounded connected set, $\partial\Omega$ of class $C^2$. Let $f \in L^1(\Omega)$ and $u \in W^{1,1}(\Omega)$. For almost every $t \in (essinf(u), esssup(u))$ (where essinf and esssup are the essential infimum and supremum) we have*

$$\frac{d}{dt}\int_{\{x\in\Omega: u(x)>t\}} f(x)dx = \mu'(t) f_u^*(\mu(t)). \tag{2.102}$$

Another application of such kind of rearrangement is the obtention of a $L^\infty$ bound for semilinear equations given by Leray-Lions type operators (see [Rak87]). The results can be generalized to weighted spaces, and similar results are obtained (see [RS97; RS93a; RS93b]).

This kind of technique has also been used in Chapter 4 to obtain estimates of very weak solutions.

We complete this chapter by saying a few words on weighted rearrangements.

## Gaussian rearrangement

The previous approaches to rearrangement yields good results in bounded domains. However, unbounded domains are not covered, since $\Omega^\star$ would be the whole of $\mathbb{R}^n$. A rearrangement based on the the Gaussian distribution

$$\varphi(x) = (2\pi)^{-\frac{n}{2}} \exp\left(-\frac{|x|}{2}\right) \tag{2.103}$$

can be applied. It has been used in several papers [BBMP02; Di 03; CO04] with good results.

The idea is to define a weighted measure to substitute the Lebesgue measure

$$|\Omega|_\varphi = \int_\Omega \varphi \in [0,1]. \tag{2.104}$$

We this definition in mind the can the rearrangement of a set as

$$\Omega_\varphi^\sharp = \{x = (x_1, \cdots, x_n) : x_1 > a\} \text{ such that } |\Omega_\varphi^\sharp| = |\Omega|. \tag{2.105}$$

For a function $u : \Omega \to \mathbb{R}$ we define

$$\mu_\varphi(t) = |\{x \in \Omega : |u(x)| > t\}|_\varphi, \tag{2.106}$$

$$u_\varphi^*(s) = \inf\{t \geq 0 : \mu(t) \leq s\}, s \in [0,1] \tag{2.107}$$

and $u_\varphi^\sharp : \Omega_\varphi^\sharp \to \mathbb{R}$ as the only function such that

$$|\{x : \Omega_\varphi^\sharp : |u_\varphi^\sharp(x)| > t\}|_\varphi = |\{x \in \Omega : |u(x)| > t\}|_\varphi. \tag{2.108}$$

Naturally $u_\varphi^\sharp(x) = u_\varphi^\sharp(x_1) = u^*(k(x_1))$. This rearrangement is well suited to compare elliptic problem in which the coefficient are of Gaussian-type. Let us state here the result for the parabolic problem

**Theorem 2.18** ([CO04]). *Let u be an analytical solution of*

$$\begin{cases} \frac{\partial u}{\partial t} - \frac{1}{\varphi}\operatorname{div}(\varphi \nabla u) = f & Q = \Omega \times (0,T), \\ u = u_0 & t = 0, \\ u = 0 & \partial\Omega. \end{cases} \tag{2.109}$$

*where $f \in L^2(Q, \varphi)$ and $u_0 \in H_0^1(\Omega, \varphi)$. Let v be the symmetrized solution*

$$\begin{cases} \frac{\partial v}{\partial t} - \frac{1}{\varphi}\operatorname{div}(\varphi \nabla v) = f_\varphi^\sharp & \Omega_\varphi^\sharp \times (0,T), \\ v = (u_0)_\varphi^\sharp & t = 0, \\ v = 0 & \partial\Omega. \end{cases} \tag{2.110}$$

*Then,*

$$\int_0^s u_\varphi^*(t,\sigma)d\sigma \leq \int_0^s v_\varphi^*(t,\sigma)d\sigma \quad \forall s \in [0,1]. \tag{2.111}$$

# Chapter 3

# Shape optimization

## 3.1 Shape differentiation

The main goal of this section is to analyze the differentiability, with respect to the domain $\Omega$, of the *effectiveness factor* (2.3)

$$\mathscr{E}(\Omega) = \frac{1}{|\Omega|} \int_{\Omega} \beta(w_{\Omega}) dx. \tag{3.1}$$

For convenience, we will sometimes refer to the ineffectiveness (2.4). This will pass by the differentiation of functions $w_{\Omega}$ (resp. $u_{\Omega}$) defined by (2.1) (resp. (2.2)).

This kind of problem falls within the family of problems already considered by Hadamard [Had08] and it has been studied by several authors in the literature (see, e.g., [MS76; Pir12; Sim80] and the references therein). In the most general formulation this family of problems may be associated to the general boundary value problem:

$$\begin{cases} A(u(D)) = f, & \text{in } D, \\ B(u(D)) = g, & \text{on } \partial D \end{cases} \tag{3.2}$$

and the question is to study the differentiability with respect to $D$ of a functional which can by given generally as

$$J(D) = \int_{D} C(u_D) dx,$$

where $A, B, C$ are operators or functions that may contain also some derivatives of $u_D$ and $D$ is a domain belonging to a certain class.

As mentioned before, our aim is to study the differentiability of functional (2.4). We consider a fixed open bounded regular domain of $\mathbb{R}^n$, $\Omega_0$, and study its deformations given by a "small" function $\theta : \mathbb{R}^n \to \mathbb{R}^n$, so that the new domain is $\Omega = (I+\theta)\Omega_0$, where $I$ is the identity function

$$I(x) = x. \tag{3.3}$$

We consider, as it is the case in chemistry catalysis, $g$ and $f$ such that $0 \leq w_\Omega, u_\Omega \leq 1$. Besides the above mentioned references we recall here the articles [Der80] for a linear problem with a Dirichlet constant boundary condition and [MPM79] where a semilinear equation arising in combustion was considered (corresponding, in our formulation to take $g(u) = -e^u$).

First we studied in [DGC15a] the case in which $g$ (and $\beta$) are smooth (in $W^{2,\infty}(\mathbb{R})$). Then, in [GC17] we studied some non smooth cases. In particular, due to its importance in Chemical Engineering, we will discuss the case of root type non linearities. This is much more difficult, since for this kind on nonlinearity a dead core appears, and therefore the non-differentiable point of the nonlinearity might be in the range of the solution.

### 3.1.1 Fréchet derivative when $\beta \in W^{2,\infty}$

In order to obtain properties in the sense of derivatives, we consider two approaches, mimicking the approach in Differential Geometry. We first consider the global differentiability of solutions (as it was done in the linear cases in [HP05; All07] and for abstract problems in [Sim80]), which unfortunately requires derivatives in spaces of very regular functions, and then we take advantage of the differentiation along curves (the approach followed in [SZZ91]).

Let us recall that, for $\Omega \subset \mathbb{R}^n$, $u_\Omega$ is the unique solution of (2.2) (we assume that the formulation of the problem leads to the uniqueness of solutions). Let $\Omega_0$ be a fixed smooth domain. We will work in the family of deformations

$$\Omega_\theta = (I+\theta)\Omega_0 \tag{3.4}$$

where $\theta \in W^{1,\infty}(\mathbb{R}^n, \mathbb{R}^n)$. We will consider the Lagrangian representation of $u_{\Omega_\theta}$ as

$$u_\theta = u_{\Omega_\theta}, \tag{3.5}$$

and the Eulerian representations

$$\widehat{u}_\theta = u_\theta \circ (I+\theta). \tag{3.6}$$

Notice that

$$u_\theta : \Omega_\theta \to \mathbb{R} \qquad \widehat{u}_\theta : \Omega_0 \to \mathbb{R}.$$

It turns out that $\widehat{u}_\theta$ simplifies the study of the differentiability of $u_\Omega$ and the functional $\eta(\Omega)$ with respect to $\Omega$.

Our proof relies heavily on the Implicit Function Theorem. The application of this theorem requires an uniform choice of functional space, which would require some additional informations on $u$. This kind of difficulties in the functional setting is well portrayed in [Bre99].

For the nonlinearity $g$ we shall consider the following assumptions:

**Assumption 3.1.** $g$ is nondecreasing.

**Assumption 3.2.** The Nemitskij operator for $g$ (which we will denote again by $g$ in some circumstances, as a widely accepted abuse of notation)

$$G : H^1(\Omega) \rightarrow L^2(\Omega) \tag{3.7}$$

$$u \mapsto g \circ u \tag{3.8}$$

is well defined and is of class $\mathscr{C}^m$ for some $m \geq 1$.

We recall that Assumption 3.2 immediately implies that $[DG](v)\varphi = g'(v)\varphi$ for $\varphi, v \in H^1(\Omega)$ and that, if $G$ is of class $\mathscr{C}^k$, with $k > 1$ then necessarily $g(s) = as + b$ for some $a, b \in \mathbb{R}$ (see, e.g., [Hen93]).

Our first result collects some general results on the differentiability of the solution $u_\Omega$ with respect to $\Omega$:

**Theorem 3.1** ([DGC15a])**.** *Let g satisfy Assumption 3.1 and 3.2. Then, the map*

$$\begin{aligned} W^{1,\infty}(\mathbb{R}^n, \mathbb{R}^n) &\rightarrow H^1_0(\Omega_0) \\ \theta &\mapsto \widehat{u}_\theta \end{aligned}$$

*(where $\widehat{u}_\theta$ is defined by* (3.6)*) is of class $\mathscr{C}^l$ in a neighbourhood of 0 if $f \in H^k(\mathbb{R}^n)$ where $l = \min\{k,l\}$. Furthermore, the application*

$$\begin{array}{rcl} u : W^{1,\infty}(\mathbb{R}^n,\mathbb{R}^n) & \to & L^2(\mathbb{R}^n) \\ \theta & \mapsto & u_\theta \end{array}$$

*(where $u_\theta$ is given by* (3.5) *and extended by zero outside $\Omega_\theta$) is differentiable at* 0*. In fact $u'(0) : W^{1,\infty}(\mathbb{R}^n,\mathbb{R}^n) \to H^1(\Omega_0)$ and*

$$u'(0)\theta + \nabla u_{\Omega_0} \cdot \theta \in H_0^1(\Omega_0).$$

**Remark 3.1.** Since the function is only differentiable at 0 we will simplify write $u'$ to represent $u'(0)$.

One of the easiests ways to characterize the global derivative is, as usually, to compute the directional derivatives.

**Definition 3.1.** We will say that $\Phi$ is a curve of deformations of $\Omega_0$ if

$$\Phi : [0,T) \to W^{1,\infty}(\mathbb{R}^n,\mathbb{R}^n)$$

is such that $\det\Phi(\tau) > 0$ and $\Phi(0) = I$.

**Assumption 3.3.** We will say that $\theta$ is a curve of small perturbations of the identity if $\Phi(\tau) = I + \theta(\tau)$ is a curve of deformations and

i) $\theta : [0,T) \to W^{1,\infty}(\mathbb{R}^n,\mathbb{R}^n)$ is differentiable at 0 (from the right)

ii) $\theta(0) = 0$.

Sometimes we will consider higher order derivatives too. We will refer to $\theta$ or $\Phi$ indistinctively, since they relate by $\Phi(\tau) = I + \theta(\tau)$. It will be common that we consider the curve of deformations

$$\Phi(\tau) = I + \tau\theta, \tag{3.9}$$

for a fixed deformation $\theta \in W^{1,\infty}(\mathbb{R}^n,\mathbb{R}^n)$. In this notation we will admit the abuse of notation $\theta(\tau) = \tau\theta$, where naturally both elements are different, but this should not lead to confusion.

In this terms, the above theorem leads to:

**Corollary 3.1.** *Let $\Phi$ be a a curve of deformations of class $\mathscr{C}^k$. Then $\tau \mapsto v_{\theta(\tau)}$ is of class $\mathscr{C}^l$ with $l = \min\{m,k\}$.*

Our second result concerns the characterization of $u'$. We have:

**Theorem 3.2.** *Let g satisfy Assumption 3.1 and 3.2. Let $\theta$ be a curve satisfying assumptions 3.3. Then u is differentiable along $\Phi$ at least at 0. That is, the directional derivative $\frac{d}{d\tau}(u \circ \Phi)$ exists, and it is the solution $u'$ of the linear Dirichlet problem*

$$\begin{cases} -\Delta u' + \lambda g'(u_{\Omega_0})u' = 0 & \text{in } \Omega_0, \\ u' = -\nabla u_{\Omega_0} \cdot \theta'(0) & \text{on } \partial\Omega_0. \end{cases} \tag{3.10}$$

We point out that the above result shows, in other terms, for $\theta \in W^{1,\infty}(\mathbb{R}^n, \mathbb{R}^n)$, that $u'(0)\theta$ is the unique weak solution of the Dirichlet problem

$$\begin{cases} -\Delta u' + \lambda g'(u_{\Omega_0})u' = 0, & \text{in } \Omega_0, \\ u' = -\nabla u_{\Omega_0} \cdot \theta, & \text{on } \partial\Omega_0. \end{cases} \tag{3.11}$$

As consequence we have:

**Corollary 3.2.** *The function $u' : W^{1,\infty}(\mathbb{R}^n, \mathbb{R}^n) \to H^1(\Omega_0)$ is continuous. In fact, since due to Assumption 3.2, the solution u of (2.2) verifies $u \in W^{2,p}(\Omega_0)$ for any $p \in [1, +\infty)$, then for any $q \in [1, +\infty)$*

$$|u'(0)(\theta)|_q \leq c|\nabla u \cdot \theta|_{L^p(\partial\Omega_0)} \leq c|\theta|_\infty |\nabla u_{\Omega_0}|_{L^p(\partial\Omega_0)} \tag{3.12}$$

$$\leq c(p)|\theta|_\infty |u_{\Omega_0}|_{W^{2,p}(\Omega_0)}. \tag{3.13}$$

Concerning the differentiability of the effectiveness factor functional we have:

**Theorem 3.3.** *On the assumptions of Theorem 3.1, let*

$$\hat{\eta}(\theta) = \int_{(I+\theta)\Omega_0} g(u_{(I+\theta)\Omega_0})dx. \tag{3.14}$$

*Then $\eta$ is of class $\mathscr{C}^m$ in a neighbourhood of 0. It holds that*

$$\hat{\eta}^{(m)}(0)(\theta_1, \cdots, \theta_m) = \int_{\Omega_0} \frac{d^n}{d\theta_n \cdots d\theta_1}(g(\widehat{u}_\theta)J_\theta)\, dx. \tag{3.15}$$

*Its first derivative can be expressed in terms of u*

$$\hat{\eta}'(0)(\theta) = \int_{\Omega_0} \left(g'(u_{\Omega_0})u' + \operatorname{div}(g(u_{\Omega_0})\theta)\right) dx, \tag{3.16}$$

*and, if $\partial\Omega_0$ is Lipschitz,*

$$\hat{\eta}'(0)(\theta) = \int_{\Omega_0} g'(u_{\Omega_0})u'\, dx + g(0)\int_{\partial\Omega_0} \theta\cdot n\, dS, \tag{3.17}$$

*where $u' = u'(0)(\theta)$.*

As a direct consequence we get:

**Corollary 3.3.** *On the assumptions of Theorem 3.1, it holds that*

$$\eta'(\theta) = \frac{1}{|\Omega_0|}\left(\int_{\Omega_0} g'(u_{\Omega_0})u'\, dx - \eta(0)\int_{\partial\Omega_0} \theta\cdot n\, dS\right).$$

**Corollary 3.4.** *On the assumptions of Theorem 3.1, if $\Phi$ is a volume preserving curve of deformations then*

$$\eta'(\theta) = \frac{1}{|\Omega_0|}\int_{\Omega_0} g'(u_{\Omega_0})u'\, dx.$$

We point out that if $g$ is Lipschitz (i.e. $g \in W^{1,\infty}(\mathbb{R})$) then we get that

$$|\eta(\theta) - \eta(0)| = |\eta'(0)(\lambda\theta)| \le c|g'|_\infty |u|_{W^{2,p}} |\theta|_\infty.$$

The details of the proof of the results in this section can be found in [DGC15a]. They will be ommited here to speed-up the presentation of results.

#### 3.1.1.1 Functional setting: Nemitskij operators and the implicit function theorem.

Let us formalize what we mean by a *shape functional*. At the most fundamental level it should be a function defined over a set of domain, that is defined over a subset $\mathfrak{C}$ of $\mathscr{P}(\mathbb{R}^n)$. Since we want to differentiate this functional we, at the very least, need to define proximity, that is a way to define the *neighbourhood of a set*. As it is usual in the literature of shape optimization we work over the set of weakly differentiable bounded deformations with bounded derivative, i.e. over the Sobolev space $W^{1,\infty}(\mathbb{R}^n,\mathbb{R}^n)$.

**Definition 3.2.** We say that a functional $I : \mathfrak{C} \subset \mathscr{P}(\mathbb{R}^n) \to \mathbb{R}$ is defined on a neighbourhood of $\Omega_0 \subset \mathbb{R}^n$ if there exists $U$ a neighbourhood of 0 on $W^{1,\infty}(\mathbb{R}^n,\mathbb{R}^n)$ such that $I$ is defined over $\{(Id+\theta)(\Omega_0) : \theta \in U\}$. We say that $I$ is differentiable at $\Omega_0$ if the application

$$\begin{array}{rcl} W^{1,\infty}(\mathbb{R}^n;\mathbb{R}^n) & \to & \mathbb{R} \\ \theta & \mapsto & I((Id+\theta)(\Omega_0)) \end{array}$$

is differentiable at 0.

We present a sufficient condition so that Assumption 3.2 holds. This is widely used in the context of partial differential equations, but as far as we know no reference is known besides it being an exercise in [Hen93]. That being the case we provide a proof[1]. Other conditions, mainly on the growth of $g$ can be considered so that Assumption 3.1 and 3.2 hold.

**Lemma 3.1.1.** *Let $g \in W^{2,\infty}(\mathbb{R})$. Then the Nemitskij operator* (3.8) *(in the sense $L^p(\Omega) \to L^2(\Omega)$) is of class $\mathscr{C}^1$ for all $p > 2$. In particular, Assumption 3.2 holds.*

*Proof.* Let us define $G$ the Nemitskij operator defined in (3.8). Consider it $G : L^p(\Omega) \to L^2(\Omega)$ for $p \geq 2$. We first have that, for $L = \max\{\|g\|_\infty, \|g'\|_\infty, \|g''\|_\infty\}$

$$\|G(u) - G(v)\|_{L^2}^2 = \int_\Omega |g(u) - g(v)|^2 dx \leq L \int_\Omega |u - v|^2 dx$$

so that $G$ is continuous. For $p > 2$ let $\varphi \in \mathscr{C}^\infty(\Omega)$ we compute

$$\|g(u+\varphi) - g(u) - g'(u)\varphi\|_{L^2}^2 = \int_\Omega |g'(\xi(x)) - g'(u(x))|^2 |\varphi(x)|^2 dx$$

for some function $\xi(x)$, between $u(x)$ and $u(x)+\varphi(x)$, due to the intermediate value theorem. We have that

$$\begin{aligned} |g'(\xi(x)) - g'(u((x))| &\leq L|\xi(x) - u(x)| \leq L|\varphi(x)| \\ |g'(\xi(x)) - g'(u(x))| &\leq 2L \\ |g'(\xi(x)) - g'(u(x))| &\leq L2^{1-\alpha}|\varphi(x)|^\alpha, \quad \forall \alpha \in (0,1). \end{aligned}$$

Therefore,

$$\|g(u+\varphi) - g(u) - g'(u)\varphi\|_{L^2}^2 \leq L^2 2^{2-2\alpha} \int_\Omega |\varphi(x)|^{2+2\alpha} dx.$$

Let $2 < p < 4$ then we have that $p = 2 + 2\alpha$ with $0 < \alpha < 1$. We then have that

$$\|g(u+\varphi) - g(u) - g'(u)\varphi\|_{L^2} \leq L2^{1-\alpha} \|\varphi(x)\|_{L^p}^{1+\alpha},$$

[1] This candidate is thankful to Prof. J.M. Arrieta for the details of the proof.

which proves the Fréchet differenciability. For $p > 4$ we have that $L^p(\Omega) \hookrightarrow L^3(\Omega)$. Furthermore, for any given dimension $n$ we can use the Sobolev inclusions $H^1(\Omega) \hookrightarrow L^p(\Omega)$ with $p > 2$, proving the desired differenciability. □

Some other well-known results are quoted now:

**Theorem 3.4.** *Let $g \in W^{1,p}(\mathbb{R}^n)$. Then the map*

$$\begin{aligned} \mathfrak{G} : W^{1,\infty}(\mathbb{R}^n, \mathbb{R}^n) &\to L^p(\mathbb{R}^n) \qquad (3.18)\\ \theta &\mapsto g \circ (I + \theta) \qquad (3.19) \end{aligned}$$

*is differentiable in a neighbourhood of* 0 *and*

$$\mathfrak{G}'(0) = (\nabla g) \circ (I + \theta).$$

**Theorem 3.5** ([HP05, Lemme 5.3.3.])**.** *Let*

$$\begin{aligned} g : W^{1,\infty}(\mathbb{R}^n, \mathbb{R}^n) &\to L^p(\mathbb{R}^n),\\ \Psi : W^{1,\infty}(\mathbb{R}^n, \mathbb{R}^n) &\to W^{1,\infty}(\mathbb{R}^n, \mathbb{R}^n) \end{aligned}$$

*continuous at* 0 *with* $\Psi(0) = I$*,*

$$\begin{aligned} W^{1,\infty}(\mathbb{R}^n, \mathbb{R}^n) &\to L^p(\mathbb{R}^n) \times L^\infty(\mathbb{R}^n; \mathbb{R}^n) \qquad (3.20)\\ \theta &\mapsto (g(\theta), \Psi(\theta)) \qquad (3.21) \end{aligned}$$

*differentiable at* 0*, with* $g(0) \in W^{1,p}(\mathbb{R}^n)$ *and*

$$g'(0) : W^{1,\infty}(\mathbb{R}^n, \mathbb{R}^n) \to W^{1,p}(\mathbb{R}^n)$$

*is continuous. Then the application*

$$\begin{aligned} \mathfrak{G} : W^{1,\infty}(\mathbb{R}^n, \mathbb{R}^n) &\to L^p(\mathbb{R}^n) \qquad (3.22)\\ \theta &\mapsto g(\theta) \circ \Psi(\theta) \qquad (3.23) \end{aligned}$$

*is differentiable at* 0 *and*

$$\mathfrak{G}'(0) = g'(0) + \nabla g(0) \cdot \Psi'(0).$$

To conclude this section we state a classical result, the Implicit Function Theorem. This result is typically a direct consequence of the Inverse Function Theorem. In the Banach space setting this result is originally due to Nash and Moser (see [Nas56; Mos66]). In the detailed survey [Ham82] the author points towards Zehnder [Zeh76] as one of the first presentations as implicit function theorem.

**Theorem 3.6** (Implicit Function Theorem)**.** *Let $X, Y$ and $Z$ be Banach spaces and let $U, V$ be neighbourhoods on $X$ and $Y$, respectively. Let $F : U \times V \to Z$ be continuous and differentiable, and assume that $D_yF(0,0) \in \mathscr{L}(Y,Z)$ is bijective. Let us assume, further, that $F(0,0) = 0$. Then there exists $W$ neighbourhood of $0$ on $X$ and a differentiable map $f : W \to Y$ such that $F(x, f(x)) = 0$. Furthermore, for $x$ and $y$ small, $f(x)$ is the only solution $y$ of the equation $F(x,y) = 0$. If $F$ is of class $\mathscr{C}^m$ then so is $f$.*

#### 3.1.1.2 Differentiation of solutions

For the reader convenience we repeat here the general result in [Sim80]:

**Theorem 3.7.** *Let $D$ be a bounded domain such that $\partial D$ be a piecewise $\mathscr{C}^1$ and assume that $D$ is locally on one side of $\partial D$. Let $u_0$ be the solution of (3.2). Let us use the notation $\mathscr{C}^k = \mathscr{C}^k(\mathbb{R}^n, \mathbb{R}^n)$ and $k \geq 1$. Assume that*

$$u(\theta) \in W^{m,p}((I+\theta)D) \tag{3.24}$$

*and that for every open set $D'$ close to $D$ (for example $D' = (I+\theta)D$ for small $\theta$ in the norm of $\mathscr{C}^k$), $A, B, C : W^{m-1,p}(D') \to \mathscr{D}'(D)$ are differentiable and*

$$\begin{cases} A : W^{m,p}(D') \to \mathscr{D}'(D') \\ B : W^{m,p}(D') \to W^{1,1}(D') \\ C : W^{m,p}(D') \to L^1(D') \end{cases} \tag{3.25}$$

*and*

$$\mathscr{C}^k \to W^{m,p} \tag{3.26}$$

$$\theta \mapsto u(\theta) \circ (I+\theta) \tag{3.27}$$

*is differentiable at $0$. Then:*

*i) The solution $u$ is differentiable in the sense that*

$$u : \mathscr{C}^k \to W^{m-1,p}_{loc}(D) \text{ is differentiable}$$

*and the derivative (i.e. the local derivative $u'$ in the direction of $\tau$) satisfies*

$$\frac{\partial A}{\partial u}(u_0)u' = 0, \text{ in } D. \tag{3.28}$$

*ii) If*

$$\begin{cases} \theta \mapsto B(u(\theta)) \circ (I+\theta) \text{ is differentiable at 0 into } W^{1,1}(D), \\ \qquad \text{(i.e. with the } W^{1,1}(D) \text{ topology in the image set)} \\ B(u_0) \in W^{2,1}(D), \\ g \in W^{2,1}(\mathbb{R}^n) \end{cases} \tag{3.29}$$

*then $u'$ satisfies*

$$\frac{\partial B}{\partial u}(u_0)u' = -\tau \cdot n \frac{\partial}{\partial n}(B(u_0) - g). \tag{3.30}$$

*iii) If*

$$\begin{cases} \theta \mapsto C(u(\theta)) \circ (I+\theta) \text{ is differentiable at 0 into } L^1(D), \\ C(u_0) \in W^{1,1}(D), \end{cases} \tag{3.31}$$

*then $\theta \mapsto J(\theta)$ is differentiable and its directional derivative in the direction of $\tau$ is:*

$$\frac{\partial J}{\partial \theta}(0)\tau = \int_D \frac{\partial C}{\partial u} u' \, dx + \int_{\partial D} \tau \cdot n C(u_0) \, dS. \tag{3.32}$$

#### 3.1.1.3 Differentiation under the integral sign

We shall follow some reasonings similar to the ones presented in [HP05]. Let us define $\Omega_\tau = \Phi(\tau, \Omega_0)$ and consider a function $f$ such that $f(\tau) \in L^1(\Omega_\tau)$. We take interest on the map

$$I : \mathbb{R} \rightarrow \mathbb{R} \tag{3.33}$$

$$\tau \mapsto \int_{\Omega_\tau} f(\tau, x) \, dx = \int_{\Omega_0} f(\tau, \Phi(\tau, y)) J(\tau, y) \, dy \tag{3.34}$$

where $f(\tau, x) = f(\tau)(x)$ and the Jacobian

$$J(\tau, y) = \det(D_y \Phi(\tau, y)).$$

**Theorem 3.8.** *Let $\Phi$ satisfy Assumption 3.3, $f$ such that*

$$\begin{array}{rcl} f:[0,T) & \to & L^1(\mathbb{R}^n) \\ \tau & \mapsto & f(\tau) \end{array}$$

*is differentiable at* 0 *and, besides, it satisfies the spatial regularity at* $\tau = 0$

$$f(0) \in W^{1,1}(\mathbb{R}^N).$$

*Then,* $\tau \mapsto I(\tau) = \int_{\Omega_\tau} f(\tau)$ *is differentiable at* 0 *and*

$$I'(0) = \int_{\Omega_0} f'(0) + \operatorname{div}\left( f(0)\frac{\partial \Phi}{\partial \tau}(0)\right).$$

*If* $\Omega_0$ *is an open set with Lipschitz boundary then*

$$I'(0) = \int_{\Omega_0} f'(0) + \int_{\partial\Omega_0} f(0) n \cdot \frac{\partial \Phi}{\partial \tau}(0).$$

In simpler terms, under regularity it holds that

$$\frac{\partial}{\partial \tau}\bigg|_{\tau=0} \left( \int_{G_\tau} f(\tau,x)dx \right) = \int_{\Omega_0} \left\{ \frac{\partial f}{\partial \tau}(0,x) + \operatorname{div}\left( f(0,x)\frac{\partial \Phi}{\partial \tau}(0,x)\right)\right\} dx. \tag{3.35}$$

We have some immediate consequences of Theorem 3.8

**Lemma 3.1.2.** *Let $g \in W^{1,1}(\mathbb{R}^N)$ and $\Psi : [0,T) \to W^{1,\infty}$ be continuous at* 0 *such that $\Psi : [0,T) \to L^\infty$ is differentiable at* 0, *and let $Z$ be its derivative. Then*

$$G : [0,T) \to L^1(\mathbb{R}^n) \tag{3.36}$$

$$\tau \mapsto g \circ \Psi(\tau) \tag{3.37}$$

*is differentiable at* 0 *and* $G'(0) = \nabla g \cdot Z$.

**Lemma 3.1.3** (Differentiation under the integral sign)**.** *Let $E$ be a Banach space and*

$$f : E \times \Omega \to \mathbb{R} \tag{3.38}$$

$$(v,y) \mapsto f(v,y) \tag{3.39}$$

*such that*

$$\tilde{f} : E \to L^1(\Omega) \tag{3.40}$$
$$v \mapsto f(v,\cdot) \tag{3.41}$$

*is differentiable at $v_0$. Let*

$$F : E \to \mathbb{R} \tag{3.42}$$
$$v \mapsto \int_\Omega f(v,y)dy \tag{3.43}$$

*Then $F$ is differentiable at $v_0$ and*

$$DF(v) = \int_\Omega (D_v \tilde{f})(v)(y).$$

### 3.1.2 Gateaux derivative when $\beta \in W^{1,\infty}$

Once the case $\beta \in W^{2,\infty}(\mathbb{R})$ is understood, let us focus on the less smooth case $\beta \in W^{1,\infty}(\mathbb{R})$. In this case, we can only prove that the shape derivative exists in the Gateaux sense (which is weaker than the Fréchet sense).

**Theorem 3.9.** *Let $\theta \in W^{1,\infty}(\mathbb{R}^n,\mathbb{R}^n)$, $\beta \in W^{1,\infty}(\mathbb{R})$ be nondecreasing such that $\beta(0) = 0$ and $f \in H^1(\mathbb{R}^n)$. Then, the applications*

$$\begin{aligned} \mathbb{R} &\to L^2(\Omega_0) \\ \tau &\mapsto u_{(I+\tau\theta)\Omega_0} \circ (I+\tau\theta), \end{aligned}$$

*and*

$$\begin{aligned} \mathbb{R} &\to L^2(\mathbb{R}^n) \\ \tau &\mapsto u_{(I+\tau\theta)\Omega_0} \end{aligned}$$

*are differentiable at* 0*. Furthermore, $\frac{du_\tau}{d\tau}|_{\tau=0}$ is the unique solution of* (3.11).

In this chapter we will be particularly interested in the case in which $\beta'$ only has blow-up at $w = 0$. Let us define

$$v = \left.\frac{dw_\tau}{d\tau}\right|_{\tau=0}. \tag{3.44}$$

We can rewrite (3.11) in terms of $w$

$$\begin{cases} -\Delta v + \beta'(w_{\Omega_0}) v = 0 & \Omega, \\ v + \nabla w_{\Omega_0} \cdot \theta = 0 & \partial\Omega. \end{cases} \tag{3.45}$$

**Remark 3.2.** In most cases, the process of homogenization developed in Chapter 1 leads to an homogeneous equation (2.1) in which $\beta$ is the same as the function in the microscopic problem, and thus it is natural that $\beta$ be singular at 0. However, it sometimes happens that the limit kinetic is different. In the homogenization of problems with particles of critical size (see [DGCPS17c]) it turns out that the resulting kinetic in the macroscopic homogeneous equation (2.1) satisfies $\beta \in W^{1,\infty}$, even when the original kinetic of the microscopic problem was a general maximal monotone graph.

#### 3.1.2.1 From $W^{2,\infty}$ to $W^{1,\infty} \cap \mathscr{C}^1$

Let us show that the shape derivative is continuously dependent on the nonlinearity, and thus that we can make a smooth transition from the Fréchet scenario presented in [DGC15a] to our current case. For the rest of the paper we will use the notation:

**Lemma 3.1.4.** *Let $f \in L^2(\mathbb{R}^n)$, $\beta \in W^{1,\infty}(\mathbb{R})$ be a nondecreasing function such that $\beta(0) = 0$ and let $\beta_n \in W^{2,\infty}(\mathbb{R})$ nondecreasing such that $\beta_n(0) = 0$. Let $w_n$ be the unique solution of*

$$\begin{cases} -\Delta w_n + \beta_n(w_n) = f & \Omega_0, \\ w_n = 1 & \partial\Omega_0. \end{cases} \tag{3.46}$$

*Then*

$$\|w_n - w\|_{H^1(\Omega_0)} \le C \|\beta_n - \beta\|_{L^\infty(\mathbb{R})} \tag{3.47}$$

$$\|w_n - w\|_{H^2(\Omega_0)} \le C(1 + \|\beta'\|_{L^\infty(\mathbb{R})}) \|\beta_n - \beta\|_{L^\infty(\mathbb{R})}. \tag{3.48}$$

*Furthermore, let $\beta \in C^1(\mathbb{R}) \cap W^{1,\infty}(\mathbb{R})$ and $v_n$ be the unique solution of*

$$\begin{cases} -\Delta v_n + \beta_n'(w_n) v_n = 0 & \Omega_0, \\ v_n + \nabla w_n \cdot \theta = 0 & \partial\Omega_0. \end{cases} \tag{3.49}$$

*Then, if $\beta_n \to \beta$ in $W^{1,\infty}(\mathbb{R})$,*

$$v_n \rightharpoonup v \text{ in } H^1(\Omega_0). \tag{3.50}$$

**Remark 3.3.** In (3.47) we used the notation

$$||\beta_n - \beta||_{L^\infty} = \sup_{x\in\mathbb{R}} |\beta_n(x) - \beta(x)|.$$

It doesn't mean that either $\beta_n$ or $\beta$ are $L^\infty(\mathbb{R})$ functions themselves, but rather that their difference is pointwise bounded. In fact, this bound is destined to go 0 as $n \to +\infty$.

### 3.1.3 Shape derivative with a dead core

We can prove that the shape derivative in the smooth case has, under some assumptions, a natural limit when $\beta$ is not smooth.

In some cases in the applications (see, e.g., [Día85]) we can take $\beta$ so that $\beta'(w_{\Omega_0})$ has a blow up. It is common, specially in Chemical Engineering, that $\beta'(0) = +\infty$ and

$$N_{\Omega_0} = \{x \in \Omega_0 : w_{\Omega_0}(x) = 0\}$$

exists and has positive measure (see [Día85]). This is region is known as a dead core. In this case $\beta'(w_{\Omega_0}) = +\infty$ in $N_{\Omega_0}$. Due to this fact, the natural behaviour of the weak solutions of (3.45) is $v = 0$ in $N_{\Omega_0}$. We have the following result

**Theorem 3.10.** *Let $\beta$ be nondecreasing, $\beta(0) = 0$, $\beta'(0) = +\infty$,*

$$\beta \in \mathscr{C}(\mathbb{R}) \cap \mathscr{C}^1(\mathbb{R}\setminus\{0\}),$$

*and assume that $|N_{\Omega_0}| > 0$, $\theta \in W^{1,\infty}(\mathbb{R}^n,\mathbb{R}^n)$ and $0 \le f \le \beta(1)$. Then, there exists $v$ a solution of*

$$\begin{cases} -\Delta v + \beta'(w_{\Omega_0})v = 0 & \Omega_0 \setminus N_{\Omega_0}, \\ v = 0 & \partial N_{\Omega_0}, \\ v = -\nabla w_{\Omega_0} \cdot \theta & \partial\Omega_0, \end{cases} \tag{3.51}$$

*in the sense that $v \in H^1(\Omega_0)$, $v = 0$ in $N_{\Omega_0}$, $v = -\nabla w_{\Omega_0} \cdot \theta$ in $L^2(\partial\Omega_0)$, $\beta'(w_{\Omega_0})v^2 \in L^1(\Omega_0)$ and*

$$\int_{\Omega_0\setminus N_{\Omega_0}} \nabla v \nabla \varphi + \int_{\Omega_0\setminus N_{\Omega_0}} \beta'(w) v\varphi = 0 \tag{3.52}$$

*for every $\varphi \in W_c^{1,\infty}(\Omega_0 \setminus N_{\Omega_0})$. Furthermore, for $m \in \mathbb{N}$, consider $\beta_m$ defined by*

$$\beta_m'(s) = \min\{m, \beta'(s)\}, \qquad \beta_m(0) = \beta(0) = 0,$$

*and let $w_m, v_m$ be the unique solutions of* (3.46) *and* (3.49). *Then,*

$$v_m \rightharpoonup v, \quad in\ H^1(\Omega_0), \tag{3.53}$$

*where $v$ is a solution of* (3.51).

The uniqueness of solutions of (3.51) when $\beta'(w_{\Omega_0})$ blows up is by no means trivial. Problem (3.51) can be written in the following way:

$$-\Delta v + V(x)v = f \tag{3.54}$$

where $V(x) = \beta'(w_{\Omega_0}(x))$ may blow up as a power of the distance to a piece of the boundary. This kind of problems are common in Quantum Physics, although their mathematical treatment is not always rigorous (cf. [Día15; Día17]).

In the next section we will show some estimates on $\beta'(w_{\Omega_0})$. Let us state here some uniqueness results depending on the different blow-up rates.

When the blow-up is subquadratic (i.e. not *too* rapid), by applying Hardy's inequality and the Lax-Migram theorem, we have the following result (see [Día15; Día17]).

**Corollary 3.5.** *Let $N_{\Omega_0}$ have positive measure and $\beta'(u(x)) \leq Cd(x, N_{\Omega_0})^{-2}$ for a.e. $x \in \Omega_0 \setminus N_{\Omega_0}$. Then the solution $v$ is unique.*

The study of solutions of problem (3.54) in $\Omega_0$ when $V \in L^1_{loc}(\Omega_0)$ was carried out by many authors (see [DR10; DGCRT17] and the references therein). Existence and uniqueness of this problem in the case $V(x) \geq Cd(x, \partial\Omega_0)^{-r}$ with $r > 2$ was proved in [DGCRT17]. Applying these techniques one can show that

**Corollary 3.6.** *Let $N_{\Omega_0}$ have positive measure and $\beta'(w(x)) \geq Cd(x, N_{\Omega_0})^{-r}, r > 2$ for a.e. $x \in \Omega_0 \setminus N_{\Omega_0}$. Then the solution $v$ is unique.*

Similar techniques can be applied to the case $\beta'(w(x)) \geq Cd(x, N_{\Omega_0})^{-2}$. This will be the subject of a further paper[2].

[2] At the time of writing this thesis, there is a draft of such a paper by this candidate jointly with J.I. Díaz and J.M. Rakotoson

### 3.1.4 Estimates of $w_{\Omega_0}$ close to $N_{\Omega_0}$

Let us study the solution $w_{\Omega_0}$ on the proximity of the dead core and the blow up behaviour of $\beta'(w_{\Omega_0})$. First, we recall a well-known example

**Example 3.1.** Explicit radial solutions with dead core are known when $\beta(w) = |w|^{q-1}w$ $(0 < q < 1)$, $\Omega_0$ is a ball of large enough radius and $f$ is radially symmetric. In this case it is known that $N_{\Omega_0}$ exists, has positive measure and

$$\frac{1}{C}d(x,N_{\Omega_0})^{-2} \leq \beta'(w_{\Omega_0}) \leq Cd(x,N_{\Omega_0})^{-2}.$$

For the details see [Día85].

In fact, we present here a more general result to study the behaviour in the proximity of the dead core, based on estimates from [Día85].

**Proposition 3.1.** *Let $f = 0$, $\beta$ be continuous, monotone increasing such that $\beta(0) = 0$, $w$ be a solution of* (2.1) *that develops a dead core $N_{\Omega_0}$ of positive measure and assume that $\partial N_{\Omega_0} \in \mathscr{C}^1$. Define*

$$G(t) = \sqrt{2}\left(\int_0^t \beta(\tau)d\tau + \alpha t\right)^{\frac{1}{2}}, \quad \textit{where } \alpha = \max\left\{0, \min_{x\in\partial\Omega_0} H(x)\frac{\partial w}{\partial n}(x)\right\}, \tag{3.55}$$

*and assume that $\frac{1}{G} \in L^1(\mathbb{R})$. Then*

$$w_{\Omega_0}(x) \leq \Psi^{-1}(d(x,N_{\Omega_0})), \quad \textit{in a neighbournood of } N_{\Omega_0}, \tag{3.56}$$

*where $\Psi(s) = \int_0^s \frac{dt}{G(t)}$.*

**Example 3.2** (Root type reactions)**.** Let $f = 0$, $\beta(s) = \lambda|s|^{q-1}s$ with $0 < q < 1$ and let $\Omega_0$ be a convex set such that $N_{\Omega_0}$ exists and satisfies that $\partial N_{\Omega_0} \in \mathscr{C}^1$. Then

$$w_{\Omega_0}(x) \leq Cd(x,N_{\Omega_0})^{\frac{2}{1-q}}. \tag{3.57}$$

Furthermore

$$\beta'(w_{\Omega_0}(x)) \geq Cd(x,N_{\Omega_0})^{-2}. \tag{3.58}$$

**Remark 3.4.** The regularity assumptions on $\partial N_{\Omega_0}$ are by no means trivials. Examples can be constructed in which this does not hold. However, there are many cases of relevance the applications in which this regularity holds.

## 3.2 Convex optimization of the homogenized solutions

This section includes results published in [DGCT15] [DGCT16].

For the homogenized problem, we have the following optimality result:

**Theorem 3.11.** *Let* $1 \leq \alpha < \frac{n}{n-2}$, $0 < \theta < |Y|$, $C, D$ *be fixed proper subsets of* $Y$ *and* $\tilde{\varepsilon} > 0$. *Let us assume that*

$$G_0 \text{ satisfies the uniform } \tilde{\varepsilon}\text{-cone property.} \tag{3.59}$$

*We define*

$$U_{adm} = \{\overline{C} \subset G_0 \subset \overline{D} : G_0 \text{ satisfies } (3.59) \text{ and } |G_0| = \theta\},$$
$$C_\theta(D) = \{G_0 \subset D : G_0 \text{ is open, convex and } |G_0| = \theta\}.$$

*Then, at fixed volume* $\theta \in (0, |Y|)$, *there exists a domain of maximal (and minimal) effectiveness for the homogenized problem (see Chapter 1) in the class of* $G_0 \in U_{adm} \cap C_\theta(D)$.

For small (non-critical) holes, we can characterize the optimizer shape in the class of fixed volume.

**Theorem 3.12.** *For the case* $1 < \alpha < \frac{n}{n-2}$, *the ball is the domain* $G_0$ *of maximal effectiveness for a set volume in the class of star-shaped* $C^2$ *domains with fixed volume.*

**Remark 3.5.** It is a curious fact that Theorem 3.12 is opposed to the homogenization with respect to the exterior domain $\Omega$. In this context, when $\Omega$ is a ball has least effectivity, as can be shown by rearrangement techniques (see [Día85]). In the context of product domains, $\Omega = B \times \Omega''$ is the least effective on the class $\Omega = \Omega' \times \Omega''$ for set volume, at least for convex or concave kinetics as presented in Chapter 2 (see [DGC15b; DGC16; KS80]).

Through standard procedures in weak solution theory, one easily gets several results (see, e.g., [Bré71a]).

### 3.2.1 Some auxiliary results for convex domains

For the optimization, we will restrict ourselves to a general enough family of domains, but in which we can define a topology which makes the family to be compact. It is well known (see, for example, [Pir84]) that the following result holds true.

**Theorem 3.13** ([Pir84]). *The class of closed subsets of a compact set $D$ is compact in $\mathscr{P}(\mathbb{R}^n)$ for the Hausdorff convergence.*

A proof for the continuity of the effective diffusion $a_0(G_0)$ (given by (1.77)) under the Hausdorff distance in $U_{adm}$ can be found in [HD95].

**Lemma 3.2.1** ([HD95]). *If $U_{adm}$ is compact with respect to the Hausdorff metric and if $(G_0^n) \subset U_{adm}$, $G_0^m \to G_0$ as $m \to \infty$, $G_0 \in U_{adm}$, then $a_0(G_0^m) \to a_0(G_0)$ in $\mathscr{M}_n(\mathbb{R})$, where $a_0$ is the effective diffusion matrix given by* (1.77).

The behaviour of the measure $|Y \setminus G_0|$ is slightly more delicate (we include a commentary even though, in our case, this will be constant). A distance with a definition similar to Hausdorff metric is the Hausdorff complementary distance

$$d_{H^c}(\Omega_1,\Omega_2) = \sup_{x\in\mathbb{R}^n} |d(x,\Omega_1^c) - d(x,\Omega_2^c)|.$$

It has the following property: given open domains $(\Omega_m)_m, \Omega$, such that $d_{H^c}(\Omega_m,\Omega) \to 0$ as $m \to \infty$ then $\liminf_m |\Omega_m| \geq |\Omega|$. However, lower semicontinuity of the measure of the boundary ($|\partial G_0|$) is, in general, false (see [HD95] for some counterexamples). Nevertheless, the set of **convex** domains has a number of very interesting properties (see [Van04]).

**Lemma 3.2.2** ([Van04]). *The topological spaces $(C_\theta(D), d_H)$ and $(C_\theta(D), d_{H^c})$ are equivalent.*

The continuity of the boundary measure is provided by the following result, proved in [BG97].

**Lemma 3.2.3** ([BG97]). *Let $(\Omega_m), \Omega \in C_\theta(D)$. If $\Omega_1 \subset \Omega_2$, then $|\partial\Omega_1| \leq |\partial\Omega_2|$. Moreover, if $\Omega_m \overset{d_H}{\to} \Omega$, then $|\Omega_m| \to |\Omega|$ and $|\partial\Omega_m| \to |\partial\Omega|$, as $m \to \infty$.*

For the continuity of solutions with respect to $G_0$, we need the following theorem on the continuity of the associated Nemitskij operators of $g$ (see, for example, [Dal93] and [Lio69]).

**Lemma 3.2.4** ([Lio69]). *Let $g : \Omega \times \mathbb{R} \to \mathbb{R}$ be a Carathéodory function such that*

$$|g(x,v)| \leq C(1+|v|^q) \tag{3.60}$$

*holds true for $q = \frac{r}{t}$ with $r \geq 1$ and $t < \infty$. Then, the map*

$$L^r(\Omega) \to L^t(\Omega) \qquad v \mapsto g(x,v(x))$$

*is continuous in the strong topologies.*

**Lemma 3.2.5.** *Let $\mathscr{A}$ be the set of elliptic matrices and let g satisfy* (3.60) *for some* $0 \leq q \leq \frac{n}{n-2}$. *Let $u(A,\lambda)$ be the unique solution of*

$$\begin{cases} -\operatorname{div}(A\nabla u) + \lambda g(u) = f, & \text{in } \Omega, \\ u = 1, & \text{on } \partial\Omega, \end{cases}$$

*Then, the application*

$$\mathscr{A} \times \mathbb{R}_+ \to H^1(\Omega) \qquad (A,\lambda) \mapsto u(A,\lambda),$$

*is continuous in the weak topology.*

*Proof.* Let us define $G(u) = \int_0^u g(s)ds$ and

$$J_{A,\lambda}(v) = \frac{1}{2}\int_\Omega (A\nabla v)\cdot\nabla v + \int_\Omega \lambda G(v) - \int_\Omega fv.$$

We know that $u(A,\lambda)$ is the unique minimizer of this functional. Let $A_m \to A$ and $\lambda_m \to \lambda$ be two converging sequences. It is easy to prove that $u_m = u(A_m,\lambda_m)$ is bounded in $H^1(\Omega)$ and, up to a subsequence, $u_m \rightharpoonup u$ in $H^1$ as $m \to \infty$. Therefore, $\int_\Omega (A\nabla u)\cdot\nabla u \leq \liminf_m \int_\Omega (A_m \nabla u_m)\cdot\nabla u_m$. We can apply Theorem 3.2.4 to show that $G(u_m) \to G(u)$ in $L^1$ as $n \to \infty$ (see details for a similar proof, for example, in [CDLT04]) and we have that $u = u(A,\lambda)$. □

**Corollary 3.7.** *The map $(I,\lambda) \mapsto u$, where I is the identity matrix, is continuous in the weak topology of $H^1$.*

**Corollary 3.8.** *In the hypotheses of Lemma 3.2.5, the maps $(A,\lambda) \mapsto \int_\Omega g(u(A,\lambda))$ and $(I,\lambda) \mapsto \int_\Omega g(u(I,\lambda))$ are continuous.*

## 3.3 Some numerical work for the case $\alpha = 1$

The following work is part of [DGCT15].

There exists a large literature on the computation and behaviour of the homogenized coefficient $a_0(G_0)$, both from the mathematics and the engineering part (see, e.g., [ABG09], [HD95], [Kri03]). In these papers, one can find power series techniques and numerical analysis, generally for spherical obstacles. As it is common in the literature (e.g. [ABG09]), we use the commercial software COMSOL. As said on the introduction, in Nanotechnology, however, it is a common misconception that the measure of the surface alone, $|\partial G_0|$, is a good indicator of the effectiveness of the obstacle.

Considering obstacles with some symmetries (for $N = 2$ it is sufficient that they are invariant under a $90^0$ rotation) in general, it is well known that

$$a_0(G_0) = \alpha(G_0)I, \tag{3.61}$$

where $\alpha(G_0)$ is a scalar (see, for example, [ABG09], [Kri03]) and $I$ is the identity matrix in $\mathscr{M}_N(\mathbb{R})$. In this case, it can be easily proved that the effectiveness is an decreasing function of

$$\lambda(G_0) = \frac{|\partial G_0|}{\alpha(G_0)|Y \setminus G_0|} \tag{3.62}$$

(it is a direct consequence of the comparison principle, see [Día85]). In fact, this is the only relevant parameter (once $g(u)$ is fixed) of the equation (1.76). The behaviour of the effectiveness with respect to the coefficient $\lambda$ can also be numerically computed:

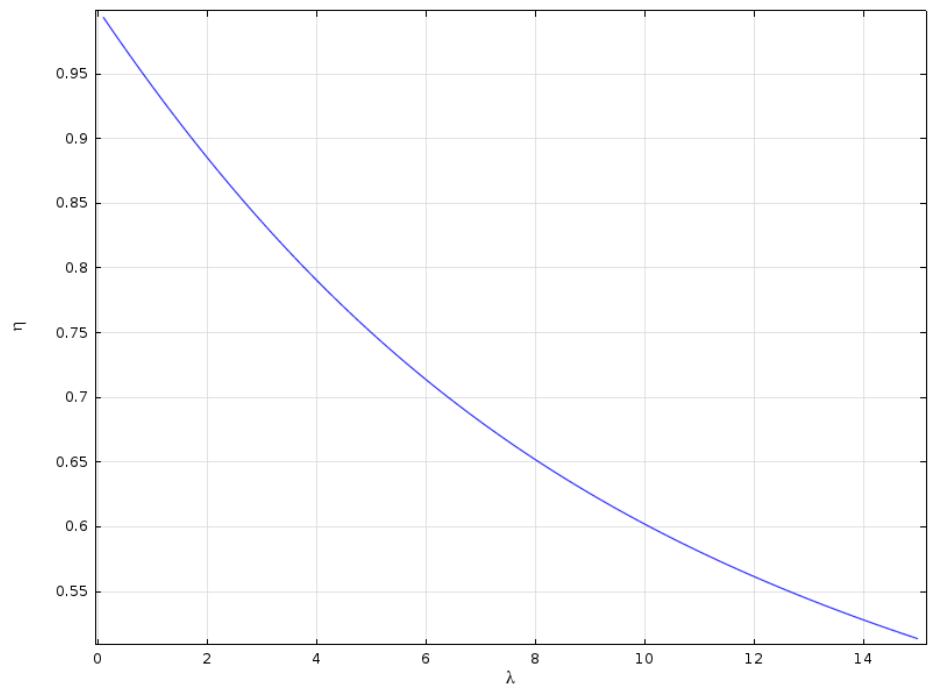

Fig. 3.1 Plot of $\eta$ as a function of $\lambda$ when $\Omega$ is a 2D circle.

Let us consider, in dimension two for simplicity, the following obstacles:

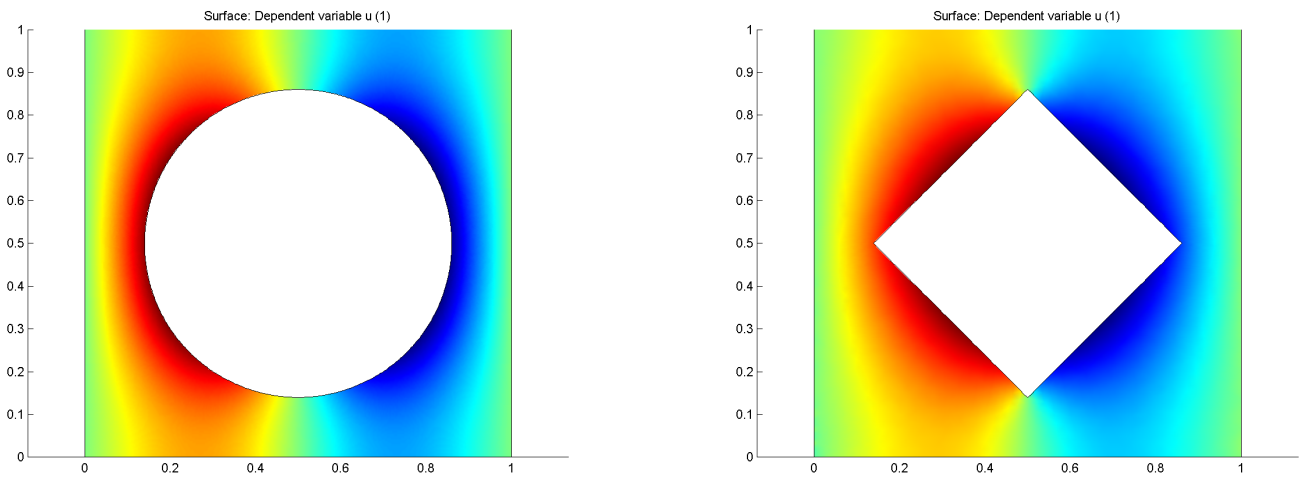

(a) Circular particle

(b) Square particle

Fig. 3.2 Two types of particle $G_0$, and the level sets of the solution of the cell problem (1.72)

We can numerically compute the homogenized diffusion coefficient $a_0(G_0)$ via a parametric sweep on the size of the particle.

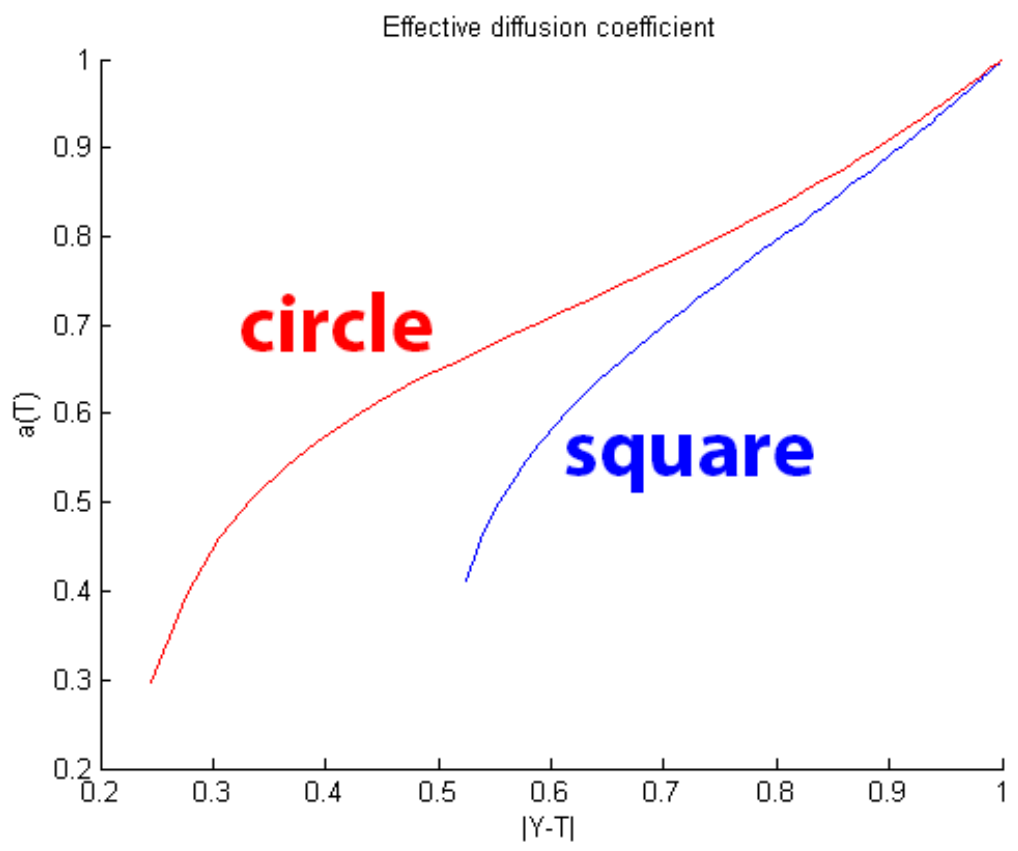

Fig. 3.3 The effective diffusion coefficient $\alpha(G_0)$ as a function of $|Y \setminus G_0|$.

Now, we can couple this with direct computations of $|\partial G_0|$ and compare the behaviour of both indicators.

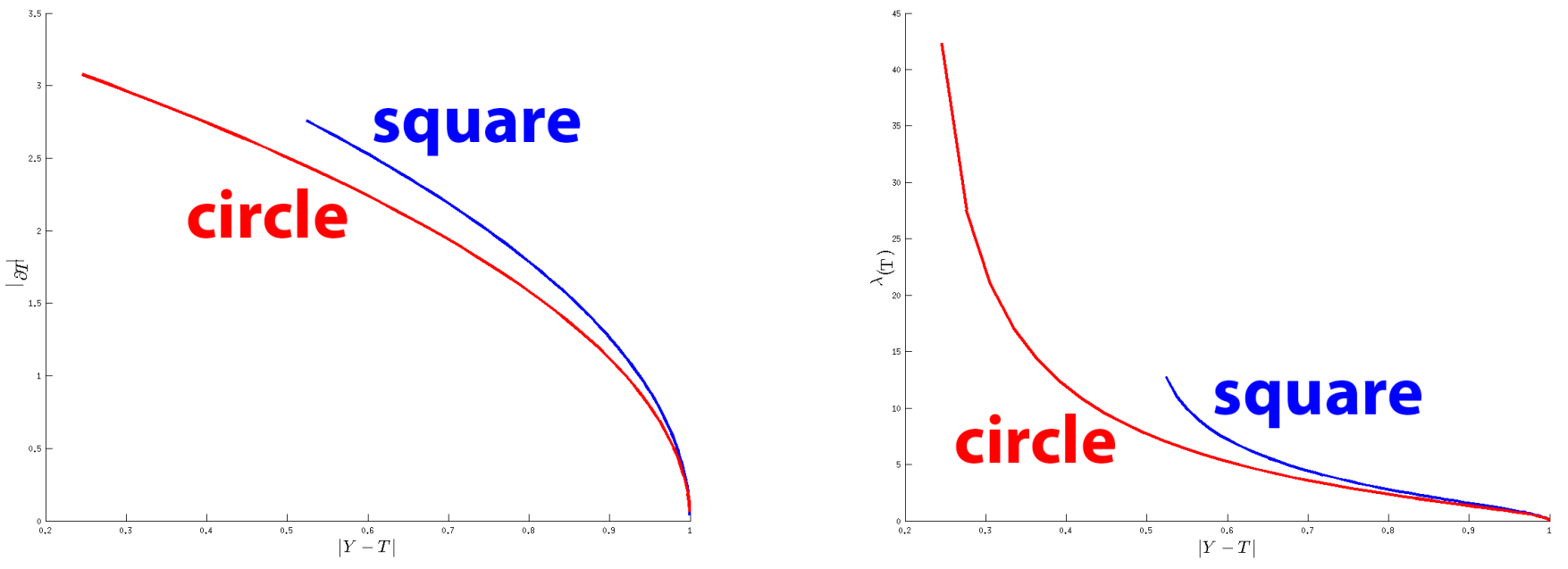

(a) Classical coefficient $|\partial G_0|$.

(b) New coefficient $\lambda(G_0)$

Fig. 3.4 Coefficients $|\partial G_0|$ and $\lambda(G_0)$ as a function of $|Y \setminus G_0|$.

# Chapter 4

# Very weak solutions of problems with transport and reaction

In studying shape differentiation when the nonlinear kinetic term $\beta(u)$ is non smooth we find that we need to understand problems of the form

$$-\Delta u + \beta'(u_0)u = f, \tag{4.1}$$

where $\beta'(u_0)$ blows up in the proximity of the boundary of the dead core. Since $\beta'(u_0)$ is a priori known (before we study $u$) we can define $V(x) = \beta'(u_0)$. The expected behaviour is, in the blow up case that $V(x) \sim d(x, \partial\Omega)^{-\alpha}$ where $\alpha > 0$. Thus, we become interested in the study of the problem

$$-\Delta u + V(x)u = f. \tag{4.2}$$

## 4.1 The origin of very weak solutions

The notion of very weak solution with data $f$ such that $fd(\cdot, \partial\Omega) \in L^1(\Omega)$ first appears in an unpublished paper by Brézis [Bré71b], and was later presented in [BCMR96]. Let

$$\delta(x) = d(x, \partial\Omega), \qquad x \in \Omega. \tag{4.3}$$

If $u \in \mathscr{C}^2(\overline{\Omega})$ is a solution of the following problem

$$\begin{cases} -\Delta u = f & \Omega, \\ u = u_0 & \partial\Omega, \end{cases} \tag{4.4}$$

then, integrating twice by parts we obtain that

$$-\int_{\Omega} u\Delta\varphi = \int_{\Omega} f\varphi - \int_{\partial\Omega} u_0 \frac{\partial\varphi}{\partial n} \tag{4.5}$$

for every $\varphi \in W^{2,\infty}(\Omega)\cap W_0^{1,\infty}(\Omega)$. A fortiori, even if $u$ is not of class $\mathscr{C}^2$, since, for this test functions, $\frac{\varphi}{\delta} \in L^\infty(\Omega)$, the problem (4.5) is well formulated for data $f$ and $u_0$ such that $f\delta \in L^1(\Omega)$ and $u_0 \in L^1(\partial\Omega)$. Equation (4.5) is known as the very weak formulation of problem (4.4).

The surprising new result introduced by Brézis in 1971 is

**Theorem 4.1** ([Bré71b])**.** *Let $f$ be measurable such that $\delta f \in L^1(\Omega)$ and let $u_0 \in L^1(\partial\Omega)$. There exists a unique $u \in L^1(\Omega)$ such that* (4.5) *for all $\varphi \in W^{2,\infty}(\Omega)\cap W_0^{1,\infty}(\Omega)$. Furthermore, there exists a constant $C > 0$ such that*

$$\|u\|_{L^1(\Omega)} \le C(\|\delta f\|_{L^1(\Omega)} + \|u_0\|_{L^1(\partial\Omega)}) \tag{4.6}$$

*Moreover, u satisfies that*

$$-\int_{\Omega} |u|\Delta\rho + \int_{\partial\Omega} |u_0| \frac{\partial\rho}{\partial n} \le \int_{\Omega} f\rho\, sign(u) \tag{4.7}$$

*for all $\rho \in W^{2,\infty}(\Omega)\cap W_0^{1,\infty}(\Omega)$, where*

$$\operatorname{sign}(s) = \begin{cases} 1 & s > 0, \\ 0 & s = 0, \\ -1 & s < 0, \end{cases} \tag{4.8}$$

He goes further in a second result that states the following:

**Theorem 4.2.** *Let $f$ be measurable in $\Omega$ such that $\delta f \in L^1(\Omega)$, $u_0 \in L^1(\partial\Omega)$ and $\beta$ monotone nondecreasing and continuous. Then there exists a unique $u \in L^1(\Omega)$ such that $\delta\beta(u) \in L^1(\Omega)$ that satisfies*

$$-\int_{\Omega} \delta\varphi + \int_{\Omega} \beta(u)\varphi = \int_{\Omega} f\varphi - \int_{\partial\Omega} u_0 \frac{\partial\varphi}{\partial n} \tag{4.9}$$

*for all* $\varphi \in W^{2,\infty}(\Omega) \cap W_0^{1,\infty}(\Omega)$*. Furthermore, if* $u$ *and* $\hat{u}$ *are two solutions corresponding to* $f, \hat{f}, u_0, \hat{u}_0$ *then*

$$\|u-\hat{u}\|_{L^1(\Omega)} + \|\delta\beta(u) - \delta\beta(\hat{u})\|_{L^1(\Omega)} \le C(\|\delta f - \delta\hat{f}\|_{L^1(\Omega)} + \|u_0 - \hat{u}_0\|_{L^1(\partial\Omega)}) \quad (4.10)$$

*where* $C$ *depends only on* $\Omega$*.*

The theory of very weak solutions developed in the $20^{\text{th}}$-century focused on the use of weighted Lebesgue spaces. Let $L^0(\Omega)$ be the space of measurable functions in $\Omega$, $\mu \in L^0(\Omega)$ and $1 \le p \le +\infty$. We define the weighted $L^p$ space as

$$L^p(\Omega,\mu) = \left\{ f \in L^0(\Omega) : \int_\Omega |f|^p \mu < +\infty \right\}. \quad (4.11)$$

However, a more modern theory will require the definition of some interpolation spaces, known as Lorentz spaces, which allow for sharp regularity results, and have nice embedding and duality properties.

## 4.2 Lorentz spaces

In order to get sharper results of regularity we introduce some interpolation spaces. Lorentz defined the following spaces in [Lor50; Lor51].

**Definition 4.1.** Given $0 < p, q \le \infty$ define

$$\|f\|_{(p,q)} = \begin{cases} \left( \int_0^\infty \left( t^{\frac{1}{p}} f^*(t) \right)^q \frac{dt}{t} \right)^{\frac{1}{q}} & q < +\infty, \\ \sup\limits_{t>0} t^{\frac{1}{p}} f^*(t) & q = +\infty, \end{cases}$$

and $L^{(p,q)}(\Omega) = \{f \text{ measurable in } \Omega : \|f\|_{(p,q)} < +\infty\}$.

There is an alternative definition of the Lorentz spaces, which is the one we have considered.

**Definition 4.2.** Let $1 \le p \le +\infty$, $1 \le q \le +\infty$ . Let $u \in L^0(\Omega)$. We define

$$||u||_{p,q} = \begin{cases} \left[ \int_{\Omega_*} \left[ t^{\frac{1}{p}} |u|_{**}(t) \right]^q \frac{dt}{t} \right]^{\frac{1}{q}} & q < +\infty, \\ \sup\limits_{0<t\le|\Omega|} t^{\frac{1}{p}} |u|_{**}(t) & q = +\infty, \end{cases} \qquad \text{where } |u|_{**}(t) = \frac{1}{t} \int_0^t |u|_*(s)ds.$$

We define

$$L^{p,q}(\Omega) = \{f \text{ measurable in } \Omega : \|u\|_{p,q} < +\infty\}.$$

These spaces are equal, and their norms equivalent, to the previously defined Lorentz spaces.

**Proposition 4.1** (Corollary 1.4.1 in [Rak08]). *Let* $1 < p \le +\infty, 1 \le q \le +\infty$. *Then*

$$L^{p,q}(\Omega) = L^{(p,q)}(\Omega)$$

*with equivalent quasi-norms.*

The functionals $\|\cdot\|_{(p,q)}$ do not, in general, satisfy the triangle inequality. However, $L^{p,q}$ is a quasi-Banach space. The following properties are known

**Proposition 4.2** ([Gra09]). *i) If* $0 < p \le \infty$ *and* $0 < q < r \le +\infty$ *then* $L^{(p,q)} \subset L^{(p,r)}$.

*ii)* $L^{(p,p)} = L^p$ *for all* $p \ge 1$.

*iii) Let* $1 \le p,q < \infty$. *Then* $(L^{(p,q)}(\Omega))' = L^{(p',q')}(\Omega)$.

*iv) If* $q < r < p$ *then* $L^{(p,\infty)}(\Omega) \cap L^{(q,\infty)}(\Omega) \subset L^r(\Omega)$ *(even for* $\Omega$ *unbounded)*

*v) If* $\Omega$ *is bounded and* $r < p$, *then* $L^{(p,\infty)}(\Omega) \subset L^r(\Omega)$

For the convenience of the reader we include an inclusion diagram for $1 \le q \le r < p < +\infty$ and $\Omega$ bounded:

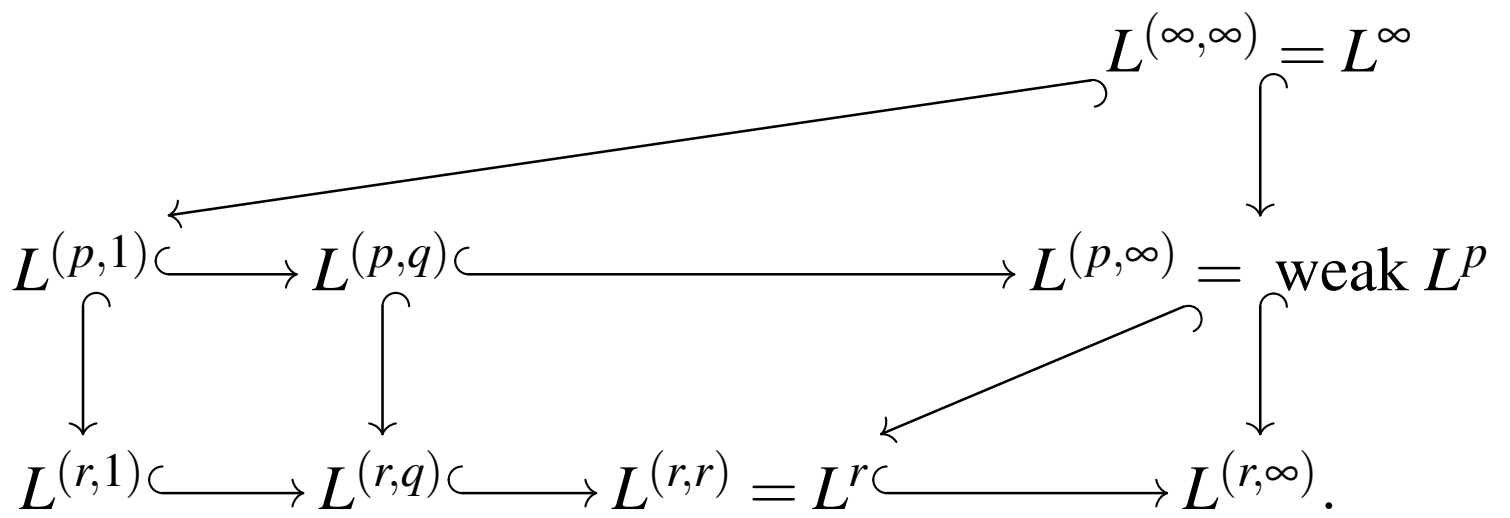

## 4.3 Modern theory of very weak solutions

Even though very weak solutions have been studied in many different contexts (see, e.g. [MV13]) the papers most linked with the research in this thesis corresponds to [DR09; DR10]. In [DR09] the regularity of very weak solutions

$$\begin{cases} -\Delta u = f & \Omega, \\ u = 0 & \partial\Omega, \end{cases} \tag{4.12}$$

is studied, which can be formulated as (4.5) when $u_0 = 0$, and solutions and its gradients are shown to be in Lorentz spaces. For the reader's convenience a brief description of this was given in Section 4.2.

Later, in [DR10] the authors tackle the problem

$$\begin{cases} -\Delta u + Vu = f & \Omega, \\ u = 0 & \partial\Omega, \end{cases} \tag{4.13}$$

where $V \geq -\lambda_1(\Omega)$, the first eigenvalue of the Laplacian, which can be written in very weak formulation as

$$\begin{cases} Vu \in L^1(\Omega, \delta), \\ -\int_\Omega u\Delta\varphi + \int_\Omega Vu\varphi = \int_\Omega f\varphi \quad \forall \varphi \in W^{2,\infty} \cap W_0^{1,\infty}(\Omega). \end{cases} \tag{4.14}$$

In [DR10] the authors prove an existence result for the most general case. As a matter of fact, they also add a nonlinear term $\beta(u)$ to (4.13).

In this chapter $\Omega \subset \mathbb{R}^n$ and

$$p' = \frac{p}{p-1}. \tag{4.15}$$

**Theorem 4.3** ([DR10]). *Let $V \in L^1_{loc}(\Omega)$ with $V \geq -\lambda \geq -\lambda_1$, where $\lambda_1$ is the first eigenvalue of $L = -\operatorname{div}(A\nabla u)$ where $A$ is symmetric, uniformly elliptic and $C^{0,1}(\bar{\Omega})$, $f \in L^1(\Omega, \delta)$. Then, there exists a solution $u \in L^{n',\infty}(\Omega) \cap W^{1,q}(\Omega, \delta)$ for every $1 \leq q < \frac{2n}{2n-1}$ satisfying*

$$\begin{cases} Vu \in L^1(\Omega, \delta), \\ \int_\Omega uL\varphi + \int_\Omega Vu\varphi = \int_\Omega f\varphi \quad \forall \varphi \in W^{2,\infty} \cap W_0^{1,\infty}(\Omega). \end{cases} \tag{4.16}$$

*Furthermore*

*i)* $\|Vu\|_{L^1(\Omega,\delta)} \leq C\|f\|_{L^1(\Omega,\delta)}$,

*ii)* $\|u\|_{L^{n',\infty}(\Omega)} C\|f\|_{L^1(\Omega,\delta)}$,

*iii) The following also holds*

$$\int_\Omega |\nabla u|^q \delta \leq C\|f\|^{\frac{q}{2}}_{L^1(\Omega,\delta)} \left(1 + \|f\|^{n'}_{L^1(\Omega,\delta)}\right)^{1-\frac{q}{2}}. \tag{4.17}$$

Nonetheless, since $\delta^{-r} \notin L^1(\Omega)$ for $r > 1$ , $\delta^{-\alpha} \notin L^1(\Omega, \delta)$ for $\alpha > 2$. Therefore, we set out to see what could be done in this case. In [DGCRT17] the author, jointly with J. I. Díaz, J. M. Rakotoson and R. Temam, solved some cases left open in [DR10].

**Theorem 4.4.** *Let $\Omega$ be bounded and $V \geq cd(x, \partial\Omega)^{-s}, s > 2$. Then there exists a unique very weak solution $u \in L^1(\Omega)$ of the problem*

$$-\Delta u + Vu = f \qquad \text{in } \Omega$$

*in the sense that*

$$\begin{cases} Vu \in L^1(\Omega, \delta), \\ -\int_\Omega u\Delta\varphi + \int_\Omega Vu\varphi = \int_\Omega f\varphi \quad \forall \varphi \in W_c^{2,\infty}(\Omega). \end{cases} \tag{4.18}$$

We will prove this result in Section 4.6.

**Remark 4.1.** Notice that the uniqueness theorem is stated without imposing any boundary conditions in a classical way (the test functions have compact support).

Later, Brezis proved the same result for $s = 2$, in personal communication to the author during his visit to Technion by an extension of the previous argument.

**Theorem 4.5.** *Let $\Omega$ be bounded and $V \geq cd(x, \partial\Omega)^{-2}$. Then there exists a unique $u \in L^1(\Omega)$ such that* (4.18)

We will include the details of the proof of this improvement[1].

### 4.3.1 Very weak solutions in problems with transport

In developing the theory, thanks to a very fruitful collaboration with J.M. Rakotoson (U. Poitiers, France) and R. Temam (U. Indiana, USA) we managed to extend the results to the problem with a transport term

$$\begin{cases} -\Delta u + \vec{b} \cdot \nabla u + Vu = f & \Omega, \\ u = 0 & \partial\Omega, \end{cases} \tag{4.19}$$

where

$$\begin{cases} \operatorname{div} \vec{b} = 0 & \Omega, \\ \vec{b} \cdot n = 0 & \partial\Omega. \end{cases} \tag{4.20}$$

[1] At the time of writing, a draft paper containing further improvements is in preparation jointly with J.I. Díaz and J.M. Rakotoson

This case is very relevant in incompressible flows. The very weak formulation of this problem can be written as

$$\int_\Omega u(-\Delta\varphi - \vec{b}\cdot\nabla\varphi + V\varphi) = \int_\Omega f\varphi \qquad \forall\varphi\in W^{2,\infty}(\Omega)\cap W_0^{1,\infty}(\Omega). \tag{4.21}$$

In order to be very clear about the definition of very weak solution, and the sense in which we define boundary conditions. We collect now some definitions:

**Definition 4.3.** Let $V, f \in L^1_{loc}$ and $\vec{b}\in L^n(\Omega)^n$, satisfy (4.20) in the sense that

$$\int_\Omega \varphi\nabla u\cdot\vec{b} = -\int_\Omega u\nabla\varphi\cdot\vec{b} \tag{4.22}$$

for all $\varphi\in W^{1,\infty}(\Omega)$ and $u\in W^{1,n'}(\Omega)$ (see an explanation of this definition in Remark 4.2).

Let us define the following types of very weak solutions

- **Local very weak solution of** (4.19) (i.e. without boundary condition). We say that $u$ is v.w.s. without b.c. if (4.21) holds for every $\varphi\in\mathscr{C}_c^2(\Omega)$.
- **Very weak solution of** (4.19) **in the sense of weights** We say that $u$ is v.w.s. with weight if it satisfies (4.21) holds for every $\varphi\in\mathscr{C}_c^2(\Omega)$ and $Vu\in L^1(\Omega,\delta)$.
- **Very weak solution of** (4.19) **in the sense of traces** We say that $u$ is v.w.s. with Dirichlet homogeneous boundary conditions if (4.21) holds for every $\varphi\in\mathscr{C}^2(\overline{\Omega})$ such that $\varphi=0$ on $\partial\Omega$ and $Vu\in L^1(\Omega,\delta)$.

We have the following result

**Theorem 4.6.** *Let $f\in L^1(\Omega,\delta), V\in L^1_{loc}(\Omega)$ and $\vec{b}\in L^{p,1}(\Omega)$ such that* $\operatorname{div}\vec{b}=0$ *in* $\Omega$ *and* $\vec{b}\cdot n=0$ *on* $\partial\Omega$ *where either*

*i)* $p>n$ *or*

*ii)* $p=n$ *and* $\vec{b}$ *is small in* $L^{n,1}$ *(in the sense that* $\|\vec{b}\|_{n,1}\le K^{s0}_{n1}$ *for a constant specified in [DGCRT17]).*

*Then,*

*i) there exists a very weak solution without boundary condition* $u\in L^{n',\infty}(\Omega)$.

*ii) If* $V\in L^1(\Omega,\delta)$*, then there exists a v.w.s. in the sense of traces in* $u\in L^{n',\infty}(\Omega)$.

*iii) If* $V\in L^{p,1}(\Omega)$*, then there exists a unique v.w.s. in the sense of traces* $u\in L^{n',\infty}(\Omega)$.

*iv) If $V \geq c\delta^{-\alpha}$ for some $\alpha > 2$, then there exists a unique v.w.s. in the sense of weights $u \in L^{n',\infty}(\Omega)$.*

We conclude this statement section by explaining our definition of (4.20).

**Remark 4.2.** Assume first that $\vec{b}, u, \varphi$ are smooth, and that $\vec{b}$ satisfies (4.20). Then

$$\begin{aligned} \operatorname{div}(u\varphi\vec{b}) &= \nabla(u\varphi)\cdot\vec{b} + u\varphi\operatorname{div}\vec{b} = \nabla(u\varphi)\cdot b \\ &= u\nabla\varphi\cdot\vec{b} + \varphi\nabla u\cdot\vec{b}. \end{aligned}$$

Integrating over $\Omega$

$$\int_{\partial\Omega} u\varphi\vec{b}\cdot\vec{n} = \int_\Omega u\nabla\varphi\cdot\vec{b} + \int_\Omega \varphi\nabla u\cdot\vec{b}.$$

Since $\vec{b}\cdot\vec{n} = 0$ on $\partial\Omega$ we have that (4.22) holds. We can pass to the limit for less smooth $\vec{b}, u, \varphi$.

## 4.4 Existence and regularity

We will construct the solution as a limit of problems with cutoff. Let us define the cut-off operator, for $k > 0$

$$T_k(s) = \begin{cases} s & |s| \leq k, \\ k\operatorname{sign}(s) & |s| > k, \end{cases} \tag{4.23}$$

and let

$$V_k = T_k \circ V. \tag{4.24}$$

### 4.4.1 Regularity of the adjoint operator $-\Delta - \vec{u}\cdot\nabla$

Given $T \in H^{-1}(\Omega)$ we focus first on the regularity of the adjoint problem

$$\int_\Omega \nabla\phi\nabla\varphi - \int_\Omega \vec{b}\cdot\nabla\phi\varphi + \int_\Omega V\phi\varphi = \langle T, \varphi\rangle \tag{4.25}$$

$\forall\varphi \in H_0^1(\Omega)$.

By applying the Lax-Milgram theorem we can show that

**Proposition 4.3.** *Let $T \in H^{-1}(\Omega)$, $V \in L^0(\Omega)$ satisfying $V \geq -\lambda > -\lambda_1$ (where $\lambda_1$ is the first eigenvalue of $-\Delta$ with Dirichlet boundary condition). Let*

$$\mathscr{W} = \{\varphi \in H_0^1(\Omega) : (V+\lambda)\phi^2 \in L^1(\Omega)\} \tag{4.26}$$

*endowed with*

$$[\varphi]_{\mathscr{W}}^2 = \|\varphi\|_{H_0^1(\Omega)}^2 + \int_\Omega (V+\lambda)\varphi^2. \tag{4.27}$$

*Then there exists a unique element $\phi \in \mathscr{W}$ such that* (4.25) *holds for every $\varphi \in \mathscr{W}$. Moreover*

$$\|\phi\|_{H_0^1(\Omega)} \leq \frac{\lambda_1}{\lambda_1 - \lambda}\|T\|_{H^{-1}}, \tag{4.28}$$

$$\left(\int_\Omega (V+\lambda)\varphi^2\right)^{\frac{1}{2}} \leq \left(\frac{\lambda_1}{\lambda_1-\lambda}\right)^{\frac{1}{2}}\|T\|_{H^{-1}}. \tag{4.29}$$

It is clear that if $T \in H^{-1}$ then there exists a unique solution $\phi_k \in H_0^1(\Omega)$ of

$$\int_\Omega \nabla\phi_k\nabla\varphi - \int_\Omega \vec{b}\nabla\phi_k\varphi + \int_\Omega V_k\phi_k\varphi = \langle T, \varphi\rangle, \qquad \forall\varphi \in H_0^1(\Omega). \tag{4.30}$$

It turns out that $\phi_k \to \phi$ strongly in $H_0^1(\Omega)$. We can show that the regularity can be improved

**Proposition 4.4.** *Let $T \in L^{\frac{n}{2},1}(\Omega) \subset H^{-1}(\Omega)$ and $V \geq 0$. Then $\phi \in L^\infty(\Omega)$ and there exists a constant $C = C(n,\Omega)$ such that*

$$\|\phi\|_{L^\infty(\Omega)} \leq C\|T\|_{L^{\frac{n}{2},1}(\Omega)}. \tag{4.31}$$

**Proposition 4.5.** *Let $V \in L^0(\Omega)$ and*

$$T = -\operatorname{div}\vec{F} \qquad \vec{F} \in L_F = \begin{cases} L^{n,1}(\Omega)^n & n \geq 3, \\ L^{2+\varepsilon}(\Omega)^2 & n = 2. \end{cases} \tag{4.32}$$

*Then $\phi \in L^\infty(\Omega)$ and there exists a constant $C = C(n,\Omega)$ such that*

$$\|\phi\|_{L^\infty(\Omega)} \leq C\|\vec{F}\|_{L_F}. \tag{4.33}$$

**Proposition 4.6.** *Let*

*i)* $\vec{b} \in L^{p,q}(\Omega)$ *with* $p > n$,

*ii)* $0 \leq V \in L^{r,q}(\Omega)$ *where* $r = \frac{np}{n+p}$,

*iii)* $T = -\operatorname{div}\vec{F}$ *where* $\vec{F} \in L^{p,q}(\Omega)$ *for* $1 \le q \le +\infty$.

*Then* $\phi \in W^1L^{p,q}(\Omega)$. *Moreover, there exists* $K_{pq} = K(p,q,n,\Omega)$ *such that*

$$\|\nabla\phi\|_{L^{p,q}} \le K_{pq}(1 + \|\vec{b}\|_{L^{p,q}} + \|V\|_{L^{r,q}})\|F\|_{L^{p,q}(\Omega)^n}. \tag{4.34}$$

**Proposition 4.7.** *If* $\vec{b}, \vec{F} \in L^{p,\infty}(\Omega)^n$ *for* $p > n$ *then* $\phi \in \mathscr{C}^{0,\alpha}(\bar{\Omega})$ *for* $\alpha = 1 - \frac{n}{p}$.

As an auxiliary space we will use the spaces of bounded mean oscillation

**Definition 4.4.** A locally integrable function $f$ on $\mathbb{R}^n$ is said to be in $\mathrm{bmo}(\mathbb{R}^N)$ if

$$\|f\|_{\mathrm{bmo}(\mathbb{R}^N)} = \sup_{0<\mathrm{diam}\,(Q)<1} \frac{1}{|Q|}\int_Q |f(x) - f_Q|\,dx + \sup_{\mathrm{diam}\,(Q)\ge 1} \frac{1}{|Q|}\int_Q |f(x)|dx < +\infty$$

where the supremum is taken over all cube $Q \subset \mathbb{R}^n$ the sides of which are parallel to the coordinates axes and

$$f_Q = \frac{1}{|Q|}\int_Q f(y)dy.$$

**Definition 4.5.** A locally integrable function $f$ on a Lipschitz bounded domain $\Omega$ is said to be in $\mathrm{bmo}_r(\Omega)$ ($r$ stands for restriction) if

$$\|f\|_{\mathrm{bmo}_r(\Omega)} = \sup_{0<\mathrm{diam}\,(Q)<1} \frac{1}{|Q|}\int_Q |f(x) - f_Q|\,dx + \int_\Omega |f(x)|dx < +\infty, \tag{4.35}$$

where the supremum is taken over all cube $Q \subset \Omega$ the sides of which are parallel to the coordinates axes.

In this case, there exists a function $\widetilde{f} \in \mathrm{bmo}(\mathbb{R}^N)$ such that

$$\widetilde{f}\Big|_\Omega = f \text{ and } \|\widetilde{f}\|_{\mathrm{bmo}(\mathbb{R}^N)} \le c_\Omega \cdot \|f\|_{\mathrm{bmo}_r(\Omega)}. \tag{4.36}$$

**Proposition 4.8.** *Let* $\vec{b}, \vec{F} \in \mathrm{bmo}_r(\Omega)^n$ *and* $V \in \mathrm{bmo}_r(\Omega)$. *Then:*

*i)* $\vec{b}\phi \in \mathrm{bmo}_r(\Omega)^n$,

*ii)* $\nabla\phi \in \mathrm{bmo}_r(\Omega)^n$.

We can even estimate some second order derivatives

**Proposition 4.9.** *Let* $\vec{b} \in L^{p,q}$, $T, V \in L^{p,q}(\Omega)$ *for some* $p > n$ *and* $1 \le q \le +\infty$. *Then* $\phi \in W^2L^{p,q}(\Omega)$ *and there exists* $K = K(p,q,n,\Omega)$ *such that*

$$\|\phi\|_{W^2L^{p,q}} \le K\frac{1 + c_{\varepsilon_0}\|\vec{b}\|_{L^{p,q}} + \|V\|_{L^{p,q}}}{1 - \varepsilon_0\|\vec{b}\|_{L^{p,q}}}\|T\|_{L^{p,q}}, \tag{4.37}$$

*where* $\varepsilon_0 > 0$ *is such that* $\varepsilon_0 \|\vec{b}\|_{L^{p,q}} < 1$ *and* $c_{\varepsilon_0} \to +\infty$ *as* $\varepsilon \to 0$.

For the proof of the existence of solutions in Theorem 4.6 the idea is to consider $u_k$ the solution of

$$-\Delta u_k + \vec{b}_k \cdot \nabla u_k + V_k u_k = f_k \tag{4.38}$$

where $f_k = T_k \circ f$ and $\vec{b}_k$ is an approximating sequence in

$$\mathscr{V} = \{\vec{b} \in \mathscr{C}_c^\infty(\Omega)^n : \operatorname{div} \vec{b} = 0 \text{ in } \Omega\}, \tag{4.39}$$

which has adherence in $L^{p,q}$ the set

V = { $\in L^{p,q}(\Omega)^n : \operatorname{div} \vec{b} = 0$ in $\Omega, \vec{b} \cdot \vec{n} = 0$ on $\partial\Omega\}$.(4.40)In order to get some uniform estimate (we want to apply the Dunford–Pettis compactness theorem) we consider the family of test functions

$$\begin{cases} -\Delta \phi_{k,E} - \vec{b}_j \cdot \nabla \phi_{k,E} = \chi_E & \Omega, \\ \phi_{k,E} = 0 & \partial\Omega. \end{cases} \tag{4.41}$$

The previously established regularity result assure that

$$\|\phi_{k,E}\|_{W^2 L^{n,1}} \le C \|\chi_E\|_{L^{n,1}} \le C|E|^{\frac{1}{n}}. \tag{4.42}$$

From this reasoning we can extract some conclusions (see [DGCRT17] for the details)

$$\int_E u_j \le C|E|^{\frac{1}{n}} \int_\Omega f_k \delta, \tag{4.43}$$

$$\|u_j\|_{L^{n',\infty}} \le C \int_\Omega f\delta, \tag{4.44}$$

$$\int_\Omega V_k u_k \delta \le C(1 + \|\vec{b}\|_{L^{n,1}}) \int_\Omega f\delta. \tag{4.45}$$

With some additional work we show that there exists $u \in L^1(\Omega)$ such that

i) $u_j \to u$ in $L^1(\Omega)$,

ii) $Vu \in L^1(\Omega)$ (by applying Fatou's lemma),

iii) $V_k u_k \delta \rightharpoonup Vu\delta$ in $L^1_{loc}(\Omega)$,

iv) if $V \in L^1(\Omega, \delta)$, then $V_k u_k \delta \rightharpoonup Vu\delta$ in $L^1(\Omega)$.

This is enough to show that $u$ is a v.w.s. without boundary condition. If $V \in L^1(\Omega,\delta)$ then $u$ is a v.w.s. in the sense of traces. If $V \geq C\delta^{-\alpha}$ with $\alpha > 2$ then

$$+\infty > \int_\Omega |u| V\delta \geq C\int_\Omega |u|\delta^{1-\alpha} \tag{4.46}$$

hence $u \in L^1(\Omega,\delta^{-r})$ for some $r > 1$. For the uniqueness of solutions in this last case we must set a suitable theory.

## 4.5 Maximum principles in some weighted spaces

The classical maximum principle states that, if $u \in \mathscr{C}^2(\Omega)$ and

$$\begin{cases} -\Delta u \leq 0 & \text{in } \Omega, \\ u \leq 0 & \text{on } \partial\Omega, \end{cases} \implies u \leq 0. \tag{4.47}$$

We will say that a space $X$ satisfies a maximum principle if

$$\begin{cases} -\Delta u \leq 0 & \text{in } \mathscr{D}'(\Omega), \\ u \in X & \text{on } \partial\Omega, \end{cases} \implies u \leq 0. \tag{4.48}$$

Since it will be used very subtly in the following sections, we recall the following definition

**Definition 4.6.** Let $u$ be an integrable function. We say that $-\Delta u = f$ in $\mathscr{D}'(\Omega)$ if

$$-\int_\Omega u\Delta\varphi = \int_\Omega f\varphi \qquad \forall \varphi \in \mathscr{C}_c^\infty(\Omega).$$

To show that some spaces satisfy the property above, let us state an approximation lemma for the space of test functions.

**Remark 4.3.** One of the useful properties of (4.47) is that it allows to prove uniqueness of solutions straightforwardly. Let $f \in \mathscr{C}(\overline{\Omega})$. Consider two solutions $u_i$ such that

$$\begin{cases} -\Delta u_i = f & \text{in } \Omega, \\ u_i = 0 & \text{on } \partial\Omega. \end{cases} \tag{4.49}$$

Then, one immediately proves that $u_1 - u_2 \leq 0$ and $u_2 - u_1 \leq 0$. Hence $u_1 = u_2$.

In order to prove the relevant results in this Chapter we will use the following maximum principle:

**Theorem 4.7.** *Let $u \in L^1(\Omega)$ be such that*

$$-\int_\Omega u\Delta\varphi \leq 0, \qquad \forall \varphi \geq 0,\ \varphi \in W_0^{1,\infty}(\Omega) \cap W^{2,\infty}(\Omega). \tag{4.50}$$

*Then $u \leq 0$.*

A very useful result to be used in conjuction with this kind of maximum principle is known as Kato's inequality (which was originally published in [Kat72]). To present it we give the definition of the *positive sign function*: for $s \in \mathbb{R}$

$$\operatorname{sign}_+(s) = \begin{cases} 1 & s > 0, \\ 0 & s \leq 0. \end{cases} \tag{4.51}$$

**Theorem 4.8** (Kato's inequality as presented in [MV13])**.** *Assume that $u, f \in L^1_{loc}(\Omega)$ and $-\Delta u \leq f$ in $\mathscr{D}'(\Omega)$. Then:*

*i) $-\Delta |u| \leq f \operatorname{sign} u$ in $\mathscr{D}'(\Omega)$.*

*ii) $-\Delta u_+ \leq f \operatorname{sign}_+ u$ in $\mathscr{D}'(\Omega)$.*

**Remark 4.4.** It is very important to compare the test functions in (4.50) with the ones of the definition of $-\Delta u \leq 0$ in $\mathscr{D}'(\Omega)$.

### 4.5.1 Some approximation lemmas

**Approximation in $W_0^{1,\infty}$ with weights** In [DGCRT17] we proved the following result, which is stated for the spaces

$$W_c^{m,\infty}(\Omega, \delta^r) = \{f \in W^{m,\infty}(\Omega, \delta^r) : \exists K \subset \Omega \text{ compact such that } f = 0 \text{ a.e. in } \Omega \setminus K\}.$$

**Theorem 4.9.** *The following density results hold:*

*i) Let $r > m$. Then $W_c^{m,\infty}(\Omega, \delta^r)$ is dense in $W^{m,\infty}(\Omega, \delta^r)$*

*ii) Let $r > m - 1$. Then $W_c^{m,\infty}(\Omega, \delta^r)$ is dense in $W_0^{1,\infty}(\Omega) \cap W^{m,\infty}(\Omega, \delta^r)$.*

**Remark 4.5.** Notice that, without the weights, the results do not hold. If a sequence $(f_n) \in W_c^{m,\infty}(\Omega)$ converges to a function $f$ in the norm of this space, due to the continuity of trace $f \in W_0^{1,\infty}(\Omega)$. Hence, the adherence of $W_c^{m,\infty}(\Omega)$ with the $W^{m,\infty}(\Omega)$ norm can not be the whole space $W^{m,\infty}(\Omega)$.

We will prove that

**Proposition 4.10.** *Let $\varphi \in W^{m,\infty}$, then for $r > m$ there exists $(\varphi_n) \subset W_c^{m,\infty}$ such that*

$$\delta^r(\partial_\alpha \varphi_n) \xrightarrow{L^\infty} \delta^r(\partial_\alpha \varphi), \qquad |\alpha| < r, m.$$

*If $\varphi \in W_0^{1,\infty} \cap W^{m,\infty}$ then*

$$\delta^r(\partial_\alpha \varphi_n) \xrightarrow{L^\infty} \delta^r(\partial_\alpha \varphi), \qquad |\alpha| < r+1, m.$$

The proof is based on the existence and bounds of the cut-off function we will define now. Let $\psi \in \mathscr{C}^\infty(\mathbb{R})$ be a non decreasing function such that $0 \le \psi \le 1$ and

$$\psi(s) = \begin{cases} 1, & s \ge 1, \\ 0, & s \le 0. \end{cases}$$

Let, for $x \in \Omega$,

$$\eta_\varepsilon(x) = \psi\left(\frac{\delta(x) - \varepsilon}{\varepsilon}\right).$$

We have constructed a function $\eta_\varepsilon$ which will be of relevance to us.

**Lemma 4.5.1.** *Let $\Omega$ be such that $\partial\Omega \in \mathscr{C}^2$. Then, there exists a sequence of function $\eta_\varepsilon$ such that*

*i)* $\operatorname{supp} \eta_\varepsilon \subset \{\delta \ge \varepsilon\}$,

*ii)* $\operatorname{supp}(1 - \eta_\varepsilon) \subset \{\delta \le 2\varepsilon\}$,

*iii)* $|D^\alpha \eta_\varepsilon(x)| \le C\varepsilon^{-|\alpha|}$.

The approximating sequence that we construct to prove Proposition 4.10 is precisely, for $\varphi \in W^{m,\infty}(\Omega)$, given by $\varphi_n = \eta_{\frac{1}{n}} \varphi$. The details of the proof (which requires several sharp estimations) can be found [DGCRT17].

**Approximation in $L^1(\Omega, \delta)'$** The mentioned improvement by Brezis is the following.

**Theorem 4.10.** *Let $u \in L^1(\Omega, \delta^{-1})$ and $\varphi \in W^{2,\infty}(\Omega) \cap W_0^{1,\infty}(\Omega)$. Then*

$$\int_\Omega u\Delta(\varphi\eta_\varepsilon) \to \int_\Omega u\Delta\varphi. \tag{4.52}$$

*Proof.* Taking into account that

$$\Delta(\eta_\varepsilon \varphi) = \eta_\varepsilon \Delta \varphi + 2\nabla \eta_\varepsilon \cdot \nabla \varphi + \varphi \Delta \eta_\varepsilon, \tag{4.53}$$

we have that

$$-\int_\Omega u\Delta(\eta_\varepsilon \varphi) = -\int_\Omega u\eta_\varepsilon \Delta\varphi - 2\int_\Omega u\nabla\eta_\varepsilon \cdot \nabla\varphi - \int_\Omega u\varphi\Delta\eta_\varepsilon. \tag{4.54}$$

Using the fact that $\frac{u}{\delta} \in L^1(\Omega)$, $\delta\eta_\varepsilon \to \delta$ in $L^\infty(\Omega)$ and $\Delta\varphi \in L^\infty(\Omega)$:

$$-\int_\Omega u\eta_\varepsilon\Delta\varphi = -\int_\Omega \frac{u}{\delta}\delta\eta_\varepsilon\Delta\varphi \to -\int_\Omega \frac{u}{\delta}\delta\Delta\varphi = -\int_\Omega u\Delta\varphi.$$

On the other hand

$$\begin{aligned}
\left|\int_\Omega u\nabla\eta_\varepsilon \cdot \nabla\varphi\right| &\le \left|\int_{\{\varepsilon<\delta<2\varepsilon\}} \frac{u}{\delta}\delta\nabla\eta_\varepsilon \cdot \nabla\varphi\right| \\
&\le \int_{\{\varepsilon<\delta<2\varepsilon\}} \frac{|u|}{\delta}\delta\|\nabla\eta_\varepsilon\|_{L^\infty(\Omega)}\|\nabla\varphi\|_{L^\infty(\Omega)} \\
&\le C\varepsilon\|\nabla\eta_\varepsilon\|_{L^\infty(\Omega)}\int_{\{\varepsilon<\delta<2\varepsilon\}} \frac{|u|}{\delta} \\
&\le C\int_{\{\varepsilon<\delta<2\varepsilon\}} \frac{|u|}{\delta}.
\end{aligned}$$

Since $u/\delta \in L^1(\Omega)$ and the Lebesgue measure $m(\{\varepsilon < \delta < 2\varepsilon\}) \to 0$ we have

$$\int_\Omega u\nabla\eta_\varepsilon \cdot \nabla\varphi \to 0 \text{ as } \varepsilon \to 0.$$

Due the Hardy inequality in $W_0^{1,\infty}(\Omega)$ we have that $\frac{\varphi}{\delta} \in L^\infty(\Omega)$. Therefore

$$\begin{aligned}
\left|\int_\Omega u\varphi\Delta\eta_\varepsilon\right| &\le \left|\int_{\{\varepsilon<\delta<2\varepsilon\}} \frac{u}{\delta}\frac{\varphi}{\delta}\delta^2\Delta\eta_\varepsilon\right| \\
&\le \left\|\frac{\varphi}{\delta}\right\|_{L^\infty(\Omega)} \varepsilon^2\|\Delta\eta_\varepsilon\|_{L^\infty(\Omega)}\int_{\{\varepsilon<\delta<2\varepsilon\}} \frac{|u|}{\delta} \\
&\le C\int_{\{\varepsilon<\delta<2\varepsilon\}} \frac{|u|}{\delta} \to 0 \text{ as } \varepsilon \to 0.
\end{aligned}$$

This concludes the proof. □

### 4.5.2 Maximum principle of $-\Delta$ in $L^1$ with weights and without boundary condition

In [DGCRT17] we proved a first result in this direction, which we will write following the definitions in [MV13].

**Theorem 4.11** ([DGCRT17])**.** *Let $u \in L^1(\Omega, \delta^{-r})$ for some $r > 1$ be such that $-\Delta u \le 0$ in $\mathscr{D}'(\Omega)$, i.e.*

$$-\int_\Omega u\Delta\varphi \le 0, \qquad \forall \varphi \in \mathscr{C}_c^\infty(\Omega), \quad \varphi \ge 0. \tag{4.55}$$

*Then $u \le 0$.*

*Proof.* Assume first that $r > 1$. Let $\varphi \in W_0^{1,\infty} \cap W^{2,\infty}$ and let $\varphi_n \in W_c^{2,\infty}$ be the approximating sequence constructed in Proposition 4.10 (e.g. $\eta_{\frac{1}{n}}\varphi$ where $\eta_\varepsilon$ is given by Lemma 4.5.1). Then

$$0 \ge -\int_\Omega u\Delta\varphi_n = -\int_\Omega u\delta^{-r}(\Delta\varphi_n)\delta^r.$$

Since $u\delta^{-r} \in L^1(\Omega)$ and $\delta^r\Delta\varphi_n \to \delta^r\Delta\varphi$ in $L^\infty$ we can pass to the limit and obtain

$$0 \ge -\int_\Omega u\delta^{-r}\Delta\varphi\delta^r = -\int_\Omega u\Delta\varphi,$$

which proves the result. □

Combining this fact with Theorem 4.8 we have the following result (without boundary condition)

**Corollary 4.1.** *Let $u \in L^1(\Omega, \delta^{-r})$ for some $r > 1$ be such that $-\Delta|u| \le 0$ then $u = 0$.*

Prof. Häim Brezis improved Theorem 4.11 in a personal communication. The proof is a refinement of the one in [DGCRT17].

**Theorem 4.12.** *Let $u \in L^1(\Omega, \delta^{-1})$ be such that $-\Delta u \le 0$. Then $u \le 0$.*

*Proof.* Let $0 \le \varphi \in W^{2,\infty}(\Omega) \cap W_0^{1,\infty}(\Omega)$. Since $0 \le \eta_\varepsilon\varphi \in W_c^{2,\infty}(\Omega)$ we can use it as a test function. We have that

$$-\int_\Omega u\Delta(\varphi\eta_\varepsilon) \le 0.$$

Therefore, due to Theorem 4.10 and the previous estimates,

$$0 \ge -\int_\Omega u\Delta(\varphi\eta_\varepsilon) \to -\int_\Omega u\Delta\varphi.$$

Finally, for any $\varphi \in W^{2,\infty}(\Omega) \cap W_0^{1,\infty}(\Omega)$, we have that

$$-\int_\Omega u\Delta\varphi \leq 0.$$

Due to Theorem 4.7, we have that $u \leq 0$. □

### 4.5.3 Maximum principle of $-\Delta + \vec{b}\cdot\nabla$ in $L^1$ with weights

**Theorem 4.13** ([DGCRT17]). *Let $\bar{u} \in W^{1,1}_{loc}(\Omega)$ and $\vec{b}\bar{u} \in L^1_{loc}(\Omega)$ and*

$$L\bar{u} = -\Delta\bar{u} + \operatorname{div}(\vec{b}\bar{u}) \in L^1_{loc}(\Omega). \tag{4.56}$$

*Define the dual operator*

$$L^*\bar{\psi} = -\Delta\bar{\psi} - \vec{b}\cdot\nabla\bar{\psi}. \tag{4.57}$$

*Then*

*i) For all $\psi \in \mathscr{D}(\Omega)$, $\psi \geq 0$ we have that*

$$\int_\Omega \bar{u}_+ L^*\psi \leq \int_\Omega \psi \operatorname{sign}_+(\bar{u})L\bar{u}. \tag{4.58}$$

*That is, $L\bar{u}_+ \leq \operatorname{sign}_+(u_+)L\bar{u}$ in $\mathscr{D}'(\Omega)$.*

*ii) $L(|\bar{u}|) \leq \operatorname{sign}(\bar{u})L\bar{u}$ in $\mathscr{D}'(\Omega)$.*

## 4.6 Uniqueness of very weak solutions of problem (4.18)

We provide here the proof of the extended result Theorem 4.5, which has not been published. For the proof of Theorem 4.4 can be found in [DGCRT17].

*Proof of Theorem 4.5.* The existence result was shown in [DGCRT17]. Since the problem is linear let us show uniqueness for $f = 0$. Since $Vu \in L^1(\Omega,\delta)$ we have that $u \in L^1(\Omega,\delta^{-1})$. We have that $-\Delta u = -Vu$ in the sense of distributions. Applying Theorem 4.8 we have that

$$-\Delta|u| \leq -(\operatorname{sign} u)Vu = -V|u| \leq 0 \tag{4.59}$$

Applying Theorem 4.12 we have $u = 0$. □

# 4.7 On weights and traces

It was noted on [Kuf85] that some power type weights in $L^1$ induce zero trace on continuous functions $L^1(\Omega,\delta^{-r})\cap\mathscr{C}(\bar\Omega)\subset\mathscr{C}_0(\bar\Omega)$. This naturally raises the question: Is $u\delta^{-\alpha}\in L^p$ for some $\alpha$ and $p$ a sufficient condition to have uniqueness in an elliptic equation even when the solution does not neccesarily have a trace? Does the weight work as a trace, even when there is no trace? A number of results in this directions are provided in [Kuf85], for $p>1$

$$W^{k,p}(\Omega,\delta^{\varepsilon}) = W_0^{k,p}(\Omega,\delta^{\varepsilon})(=\overline{C_0^\infty(\Omega)}^{W^{k,p}(\Omega,\delta^{\varepsilon})}) \tag{4.60}$$

and

$$u\in W^{1,p}(\Omega) \quad \text{and} \quad \frac{u}{\delta}\in L^p(\Omega) \iff u\in W_0^{1,p}(\Omega). \tag{4.61}$$

In this sense it is natural that something like this might be used as a boundary condition. In some cases, the fact that the solution is in such a weighted appears naturally.

In fact, due to Theorem 4.12, weights in the form of negative powers of the distance to the boundary can be used to define “Dirichlet boundary conditions” for elliptic equations and ensure uniqueness. In particular, our aim is to show that, if we assume

$$\begin{cases} -\Delta u = 0 \text{ in } \mathscr{D}'(\Omega),\\ u\in L^1(\Omega,\delta^{-1}), \end{cases} \qquad \text{or} \qquad \begin{cases} -\Delta u + Vu = 0 \text{ in } \mathscr{D}'(\Omega),\\ u\in L^1(\Omega,\delta^{-1}), \end{cases}$$

then $u=0$ in $\Omega$.

## 4.7.1 Weights and Hardy's inequalities in $L^1(\Omega)$

The equivalence (4.61), which holds for $p>1$ and is useful throughout this Chapter, is heavily linked with Hardy's inequality for $W_0^{1,p}(\Omega)$ for $p>1$:

$$\int_\Omega \left(\frac{|u|}{\delta}\right)^p \le C\int_\Omega |\nabla u|^p \qquad \forall u\in\mathscr{C}_c^\infty(\Omega). \tag{4.62}$$

(see [Har25; BM97]). Neither of these results is true for $p=1$ (see, e.g., [Psa13]).

The following facts for the case $p=1$ are known for a smooth bounded domain in $\mathbb{R}^n$:

i) If $u\in W_0^{1,1}(\Omega)$ and $-\Delta u=0$ in $\mathscr{D}'(\Omega)$ then $u=0$.

ii) If $u \in W^{1,1}(\Omega)$ and $\frac{u}{\delta} \in L^1(\Omega)$ then $u \in W_0^{1,1}(\Omega)$.

iii) $u \in W_0^{1,1}(\Omega)$ does not imply $\frac{u}{\delta} \in L^1(\Omega)$.

We proved that

iv) If $\frac{u}{\delta} \in L^1(\Omega)$ and $\Delta u = 0$ in $\mathscr{D}'(\Omega)$ then $u = 0$.

The following question[2] seems natural:

> Does it exist a weight function $p(x)$ such that:
>
> $$u \in W_0^{1,1}(\Omega) \implies \frac{u}{p} \in L^1(\Omega),$$
>
> and
>
> $$\left.\begin{array}{c} \Delta u = 0 \text{ in } \mathscr{D}'(\Omega) \\ \frac{u}{p} \in L^1(\Omega) \end{array}\right\} \implies u = 0, \tag{4.63}$$
>
> both hold?

If $p(x)$ satisfies (4.63), we will say that weight $\frac{1}{p}$ gives a Dirichlet boundary condition in a generalized sense.

We will focus in the case $\Omega = (0,1) \subset \mathbb{R}$. We consider the set of admissible weights:

$$\mathbb{X} = \{p \in \mathscr{C}([0,1]) : p(0) = 0, p(1) = 0, p > 0 \text{ in } (0,1)\}.$$

Naturally, the distance to the boundary is a function in this set.

The map $u \in W_0^{1,1}(0,1) \mapsto \frac{u}{p} \in L^1(0,1)$ is continuous if and only if there exists $C > 0$ such that the following Hardy-type inequality is satisfied:

$$\int_0^1 |u'| \geq C \int_0^1 \frac{|u|}{p}, \qquad \forall u \in \mathscr{C}_c^\infty(0,1). \tag{4.64}$$

**Remark 4.6.** In [Psa13] the author studies the possible nature of weights $p$ such that (4.64) holds.

We can answer negatively to the question above.

[2] which was raised to the candidate by H. Brezis in May 2017

**Theorem 4.14.** *There exists no $p \in \mathbb{X}$ that satisfies both* (4.63) *and* (4.64) *for* $\Omega = (0,1)$.

For the proof we will state several intermediate results.

**Lemma 4.7.1.** *Let $p \in \mathbb{X}$ satisfy* (4.64)*. Then $\frac{1}{p} \in L^1(0,1)$.*

*Proof.* For $0 < \varepsilon < 1$ define $u_\varepsilon = \chi_{[\varepsilon,1-\varepsilon]} \in BV(0,1)$, the characteristic function of the interval $[\varepsilon, 1-\varepsilon]$. We have that $u'_\varepsilon = \delta_\varepsilon - \delta_{1-\varepsilon}$ and $|u'| = \delta_\varepsilon + \delta_{1-\varepsilon}$. By passing to the limit by an approximating sequence in $\mathscr{C}_c^\infty(0,1)$ and applying the coarea formula (see Section 2.4), we write (4.64) as

$$2 \geq C \int_\varepsilon^{1-\varepsilon} \frac{1}{p}.$$

As $\varepsilon \to 0$ we deduce that

$$\int_0^1 \frac{1}{p} \leq \frac{2}{C}.$$

This proves the lemma. □

**Lemma 4.7.2.** *If $p \in \mathbb{X}$ satisfies* (4.63) *then $\frac{1}{p} \notin L^1$.*

*Proof.* If $\frac{1}{p} \in L^1$ then we can take $u = 1$ and (4.63) is not satisfied. □

We have the following extra information:

**Lemma 4.7.3.** *If $\frac{1}{p} \notin L^1(0, \frac{1}{2})$ and $\frac{1}{p} \notin L^1(\frac{1}{2}, 1)$ then* (4.63) *holds.*

*Proof.* Let $u \in \mathscr{D}'(0,1)$ be such that $u'' = 0$. Then $u(x) = a + bx$ for some $a, b \in \mathbb{R}$. Assume, towards a contradiction that $u \not\equiv 0$. There exists at most one $c \in [0,1]$ such that $u(c) = 0$. We distinguish 4 cases. If no $c$ exists then $|u(x)| \geq D > 0$. Then

$$+\infty > \frac{1}{D} \int_0^1 \frac{|u|}{p} \geq \int_0^1 \frac{1}{p}.$$

This is a contradiction. If $c = 0$ then $|u| \geq D > 0$ in $(\frac{1}{2}, 1)$. Then

$$+\infty > \frac{1}{D} \int_{\frac{1}{2}}^1 \frac{|u|}{p} \geq \int_{\frac{1}{2}}^1 \frac{1}{p}.$$

This is also a contradiction. The same happens if $c = 1$. If $c \in (0,1)$ then $|u| \geq D$ in $(0, \varepsilon) \cup (1-\varepsilon, 1)$. Then

$$+\infty > \frac{1}{D} \left( \int_0^\varepsilon \frac{|u|}{p} + \int_{1-\varepsilon}^1 \frac{1}{p} \right) \geq \int_0^\varepsilon \frac{1}{p} + \int_{1-\varepsilon}^1 \frac{|u|}{p}$$

This concludes the proof. □

### 4.7.2 A decomposition problem

The notions of trace and weighted boundary condition do not inter-relate. A question that emerges[3] is the follow:

> what happens if we know that a function is the sum of two parts, one satisfying a boundary condition in the sense of traces and the other one in the sense of weights (a generalized version of it).

We have the following result:

**Proposition 4.11.** *Let u satisfy the following:*

*i)* $\Delta u = 0$ *in* $\mathscr{D}'(\Omega)$.

*ii)* $u = u_1 + u_2$

*iii)* $u_1 \in W_0^{1,1}(\Omega)$

*iv)* $\frac{u_2}{\delta} \in L^1(\Omega)$.

*Then* $u = 0$.

*Proof.* Since $\frac{u_2}{\delta} \in L^1(\Omega)$, due to Theorem 4.10, it holds that

$$\int_\Omega u_2 \Delta(\eta_\varepsilon \varphi) \to \int_\Omega u_2 \Delta \varphi. \tag{4.65}$$

On the other hand, since $u_1 \in W_0^{1,1}(\Omega)$:

$$-\int_\Omega u_1 \Delta(\eta_\varepsilon \varphi) = \int_\Omega \nabla u_1 \nabla(\eta_\varepsilon \varphi) = \int_\Omega \eta_\varepsilon \nabla u_1 \nabla \varphi + \int_\Omega \varphi \nabla u_1 \nabla \eta_\varepsilon. \tag{4.66}$$

Therefore, applying the properties of $\eta_\varepsilon$ we have that:

$$\left| \int_\Omega u_1 \Delta(\eta_\varepsilon \varphi) - \int_\Omega u_1 \Delta \varphi \right| \leq \int_\Omega |1 - \eta_\varepsilon| |\nabla u_1| |\nabla \varphi| + \int_\Omega |\varphi| |\nabla u_1| |\nabla \eta_\varepsilon| \tag{4.67}$$

$$\leq \int\limits_{\delta < 2\varepsilon} |\nabla u_1 \nabla \varphi| + C \int\limits_{\varepsilon < \delta < 2\varepsilon} \delta |\nabla u_1| \varepsilon^{-1} \tag{4.68}$$

$$\leq C \left( \int\limits_{\delta < 2\varepsilon} |\nabla u_1| + \int\limits_{\varepsilon < \delta < 2\varepsilon} |\nabla u_1| \right) \tag{4.69}$$

$$\to 0, \tag{4.70}$$

[3] Raised by H. Brezis to the candidate (Haifa, June 2017)

since $|\nabla u_1| \in L^1(\Omega)$. Therefore

$$\int_\Omega u_1 \Delta(\eta_\varepsilon \varphi) \to \int_\Omega u_1 \Delta\varphi. \tag{4.71}$$

Hence

$$0 = \int_\Omega u_1 \Delta(\eta_\varepsilon \varphi) + \int_\Omega u_2 \Delta(\eta_\varepsilon \varphi) \to \int_\Omega u_1 \Delta\varphi + \int_\Omega u_2 \Delta\varphi = \int_\Omega u \Delta\varphi. \tag{4.72}$$

We have that

$$\int_\Omega u\Delta\varphi = 0 \qquad \forall \varphi \in W_0^{1,\infty}(\Omega) \cap W^{2,\infty}(\Omega). \tag{4.73}$$

Therefore $u = 0$. □

**Remark 4.7.** The conclusion of this result can be useful to prove the uniqueness of solutions of some suitable non-standard linear boundary value problems.

### 4.7.3 The $L^1$ weight as a trace operator in $W^{1,q}, q > 1$

Another approach to this problem is to study whether being in $L^1(\Omega, \delta^{-r})$ does imply having trace 0 at least for functions in $W^{1,p}(\Omega)$ for $p > 1$. In this direction, in a more functional presentation, we have proved the following new result[4]:

**Theorem 4.15.** *Let $\Omega$ be a bounded domain of class $\mathscr{C}^{0,1}$. Then, for all $r > 1$ and $q > 1$*

$$L^1(\Omega, \delta^{-r}) \cap W^{1,q}(\Omega) \hookrightarrow W_0^{1,q}(\Omega), \tag{4.74}$$

**Lemma 4.7.4.** *Let $1^* = \frac{n}{n-1}$, $n \geq 2$ and $\alpha > 1$. Then, there exists $c_\Omega$ such that for $u \in L^1(\Omega, \delta^{-\alpha}) \cap W^{1,1}(\Omega)$ one has*

$$\left( \int_\Omega \left| \frac{u}{\delta} \right|^p \right)^{\frac{1}{p}} \leq c_\Omega \|u\|_{L^{1^*}(\Omega)}^{1-\frac{1}{\alpha}} \left( \int_\Omega |u| \delta^{-\alpha} dy \right)^{\frac{1}{\alpha}}, \quad 1 < p < \min\left\{ \alpha, \frac{1^*\alpha}{\alpha - 1 + 1^*} \right\}. \tag{4.75}$$

*Proof.* By Hölder's inequality

$$\int_\Omega |u|^p \delta^{-p} = \int_\Omega |u|^{p(1-\frac{1}{\alpha})} \delta^{-p} |u|^{\frac{p}{\alpha}} \leq \left( \int_\Omega |u|^{\frac{p(\alpha-1)}{\alpha-p}} \right)^{1-\frac{p}{\alpha}} \left( \int_\Omega |u| \delta^{-\alpha} \right)^{\frac{p}{\alpha}}. \tag{4.76}$$

We impose that

$$\frac{p(\alpha-1)}{\alpha-p} \leq 1^* \tag{4.77}$$

[4] This candidate thanks J.M. Rakotoson and J.I. Díaz for their coversation on this topic.

which is exactly the condition on the statement. $\square$

*Proof of Theorem 4.15.* Let $q < \min\left\{\alpha, \frac{1^*\alpha}{\alpha-1+1^*}\right\}$. Then $\frac{u}{\delta} \in L^q$ and therefore $u \in W_0^{1,q}$. If $q \geq \min\left\{\alpha, \frac{1^*\alpha}{\alpha-1+1^*}\right\}$ first we observe that the result holds for $\bar{q}$ in the previous case and hence we see that $u \in W_0^{1,\bar{q}}(\Omega)$. Therefore

$$u \in W^{1,q}(\Omega) \cap W_0^{1,\bar{q}}(\Omega) = W_0^{1,q}(\Omega). \tag{4.78}$$

This proves the result. $\square$

**Remark 4.8.** Notice that we have substituted $L^p$ to $L^1$ in the known result (4.61).

# Part II

# A problem in Fourier representation

# Chapter 5

# Optimal basis in Fourier representation

This chapter presents work developed while on a visit to Prof. Häim Brezis at Technion - Israel Institute of Technology in Haifa, Israel in April-July 2017. The candidate wishes to extend to Häim Brezis his warmest thanks for the hospitality and the mentoring. The visit and the work led to the publication of [BGC17].

## 5.1 A problem in image representation

While studying compression of meshes for 3D representation Ron Kimmel and his group stumbled upon the following question, of a strict mathematical nature:

> Which is the basis of $L^2(\Omega)$ that provides the best finite dimensional projections of functions in $H_0^1(\Omega)$?

First, we need to define the term "optimal basis". It is natural to define as optimal a basis $b = (b_i)$ of $L^2(\Omega)$ such that, for all $m \geq 1$,

$$\left\| f - \sum_{i=1}^{m} (f, b_i) b_i \right\|_{L^2}^2 \leq \alpha_m \|\nabla f\|_{L^2}^2 \qquad \forall f \in H_0^1(\Omega). \tag{5.1}$$

with optimal constants $\alpha_m$. This technique led the group of Ron Kimmel to the publication of several paper in this direction, in collaboration with Häim Brezis (see [ABK15; ABBKS16]).

## 5.2 The mathematical treatment

In the works above the authors had shown that, in a bounded smooth set $\Omega \subset \mathbb{R}^n$, an optimal basis for $H_0^1(\Omega)$-representation in the sense of (5.1) was formed by the eigenfunctions $e_i$ of the Laplace operator

$$\begin{cases} -\Delta e_i = \lambda_i e_i & \text{in } \Omega, \\ e_i = 0 & \text{on } \partial\Omega, \end{cases} \tag{5.2}$$

where $0 < \lambda_1 < \lambda_2 \leq \lambda_3 \leq \cdots$ is the ordered sequence of eigenvalues repeated according to their multiplicity.

It is a classical result that

**Theorem 5.1.** *We have, for all $m \geq 1$,*

$$\left\| f - \sum_{i=1}^{m} (f, e_i) e_i \right\|_{L^2}^2 \leq \frac{\|\nabla f\|_{L^2}^2}{\lambda_{n+1}} \qquad \forall f \in H_0^1(\Omega). \tag{5.3}$$

The proof of this fact is tremendously simple, due to the orthogonality of the eigenfunctions. Indeed

$$\left\| f - \sum_{i=1}^{m} (f, e_i) e_i \right\|_{L^2}^2 = \left\| \sum_{i=m+1}^{+\infty} (f, e_i) e_i \right\|_{L^2}^2 = \sum_{i=m+1}^{+\infty} (f, e_i)^2$$

and

$$\|\nabla f\|_{L^2}^2 = \sum_{i=1}^{+\infty} \lambda_i (f, e_i)^2 \geq \sum_{i=m+1}^{+\infty} \lambda_i (f, e_i)^2 \geq \lambda_{m+1} \sum_{i=m+1}^{+\infty} (f, e_i)^2.$$

Combining these expressions yields the result.

The authors of [ABK15] and [ABBKS16] have investigated the “optimality” in various directions of the basis $(e_i)$, with respect to inequality (5.3). Here is one of their results restated in a slightly more general form:

**Theorem 5.2** (Theorem 3.1 in [ABK15])**.** *There is no integer $m \geq 1$, no constant $0 \leq \alpha < 1$ and no sequence $(\psi_i)_{1 \leq i \leq m}$ in $L^2(\Omega)$ such that*

$$\left\| f - \sum_{i=1}^{m} (f, \psi_i) \psi_i \right\|_{L^2}^2 \leq \frac{\alpha}{\lambda_{m+1}} \|\nabla f\|_{L^2}^2 \qquad \forall f \in H_0^1(\Omega). \tag{5.4}$$

The proof in [ABK15] relies in the Fischer-Courant max-min principle (see, e.g., [Lax02] or [Wei74]). For the convenience of the reader we present a very elementary proof based on a simple and efficient device originally due to H. Poincaré [Poi90, p. 249-250] (and later rediscovered by many people, e.g. H. Weyl [Wey12, p. 445] and R. Courant [Cou20, p. 17-18]; see also H. Weinberger [Wei74, p. 56] and P. Lax [Lax02, p. 319]).

Suppose not, and set

$$f = c_1 e_1 + c_2 e_2 + \cdots + c_m e_m + c_{m+1} e_{m+1} \tag{5.5}$$

where $c = (c_1, c_2, \cdots, c_m, c_{m+1}) \in \mathbb{R}^{m+1}$. The under-determined linear system

$$(f, \psi_i) = 0, \qquad \forall i = 1, \cdots, m \tag{5.6}$$

of $m$ equations with $m+1$ unknowns admits a non-trivial solution. Inserting $f$ into (5.4) yields

$$\lambda_{m+1} \sum_{i=1}^{m+1} c_i^2 \leq \alpha \sum_{i=1}^{m} \lambda_i c_i^2 \leq \alpha \lambda_{m+1} \sum_{i=1}^{m+1} c_i^2. \tag{5.7}$$

Therefore $\sum_{i=1}^{m+1} c_i^2 = 0$ and thus $c = 0$. A contradiction. This proves Theorem 5.2. □

The authors of [ABBKS16] were thus led to investigate the question of whether inequality (5.3) holds *only* for the orthonormal bases consisting of eigenfunctions corresponding to ordered eigenvalues. They established that a "discrete", i.e. finite-dimensional, version does hold; see [ABBKS16, Theorem 2.1]. But their proof of "uniqueness" could not be adapted to the infinite-dimensional case (because it relied on a "descending" induction). It was raised there as an open problem (see [ABBKS16, p. 1166]). The following result solves this problem.

**Theorem 5.3** ([BGC17]). *Let $(b_i)$ be an orthonormal basis of $L^2(\Omega)$ such that, for all $m \geq 1$,*

$$\left\| f - \sum_{i=1}^{m} (f, b_i) b_i \right\|_{L^2}^2 \leq \frac{\|\nabla f\|_{L^2}^2}{\lambda_{m+1}} \quad \forall f \in H_0^1(\Omega). \tag{5.8}$$

*Then, $(b_i)$ consists of an orthonormal basis of eigenfunctions of $-\Delta$ with corresponding eigenvalues $(\lambda_i)$.*

In fact, a more general result, which was introduced in [BGC17] as a remark, also holds:

**Theorem 5.4.** *Let $V$ and $H$ be Hilbert spaces such that $V \subset H$ with compact and dense inclusion ($\dim H \leq +\infty$). Let $a : V \times V \to \mathbb{R}$ be a continuous bilinear symmetric form for which there exist constants $C, \alpha > 0$ such that, for all $v \in V$,*

$$a(v,v) \geq 0,$$
$$a(v,v) + C|v|_H^2 \geq \alpha \|v\|_V^2.$$

*Let $0 \leq \lambda_1 \leq \lambda_2 \leq \cdots$ be the sequence of eigenvalues associated with the orthonormal (in $H$) eigenfunctions $e_1, e_2, \cdots \in V$, i.e.,*

$$a(e_i, v) = \lambda_i (e_i, v) \qquad \forall v \in V,$$

*where $(\cdot,\cdot)$ denotes the scalar product[1] in $H$. For every $m \geq 1$ and $f \in V$:*

$$\lambda_{m+1} \left| f - \sum_{i=1}^{m} (e_i, f) e_i \right|_H^2 \leq a(f,f). \tag{5.9}$$

*Let $(b_i)$ be an orthonormal basis of $H$ such that for all $m \geq 1$ and $f \in V$*

$$\lambda_{m+1} \left| f - \sum_{i=1}^{m} (b_i, f) b_i \right|_H^2 \leq a(f,f). \tag{5.10}$$

*Then, $(b_i)$ consists of an orthonormal basis of eigenfunctions of $a$ with corresponding eigenvalues $(\lambda_i)$.*

**Remark 5.1.** When $\dim H < +\infty$ and $V = H$ this result is originally due to [ABBKS16]. The proof of “rigidity” was quite different and could not be adapted to the infinite dimensional case. It was raised there as an open problem.

This more general formulation allows us to cover some of the most relevant situations in the applications:

- For the optimal representation of function in $H^1(\Omega)$ we must take

$$H = L^2(\Omega), \qquad V = H^1(\Omega), \qquad a(f,h) = \int_\Omega \nabla f \cdot \nabla h + \mu \int_\Omega fh \tag{5.11}$$

[1] We point that, in this general setting, it may happen that $\lambda_1 = 0$ (e.g. $-\Delta$ with Neumann boundary conditions); and $\lambda_1$ may have multiplicity greater than 1

where $\mu$ is a positive constant. Then, $e_i$ are solutions of

$$\begin{cases} -\Delta e_i + \mu e_i = \lambda_i e_i \\ \frac{\partial e_i}{\partial n} = 0. \end{cases} \tag{5.12}$$

Notice that, depending of the choice of bilinear product, we have a different choice of eigenfunctions.

- Let $\mathscr{M}$ be a compact Riemmanian manifold without boundary. Then one can choose

$$H = L^2(\mathscr{M}), \qquad V = H^1(\mathscr{M}), \qquad a(f,h) = \int_{\mathscr{M}} \nabla_g f \cdot \nabla_g h \tag{5.13}$$

where $g$ is the Riemmanian metric. Then, the basis are the solutions of

$$-\Delta_g e_i = \lambda_i e_i. \tag{5.14}$$

where $-\Delta_g$ is the Laplace-Beltrami operator. Since there is no boundary, there is no boundary condition.

The basic ingredient of our proof is the following lemma, the proof of which is based on Poincaré's *magic trick*:

**Lemma 5.2.1.** *Assume that* (5.8) *holds for all* $m \geq 1$ *and all* $f \in H_0^1(\Omega)$, *and that*

$$\lambda_i < \lambda_{i+1} \tag{5.15}$$

*for some* $i \geq 1$. *Then*

$$(b_j, e_k) = 0, \qquad \forall j,k \textit{ such that } 1 \leq j \leq i < k. \tag{5.16}$$

Applying this lemma we can quickly complete the proof in the case of simple eigenvalues. Since $\lambda_1 < \lambda_2$ then, by the lemma,

$$(b_1, e_k) = 0 \qquad \forall k \geq 2. \tag{5.17}$$

Thus $b_1 = \pm e_1$. Next we apply the lemma with $\lambda_2 < \lambda_3$. We have that

$$(b_2, e_k) = 0 \qquad \forall k \geq 3. \tag{5.18}$$

Also, we have that

$$(b_2, e_1) = \pm(b_2, b_1) = 0. \tag{5.19}$$

Therefore $b_2 = \pm e_2$. Similarly, we have that $b_i = \pm e_i$ for $i \geq 3$.

## 5.3 Connection to the Fischer-Courant principles

It is a very relevant part of the proof in [BGC17] that (5.3) can be understood under the light of the Fischer-Courant principles. In particular, if one considers the functions

$$0 \neq f \in \operatorname{span}(e_1, \cdots, e_m)^\perp$$

then, automatically,

$$\lambda_{m+1} \leq \frac{\|\nabla f\|_{L^2}}{\|f\|_{L^2}} \qquad \forall f \in \operatorname{span}(e_1, \cdots, e_m)^\perp, f \neq 0, \quad \forall m \geq 1. \tag{5.20}$$

Recall that the usual Fischer-Courant max-min principle asserts that for every $m \geq 1$ we have

$$\lambda_{m+1} = \max_{\substack{M \subset L^2(\Omega) \\ M \text{ linear space} \\ \dim M = m}} \min_{\substack{0 \neq f \in H_0^1(\Omega) \\ f \in M^\perp}} \frac{\|\nabla f\|_{L^2}^2}{\|f\|_{L^2}^2}, \tag{5.21}$$

(see, e.g., [Lax02] or [Wei74]). Therefore, in some sense our basis $b$ must be a maximizer of (5.21) for every $m \geq 1$.

Applying the same technique as in the proof of our Theorem 5.3, we can prove the following:

**Proposition 5.1.** *Let $(b_i)$ be an orthonormal sequence in $L^2(\Omega)$ such that, for every $m \geq 1$,*

$$\lambda_{m+1} = \min_{\substack{0 \neq f \in H_0^1(\Omega) \\ f \in M_m^\perp}} \frac{\|\nabla f\|_{L^2}^2}{\|f\|_{L^2}^2} \qquad \text{where } M_m = \operatorname{span}(b_1, b_2, \cdots, b_m). \tag{5.22}$$

*Then, each $b_i$ is an eigenfunction associated to $\lambda_i$.*

The natural way to establish eigen-decomposition is through a compact, symmetric operator $A : H \to H$. The resolvent operator of Dirichlet problem $A = (-\Delta)^{-1} : f \mapsto u$ where

$u$ is given by the weak solution of

$$\begin{cases} -\Delta u = f & \Omega, \\ u = 0 & \partial\Omega. \end{cases} \tag{5.23}$$

satisfies this properties with $H = L^2(\Omega)$, due to the compact embedding $H^1(\Omega) \to L^2(\Omega)$. For simplicity, we will consider $\mu_n$ its eigenvalues. Notice that

$$Ae_i = \lambda e_i \implies \frac{1}{\lambda_i} e_i = A^{-1} e_i. \tag{5.24}$$

Thus, we get that

$$\mu_m((-\Delta)^{-1}) = \frac{1}{\lambda_m(-\Delta)}. \tag{5.25}$$

The spectral theorem guaranties that $A = (-\Delta)^{-1}$ has a basis of eigenvalues that expand $L^2(\Omega)$, and the existence of a sequence of positive eigenvalues $\mu_m \to 0$. However, this guaranties the spectral decomposition for $-\Delta$.

The Courant-Fischer principles are usually written in the literature for Rayleigh quotient

$$R_A(x) = \frac{(Ax, x)}{\|x\|^2}. \tag{5.26}$$

of $A = (-\Delta)^{-1}$, rather than $(-\Delta)$. On the other hand, (5.21) is written in terms of $R_{-\Delta}$. Nonetheless, once the eigendecomposition of $A$ is established, Theorems 5.2 and 5.3 and (5.20) gives us a direct proof of (5.21).

**Remark 5.2.** The Rayleigh quotients $R_{-\Delta}$ and $R_{(-\Delta)^{-1}}$ do not seem to be directly related. Notice that

$$R_{-\Delta}(u) = \frac{(-\Delta u, u)}{\|u\|_{L^2}^2} = \frac{\|\nabla u\|_{L^2}^2}{\|u\|_{L^2}^2} \tag{5.27}$$

$$R_{(-\Delta)^{-1}}(f) = \frac{(Af, f)}{\|f\|_{L^2}^2} = \frac{(u, -\Delta u)}{\|\Delta u\|_{L^2}^2} \tag{5.28}$$

$$= \frac{\|\nabla u\|_{L^2}^2}{\|\Delta u\|_{L^2}^2}, \tag{5.29}$$

where (5.23). Nonetheless, notice that

$$R_{-\Delta}(e_i) = \lambda_i, \tag{5.30}$$

$$R_{(-\Delta)^{-1}}(e_i) = \frac{1}{\lambda_i}. \tag{5.31}$$

### 5.3.1 Some controversy about the Fischer-Courant principles

Principle (5.21) has several different presentations in the literature. Currently, there are two main presentations, which are due to Fischer in 1905 [Fis05] and Courant in 1920 [Cou20].

For the rest of the section we will focus on compact symmetric operators defined over the whole Hilbert space $H$. Many of the references provide sharper results, and only simplified versions are stated here.

Let us, first, state the principles as Lax does in [Lax02]. This appears to be the commonly accepted nomenclature.

**Theorem 5.5.** *Let $A$ be a compact symmetric operator in a Hilbert space $H$ and let $\mu_n$ be its eigenvalues. Then, the following statements hold:*

- *Fischer's principle:*

$$\mu_m = \max_{S_m} \min_{x \in S_m} R_A(x), \tag{5.32}$$

*where $S_m$ is any linear subspace of $H$ of dimension $m$*

- *Courant's principle:*

$$\mu_m = \min_{S_{m-1}} \max_{x \perp S_{m-1}} R_A(x). \tag{5.33}$$

**Remark 5.3.** It is important to notice that (5.33) with (5.21) are both the Courant principle even though max and min are in reverse order. This is due to (5.25). This relates strongly to Remark 5.2.

However, Weinberger in [Wei74] assigns the credit differently. Here, (5.32) is named Poincaré's principle (see [Wei74, Theorem 5.1]), due to Poincaré's seminal paper [Poi90] in 1890, in which he starts the theory of eigen-decomposition. Also, (5.33) is named Courant-Weyl's principle (see [Wei74, Theorem 5.2]) and it is written in a slightly more general way

$$\mu_m = \min_{\substack{l_1, \cdots, l_m \\ \text{linear functionals}}} \sup_{\substack{v \in H \\ l_1(v) = \cdots = l_m(v) = 0}} R_A(x). \tag{5.34}$$

Notice that, in infinite dimensional spaces, a linear functional $l_i$ need not be continuous, so it not be written $l_i(v) = (w, v)$ for some $w \in H$. Hence, there are many more functionals in this characterization. In Weinberger's text, the name of Fischer does not appear.

The inclusion of the name of Weyl is due to his paper [Wey12] in which he proves the asymptotic behaviour of eigenvalues (see also [Wey11]). Some books, e.g. [WS72], go as far as stating the following:

> An even more important variational characterization, the maximum- minimum principle, is claimed by Weyl, who used some of its consequences in his famous theory of asymptotic distribution of eigenvalues [W31, W32]. Later, Courant applied the principle contained in Weyl's fundamental inequality to a fairly general typical situation [C2].

In [WS72] (where the authors use $A$ as the operator with increasing eigenvalues, and thus in direct conflict with [Lax02]) the following is stated.

**Lemma 5.3.1** ([Wey12], as extracted from [WS72]). *Let $A$ be a symmetric, compact operator. Let $p_1, \cdots, p_{m-1}$ be any arbitrary vectors in $H$. Then*

$$\max_{\substack{f \in H \\ (f,p_1)=\cdots=(f,p_{m-1})=0}} R_A(f) \geq \mu_m. \tag{5.35}$$

The proof again passes by the use of Poincaré's *magic trick*. This, which is presented in [WS72] as "Weyl's lemma" must not be confused with what is usually called Weyl's lemma, that is the regularity of functions such that $-\Delta u = 0$ in $\mathscr{D}'(\Omega)$. From this result the author extracts the proof of Courant's principle.

## 5.4 Some follow-up questions

In [WS72] the problem of whether the equality can hold in (5.35) for a finite orthonormal set $b = (b_1, \cdots, b_N)$ is studied. This question was also answered in [BGC17] in some of the relevant cases.

**Remark 5.4.** After the publication of the paper the authors were made aware of the interest of this question by many authors. The terminology employed by the specialists in this field is *n-widths*. See, e.g., [Pin85; EBBH09; FS17] and the references therein.

Let us present merely the case in which only $N = 2$.

**Remark 5.5.** Assume that $b = b_1 \in L^2(\Omega)$ is such that $\|b\|_{L^2} = 1$ and

$$\|f - (f,b)b\|_{L^2}^2 \leq \frac{1}{\lambda_2}\|\nabla f\|_{L^2}^2 \qquad \forall f \in H_0^1(\Omega). \tag{5.36}$$

Of course, (5.36) holds with $b = e_1$. From Lemma 5.2.1 we know that (5.36) implies that

$$(e_2, b) = 0. \tag{5.37}$$

Clearly, (5.37) is not sufficient. Indeed, take $b = e_3$. Then, (5.37) holds but (5.36) fails for $f = e_1$. We do not have a simple characterization of the functions $b$ satisfying (5.36). But we can construct a large family of functions $b$ (which need not be smooth) such that (5.36) holds. Assume that $0 < \lambda_1 \leq \lambda_2 < \lambda_3$. Let $\chi \in L^2(\Omega)$ be any function such that

$$(e_1, \chi) = 0, \tag{5.38}$$

$$(e_2, \chi) = 0, \tag{5.39}$$

$$\|\chi\|_{L^2}^2 = 1. \tag{5.40}$$

Set

$$b = \alpha e_1 + \varepsilon\chi \qquad \alpha^2 + \varepsilon^2 = 1, \text{ with } 0 < \varepsilon < 1. \tag{5.41}$$

Then, there exists $\varepsilon_0 > 0$, depending on $(\lambda_i)_{1\leq i\leq 3}$, such that for every $0 < \varepsilon < \varepsilon_0$ (5.36) holds (see [BGC17]).

**Remark 5.6.** In the general setting of Theorem 5.4 it may happen that $0 = \lambda_1 < \lambda_2$. Suppose now that $b \in H$ is such that $\|b\|_H = 1$ and

$$\|f - (f,b)b\|_H^2 \leq \frac{1}{\lambda_2} a(f,f) \qquad \forall f \in V. \tag{5.42}$$

Claim: we have $b = \pm e_1$. Indeed, let $f = e_1$ in (5.42) we have that

$$\|e_1 - (e_1,b)b\|_H^2 \leq \frac{\lambda_1}{\lambda_2} = 0, \tag{5.43}$$

Therefore $b = \pm e_1$.

# Index

# Homogenization and Shape Differentiation of Quasilinear Elliptic Equations

**Publications**

**David Gómez Castro**

Advisor: Prof. Jesús Ildefonso Díaz Díaz

Dpto. de Matemática Aplicada &
Instituto de Matemática Interdisciplinar
Universidad Complutense de Madrid

Esta tesis se presenta dentro del
*Programa de Doctorado en Ingeniería Matemática, Estadística e Investigación Operativa*

Diciembre 2017

*The papers have not been attached to the arxiv version of this thesis due to copyright reasons.*